\documentclass[10pt]{article}
\usepackage{amsmath,amssymb} 

\begin{document}
\title{Toroidal Dehn fillings on hyperbolic 3-manifolds}
\author{Cameron McA.\ Gordon$^1$ and Ying-Qing Wu$^2$}
\date{}
\maketitle

\footnotetext[0]{ Mathematics subject classification (1991):  {\em Primary 
57N10.}}
\footnotetext[0]{ Keywords and phrases: Toroidal manifolds, Dehn fillings}
\footnotetext[1]{ Partially supported by NSF grant DMS 0305846.}
\footnotetext[2]{ Partially supported by NSF grant DMS 0203394.}

\begin{abstract}
  We determine all hyperbolic $3$-manifolds $M$ admitting two toroidal
  Dehn fillings at distance $4$ or $5$. We show that if $M$ is a
  hyperbolic $3$-manifold with a torus boundary component $T_0$, and
  $r,s$ are two slopes on $T_0$ with $\Delta(r,s) = 4$ or $5$ such
  that $M(r)$ and $M(s)$ both contain an essential torus, then $M$ is
  either one of $14$ specific manifolds $M_i$, or obtained from $M_1,
  M_2, M_3$ or $M_{14}$ by attaching a solid torus to $\partial M_i -
  T_0$. All the manifolds $M_i$ are hyperbolic, and we show that only
  the first three can be embedded into $S^3$. As a consequence, this
  leads to a complete classification of all hyperbolic knots in $S^3$
  admitting two toroidal surgeries with distance at least $4$.
\end{abstract}

\newcommand{\proof}{\noindent {\bf Proof.} }
\newcommand{\qed}{\quad $\Box$ \medskip}

\newtheorem{thm}{Theorem}[section]
\newtheorem{prop}[thm]{Proposition} 
\newtheorem{lemma}[thm]{Lemma} 
\newtheorem{cor}[thm]{Corollary} 
\newtheorem{defn}[thm]{Definition} 
\newtheorem{notation}[thm]{Notation} 
\newtheorem{qtn}[thm]{Question} 
\newtheorem{example}[thm]{Example} 
\newtheorem{remark}[thm]{Remark} 
\newtheorem{conj}[thm]{Conjecture} 
\newtheorem{prob}[thm]{Problem} 
\newtheorem{rem}[thm]{Remark} 

\newcommand{\bdd}{\partial}
\newcommand{\Int}{{\rm Int}}
\renewcommand{\d}{\delta}
\newcommand{\e}{\epsilon}
\newcommand{\br}{{\Bbb R}}
\renewcommand{\r}{r}
\newcommand{\ga}{\Gamma_a}
\newcommand{\gb}{\Gamma_b}
\newcommand{\rg}{\hat \Gamma}
\newcommand{\gap}{\Gamma_a^+}
\newcommand{\gbp}{\Gamma_b^+}
\newcommand{\rga}{\hat \Gamma_a}
\newcommand{\rgb}{\hat \Gamma_b}
\newcommand{\rgap}{\hat \Gamma_a^+}
\newcommand{\rgbp}{\hat \Gamma_b^+}
\newcommand{\rgaa}{\hat \Gamma_1}
\newcommand{\rgbb}{\hat \Gamma_2}
\newcommand{\rgaap}{\hat \Gamma_1^+}
\newcommand{\rgbbp}{\hat \Gamma_2^+}
\newcommand{\he}{\hat e}
\newcommand{\na}{n_{a}}
\newcommand{\nb}{n_{b}}
\newcommand{\sign}{\text{\rm sign}}
\newcommand{\BH}{\Bbb H^3}

\input epsf.tex

\section{Introduction}

Let $M$ be a {\it hyperbolic\/} 3-manifold, by which we shall mean a
compact, connected, orientable 3-manifold such that $M$ with its
boundary tori removed admits a complete hyperbolic structure with
totally geodesic boundary, and suppose that $M$ has a torus boundary
component $T_0$.  If $r$ is a slope on $T_0$ then $M(r)$ will denote
the 3-manifold obtained by $r$-Dehn filling on $M$, i.e.\ attaching a
solid torus $V_r$ to $M$ along $T_0$ in such a way that $r$ bounds a
disk in $V_r$. The Dehn filling $M(r)$ and the slope $r$ are said to
be {\it exceptional\/} if $M(r)$ is either reducible,
$\bdd$-reducible, annular, toroidal, or a small Seifert fiber space.
Modulo the Geometrization Conjecture, the manifold $M(r)$ is
hyperbolic if and only if $M(r)$ is not exceptional.

Thurston's Hyperbolic Dehn Surgery Theorem asserts that there are only
finitely many exceptional Dehn fillings on each torus boundary
component of $M$. It is known that if $r,s$ are both exceptional then
the geometric intersection number $\Delta = \Delta(r,s)$, also known
as the {\it distance\/} between $r$ and $s$, is small. In fact, the
least upper bounds for $\Delta$ have been determined for all cases
where neither $M(r)$ nor $M(s)$ is a small Seifert fiber space, by the
work of many people. See [GW2] and the references therein.

For toroidal fillings, it was shown by Gordon [Go] that if $r,s$ are
toroidal slopes then $\Delta(r,s) \leq 8$, and moreover there are
exactly two manifolds $M$ with $\Delta=8$, one with $\Delta=7$, and
one with $\Delta=6$. In this paper we classify all the hyperbolic
3-manifolds which admit two toroidal Dehn fillings with $\Delta=4$ or
$5$.

Already when $\Delta=5$ there are infinitely many such manifolds. To
see this, let $M$ be the exterior of the Whitehead sister link, also
known as the $(-2,3,8)$-pretzel link. The boundary of $M$ consists of
two tori $T_0$ and $T_1$, and there are slopes $r,s$ on $T_0$ with
$\Delta(r,s) = 5$ such that the Dehn filled manifolds $M(r) = M(r,*)$,
$M(s) = M(s,*)$ are toroidal; see for example [GW3]. Now for
infinitely many slopes $t$ on $T_1$, $M_t =M(*,t)$ will be hyperbolic
and $M_t(r) = M(r)(t)$, $M_t(s) = M(s)(t)$ will be toroidal. In this
way we get infinitely many hyperbolic 3-manifolds with boundary a
single torus having two toroidal fillings at distance $5$. We shall
show that, modulo this phenomenon, there are only finitely many $M$
with two toroidal fillings at distance $4$ or $5$, and explicitly
identify them.  Define two triples $(N_1, r_1, s_1)$ and $(N_2, r_2,
s_2)$ to be {\it equivalent}, denoted by $(N_1, r_1, s_1) \cong (N_2,
r_2, s_2)$, if there is a homeomorphism from $N_1$ to $N_2$ which
sends the boundary slopes $(r_1, s_1)$ to $(r_2, s_2)$ or $(s_2,
r_2)$.

\begin{thm} There exist $14$ $3$-manifolds $M_i$, $1\leq i\leq 14$,
such that

(1) $M_i$ is hyperbolic, $1\leq i\leq 14$;

(2) $\bdd M_i$ consists of two tori $T_0, T_1$ if $i \in
\{1,2,3,14\}$, and a single torus $T_0$ otherwise;

(3) there are slopes $r_i, s_i$ on the boundary component $T_0$ of
$M_i$ such that $M(r_i)$ and $M(s_i)$ are toroidal, where
$\Delta(r_i,s_i) = 4$ if $i \in \{1,2,4,6,9,13,14\}$, and
$\Delta(r_i,s_i) = 5$ if $i \in \{3,5,7,8,10,11,12\}$;

(4) if $M$ is a hyperbolic 3-manifold with toroidal Dehn fillings
$M(r), M(s)$ where $\Delta(r,s) = 4$ or $5$, then $(M,r,s)$ is
equivalent either to $(M_i,r_i,s_i)$ for some $1\leq i \leq 14$, or to
$(M_i(t),r_i,s_i)$ where $i \in \{1,2,3,14\}$ and $t$ is a slope on
the boundary component $T_1$ of $M_i$.
\end{thm}

\proof The manifolds $M_i$ are defined in Definition 21.3.  (1) is
Theorem 23.14. (2) follows from the definition.  (3) and (4) follow
from Theorem 21.4.  \qed

\begin{remark}
Part (4) in Theorem 1.1 is still true if
the hyperbolicity is replaced by the assumption that $M$ is compact,
connected, orientable, irreducible, atoroidal, and non Seifert
fibered, in other words, $M$ may be annular or $\bdd$-reducible but
not Seifert fibered.  
\end{remark}

\proof First assume $M$ is $\bdd$-irreducible.  Then any essential
annulus must have at least one boundary component on a non-toroidal
boundary component of $M$ as otherwise $M$ would be either toroidal or
Seifert fibered.  Attaching a hyperbolic manifold $X$ to each
non-toroidal boundary component of $M$ will produce a hyperbolic
manifold $M'$, and we may choose $X$ so that $M'$ has more than three
boundary components.  One can show that the $M'(r_i)$ are still
toroidal for $i=1,2$, which is a contradiction to Theorem 21.4.  Now
assume $M$ is $\bdd$-reducible.  Then $M$ is obtained by attaching
1-handles to the boundary of a manifold $M''$.  If a 1-handle is
attached to a toroidal component $T$ of $\bdd M''$ then either $M'' =
T \times I$, which is impossible because $M(r_i)$ would be a
handlebody and hence atoroidal, or $T$ would be an essential torus in
$M$, contradicting the assumption.  It follows that $M''$ has a higher
genus boundary component.  One can check that the $M''(r_i)$ are still
toroidal, which leads to a contradiction to Theorem 21.4 as above.
\qed

The manifolds $M_1, M_2$ and $M_3$ were discussed in [GW1]; $M_i$,
$i=1,2,3$ is the exterior of a link $L_i$ in $S^3$, where $L_1$ is the
Whitehead link, $L_2$ is the $2$-bridge link associated to the
rational number $3/10$, and $L_3$ is the Whitehead sister link. See
Figure 24.1.  The other $M_i$ can be built using intersection graphs
on tori, see Definition 21.3 for more details.  For $i \neq 4,5$, each
$M_i$ can also be described as a double branched cover of a tangle
$Q_i = (W_i, K_i)$, where $W_i$ is a 3-ball for $i=6,...,13$, and a
once punctured 3-ball for $i=1,2,3,14$.  This is done in [GW1] for
$i=1,2,3$, and in Section 22 for the other cases.  See Lemma 22.2.

Some results on the case $\Delta=5$ have been independently obtained
by Teragaito [T2]. He obtains a finite set of pairs of intersection
graphs of punctured tori at distance $5$ which must contain all the
pairs of graphs that arise from two toroidal fillings on a hyperbolic
3-manifold at distance 5.  One of his pairs produces a non-hyperbolic
manifold while the others correspond to the manifolds in our list for
$\Delta=5$.  

We remark that the related problem of determining all hyperbolic
$3$-manifolds with two Dehn fillings at distance at least $4$ that
yield manifolds containing Klein bottles has been solved by Lee [L1,
L2] (see also [MaS]).  

Since $M_i$ has more than one boundary component only when $i \in
\{1,2,3,14\}$, we have the following corollary to Theorem 1.1
(together with [Go]), which in the case $\Delta = 5$ is due to Lee
[L1].  Note that all boundary components of the manifolds are tori.

\begin{cor} Let $M$ be a hyperbolic 3-manifold with more than one
  boundary component, having toroidal Dehn fillings $M(r), M(s)$ with
  $\Delta = \Delta(r,s) \geq 4$. Then each boundary component of $M$
  is a torus, and either

  (1) $\Delta = 4$ and $(M,r,s) \cong (M_i,r_i,s_i)$ for $i \in
  \{1,2,14\}$, or

  (2) $\Delta = 5$ and $(M,r,s) \cong (M_3,r_3,s_3)$.
\end{cor}

In [GW1] and [GW3] it is shown that if $M$ is a hyperbolic 3-manifold
with fillings $M(r)$ and $M(s)$, one of which is annular and the other
either toroidal or annular, then either $(M,r,s) \cong (M_i,r_i,s_i)$ for $i
\in \{1,2,3\}$, or $\Delta(r,s) \leq 3$. It is also known that if $M(r)$
contains an essential sphere or disk, and $M(s)$ contains an essential
sphere, disk, annulus or torus, then $\Delta(r,s) \leq 3$; see [GW2] and
the references listed there. Corollary 1.3 then gives

\begin{cor} Let $M$ be a hyperbolic 3-manifold with a torus boundary
  component $T_0$ and at least one other boundary component. Let $r,s$
  be exceptional slopes on $T_0$. Then either $(M,r,s) \cong
  (M_i,r_i,s_i)$ for $i \in \{1,2,3,14\}$, or $\Delta(r,s) \leq 3$.
\end{cor}

A pair $(M, T_0)$ is called a {\it large manifold\/} if $T_0$ is a
torus on the boundary of the 3-manifold $M$ and $H_2(M, \bdd M - T_0)
\neq 0$ (see [Wu3]).  Teragaito [T2] proved that there is no large
hyperbolic manifold $M$ admitting two toroidal fillings of distance at
least 5.  The following corollary clarifies the case of distance $4$.

\bigskip
\noindent {\bf Theorem 22.3} \; {\em 
  Suppose $(M,T_0)$ is a large manifold and $M$ is hyperbolic and
  contains two toroidal slopes $r_1, r_2$ on $T_0$ with $\Delta(r_1,
  r_2) \geq 4$.  Then $M$ is the Whitehead link exterior, and
  $\Delta(r_1, r_2)=4$.
}  \bigskip

Theorem 1.1 gives information about toroidal Dehn surgeries on
hyperbolic knots in $S^3$. It follows from [Go] that the only such
knot with two toroidal surgeries at distance $> 5$ is the figure eight
knot, for which the $4$ and $-4$ surgeries are toroidal. Teragaito has
shown [T1] that the only hyperbolic knots with two toroidal surgeries
at distance $5$ are the Eudave-Mu\~noz knots $k(2,-1,n,0)$, $n \neq
1$.  We can now determine the knots with toroidal surgeries at
distance $4$.  Denote by $L_i = K'_i \cup K''_i$ the link in Figure
24.1(i), where $K'_i$ is the component on the left.  Denote by
$L_i(n)$ the knot obtained from $K''_i$ by $1/n$ surgery on $K'_i$.
One can check that $L_3(n)$ is the same as the Eudave-Mu\~noz knot
$k(3,1,-n,0)$ in [Eu, Figure 25], which is the mirror-image of
$k(2,-1,1+n,0)$ [Eu, Proposition 1.4].

\bigskip
\noindent
{\bf Theorem 24.4} \;
{\em Suppose $K$ is a hyperbolic knot in $S^3$ admitting two
toroidal surgeries $K(r_1), K(r_2)$ with $\Delta(r_1, r_2) \geq 4$.
Then $(K, r_1, r_2)$ is equivalent to one of the following, where $n$
is an integer.

(1) $K = L_1(n)$, $r_1 = 0$, $r_2 = 4$.

(2) $K = L_2(n)$, $r_1 = 2-9n$, $r_2 = -2-9n$.

(3) $K = L_3(n)$, $r_1 = -9 - 25n$, $r_2 = -(13/2) - 25n$.

(4) $K$ is the Figure 8 knot, $r_1 = 4$, $r_2 = -4$.
}
\bigskip

The only hyperbolic knots known to have more than two toroidal
surgeries are the figure eight knot and the $(-2,3,7)$-pretzel knot,
with toroidal slopes $\{-4,0,4\}$ and $\{16,37/2,20\}$ respectively.
This led Eudave-Mu\~noz [Eu] to conjecture that a hyperbolic knot in
$S^3$ has at most three toroidal surgeries.  Teragaito [T1] showed
that there can be at most five toroidal surgeries.  Theorem 1.1 and
[T1, Corollary 1.2] lead to the following improvement.

\bigskip
\noindent
{\bf Corollary 24.5} \;
{\em  A hyperbolic knot in $S^3$ has at most four toroidal surgeries. If
  there are four, then they are consecutive integers.}
\bigskip

Here is a sketch of the proof of Theorem 1.1.  A toroidal Dehn filling
$M(r)$ on a hyperbolic $3$-manifold $M$ gives rise to an essential
punctured torus $F$ in $M$ whose boundary consists of $n>0$ circles of
slope $r$ on $T_0$, where the capped-off surface $\hat F$ is an
essential torus in $M(r)$. Hence, in the usual way (see Section 2),
two toroidal fillings $M(r_1), M(r_2)$ give rise to a pair of
intersection graphs $\Gamma_1, \Gamma_2$ on the tori $\hat F_1, \hat
F_2$, with $n_1, n_2$ vertices respectively. The proof consists of a
detailed analysis of the possible pairs of intersection graphs with
$\Delta(r_1,r_2)=4$ or $5$, using Scharlemann cycles and other tools
developed in earlier works in this area. This enables us to eliminate
all but $17$ pairs of graphs.  As is usual in this kind of setting,
the permissible graphs all have small numbers of vertices. Eleven of
the pairs correspond to the manifolds $M_i$, $4\leq i\leq 14$.  We
show that any of the remaining pairs must correspond to a pair of
fillings on $M_i$ or $M_i(t)$ for $i \in \{1, 2, 3\}$.

Here is a more detailed summary of the organization of the paper.
Section 2 contains the basic definitions and some preliminary lemmas.
In Sections 3-5 we deal with the generic case $n_1, n_2 > 4$,
ultimately showing (Proposition 5.11) that this case cannot occur.
More specifically, Section 3 shows that the reduced positive graph
$\rgap$ of $\ga$ (see Section 2 for definitions) has no interior
vertices, and this is strengthened in Section 4 to showing that each
component of $\rgap$ must be one of the 11 graphs in Figure 4.2. These
are ruled out one by one in Section 5. In Sections 6-11 we consider
the case where some $n_a=4$. Section 6 discusses the situation where
the graph $\ga$ is {\it kleinian}; this arises when the torus $\hat
F_a$ is the boundary of a regular neighborhood of a Klein bottle in
$M(r_a)$. (The results here are also used in the discussion of the
case $n_1, n_2 \leq 2$.)  Sections 7, 8 and 9 show that if $n_a = 4$
and $\gb$ is non-positive then $n_b\leq 4$. Section 10 shows that if
$\Gamma_1$ and $\Gamma_2$ are both non-positive then $n_1=n_2=4$ is
impossible.  Section 11 shows (Proposition 11.9) that if $\gb$ is
positive then there are exactly two pairs of graphs, one with $n_b=2$,
the other with $n_b=1$. These give the manifolds $M_4$ and $M_5$
respectively.  If we suppose $n_a \leq n_b$, it now easily follows
(Proposition 11.10) that $n_a\leq 2$.

In Sections 12-16 we deal with the case $n_a\leq 2$, $n_b\geq 3$. The
conclusion (Proposition 16.8) is that here there are exactly six pairs
of graphs.  Two of these are the ones described in Section 11, and the
four new pairs give the manifolds $M_6, M_7, M_8$ and $M_9$. More
precisely, in Section 12 we rule out the case where $\gb$ is positive,
and in Sections 13 and 14 we consider the case where $n_b>4$ and both
graphs $\Gamma_1$ and $\Gamma_2$ are non-positive. It turns out that
here there is exactly one pair of graphs (Proposition 14.7),
corresponding to the manifold $M_6$. We may now assume that $n_b = 3$
or $4$. Section 15 establishes some notation and elementary properties
for graphs with $n_a\leq 2$. In Section 16 we show that if $\Gamma_1$
and $\Gamma_2$ are non-positive then $n_b=3$ is impossible and if
$n_b=4$ then there are exactly three examples, $M_7$, $M_8$ and $M_9$.

Sections 17-20 deal with the remaining cases where both $n_1$ and
$n_2$ are $\leq 2$. In Section 17 we introduce an equivalence
relation, {\it equidistance}, on the set of edges of a graph $\ga$,
and show that, under the natural bijection between the edges of
$\Gamma_1$ and $\Gamma_2$, the two graphs induce the same equivalence
relation. This gives a convenient way of ruling out certain pairs of
graphs. Section 18 considers the case $n_a=2$ and $n_b=1$, and shows
that here there are exactly three examples. Section 19 considers the
case $n_1=n_2=2$, $\Gamma_b$ positive, showing that there are two
examples. Finally, in Section 20 we consider the case $n_1=n_2=2$,
$\Gamma_1$ and $\Gamma_2$ both non-positive, and show that there are
exactly six pairs of graphs in this case. The final list of all 11
possible pairs of graphs with $n_1, n_2 \leq 2$ is given in
Proposition 20.4. Five of these correspond to the manifolds $M_{10},
M_{11}, M_{12}, M_{13}$ and $M_{14}$.

The remaining six pairs of graphs in Proposition 20.4 have the
property that one of the graphs has a non-disk face. In Section 21 we
show (Lemma 21.2), using the classification of toroidal/annular and
annular/annular fillings at distance $\geq 4$ given in [GW1] and
[GW3], that in this case the manifold $M$ is either $M_1, M_2$ or
$M_3$, or is obtained from one of those by Dehn filling along one of
the boundary components.

In Section 22 we show how the manifolds $M_i$, $6\leq i\leq 14$ may be
realized as double branched covers. Using this, in Section 23 we show
that the manifolds $M_i$ are hyperbolic. Finally, in Section 24 we
give the applications to toroidal surgeries on knots in $S^3$.

\section{Preliminary Lemmas}

Throughout this paper, we will fix a hyperbolic 3-manifold $M$, with a
torus $T_0$ as a boundary component.  A compact surface properly
embedded in $M$ is {\it essential\/} if it is $\pi_1$-injective, and
is not boundary parallel.  We use $a, b$ to denote the numbers $1$ or
$2$, with the convention that if they both appear in a statement then
$\{a, b\} = \{1,2\}$.

A slope on $T_0$ is a {\it toroidal slope\/} if $M(r_{a})$ is
toroidal.  Let $r_{a}$ be a toroidal slope on $T_0$.  Denote by
$\Delta = \Delta(r_1, r_2)$ the minimal geometric intersection number
between $r_1$ and $r_2$.  When $\Delta > 5$ the manifolds $M$ have been
determined in [Go].  We will always assume that $\Delta = 4$ or $5$.
Let $\hat F_a$ be an essential torus in $M(r_a)$, and let $F_a = \hat
F_a \cap M$.  If $M(r_a)$ is reducible then by [Wu1] and [Oh] we would
have $\Delta \leq 3$, which is a contradiction.  Therefore both
$M(r_a)$ are irreducible.

Let $n_{a}$ be the number of boundary components of $F_{a}$ on $T_0$.
Choose $\hat F_{a}$ in $M(r_{a})$ so that $n_{a}$ is minimal among
all essential tori in $M(r_{a})$.  Minimizing the number of
components of $F_1 \cap F_2$ by an isotopy, we may assume that
$F_1\cap F_2$ consists of arcs and circles which are essential on both
$F_{a}$.  Denote by $J_{a}$ the attached solid torus in $M(r_{a})$,
and by $u_i$ ($i=1,...,n_{a}$) the components of $\hat F_{a}\cap
J_{a}$, which are all disks, labeled successively when traveling
along $J_{a}$.  Similarly let $v_j$ be the disk components of $\hat F_b\cap
J_b$.  Let $\Gamma_{a}$ be the graph on $\hat F_{a}$ with the
$u_i$'s as (fat) vertices, and the arc components of $F_1\cap F_2$ as
edges.  Similarly for $\gb$.  The minimality of the number of
components in $F_1\cap F_2$ and the minimality of $n_{a}$ imply that
$\Gamma_{a}$ has no trivial loops, and that each disk face of
$\Gamma_{a}$ in $\hat F_{a}$ has interior disjoint from $F_{b}$.

If $e$ is an edge of $\Gamma_{a}$ with an endpoint $x$ on a fat
vertex $u_i$, then $x$ is labeled $j$ if $x$ is in $u_i\cap v_j$.  In
this case $e$ is called a {\it $j$-edge} in $\ga$, and an $i$-edge in
$\gb$.  Labels in $\Gamma_a$ are considered as mod $n_b$ integers; in
particular, $n_b +1 = 1$.  When going around $\bdd u_i$, the labels of
the endpoints of edges appear as $1, 2, \ldots, n_b$ repeated $\Delta$
times.  Label the endpoints of edges in $\Gamma_b$ similarly.

Each vertex of $\Gamma_{a}$ is given a sign according to whether
$J_{a}$ passes $\hat F_{a}$ from the positive side or negative side
at this vertex.  Two vertices of $\Gamma_{a}$ are {\it parallel\/} if
they have the same sign, otherwise they are {\it antiparallel.}  Note
that if $\hat F_{a}$ is a separating surface, then $n_{a}$ is even,
and $v_i, v_j$ are parallel if and only if $i, j$ have the same
parity.  We use $val(v, G)$ to denote the valence of a vertex $v$ in a
graph $G$.  If $G$ is clear from the context, we simply denote it by
$val(v)$.

When considering each family of parallel edges of $\Gamma_{a}$ as a
single edge $\hat e$, we get the {\it reduced graph\/}
$\hat{\Gamma}_{a}$ on $\hat F_{a}$.  It has the same vertices as
$\Gamma_{a}$.  Each edge of $\rga$ represents a family of parallel
edges in $\ga$.  We shall often refer to a family of parallel edges as
simply a {\it family}.

\begin{defn}
(1) An edge of $\Gamma_{a}$ is a {\it positive edge\/} if it
connects parallel vertices.  Otherwise it is a {\it negative edge}.

(2) The graph $\ga$ is {\it positive\/} if all its vertices are
parallel, otherwise it is {\it non-positive}.  
p\end{defn}

We use $\Gamma_{a}^+$ (resp.\ $\Gamma_{a}^-$) to denote the
subgraph of $\Gamma_{a}$ whose vertices are the vertices of
$\Gamma_{a}$ and whose edges are the positive (resp.\ negative)
edges of $\Gamma_{a}$.  Similarly for $\hat{\Gamma}_{a}^+$
and $\hat{\Gamma}_{a}^-$.

A cycle in $\Gamma_{a}$ consisting of positive edges is a {\it
  Scharlemann cycle\/} if it bounds a disk with interior disjoint from
the graph, and all the edges in the cycle have the same pair of labels
$\{i, i+1\}$ at their two endpoints, called the {\it label pair\/} of
the Scharlemann cycle.  A Scharlemann cycle containing only two edges
is called a {\it Scharlemann bigon.}  A Scharlemann cycle with label
pair, say, $\{1,2\}$ will also be called a $(12)$-Scharlemann cycle.
If $\gb$ contains a Scharlemann cycle with label pair $\{i,i\pm 1\}$,
we shall sometimes abuse terminology and say that the vertex $u_i$ of
$\ga$ is a {\it label of a Scharlemann cycle}.  An {\it extended
  Scharlemann cycle\/} is a cycle of edges $\{e_1, ..., e_k\}$ such
that there is a Scharlemann cycle $\{e'_1, ..., e'_k\}$ with $e_i$
parallel and adjacent to $e'_i$ and $e_i \neq e'_j$, $1\leq i,j \leq
k$.  If $\{e_1, ..., e_k\}$ is a Scharlemann cycle in $\ga$ then the
subgraph of $\gb$ consisting of these edges and their vertices is
called a {\it Scharlemann cocycle}.

A subgraph $G$ of a graph $\Gamma$ on a surface $F$ is {\it
essential\/} if it is not contained in a disk in $F$.  The following
lemma contains some common properties of the graphs $\Gamma_{a}$.
It can be found in [GW1, Lemma 2.2].

\begin{lemma}
(1) {\rm (The Parity Rule)} An edge $e$ is a positive edge in
$\Gamma_1$ if and only if it is a negative edge in $\Gamma_2$.

(2) A pair of edges cannot be parallel on both $\Gamma_1$ and
$\Gamma_2$. 

(3) If $\Gamma_{a}$ has a set of $n_{b}$ parallel negative
edges, then on $\Gamma_{b}$ they form mutually disjoint essential
cycles of equal length.

(4) If $\Gamma_{a}$ has a Scharlemann cycle, then $\hat F_{b}$ is
separating.  In particular, $\gb$ has the same number of positive and
negative vertices, so $n_{b}$ is even, and two vertices $v_i, v_j$
of $\gb$ are parallel if and only if $i, j$ have the same parity.

(5) If $\Gamma_{a}$ has a Scharlemann cycle $\{e_1,\ldots, e_k\}$,
then the corresponding Scharlemann cocycle on $\Gamma_{b}$ is
essential.

(6)  If $n_{b} > 2$, then $\Gamma_{a}$ contains no extended
Scharlemann cycle. 
\end{lemma}

Let $\hat e$ be a collection of parallel negative edges on
$\Gamma_{b}$, oriented from $v_1$ to $v_2$.  Then $\hat e$ defines a
permutation $\varphi: \{1, \ldots, n_{a}\} \to \{1, \ldots,
n_{a}\}$, such that an edge $e$ in $\hat e$ has label $k$ at $v_1$ if
and only if it has label $\varphi(k)$ at $v_2$.  Call $\varphi$ the
{\it transition function associated to $\hat e$}.  Define the {\it
transition number\/} to be the mod $n_a$ integer $s = s(\hat e)$
such that $\varphi(k) = k + s$.  If we reverse the orientation of
$\hat e$ then the transition function is $\varphi^{-1}$, and the
transition number is $-s$; hence if $\hat e$ is unoriented then
$\varphi$ is well defined up to inversion, and $s(\hat e)$ is well
defined up to sign.

\begin{lemma}
(1) If a family of parallel negative edges in
$\Gamma_a$ contains more than $n_b$ edges (in particular, if the
family contains 3 edge endpoints with the same label), then $\Gamma_b$
is positive, and the transition function associated to this family is
transitive.

(2) If $\Gamma_a$ contains two Scharlemann cycles with disjoint label
pairs $\{i, i+1\}$ and $\{j, j+1\}$, then $i\equiv j$ mod 2.

(3) If $n_b > 2$ then a family of parallel positive edges in
$\Gamma_a$ contains at most $n_b/2 + 2$ edges, and if it does contain
$n_b/2 +2$ edges, then $n_b \equiv 0$ mod 4.

(4) $\ga$ has at most four labels of Scharlemann cycles, at most two
for each sign.

(5) A loop edge $e$ and a non-loop edge $e'$ on $\ga$ cannot be
parallel on $\gb$.

(6) If $n_b \geq 4$ then $\ga$ contains at most $2n_b$ parallel
negative edges.
\end{lemma}

\proof 
(1)  This is obvious if $n_b \leq 2$, and it can be found in [GW1,
Lemma 2.3] if $n_b > 2$.

(2) and (3) are basically Lemmas 1.7 and 1.4 of [Wu1].  If $\ga$ has
$n_b/2 + 2$ parallel positive edges then the two outermost pairs form
two Scharlemann bigons.  One can then check the labels of these
Scharlemann bigons and use (2) to show that $n_b \equiv 0$ mod 4.

(4) If $\ga$ has more than four labels of Scharlemann cycles, then
either one can find two Scharlemann cycles with disjoint label pairs
$\{i, i+1\}$ and $\{j, j+1\}$ such that $i - j \equiv 1$ mod 2, which
is a contradiction to (2), or one can find three Scharlemann cycles
with mutually disjoint label pairs, in which case one can replace
$\hat F_a$ by another essential torus to reduce $n_a$ and get a
contradiction.  See [Wu1, Lemma 1.10]. 

If $\ga$ has three positive labels of Scharlemann cycles $u_{i_j}$
then it has negative labels of Scharlemann cycles $u_{i_j +
\epsilon_j}$ for some $\epsilon_j = \pm 1$, which cannot all be the
same, hence $\ga$ has at least 5 labels of Scharlemann cycles,
contradicting the above.  

(5) Since $e$ is positive on $\ga$, it is negative in $\gb$.  If $e$
has endpoints on $u_i$ in $\ga$ then on $\gb$ its two endpoints are
both labeled $i$, hence the corresponding transition number is $0$, so
any edge $e'$ parallel to $e$ on $\gb$ must also have the same label
at its two endpoints, which implies that $e'$ is a loop on $\ga$.

(6) This is [Go, Corollary 5.5].  
\qed

\begin{lemma}
If a label $i$ appears twice among the endpoints
of a family $\hat e$ of parallel positive edges in $\Gamma_{a}$, then
$i$ is a label of a Scharlemann bigon in $\hat e$.  In particular, if
$\hat e$ has more than $n_{b}/2$ edges, then it contains a
Scharlemann bigon.  
\end{lemma}

\proof Since the edges are positive, by the parity rule $i$ cannot
appear at both endpoints of a single edge in this family.  Let $e_1,
e_2, \ldots, e_k$ be consecutive edges of $\hat e$ such that $e_1$ and
$e_k$ have $i$ as a label.  Now $k$ must be even, otherwise the edge
$e_{(k+1)/2}$ would have the same label at its two endpoints.  If
$k\geq 4$ and $n_b>2$ one can see that these edges contain an extended
Scharlemann cycle, which contradicts Lemma 2.2(6).  Therefore $k=2$ or
$n_b=2$, in which case $e_1, e_2$ form a Scharlemann bigon with $i$ as
a label.

If $\hat e$ has more than $n_{b}/2$ edges, then it has more than
$n_{b}$ endpoints, so some label must appear twice.
\qed

\begin{lemma}
$\rga$ contains at most $3n_a$ edges.
\end{lemma}

\proof Let $V,E,F$ be the number of vertices, edges and disk faces of
$\rga$.  Then $V - E + F \geq 0$ (the inequality may be strict if
there are some non-disk faces.)  Each face of $\rga$ has at least
three edges, hence we have $ 3F \leq 2E$.  Solving those two
inequalities gives $E \leq 3V$.  \qed

\begin{lemma} If $n_b> 4$ then the vertices of $\ga$ cannot all
be parallel.  \end{lemma}

\proof By Lemma 2.5 the reduced graph $\rga$ has at most $3n_{a}$
edges.  For any $i$, since $v_i$ on $\gb$ has valence at least
$4n_{a}$, there are $4n_{a}$ $i$-edges on $\ga$, hence two of them
must be parallel, so $i$ is a label of a Scharlemann cycle.  Since
there are at most $4$ such labels (Lemma 2.3(4)), we would have
$n_{b} \leq 4$, contradicting the assumption.  \qed

A vertex $v$ of a graph is a {\it full vertex\/} if all edges incident
to it are positive.

\begin{lemma} Suppose $n_b > 4$.  Then

(1) a family of parallel negative edges in $\gb$ contains at most
$n_{a}$ edges, hence any label $i$ appears at most twice among the
endpoints of such a family;

(2) two families of positive edges in $\ga$ adjacent at a vertex
contain at most $n_b + 2$ edges; and

(3) three families of positive edges in $\ga$ adjacent at a vertex
contain at most $2n_b$ edges, and if there are $2n_b$ then $n_b = 6$.
\end{lemma}

\proof (1) If a family of parallel negative edges on $\gb$ contains
more than $n_a$ edges then by Lemma 2.3(1) all vertices of $\ga$ are
parallel, which contradicts Lemma 2.6.

If $i$ appears three times among the endpoints of a family of parallel
negative edges in $\gb$ then this family would contain more than
$n_{a}$ edges, which is a contradiction.

(2) By Lemma 2.3(3) a family of parallel positive edges contains $r
\leq n_b/2 + 2$ edges.  If two adjacent families $\hat e_1, \hat e_2$
contain more than $n_b + 2$ edges, then one of them, say $\hat e_1$,
has $n_b/2+2$ edges while the other one has either $n_b/2+1$ or
$n_b/2+2$ edges.  Now $\hat e_1$ contains two Scharlemann bigons,
which must appear on the two sides of the family because there is no
extended Scharlemann cycle.  There is also at least one Scharlemann
bigon in $\hat e_2$.  Examining the labels of these Scharlemann bigons
we can see that they contain at least $5$ labels, which contradicts
Lemma 2.3(4).

(3) Assume the three families contain $r \geq 2n_b$ edges.  Then one of
the families contains more than $n_b/2$ edges, so by Lemma 2.2(4) $n_b$
is even.  By (2) two adjacent families of parallel edges contain at
most $n_b + 2$ edges, while by Lemma 2.3(3) the other family has at
most $n_b/2+2$ edges, so we have $2n_b \leq r \leq (n_b+2) + (n_b/2+2)$,
which gives $n_b \leq 8$.

If $n_b = 8$ then the above inequalities force the three families to
have $6,4,6$ edges, and we see that all 8 labels appear as labels of
Scharlemann bigons, which contradicts Lemma 2.3(4).  So we must have
$n_b = 6$.  By Lemma 2.3(3) we have $2n_b \leq r \leq 3(n_b/2 + 1) =
12 = 2n_b$.  Hence $r = 2n_b$.  \qed

\begin{lemma} If a vertex $u_i$ of $\ga$ is incident to more
than $n_b$ negative edges, then $\gb$ has a Scharlemann cycle.
\end{lemma}

\proof In this case there are $n_b + 1$ positive $i$-edges in $\gb$,
which cut the surface $F_b$ into faces, at least one of which is a disk
face in the sense that it is a topological disk whose interior
contains no vertices of $\gb$.  Hence the subgraph of $\gb$ consisting
of these edges is a $x$-edge cycle in the sense of Hayashi-Motegi [HM,
Page 4468].  By [HM, Proposition 5.1] a disk face of this $x$-edge
cycle contains a disk face of a Scharlemann cycle.  \qed

Consider a graph $G$ on a closed surface $F$, and assume that $G$ has
no isolated vertex.  If the vertices of $G$ have been assigned $\pm$
signs (for example $\rgap$), let $X$ be the union of $G$ and all its
faces $\sigma$ such that all vertices on $\bdd \sigma$ have the same
sign, otherwise let $X$ be the union of $G$ and all its disk faces.  A
vertex $v$ of $G$ is an {\it interior vertex\/} if it lies in the
interior of $X$.  A vertex $v$ of $G$ is a {\it cut vertex\/} if a
regular neighborhood of $v$ in $X$ with $v$ removed is not connected.
A vertex $v$ of $G$ is a {\it boundary vertex\/} if it is not an
interior or cut vertex.  Note that if $G = \rgap$ then an interior
vertex is a full vertex.  Alternatively, let $\delta(v)$ be the number
of corners around $v$ which lie in $X$.  Then $v$ is an interior
vertex if $\delta(v) = val(v, G)$, a boundary vertex if $\delta(v) =
val(v, G) - 1$, and a cut vertex if $\delta(v) \leq val(v,G) - 2$.

Given a graph $G$ on a surface $D$, let $c_i(G)$ be the number of
boundary vertices of $G$ with valence $i$.  Define 
\begin{eqnarray*}
\varphi(G) & = & 6c_0(G) + 3 c_1(G) + 2 c_2(G) + c_3(G) \\
\psi(G) & = & c_0(G) + c_1(G) + c_2(G) + c_3(G).
\end{eqnarray*}
Note that $\psi(G)$ is the number of boundary vertices of $G$ with
valence at most 3.

\begin{lemma} Let $G$ be a connected reduced graph in a disk
$D$ such that any interior vertex of $G$ has valence at least $6$.
Then $\varphi(G) \geq 6$.  Moreover, if $G$ is not homeomorphic to an
arc or a single point then $\psi(G) \geq 3$.  \end{lemma}

\proof Let $X$ be the union of $G$ and all its disk faces.  The result
is obviously true if $G$ is a tree.  So we assume that $G$ has some
disk faces.

First assume that $X$ has no cut vertex, so it is a disk, and $c_0(G)
= c_1(G) = 0$.  The double of $G$ along $\bdd X$ is then a graph
$\tilde G$ on the double of $X$, which is a sphere.  Note that the
valence of a vertex $v$ of $\tilde G$ is either at least $6$, or it is
$2$ or $4$ when $v$ is a boundary vertex of $G$ with valence $2$ or
$3$, respectively.  Since each face has at least three edges, an
Euler characteristic argument gives $$
2 = V - E + F \leq V - \frac 13
E = \sum_i (1 - \frac 16 val(v_i, \tilde G)) \leq \frac 23 c_2(G) +
\frac 13 c_3(G)$$
Therefore $\varphi(G) = 2c_2(G) + c_3(G) \geq 6$.
Since $c_0(G) = c_1(G) = 0$, we also have $\psi(G) = c_2(G) + c_3(G)
\geq \frac 12 \varphi(G) \geq 3$.

Now assume that $X$ has a cut vertex $v$.  Since $G$ is connected and
contained in a disk, $X$ is simply connected, so we can write $X = X_1
\cup X_2$, where $X_i$ are subcomplexes of $X$ such that $X_1 \cap X_2
= v$, and $G_i = G \cap X_i$ are nontrivial connected subgraphs of
$G$.  The valence of $v$ in $G_i$ is at least $1$, so its
contribution to $\varphi(G_i)$ is at most 3.  Hence by induction we
have
$$\varphi(G) \geq (\varphi(G_1)-3) + (\varphi(G_2)-3) \geq 6.$$
By assumption $X$ is not homeomorphic to an arc, so at least one of
the $X_i$, say $X_1$, is not homeomorphic to an arc, and the other one
has at least 2 boundary vertices of valence at most 3, whether it is
homeomorphic to an arc or not.  Hence 
$$
\psi(G) \geq \psi(G_1) + \psi(G_2) - 2 \geq 3 + 2 - 2 = 3.
$$
\qed

\begin{lemma} Let $G$ be a reduced graph on a torus $T$ with
no interior or isolated vertex.  Let $V$ and $E$ be the number of
vertices and edges of $G$, and let $k$ be the number of boundary
vertices of $G$.

(1) $k \geq E-V$, and equality holds if and only if all disk face of
$G$ are triangles, all non-disk faces are annuli, and each cut vertex
has exactly two corners on annular faces.

(2) $G$ has at most $2V$ edges.
\end{lemma}

\proof (1) Let $D$ be the number of disk faces of $G$.  Then $0 =
\chi(T) \leq V - E + D$, and equality holds if and only if all
non-disk faces are annuli.  Thus $D \geq E - V$.  For each vertex $u$
of $G$, let $\delta(u)$ be the number of corners of disk faces
incident to $u$.  Then $\sum_u val(u) = 2E$, and $\sum_u \delta(u)
\geq 3D \geq 3(E -V)$.  Since there is no isolated or interior vertex,
we have $val(u) - \delta(u) \geq 1$, and equality holds if and only if
$u$ is a boundary vertex.  Let $p$ be the number of non-boundary
vertices.  Then $$
p \leq \sum_u(val(u) - \delta(u) - 1) \leq 2E -
3(E-V) - V = 2V - E.$$
It follows that the number of boundary vertices
is $k = V - p \geq E - V$, and equality holds if and only if (i)
$V-E+D=0$, i.e.\ non-disk faces are annuli, (ii) $\sum_u \delta(u) =
3D$, so all disk faces are triangles, and (iii) $val(u) - \delta(u) -
1 = 1$ for any cut vertex, i.e.\ each cut vertex has exactly two
corners not on disk faces.

(2) Since the number of boundary vertices is at most $V$, by (1) we
have $V \geq k \geq E - V$, hence $E \leq 2V$.
\qed

\begin{lemma} Suppose all interior vertices of $\rgap$ have
valence at least 6, and all boundary vertices of $\rgap$ have valence
at least 4.  Let $G$ be a component of $\rgap$.  Then either (i) $G$
is topologically an essential circle on the torus $\hat F_a$, or (ii)
$G$ has no cut vertex, all interior vertices of $G$ are of valence
exactly 6, and all boundary vertices of $G$ are of valence exactly
4.  \end{lemma}

\proof Let $X$ be the union of $G$ and all its disk faces.  If $X$ is
the whole torus then all vertices are interior vertices, and an easy
Euler characteristic argument shows that all vertices must be of
valence 6, so (ii) follows.  Also, by Lemma 2.9 $X$ is not in a disk
in $\hat F_a$ as otherwise $G$ would have a boundary vertex of valence at
most 3.  Therefore we may assume that $X$ has the homotopy type of
a circle.

First assume that $X$ has a cut vertex $v$.  Recall that $X$ is
homotopy equivalent to a circle, so if $X - v$ is not connected, then
$v$ cuts off a subcomplex $W$ of $X$ which lies in a disk in $\hat
F_a$.  By Lemma 2.9 the graph $G \cap W$ has at least two boundary
vertices of valence at most 3, hence at least one such vertex $v'$
other than $v$, which contradicts the assumption because $v'$ is then
a boundary vertex of $\rgap$ of valence at most 3.  Therefore we
may assume that $X - v$ is connected.  Since $X$ has the homotopy type
of a circle, $X$ cut at $v$ is a simply connected planar complex $W$,
and $X$ is obtained by identifying exactly two points of $W$.  Let
$G'$ be the corresponding graph on $W$.  We may assume that $X$ is not
a circle as otherwise (i) is true.  Thus $W$ is not homeomorphic to an
arc.  Therefore by Lemma 2.9 we have $\psi(G') \geq 3$, hence $G'$ has
at least one boundary vertex $v'$ of valence at most 3 which is not
identified to $v$ in $G$.  By definition $v'$ is a boundary vertex of
$\rgap$ of valence at most 3, which is a contradiction.  This
completes the proof that $X$ has no cut vertex.

We may now assume that $X$ is an annulus, so all vertices of $G$ are
either interior vertices of valence at least 6 in the interior of $X$,
or boundary vertices of valence at least 4 on $\bdd X$.
Consider the double $G''$ of $G$ on the double of $X$ along $\bdd X$.
Since each boundary vertex of $G$ of valence $k$ gives rise to a
vertex of valence $2k - 2$ in $G''$, we see that $G''$ is a reduced
graph on a torus such that all of its vertices have valence at
least 6.  An Euler characteristic argument shows that all vertices of
$G''$ must have valence exactly 6, hence (ii) follows.
\qed

\begin{lemma} If $M(r_a)$ contains a Klein bottle $K$, then

(1) $T = \bdd N(K)$ is an essential torus in $M(r_a)$; and

(2) $K$ intersects the core $K_a$ of the Dehn filling solid torus
at no less than $n_a/2$ points.  \end{lemma}

\proof $T$ bounds a twisted $I$-bundle over the Klein bottle $N(K)$ on
one side.  Since $M(r_a)$ is assumed irreducible, if $T$ is
compressible on the other side then $M(r_a)$ is a Seifert fiber space
over a sphere with (at most) three singular fibers of indices
$(2,2,p)$ for some $p$, and if $T$ is boundary parallel then $M(r_a)$
is a twisted $I$-bundle over the Klein bottle.  Either case
contradicts the assumption that $M(r_a)$ is toroidal.  Therefore $T$
is an essential torus.  If $|K \cap K_a| < n_a/2$ then $T$ would
intersect $K_a$ in less than $n_a$ points, contradicting the choice of
$n_a$.  \qed

\begin{lemma} Suppose $n_a>2$, and $\gb$ has both a
$12$-Scharlemann bigon $e_1 \cup e_2$ and a $23$-Scharlemann bigons
$e_3 \cup e_4$. If $e_1 \cup e_2$ and $e_3 \cup e_4$ are isotopic on
$\hat F_a$, then the disk face $D$ they bound on $\hat F_a$ contains
at least $(n_a/2)-1$ vertices in its interior.  \end{lemma}

\proof Let $m$ be the number of vertices in the interior of $D$.  Let
$D_1, D_2$ be the disk faces of $(12)$- and $(23)$-Scharlemann bigons
in $\gb$.  Shrinking the Dehn filling solid torus of $M(r_a)$ to its
core $K_a$, the union $D_1 \cup D_2 \cup D$ is a Klein bottle $Q$ in
$M(r_a)$.  A regular neighborhood of $Q$ intersects $K_a$ at an arc
from $u_1$ to $u_2$ then to $u_3$, and one arc for each vertex of
$\ga$ in the interior of $D$.  Hence $Q$ can be perturbed to intersect
$K_a$ at $1+m$ points.  By Lemma 2.12(2) we have $m+1 \geq n_a/2$,
hence the result follows.  \qed

An edge $e$ of $\ga$ is a {\it co-loop\/} edge if it has the same
label on its two endpoints, in other words, it is a loop on the other
graph $\gb$.  Given a codimension 1 manifold $X$ in a manifold $Y$,
use $Y|X$ to denote the manifold obtained by cutting $Y$ along
$X$.  

\begin{lemma} Let $\hat e$ be a family of negative edges in
$\ga$.  Let $G$ be the subgraph of $\gb$ consisting of the edges of
$\hat e$ and their vertices.

(1) Each cycle component of $G$ is an essential loop on $\hat F_b$.

(2) (The 3-Cycle Lemma.) $G$ cannot contain three disjoint cycles; in
particular, $\ga$ cannot have three parallel co-loop edges.

(3) (The 2-Cycle Lemma.)  If $\gb$ is positive then $G$ cannot contain
two disjoint cycles; in particular, $\ga$ cannot have two parallel
co-loop edges.  
\end{lemma}

\proof (1) Assume to the contrary that some cycle component of $G$ is
inessential on $\hat F_b$.  Let $D$ be a disk bounded by an innermost
cycle component of $G$, and let $D'$ be the bigon disks on $F_a$
between edges of $\hat e$.  Let $V_b$ be the Dehn filling solid torus
in $M(r_b)$.  Then a regular neighborhood $W$ of $D \cup V_b \cup D'$
is a solid torus containing the core of $V_b$ as a cable knot winding
along the longitude at least twice.  See the proof of [GLi,
Proposition 1.3].  In this case $W \cap M$ is a cable space, which is
a contradiction to the assumption that $M$ is a hyperbolic manifold.

(2) Let $\hat e = e_1 \cup ... \cup e_k$, oriented
consistently, with tails at $u'$ and heads at $u''$ on $\ga$.  Let $s$
be the transition number of $\hat e$.  We may assume that $e_i$ has
label $i$ at its tail, so it has label $i+s$ at its head.  Let $D_j$
be the bigon on $\hat F_a$ between $e_j$ and $e_{j+1}$.

If $k>n_b$ then by Lemma 2.3(1) the transition function associated
with $\hat e$ has only one orbit, hence we may assume $k\leq n_b$.  On
$\gb$ these edges form disjoint cycles and chains.  Assume there are
at least three cycles.  Then $e_1, e_2, e_3$ belong to three distinct
cycles $C_1, C_2, C_3$.  Thus for $i=1,2,3$, $$C_i = e_i \cup e_{i+s}
\cup ... \cup e_{i+(p-1)s}$$
is an oriented cycle on $\hat F_b$ for some fixed $p$.  By (1) these
are essential loops on $\hat F_b$, so they are parallel as unoriented
loops.  

Each bigon $D_{1+js}$ gives a parallelism between an edge of $C_1$ and
an edge of $C_2$, hence when shrinking the Dehn filling solid torus
$V_b$ to its core knot $K_b$, the union $A_1 = \cup D_{1+js}$ is an
annulus in $M(r_b)$ with $\bdd A_1 = C_1 \cup C_2$.  Similarly, $A_2 =
\cup D_{2+js}$ is an annulus in $M(r_b)$ with $\bdd A_2 = C_2 \cup
C_3$.  These $A_i$ are essential in $M(r_b)|\hat F_b$, the manifold
obtained from $M(r_b)$ by cutting along $\hat F_b$, otherwise $K_b$
would be isotopic to a curve having fewer intersections with $\hat
F_b$.

Let $A'_1, A'_2, A'_3$ be the annuli $\hat F_b | (C_1 \cup C_2 \cup
C_3)$, with $\bdd A'_i = C_i \cup C_{i+1}$ (subscripts mod 3.)  Let
$m_i$ be the number of times that $K_b$ intersects the interior of
$A'_i$.  Then
$$\sum m_i + 3p = n_b$$

The annulus $A_i$ is said to be of type I if a regular neighborhood of
$\bdd A_i$ lies on the same side of $\hat F_b$, otherwise it is of
type II.  Note that if $\hat F_b$ is separating then $A_i$ must be of
type I.  There are several possibilities.  In each case one can find
an essential torus $T'$ in $M(r_b)$ which has fewer intersections with
$K_b$.  This will contradict the choice of $\hat F_b$ and complete
the proof of (1).

\smallskip
Case 1.  {\it $C_2$ is anti-parallel to both $C_1$ and $C_3$. }

In this case each $T_i = A_i \cup A'_i$ is a Klein bottle for $i=1,2$,
which can be perturbed to intersect $K_b$ at $p + m_i$ points.  Since
$\sum m_i + 3p = n_b$, either $T_1$ or $T_2$ can be perturbed to
intersect $K_b$ at fewer than $n_b/2$ points, contradicting Lemma 2.12.

\smallskip
Case 2.  {\it $C_2$ is anti-parallel to $C_1$, say, and parallel to the
other cycle $C_3$.}

Let $T_1 = A_1 \cup A'_1$ and $T_2 = A_1 \cup A_2 \cup A'_3$.  Then
$T_i$ are Klein bottles, and they can be perturbed to intersect $K_b$
at $p+m_1$ and $m_3$ points, respectively.  One of these contradicts
Lemma 2.12.

\smallskip
Case 3.  {\it $C_2$ is parallel to both $C_1$ and $C_3$.}

If one of the $A_i$, say $A_1$, is of type II, then $T_1 = A_1 \cup
A'_1$ is a non-separating torus (because it can be perturbed to
intersect $\hat F_b$ transversely at a single circle), and it
intersects $K_b$ at $p+m_i < n_b$ points.  Since $M(r_b)$ is
irreducible, $T_1$ is incompressible and hence essential, which
contradicts the choice of $\hat F_b$.

If both $A_i$ are of type I then one can show that $A_1 \cup A_2 \cup
A'_3$ is an essential torus $T$ which can be perturbed to intersect
$K_b$ in $m_3 + p < n_b$ points.  The proof is standard: The torus
$\hat F_b$ and the annuli $A_1, A_2$ cut $M(r_b)$ into a manifold
whose boundary contains four tori $T_1 = A_1 \cup A'_1$, $T_2 = A_1
\cup A'_2 \cup A'_3$, $T_3 = A_2 \cup A'_2$, and $T_4 = A_2 \cup A'_1
\cup A'_3$.  Each of these tori $T_i$ can be perturbed to have fewer
than $n_b$ intersections with the knot $K_b$, and hence bounds a
manifold $W_i$ which is either a solid torus or a $T^2 \times I$
between $T_i$ and a component of $\bdd M_(r_b)$.  Moreover, if $W_i$
is a solid torus $W_i$ then the annulus $T_i \cap \hat F_b$ is
essential on $\bdd W_i$ in the sense that it is neither meridional nor
longitudinal (otherwise $\hat F_b$ would be compressible or could be
isotoped to have fewer intersections with $K_b$).  Now we have $M(r_b)
= (W_1 \cup W_4) \cup (W_2 \cup W_3) = W' \cup W''$, with $W' \cap
W''$ a torus $T = A_1 \cup A_2 \cup A'_3$ which can be perturbed to
intersect $K_b$ at $m_3 + p < n_b$ points.  Since $W' = W_1 \cup
_{A'_1} W_4$ and $A'_1$ is essential in both $W_1$ and $W_4$, $T$ is
incompressible and not boundary parallel in $W'$; similarly for $W''$.
It follows that $T$ is a contradiction to the choice of $\hat F_b$.

(3) The proof of this part is much simpler.  Let $A_1, A'_1$ be as
above, and let $A''_1$ be the complement of $A'_1$ on $\hat F_b$.  If
$C_1, C_2$ are parallel then $A_1 \cup A'_1$ is a nonseparating torus
in $M(r_b)$ which can be perturbed to intersect $K_b$ less than $n_b$
times, contradicting the choice of $\hat F_b$.  If $C_1, C_2$ are
anti-parallel then $A_1 \cup A'_1$ and $A_1 \cup A''_1$ are Klein
bottles, which can be perturbed to intersect $K_b$ at a total of
$n-2p$ points, where $p$ is the number of vertices in $C_i$; hence one
of those will intersect $K_b$ less than $n_b/2$ times, which
contradicts Lemma 2.12.  \qed

When studying Dehn surgery via intersection graphs, we usually fix the
surfaces $F_1, F_2$, and hence the graphs $\Gamma_1, \Gamma_2$ are
also fixed.  The following technique will allow us to modify the
surfaces and hence the graphs in certain situation.  Lemma 2.15 will
be used in the proofs of Lemmas 12.16 and 19.6.

Consider two surfaces $F_1, F_2$ in a 3-manifold $M$ with boundary
slopes $r_1, r_2$ respectively and suppose they intersect minimally.
Let $\ga, \gb$ be the intersection graphs on $\hat F_1, \hat F_2$,
respectively.  Let $\alpha$ be a proper arc on a disk face $D$ of
$\ga$ with boundary on edges of $\ga$.  Then one can replace two small
arcs of $\ga$ centered at $\bdd \alpha$ by two parallel copies of
$\alpha$ to obtain a new graph $\ga'$, called the graph obtained from
$\ga$ by {\it surgery along $\alpha$}.  

\bigskip
\leavevmode

\centerline{\epsfbox{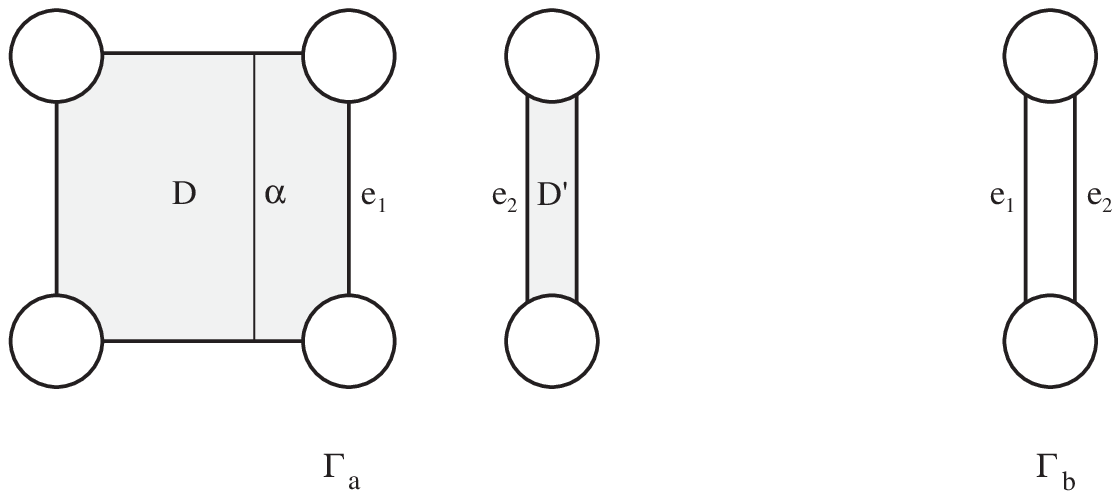}}
\bigskip
\centerline{Figure 2.1}
\bigskip

A face $D'$ of $\ga$ is called a {\it coupling face\/} to another face
$D$ of $\ga$ along an edge $e_1$ of $D$ if $D'$ has an edge $e_2$ such
that $e_1, e_2$ are adjacent parallel edges on $\gb$, and the
neighborhoods in $D$ and $D'$ of the $e_i$'s lie (locally) on the same
side of $\hat F_b$.  Note that this is independent of whether $\hat
F_b$ is orientable or separating in $M$.  See Figure 2.1.  By
definition $D$ has no coupling face along $e_1$ if $e_1$ has no
parallel edge on $\gb$, one coupling face along $e_1$ if $e_1$ has
some parallel edges and is a border edge of the family, and two
otherwise.  A 4-gon face $D$ of $\ga$ looks like a ``saddle surface''
in $M | F_b$.  In general it is not possible to push the saddle up or
down to change the intersection graph.  However, if some coupling face
to an edge of $D$ is a bigon then this is possible.  See Figure 2.2.
More explicitly, we have the following lemma.

\bigskip
\leavevmode

\centerline{\epsfbox{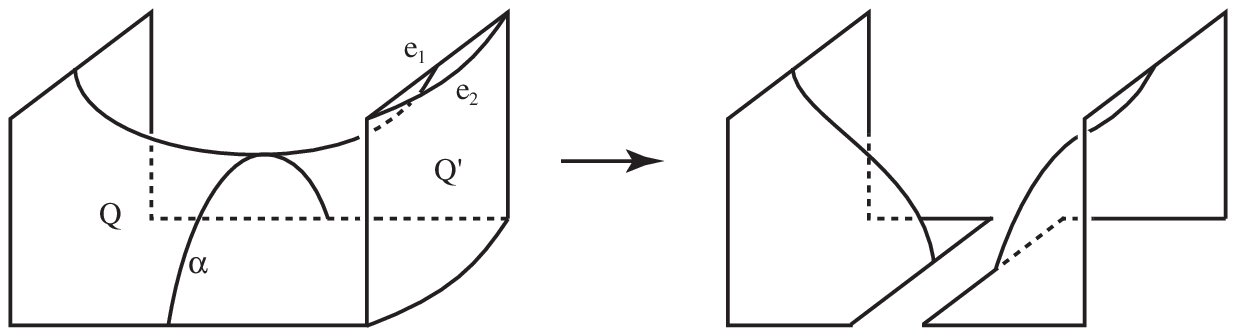}}
\bigskip
\centerline{Figure 2.2}
\bigskip

\begin{lemma} Let $\Gamma_1, \Gamma_2$ be a pair of
intersection graphs.  Let $Q$ be a face of $\ga$, and let $e$ be an
edge on $\bdd Q$.  Let $\alpha$ be an arc on $Q$ with boundary in the
interior of edges of $\ga$, cutting off a disk $B_1$ containing $e$
and exactly two corners of $Q$.  If some coupling face $Q'$ of $Q$
along $e$ is a bigon, then $F_a$ can be isotoped so that the new
intersection graph $\ga'$ is obtained from $\ga$ by surgery along
$\alpha$.  \end{lemma}

\proof Cut $M$ along $F_b$.  Then the face $Q$ is as shown in Figure
2.2.  Let $B_2$ be the bigon in $\gb$ between $e$ and the edge $e'$ on
$Q'$.  After shrinking the Dehn filling solid torus $V_b$ to its core
knot $K_b$, the union $B_1 \cup B_2 \cup Q'$ is a disk $Q''$ with
boundary the union of $\alpha$ and an arc on $\hat F_b$.  Pushing
$Q''$ off $K_b$ gives a disk $P$ in $M$ which has boundary the union
of $\alpha$ and an arc on $F_b$, and has interior disjoint from $F_b
\cup F_a$.  Therefore we can isotope $F_a$ through this disk $P$ to
get a new surface $F'_a$.  It is clear that the new intersection graph
$\ga'$ is obtained from $\ga$ by surgery along $\alpha$.  \qed

Let $u$ be a vertex of $\Gamma_a$, and $P, Q$ two edge endpoints on
$\bdd u$.  Let $I$ be the interval on $\bdd u$ from $P$ to $Q$ along
the direction induced by the orientation of $u$.  The edge endpoints
of $\ga$ cut $I$ into $k$ subintervals for some $k$.  Then the {\it
  distance\/} from $P$ to $Q$ on $\bdd u$ is defined as $d_u(P, Q) =
k$.  Some times we also use $d_{\Gamma_a}(P,Q)$ to denote $d_u(P,Q)$.
If $P, Q$ are the only edge endpoints of $e_1, e_2$ on $\bdd u$,
respectively, then we define $d_u(e_1, e_2) = d_u(P,Q)$.  Notice that
if the valence of $u$ is $m$, then $d_u(Q, P) = m - d_u(P, Q)$.
The following lemma can be found in [Go].

\begin{lemma} {\rm [Go, Lemma 2.4]} \quad
(i) Suppose $P, Q \in \bdd u_i \cap \bdd v_k$ and $R, S \in \bdd u_j
\cap \bdd v_l$.  If $d_{u_i}(P, Q) = d_{u_j}(R,S)$ then $d_{v_k}(P,Q)
= d_{v_l}(R,S)$.

(ii) Suppose that $P \in u_i \cap v_k$, $Q \in u_i \cap v_l$, $R \in
u_j \cap v_k$, and $S \in u_j\cap v_l$.  If $d_{u_i}(P, Q) =
d_{u_j}(R, S)$, then $e_{v_k}(P,R) = d_{v_l}(Q,S)$.  
\end{lemma}

Suppose two edges $e_1, e_2$ of $\ga$ connect the same pair of
vertices $u_i, u_j$.  Let $p_k, q_k$ be the endpoints of $e_k$ on
$u_i, u_j$, respectively, $k=1,2$.  Then $e_1, e_2$ are {\it
  equidistant\/} if $d_{u_i}(p_1, p_2) =d_{u_j}(q_2, q_1)$.  (Note
that the orders of the edge endpoints have been reversed.)  Thus for
example a pair of parallel positive edges is always equidistant, but a
pair of parallel negative edges is not unless their distance is
exactly half of the valence of the vertices.

Note that when $u_i \neq u_j$ the above equation can be written as
$d_{u_i}(e_1, e_2) = d_{u_j}(e_2, e_1)$.  When $u_i = u_j$,
$d_{u_i}(e_1, e_2)$ is not defined, and there are two choices for the
pair $p_k, q_k$, but one can check that whether the equality
$d_{u_i}(p_1, p_2) =d_{u_j}(q_2, q_1)$ holds is independent of the
choice of $p_i, q_i$.

The following lemma is called the {\it Equidistance Lemma}.  It
follows from Lemma 2.16, and can also be found in [GW1].

\begin{lemma} {\rm [GW1, Lemma 2.8]} \quad Let $e_1, e_2$ be a
pair of edges with $\bdd e_1 = \bdd e_2$ in both $\Gamma_1$ and
$\Gamma_2$.  Then $e_1, e_2$ are equidistant in $\Gamma_1$ if and only
if they are equidistant in $\Gamma_2$.  \end{lemma}

Given two oriented slopes $r_1, r_2$ on $T_0$, choose an oriented
meridian-longitude pair $m,l$ on the torus $T_0$ so that $r_1 = m$,
then the slope $r_2$ is homologous to $Jm+\Delta l$ for some mod
$\Delta$ integer $J = J(r_1, r_2)$, called the {\it jumping number\/}
between $r_1, r_2$.  Note that if $\Delta = 4$, then $J =\pm 1$, and
if $\Delta = 5$, then $J = \pm 1$ or $\pm 2$.  The following lemma is
call the {\it Jumping Lemma} and can be found in [GW1].

\begin{lemma} {\rm [GW1, Lemma 2.10]} \quad Let $P_1, \ldots,
P_{\Delta}$ be the points of $\bdd u_i \cap \bdd v_j$, labeled
successively on $\bdd u_i$.  Let $J = J(r_1, r_2)$ be the jumping
number of $r_1, r_2$.  Then on $v_j$ these points appear in the order
of $P_J, P_{2J},\ldots, P_{\Delta J}$.  In particular, they appear
successively as $P_1,...,P_{\Delta}$ along some direction of $\bdd
v_j$ if and only if $J = \pm 1$.
\end{lemma}

\begin{lemma} Let $e_1 \cup ... \cup e_p$ and $e'_1 \cup ...
\cup e'_q$ be two sets of parallel edges on $\ga$.  Suppose $e_1$ is
parallel to $e'_1$ and $e_p$ parallel to $e'_q$ on $\gb$.  Then $p =
q$.  \end{lemma}

\proof Let $D_1, D_2, D_3, D_4$ be the disks realizing the
parallelisms of $e_1 \cup e_p$ and $e'_1 \cup e'_q$ on $\ga$, and $e_1
\cup e'_1$ and $e_p \cup e'_q$ on $\gb$.  Then the union $A = D_1 \cup
... \cup D_4$ is a M\"obius band or annulus in $M$ with boundary on
$T_0$.  (It is embedded in $M$, otherwise there is a pair of edges
parallel in both graphs, contradicting Lemma 2.2(2).)  If $A$ is a
M\"obius band then it is already a contradiction to the hyperbolicity of
$M$.  If $A$ is an annulus and $p \neq q$ then a boundary component
$c$ of $A$ has intersection number $p-q \neq 0$ with $\cup \bdd v_i$
and hence is an essential curve on $T_0$.  Since $e_1$ is an essential
arc on both $A$ and $F_1$ and $F_1$ is boundary incompressible, $A$
cannot be boundary parallel.  It follows that $A$ is an essential
annulus in $M$, which again contradicts the assumption that $M$ is
hyperbolic.  \qed

\begin{lemma} Suppose $\gb$ is positive, $n_b \geq 3$, and
$\ga$ contains bigons $e_1 \cup e_2$ and $e'_1 \cup e'_2$, such that
$e_1,e'_1$ have label pair $\{i, j\}$ and $e_2,e'_2$ have label pair
$\{i+1, j+1\}$, where $j \neq i$.  Let $C_1 = e_1 \cup e'_1$ and $C_2
= e_2 \cup e'_2$ be the loops on $\hat F_b$.  If $C_1$ is essential on
$\hat F_b$ then $C_2$ is essential on $\hat F_b$ and not homotopic to
$C_1$.  \end{lemma}

\proof
Let $B$ and $B'$ be the bigon faces bounded by $e_1 \cup e_2$ and
$e'_1 \cup e'_2$, respectively.  Shrinking the Dehn filling solid
torus to the core knot $K_b$, the union $B \cup B'$ becomes an annulus
$A_1$ in $M(r_b)$ with boundary $C_1 \cup C_2$.  Since $\hat F_b$ is
incompressible and $C_1$ is essential on $\hat F_b$, it follows that
$C_2$ must also be essential on $\hat F_b$.

Now assume $C_i$ are essential and homotopic on $\hat F_b$.  
Since $i\neq j$ and $n_b>2$, $C_1, C_2$ have at most one
vertex in common.  If $C_1, C_2$ are disjoint, let $A_2$ be an annulus
on $\hat F_b$ bounded by $C_1 \cup C_2$.  If $C_1, C_2$ has a common
vertex $v_{i+1} = v_j$, let $A_2$ be the disk face of $C_1 \cup C_2$
in $\hat F_b$, which will be considered as a degenerate annulus as it
can be obtained from an annulus by pinching an essential arc to a
point.  Let $A'_2$ be the closure of $\hat F_b - A_2$.  Let $m$ and
$m'$ be the number of vertices in the interior of $A_2$ and $A'_2$,
respectively.  Then $n_b = m+m'+k$, where $k$ is the number of vertices
on $C_1 \cup C_2$, i.e., $k=4$ if $C_1 \cap C_2 = \emptyset$, and
$k=3$ otherwise.

First consider the case that $C_1 \cap C_2 = \emptyset$.  Orient $C_1,
C_2$ so that they are parallel on the annulus $A_1$.  If they are also
parallel on $\hat F_b$ then $A_1 \cup A_2$ is a nonseparating torus
which can be perturbed to intersect $K_b$ at $m+2<n$ points, which is
a contradiction.  If they are anti-parallel then $A_1 \cup A_2$ and
$A_1 \cup A'_2$ are Klein bottles which can be perturbed to intersect
$K_b$ at $m$ and $m'$ points, respectively.  Since at least one of $m,
m'$ is less than $n_b/2$, this contradicts Lemma 2.12.

The case that $C_1 \cap C_2 \neq \emptyset$ is similar.  If $C_1, C_2$
are parallel then $A_1 \cup A_2$ is a torus and can be perturbed to
intersect $K_b$ at $m+1<n$ points; if they are anti-parallel then $A_1
\cup A_2$ and $A_1 \cup A'_2$ can be perturbed to be Klein bottles
intersecting $K_b$ at $m+1$ and $m'$ points, respectively, which leads
to contradictions as above because $m+1+m'<n_b$ implies either
$2(m+1)<n_b$ or $2m' < n_b$.
\qed

A triple of edge endpoints $(p_1, p_2, p_3)$ on $\gb$ is {\it
  positive\/} if they appear on the boundary of the same vertex $v_i$,
and in this order on $\bdd v_i$ along the orientation of $\bdd v_i$.
Note that this is true if and only if $d_{v_i}(p_1, p_2) +
d_{v_i}(p_2, p_3) = d_{v_i}(p_1, p_3)$.

\begin{lemma} (1) Suppose $(p_1,p_2, p_3)$ is a positive
triple on $\gb$.  Let $k$ be a fixed integer and let $p'_i$ be edge
endpoints such that $d_{\ga}(p_i, p'_i) = k$ for all $i$.  Then
$(p'_1, p'_2, p'_3)$ is also a positive triple on $\gb$.

(2) Let $e_1 \cup ... \cup e_r$ be a set of parallel negative edges
with end vertices $u_1, u_2$ in $\ga$.  Let $u(p) \in \{u_1,u_2\}$ for
$p=1,2,3$, and let $e_j(u(p))$ be the endpoint of $e_j$ at $u(p)$.  If
$(e_i(u(1)), e_j(u(2)), e_k(u(3)))$ is a positive triple and $i, j, k
\leq r-t$, then $(e_{i+t}(u(1)), e_{j+t}(u(2)), e_{k+t}(u(3)))$ is
also a positive triple.  \end{lemma}

\proof  (1) Geometrically this is obvious:  Flowing on $T_0$ along $\bdd
F_a$ moves the first triple to the second triple, hence the
orientations of the components of $\bdd F_b$ containing these triples
are the same on $T_0$.  

Alternatively one may use Lemma 2.16(ii) to prove the result.  Since
$d_{\ga}(p_i, p'_i) = d_{\ga}(p_j, p'_j) = k$ for all $i,j$, by Lemma
2.16(ii) we have $d_{\gb}(p_i, p_j) = d_{\gb}(p'_i, p'_j)$ for all
$i,j$.  Therefore $d_{\gb}(p_1, p_2) + d_{\gb}(p_2, p_3) =
d_{\gb}(p_1, p_3)$ if and only if $d_{\gb}(p'_1, p'_2) + d_{\gb}(p'_2,
p'_3) = d_{\gb}(p'_1, p'_3)$.

(2) This is a special case of (1) because 
\begin{eqnarray*}
d_{\ga}(e_i(u(1)), e_{i+t}(u(1))) & = & d_{\ga}(e_j(u(2)), e_{j+t}(u(2))) \\
& = & d_{\ga}(e_k(u(3)), e_{k+t}(u(3))) = t \qquad 
\end{eqnarray*}
\qed

\begin{lemma} Suppose $\gb$ is positive and $n = n_b \geq 3$. 

(1) Suppose $\hat e \supset e_1 \cup ... \cup e_{n+2}$, and the
transition number $s = 1$.  Let $A$ be the annulus obtained by cutting
$\hat F_b$ along the cycle $e_1 \cup ... \cup e_n$.  Then the
edges $e_{n+1}, e_{n+2}$ lie in $A$ as shown in Figure 2.3, up to
reflection along the center circle of the annulus.

(2) If $s_1 = 1$, $n=3$ and $\hat e_1$ contains $6$ edges $e_1 \cup
... \cup e_6$ then the edges are as shown in Figure 2.4.

(3) Any family of parallel negative edges in $\ga$ contains at most
$2n$ edges.
\end{lemma}

\proof (1) Let $u, u'$ be the end vertices of $\hat e$ in $\ga$.
Orient $e_i$ from $u$ to $u'$ and assume without loss of generality
that $e_i$ has label $i$ at its tail $e_i(t)$ in $\ga$.  Since $s =
1$, the head of $e_i$, denoted by $e_i(h)$, has label $i+1$ in $\ga$.
The edges $e_1, ..., e_n$ form an essential loop on the torus $\hat
F_b$.  Cutting $\hat F_b$ along this loop produces an annulus, as
shown in Figure 2.3.

Up to reflection along the center circle of the annulus we may
assume that the edge $e_{n+1}$ appears in this annulus as shown in
Figure 2.3.  We need to prove that $e_{n+2}$ appears in $\gb$ as shown
in the figure.

\bigskip
\leavevmode

\centerline{\epsfbox{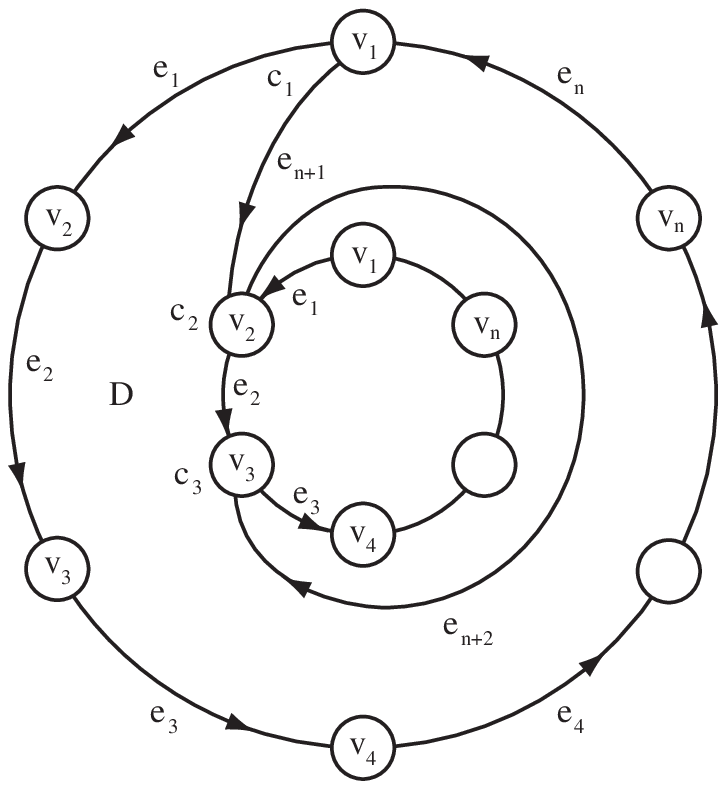}}
\bigskip
\centerline{Figure 2.3}
\bigskip

Since $\gb$ is positive, we may assume that all vertices on Figure 2.3
are oriented counterclockwise.  Note that $(e_1(h), e_{n+1}(h),
e_2(t))$ is a positive triple on $\gb$.  By Lemma 2.21(2) the triple
$(e_2(h), e_{n+2}(h), e_3(t))$ is also a positive triple.  This
determines the location of the head of $e_{n+2}$, as shown in Figure
2.3.  Applying Lemma 2.20 to $e_1 \cup e_2$ and $e_{n+1} \cup
e_{n+2}$, we see that the loop $e_2 \cup e_{n+2}$ is essential and not
homotopic to $e_1 \cup e_{n+1}$ on $\hat F_b$, so these two loops must
intersect transversely at the common vertex $v_2$ on $\hat F_b$.
Hence the edge $e_{n+2}$ must appear as shown in Figure 2.3.

(2) By (1) the first 5 edges must be as shown in Figure 2.4.  These
cut the torus into a 3-gon and a 7-gon.  The edge $e_6$ is not
parallel to the other $e_i$'s on $\gb$ and hence must lie in the
7-gon, connecting $v_3$ to $v_1$.  For the same reason as above,
$(e_3(h), e_6(h), e_4(t))$ is a positive triple on $\bdd v_1$, hence
the head of $e_6$ must be in the corner on $\bdd v_1$ from $e_1(t)$ to
$e_4(t)$ because the corner from $e_3(h)$ to $e_1(t)$ lies in the
3-gon.  Similarly, since $(e_1(h), e_5(t), e_4(h))$ is a positive
triple on Figure 2.4, by Lemma 2.21(2) $(e_2(h), e_6(t), e_5(h))$ is
also a positive triple, which determines the position of the tail of
$e_6$.  Therefore $e_6$ must be as shown in Figure 2.4.

(3) This follows from [Go, Corollary 5.5] when $n\geq 4$.  Now assume
$n=3$ and suppose there exist $2n+1 = 7$ parallel edges $e_1 \cup ...
\cup e_7$ on $\ga$.  By the 3-Cycle Lemma 2.14(2) we may assume that the
transition number $s \neq 0$.  Since $n= 3$, we may assume without
loss of generality that $s = 1$, hence by (2) the subgraph of $\gb$
consisting of the edges $e_1 \cup ... \cup e_6$ is as shown in Figure
2.4.  By the same argument as above, $(e_4(h), e_7(h), e_5(t))$ and
$(e_3(h), e_7(t), e_6(h))$ are positive triples on $\bdd v_2$ and
$\bdd v_1$, respectively.  Since $e_7$ must lie in the $6$-gon face
$D$ in Figure 2.4, this is possible only if $e_7$ is parallel to
$e_1$, which is a contradiction to Lemma 2.2(2).  \qed

\bigskip
\leavevmode

\centerline{\epsfbox{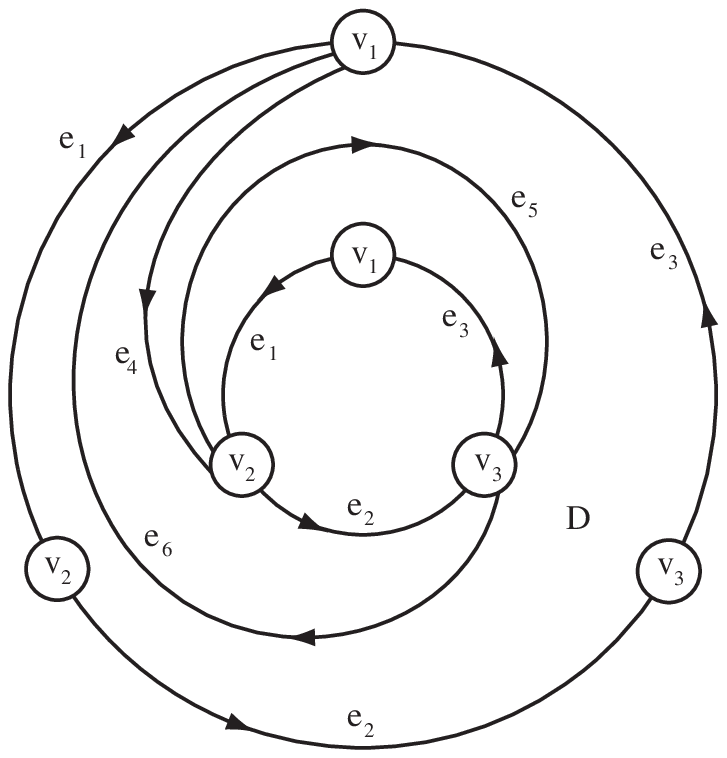}}
\bigskip
\centerline{Figure 2.4}
\bigskip

\begin{lemma}  Let $sign(v)$ be the sign of a vertex $v$, and
define $p_a$ to be the sum of the signs of the vertices of $\ga$.
Then either $p_1=0$ or $p_2=0$.  In particular, $n_1, n_2$ cannot both
be odd. 
\end{lemma}

\proof
For each edge endpoint $c$ on $u_i \cap v_j$, define $\sign(c) =
\sign(u_i)\; \sign(v_j)$.  Then the parity rule says that the 
two endpoints of an edge $e$ have different sign.  Summing over all
edge endpoints on $\ga$ gives 
$$0 = \sum_{i,j}  \Delta (\sign(u_i)\; \sign(v_j)) =
\Delta \sum_i \sign(u_i)\; \sum_j \sign(v_j) = \Delta p_1 p_2
$$
hence either $p_1 = 0$ or $p_2 = 0$.
\qed

\section{$\rgap$ has no interior vertex}

In this section we will show that if $\nb > 4$ then the graph $\ga$
does not have interior vertices; in particular the vertices of $\ga$
cannot all be parallel.  Recall that we have assumed that $\Delta \geq 4$.

\begin{lemma} If $\nb > 4$ then $\rga$ has no full vertex of
valence at most 6.  \end{lemma}

\proof 
First assume $n_b \neq 6$.  Then by Lemma 2.7(3), three adjacent
families of positive edges in $\ga$ contain at most $2n_b -1$ edges,
hence if $\ga$ has a full vertex of valence at most 6 then 
$$ 4n_b \leq \Delta n_b \leq 2(2n_b - 1) = 4n_b - 2,$$
a contradiction.

So suppose $n_b = 6$, and let $u_1$ be a full vertex of $\rga$ of
valence at most $6$.  By Lemma 2.3(3), each family of edges incident
to $u_1$ contains at most $4$ edges.  Since there are at least $24$
edges in at most $6$ families, there must be exactly 6 families, each
containing exactly 4 edges.  Each family contains a Scharlemann bigon,
so there are six Scharlemann bigons at the vertex $u_1$.  Since there
are no extended Scharlemann cycles, each Scharlemann bigon appears at
one end of a family of parallel edges.  Thus by examining the labels
around the vertex $u_1$, one can see that if one Scharlemann bigon has
label pair $\{1,2\}$ then the others must have label pair $\{1,2\}$,
$\{3,4\}$ or $\{5,6\}$, and that at least two pairs do occur as label
pairs of Scharlemann bigons.  On the other hand, by Lemma 2.3(4) all
three pairs cannot appear as label pairs of Scharlemann bigons.
Hence, without loss of generality, there are incident to $u_1$ at
least three $(12)$-Scharlemann bigons and a $(34)$-Scharlemann bigon.
Since on $\gb$ the edges of the $(34)$-Scharlemann bigon form an
essential loop on $\rgb$, there are at most two edges of $\hat
F_b$ joining $v_1$ to $v_2$. Since the three $(12)$-Scharlemann bigons
give rise to six negative $1$-edges of $\gb$ joining $v_1$ to $v_2$,
three of these must be parallel, contradicting Lemma 2.7(1).
\qed

\begin{lemma} If $\nb > 4$ then $\rgap$ has no interior
vertices.  \end{lemma}

\proof This follows from Lemma 3.1 if $n_a \leq 2$ because in this
case either $\rgap = \rga$ and there is a full vertex of valence at
most $6$, or $n_a = 2$ and there is no interior vertex.  Therefore we
may assume that $n_a \geq 3$.

Suppose to the contrary that $\rgap$ has an interior vertex $u_i$.  By
Lemma 3.1 all interior vertices of $\rgap$ have valence at least 7,
hence we can apply Lemma 2.11 to conclude that $\rgap$ has a boundary
vertex $u_1$ of valence at most 3.

By Lemma 2.7(3) the three families of adjacent positive edges at $u_1$
contain at most $2n_b$ edges, hence there are $2n_b$ adjacent negative
edges.  On $\gb$ this implies that each vertex $v_j$ is incident to
two positive edges with label $1$ at $v_j$, which cannot be parallel
as otherwise there would be at least $n_a+1 > n_a/2 + 2$ parallel
positive edges, contradicting Lemma 2.3(3).  Therefore the reduced
graph $\rgb$ contains at least $n_b$ positive edges.  On the other
hand, the existence of an interior vertex in $\rgap$ implies that
$\rgb$ contains at least $2n_b$ negative edges, as shown in the proof
of Lemma 3.1.  Since $\rgb$ has at most $3n_b$ edges (Lemma 2.5), it
must have exactly $n_b$ positive edges and $2n_b$ negative edges.
Since we have shown above that each vertex in $\rgb$ is incident to at
least two positive edges, it follows that it is incident to exactly
two positive edges.

We claim that a family of parallel positive edges in $\gb$ contains at
most $n_a/2$ edges.  If such a family contains more than $n_a/2$
edges, then there is a Scharlemann bigon on one side of the family,
and by looking at the labels one can see that all labels appear among
the endpoints of this family, which is impossible because $u_i$ being
an interior vertex in $\rgap$ implies that all edges in $\gb$ with $i$
as a label are negative.

Since each vertex $v_j$ is incident to two families of positive edges,
each containing at most $n_a/2$ edges, we see that $v_j$ is incident
to at least $3n_a$ negative edges.  By Lemmas 3.1 and 2.5 we see that
$\rga$ has less than $3n_a$ positive edges, hence two of the negative
edges incident to $v_j$ are parallel in $\ga$, so $j$ is a label of a
Scharlemann bigon in $\ga$.  Since this is true for all vertices in
$\gb$, by Lemma 2.3(4) we have $n_b \leq 4$, which is a contradiction.
\qed

\section{Possible components of $\rgap$ }

\begin{lemma} Suppose $\rgap$ has no isolated vertex or
interior vertex.  If some $v_i$ of $\gb$ is incident to more than
$2n_{a}$ negative edges in $\gb$, or if $n_a > 4$ and $v_i$ is
incident to at most two families of positive edges in $\gb$, then $i$
is a label of a Scharlemann bigon in $\ga$.  \end{lemma}

\proof If $v_i$ is a vertex of $\gb$ incident to more than $2n_{a}$
negative edges then by Lemma 2.10(2) two of them are parallel in
$\ga$, so by Lemma 2.4 they form a Scharlemann bigon, hence $i$ is the
label of a Scharlemann bigon in $\ga$.  If $n_a>4$ and $v_i$ is
incident to two families of positive edges in $\gb$ then by Lemma
2.3(3) each family contains less than $n_{a}$ edges, hence $v_i$ is
incident to more than $2n_{a}$ negative edges and the result follows
from the above.  \qed

In the rest of this section we assume $n_{a} > 4$ for $a = 1,2$, and
$\Delta \geq 4$.  By Lemma 3.2 $\rgap$ has no interior vertices.  We
will show that each component of $\rgap$ must be one of the 11 graphs
in Figure 4.2.

\begin{lemma} No vertex $u$ of $\rga$ is incident to at most
four positive edges and at most one negative edge.  \end{lemma}

\proof By Lemmas 2.7(1) and 2.7(2) a family of negative edges contains
at most $n_{b}$ edges, and four adjacent families of positive edges
contain at most $2(n_{b} + 2) = 2n_{b} + 4$ edges.  Since $\ga$ has
at least $4n_{b}$ edges incident to $u$, we would have $n_{b} \leq
4$, which is a contradiction to our assumption.  \qed

\begin{lemma} Suppose $u_i$ is incident to at most three
positive edges in $\rga$, and if there are three then two of them are
adjacent.  Then $i$ is a label of a Scharlemann bigon in $\gb$.
\end{lemma}

\proof In this case each label appears at the endpoint of some
negative edge at $u_i$, so $\rgbp$ has no isolated vertex.  By Lemma
4.1 the result is true if $u_i$ is incident to more than $2\nb$
negative edges.  So we assume that $u_i$ is incident to no more than
$2n_{b}$ negative edges, and hence at least $2n_{b}$ positive edges.
By Lemma 2.7(2) the two adjacent families of positive edges contain at
most $\nb+2$ edges, while the other positive family contains no more
than $\nb/2+2$ edges.  Thus $(\nb+2) + (\nb/2+2) \geq 2\nb$, which
gives $\nb \leq 8$.  Since one of the positive families contains more
than $\nb/2$ edges, it contains a Scharlemann bigon; by Lemma 2.2(4)
$\nb$ must be even, so $\nb = 8$ or $6$.  Using the above inequality
and the fact that when $\nb = 6$ each positive family contains at most
4 edges (Lemma 2.3(3)), we see that $u_i$ is incident to exactly
$2\nb$ positive edges and $2\nb$ negative edges.  Dually, this implies
that in $\gb$ there are exactly $2\nb$ positive $i$-edges and $2\nb$
negative $i$-edges.  (As always, an edge with both endpoints labeled
$i$ is counted twice.)

If $i$ is not a label of a Scharlemann bigon in $\gb$ then the $2\nb$
positive $i$-edges in $\gb$ are mutually nonparallel, so $\rgb$ has at
least $2n_{b}$ positive edges.  By Lemma 2.5 the reduced graph $\rgb$
has no more than $3\nb$ edges, so it has at most $\nb$ negative edges.
On the other hand, by Lemma 2.7(1) each family of parallel negative
edges in $\gb$ has at most two endpoints labeled $i$; since there are
$2\nb$ such endpoints, $\gb$ must have at least $\nb$ families of
negative edges.  It follows that $\gb$ has exactly $\nb$ families of
negative edges, each having exactly two endpoints labeled $i$.

Suppose $\nb = 6$.  Then there are 12 edges in the three families
incident to $u_i$, and by Lemma 2.3(3) each family contains at most
four edges, hence each family contains exactly four edges.  If some of
these edges are loops, then there are four loops and four non-loop
edges.  No loop can be parallel to a non-loop edge in $\gb$
since otherwise the label $i$ would appear three times among a set of
parallel edges in $\gb$.  It follows that all the 8 positive edges
incident to $u_i$ are mutually nonparallel in $\gb$, so the reduced
graph $\rgb$ would have at least 8 negative edges, which is a
contradiction as we have shown above that $\rgb$ has exactly $n_b = 6$
negative edges.  Hence we can assume there is no loop based at $u_i$.
Note that a family of four parallel edges in $\ga$ contains a
Scharlemann bigon.  If the label pair of the Scharlemann bigon is
$\{j, j+1\}$, then these two labels appear twice among the endpoints
of this family, and each of the other four labels appears exactly
once.  By Lemma 2.3(4) at most four labels are the labels of some
Scharlemann bigons in $\gb$, so there is some $k$ which is not a label
of a Scharlemann bigon and hence appears exactly three times among the
endpoints of the positive edges incident to $u_i$.  Dually, this
implies that some negative edge in $\rgb$ contains only one $i$-edge,
which is a contradiction as we have shown above that each negative
edge in $\rgb$ must contain exactly two negative $i$-edges.

The proof for $\nb=8$ is similar.  In this case the numbers of edges
in the three positive families incident to $u_i$ are either $(6, 5,
5)$ or $(6, 6, 4)$.  Using the fact that there are at most four labels
of Scharlemann cycles one can show that in either case some label
appears three times among the endpoints of these edges, which would
lead to a contradiction as above.    
\qed

\begin{lemma} No vertex $u_i$ is incident to at most
one edge in $\rgap$.
\end{lemma}

\proof By Lemma 3.2 there are no interior vertices, hence by Lemma 2.11
either (i) $\rgbp$ has a circle component, or (ii) $\rgbp$ has a
boundary vertex of valence at most 3, or (iii) all vertices of $\rgbp$
are boundary vertices of valence 4.

In case (i) a vertex $v_j$ on the circle component is incident to at
most two positive edges with label $i$ at $v_j$, hence dually there
are at most two negative edges with label $j$ at $u_i$, and hence at
least $\Delta - 2 \geq 2$ positive edges with label $j$ at $u_i$,
which is impossible because $u_i$ is incident to at most one family of
positive edges and by Lemma 2.3(3) such a family contains at most one
edge with label $j$ at $u_i$.

The proof for case (ii) is similar because by Lemma 2.7(3) a valence 3
boundary vertex $v_j$ of $\rgbp$ is incident to at most $2n_a$
positive edges of $\gb$ and hence at most two positive edges with
label $i$ at $v_j$.

In case (iii), since $u_i$ is incident to at most $n_b/2+2 < n_b$
positive edges, there is a label $j$ such that all four edges with
label $j$ at $u_i$ are negative.  Dually $v_j$ has four positive
$i$-edges.  Since it is a boundary vertex, it is incident to at least
$3n_a + 1$ positive edges.  On the other hand, since $v_j$ has valence
4 in $\rgbp$, by Lemma 2.7(2) it has at most $2(n_a + 2) < 3n_a$
positive edges, a contradiction.
\qed

\begin{cor}  Each component of $\rgap$ is contained in an
essential annulus but not a disk on $\hat F_a$.
\end{cor}

\proof By Lemma 2.6 $\rgap$ has at least two components, so if the
result is not true then one can find a disk $D$ on $\hat F_a$ such
that $D \cap \rgap$ is a component $G$ of $\rgap$.  By Lemma 4.4 $G$
is not an arc, so by Lemma 2.9 it has at least three boundary vertices
of valence at most 3.  By Lemma 4.3 these vertices are labels of
Scharlemann cycles in $\gb$, which is a contradiction because by Lemma
2.3(4) $\ga$ contains at most two labels of Scharlemann cycles of each
sign.  \qed

Let $G$ be a component of $\rgap$ contained in the interior of an
essential annulus $A$ on $\hat F_a$.  By Corollary 4.5, $G$ is not
contained in a disk, hence it contains some cycles which are
topologically essential simple closed curves on $\hat F_a$, and all
such cycles are isotopic to the core of $A$.  We call such a cycle an
{\it essential cycle\/} on $G$.  Note that a cycle may have more than
two edges incident to a vertex, but an essential cycle does not.  An
essential cycle $C$ of $G$ is {\it outermost on $A$\/} if all
essential cycles of $G$ lie in one component of $A|C$.  By cutting and
pasting one can see that outermost essential cycles always exist, and
there are at most two of them, which we denote by $C_l$ and $C_r$,
called the {\it leftmost cycle\/} and the {\it rightmost cycle},
respectively.  Let $A_l^l$ and $A_l^r$ be the components of $A|C_l$,
called the {\it left annulus\/} and the {\it right annulus\/} of
$C_l$, respectively, labeled so that $A_l^l$ contains no essential
cycles of $G$ other than $C_l$.  Similarly for $A_r^l$ and $A_r^r$,
where the right annulus $A_r^r$ of $C_r$ is the one that contains no
essential cycles other than $C_r$.

\begin{lemma}  The interiors of $A_l^l$ and $A_r^r$ do not
intersect $G$.
\end{lemma}

\proof Assuming the contrary, let $G'$ be the closure of a component
of $G \cap A_l^l$.  Since $G$ is connected, $G'$ must intersect $C_l$
at some vertex $v$, but it cannot intersect $C_l$ at more than one
vertex, as otherwise the union of an arc in $G'$ and an arc on $C_l$
would be an essential cycle in $A_l^l$ other than $C_l$, contradicting
the definitions of leftmost cycle and its left annulus.  For the same
reason, $G'$ contains no essential cycles, hence it lies on a disk $D$
in $A_l^l$.  By Lemma 4.4 $G$ has no vertex of valence 1, so $G'$ is
not homeomorphic to an arc.  By Lemma 2.9 $G'$ has at least three
boundary vertices of valence at most 3.  Let $v^1$ and $v^2$ be
such vertices other than $v$.  They are boundary vertices of $G$ lying
in the interior of $A_l^l$ with valence at most $3$, and $v^i \neq v$.

By Lemma 4.3, for $i=1,2$ there is a Scharlemann bigon $\{e^i_1,
e^i_2\}$ on $\gb$ with $v^i$ as a label, and by Lemma 2.2(5) $C_i =
e^i_1 \cup e^i_2$ is an essential curve on $\hat F_{a}$ containing
$v^i$.  Since $v^2$ is a boundary vertex of $G'$, it is not a cut
vertex, hence there is an arc $C'$ on $G'$ connecting $v^1$ to $v$
which is disjoint from $v^2$.  Now the union $C_1 \cup C' \cup C_l$
cuts $\hat F_{a}$ into an annulus and a disk $D$ containing $v^2$ in its
interior, so the cycle $C_2$ is also contained in the disk $D$, which
is a contradiction to the fact that $C_2$ is topologically an
essential curve on $\hat F_{a}$.  \qed

Lemma 4.6 shows that $G$ is contained in the region $R$ between $C_l$
and $C_r$.  Since $G$ has no interior vertices, all its vertices are on
$C_l \cup C_r$.  If $C_l$ is disjoint from $C_r$ then $R$ is an
annulus, and if $C_l = C_r$ then $R = C_l = C_r$ is a circle.  In the
generic case we have $C_l \cap C_r = E_1 \cup ... \cup E_k$, where
each $E_i$ is either a vertex or an arc.  The region $R$ is then a
union of these $E_i$ and some disks $D_1, ..., D_k$, such that $\bdd
D_i$ is the union of two arcs, one in each of $C_r$ and $C_l$.  When
$k=1$ and $E_1 = v$ is a vertex, $D_1$ is a disk with a pair of
boundary points identified to the single point $v$.  Note that a
vertex of $G$ is a boundary vertex if and only if it is on $C_l \cup
C_r - C_l \cap C_r$.

\begin{lemma} Let $C = C_l$ or $C_r$.  

(1) If $C$ has a boundary vertex $u_i$ of valence at most $3$ then it
has no other boundary vertex of valence at most $4$.

(2) If $C$ has a boundary vertex $u_i$ of valence $2$ then it has no
other boundary vertex.  \end{lemma}

\proof (1) By Lemma 4.3, $i$ is a label of a Scharlemann bigon in
$\gb$.  On $\ga$ the edges of this Scharlemann bigon form a cycle $C'$
containing $u_i$ and another vertex $u_k$.  By Lemma 2.2(5) $C'$ is
topologically an essential circle on the torus $\hat F_a$.  Since
$u_i$ is a boundary vertex, one can see that $C'$ is topologically
isotopic to $C$.  By Lemma 4.6 applied to $G$ and to the component of
$\rgap$ containing $u_k$, there are no other vertices of $\ga$ between
$C'$ and $C$.  Hence any boundary vertex $u_j \neq u_i$ on $C$ is
incident to at at most one family of parallel negative edges,
connecting it to $u_k$.  The result now follows from Lemma 4.2.

(2) Note that since $u_i$ is a boundary vertex, the edges of any
Scharlemann bigon on $\gb$ with $i$ as a label must connect $u_i$ to
the same vertex $u_k$ on $\ga$, so there are at most $n_{b}$ such
bigons because there are only two edges on $\rga$ connecting $u_i$ to
$u_k$, each representing a family of at most $n_{b}$ edges.  By Lemma
2.7(2) $u_i$ is incident to at most $n_{b}+2$ positive edges, hence
at least $3n_{b} - 2$ negative edges.  If a pair of these edges are
parallel on $\gb$ then they form a Scharlemann bigon.  Hence by the
above we see that there are at most $n_{b}$ pairs of such edges.  It
follows that $\rgb$ has at least $3n_{b} -2 - n_{b} = 2n_{b} - 2$
positive edges.

If $u_j$ is a boundary vertex of $C$ other than $u_i$ then as in the
proof of (1) it is incident to at most one family of negative edges,
so it has at least $3n_{b}$ positive edges.  Since no three of those
are parallel on $\gb$, we see that $\rgb$ has at least $3n_{b}/2$
negative edges, so $\rgb$ would have a total of at least $2n_{b} -
2 + 3n_{b}/2 > 3n_{b}$ edges, contradicting Lemma 2.5.  \qed

Now suppose $C_l \cap C_r \neq \emptyset$, and $C_l \neq C_r$.  Then
the region $R$ between $C_l$ and $C_r$ can be cut along vertices of
$C_l \cap C_r$ to obtain a set of disks, and possibly some arcs.  Let
$D$ be such a disk.  If $C_l \cap C_r$ is a single vertex $v$ then $D$
is obtained by cutting $R$ along $v$, in which case we use $D \cap G$
to denote the graph on $D$ obtained by cutting $G$ along $v$.

\bigskip
\leavevmode

\centerline{\epsfbox{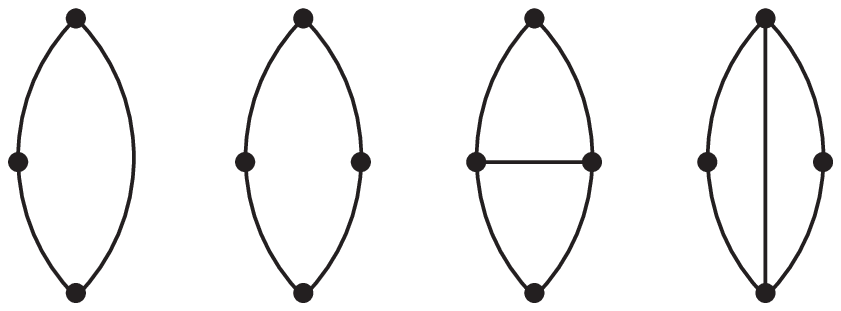}}
\bigskip
\centerline{Figure 4.1}
\bigskip

\begin{lemma} 
$G' = D \cap G$ is one of the four graphs in Figure 4.1.  \end{lemma}

\proof Let $v', v''$ be the vertices of $G'$ lying on both $C_l$ and
$C_r$.  (Note that they are distinct vertices on $G'$ but may be
identified to a single vertex on $G$.)  These vertices divide $\bdd
D$ into two arcs $E_1$ and $E_2$, with $E_1 \subset C_l$ and $E_2
\subset C_r$.

By Lemma 4.7, each $E_j$ contains at most one vertex of valence at most
3 in its interior.  Therefore, if $D$ contains at most one
interior edge then $G'$ has at most four vertices, so it is one of the
four graphs in Figure 4.1.  We need to show that $D$ cannot have more
than one interior edge.

First suppose there is an interior edge $e$ of $G'$ which has both
endpoints on $E_1$.  We may choose $e$ to be outermost in the sense
that there is an arc $E'$ on $E_1$ with $\bdd E' = \bdd e$, and there
is no edge of $G'$ inside the disk bounded by $E' \cup e$.  Since $G'$
has no parallel edges, there must be a vertex $v$ in the interior of
$E'$, which has valence 2.  By Lemma 4.7(2), in this case $C_l$ has no
other boundary vertices, so $E_1$ has no vertex other than $v$ in its
interior; in particular, $e$ must have its endpoints on $v'$ and
$v''$.  This implies that all interior edges have both endpoints on
$E_2$, and by the same argument as above we see that $E_2$ has
exactly one vertex in its interior, and all edges must have endpoints
on $v'$ and $v''$.  Since $G'$ has no parallel edges, it can have at
most one edge connecting $v'$ to $v''$, and we are done.

We can now assume that every interior edge of $G'$ has one endpoint in
the interior of each $E_i$.  Let $G''$ be the union of the interior
edges.  The above implies that $G''$ cannot have a cycle, so it is a
union of several trees with endpoints in the interiors of $E_1$ and
$E_2$.  A vertex of valence 1 in $G''$ is a vertex of valence 3 in
$G'$, and by Lemma 4.7(1) there is at most one such for each $E_i$.
Therefore $G''$ is a chain, with two vertices of valence 1 and $k\geq
0$ vertices of valence 2, so $G'$ has one vertex of valence 3 on each
$E_i$, and $k$ vertices of valence 4.  Note that these are boundary
vertices.  However, by Lemma 4.7(1), if $G$ has a vertex of valence
$3$ on $C_l$ then it has no boundary vertex of valence at most $4$ on
$C_l$, and similarly for $C_r$.  It follows that $k=0$, which again
implies that $G'$ has only one interior edge.  \qed

\bigskip
\leavevmode

\centerline{\epsfbox{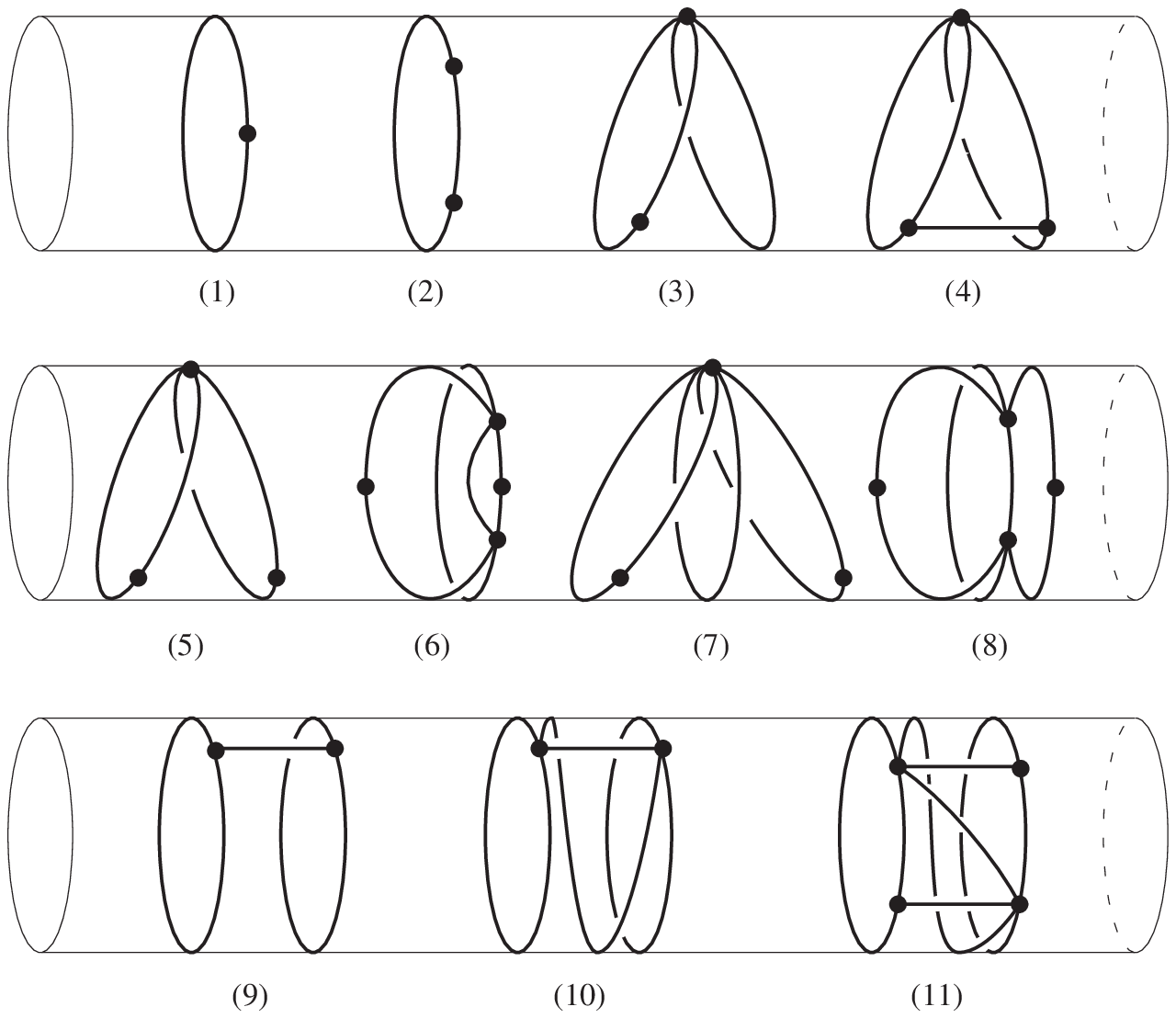}}
\bigskip
\centerline{Figure 4.2}
\bigskip

\begin{lemma}  If $G$ is a component of $\rga$
and $C_l \cap C_r \neq \emptyset$, then $G$ is one of the graphs in
Figure 4.2 (1) -- (8).
\end{lemma}

\proof If $C_l = C_r$ then $G$ is a simple cycle, in which case each
vertex has valence 2 and hence is a label of a Scharlemann bigon by
Lemma 4.3.  By Lemma 2.3(4), $G$ has at most two such vertices, hence
$G$ is the graph in Figure 4.2(1) or (2).

Suppose $C_l \neq C_r$ and $C_l \cap C_r \neq \emptyset$.  We call the
endpoints of $C_l \cap C_r$ {\it breaking points} of $G$, which cut the
region $R$ between $C_l$ and $C_r$ into several disks $D_1, ..., D_k$
and possibly some arcs.  By Lemma 4.8 each $G_i = D_i \cap G$ is one
of the graphs in Figure 4.1.  We say that $G_i$ is of type ($j$) if it
is the graph in Figure 4.1(j).  Since $G$ can have at most two
boundary vertices of valence at most three, we see that either $k=1$,
or $k=2$ and both $G_i$ are of type (1).

First assume that $k=2$ and $G_1, G_2$ are of type (1).  By Lemma 4.7
the two boundary vertices of $G_i$ must be one on each of $C_l,
C_r$.  If the component of $C_l \cap C_r$ containing a breaking point
$v'$ on $G_i$ is an arc instead of a vertex, then $v'$ would be a
vertex of $\rgb$ which is incident to three positive edges, two of
which are adjacent, in which case by Lemma 4.3 $v'$ is a label of a
Scharlemann bigon in $\gb$.  Since $G$ contains no more than two
Scharlemann bigon labels, this cannot happen.  It follows that $G$ is
the graph shown in Figure 4.2(6).

We can now assume $k=1$.  For the same reason as above, we see that if
$G_1$ is of type (1), (2) or (3), then $G$ is as shown in Figure 4.2(3),
(5) or (4), respectively. If $G_1$ is of type (4), the breaking vertices
may be incident to an edge in $C_l \cap C_r$, so $G$ is the graph in
Figure 4.2(7) or (8).  \qed

\begin{lemma} If $G$ is a component of $\rgap$ and $C_l \cap
C_r = \emptyset$, then $G$ is one of the three graphs in Figure 4.2
(9), (10), or (11).  \end{lemma}

\proof Note that in this case all vertices on $C_l$ and $C_r$ are
boundary vertices.  If $C_l$ has a vertex of valence 2 then by Lemma
4.7 it has no other vertices, in which case $C_l$ is a loop and we have
$G = C_l$, so $C_l = C_r$, a contradiction.  Therefore $C_l$ and $C_r$
have no vertices of valence 2, hence all vertices of $G$ have valence
at least 3.

Doubling the annulus and calculating Euler characteristic, we see that
$$\sum (4 - {\text val}(v_i)) \geq 0.$$
By Lemma 4.7 $G$ has at most two vertices of valence 3.  

First assume that $G$ has two vertices $v_1, v_2$ of valence 3.  By
Lemma 4.7 $v_1, v_2$ cannot both be on $C_l$ or $C_r$, hence each of
$C_l$ and $C_r$ contains exactly one vertex of valence 3.  By Lemma
4.7 they cannot contain vertices of valence 4.  By the above formula
$G$ has either (i) no other vertex, or (ii) one other vertex with
valence 6 or 5, or (iii) two other vertices, both having valence 5.
One can check that in Case (i) the graph is that of Figure 4.2(9),
Case (ii) does not happen, and in Case (iii) the graph is the one in
Figure 4.2(11).

If $G$ contains only one vertex $v$ of valence 3, then by the above
formula it contains at least one vertex of valence 5, and all other
vertices are of valence 4.  If $C_l$ contains $v$ then by Lemma 4.7 it
contains no vertices of valence 4.  Since each edge of $G$ must have
one endpoint on each of $C_l$ and $C_r$, we see that $C_l$ must
contain a vertex of valence 5, and $C_r$ contains exactly two
vertices, each of valence 4.  One can check that there is no reduced
graph satisfying these conditions.

Now assume that $G$ has no vertices of valence 3.  Then by the above
formula all vertices of $G$ are of valence 4.  Since $G$ has no
parallel edges, all edges of $G - C_l \cup C_r$ must connect $C_r$ to
$C_l$, so the graph $G$ is completely determined by the number of
vertices $k$ on $C_r$, which must be the same as that on $C_l$.
Denote such a graph by $G_k$.  When $k = 1$, the graph $G = G_1$ is
shown in Figure 4.2(10).  We need to show that $k>1$ does not happen.

Suppose $k>1$.  By Lemma 4.2 each vertex on $C_l$ is incident to at
least two negative edges in $\rga$.  Let $G'$ be the component of
$\rgap$ adjacent to $C_l$, and let $C'_r$ be the outermost cycle of
$G'$ adjacent to $C_l$.  If $C'_r$ has only one vertex then two
negative edges based at some vertex $v_i$ on $C_l$ form an essential
loop on the annulus between $C'_r$ and $C_l$, so there is only one
negative edge of $\rga$ incident to any other vertex on $C_l$, which
is a contradiction.  Similarly if some vertex on $C'_r$ is a boundary
vertex of valence at most 3 then by Lemma 4.3 it is a label of a
Scharlemann bigon, which is again a contradiction because the two
edges of the Scharlemann bigon would form an essential loop as above.
This rules out the possibility of $G'$ being a graph in Figure 4.2 (1)
or (3) -- (11).  If $G'$ is the one in Figure 4.2(2) then by Lemma 4.3
both of its vertices are labels of Scharlemann bigons.  The edges of
these two Scharlemann bigons form two cycles, which cannot be on the
same side of $G'$ as otherwise one of them would lie on a disk, which
contradicts Lemma 2.2(5).  Hence one of the pairs of edges connect a
vertex of $G'$ to a vertex of $C_l$, which is again a contradiction.

It now follows that if some component of $\rgap$ is a $G_k$ for $k\geq
2$, then so are all the other components.  Moreover, none of the
vertices is a label of a Scharlemann cycle as otherwise some vertex
would be incident to a single negative edge in $\rga$, which would
contradict Lemma 4.2.  Hence $\gb$ has no Scharlemann cycles.  On the
other hand, by Lemma 2.7(2) the four families of positive edges at a
vertex $v_i$ of $G_k$ contain at most $2n_{b}+4$ edges, so $v_i$ is
incident to at least $\Delta n_{b} - (2n_{b}+2) > n_{b}$ negative
edges.  By Lemma 2.8 this implies that $\gb$ does have a Scharlemann
cycle, which is a contradiction.  \qed

\begin{cor}  Suppose $\Delta \geq 4$, and $n_a > 4$
for $a=1,2$.  Then

(1) each component of $\rgap$ is one of the 11 graphs in Figure 4.2;
and

(2) each $\ga$ contains a Scharlemann cycle, hence $\hat F_{b}$ is
separating, and $\nb$ is even for $b = 1,2$.

\end{cor}

\proof (1) This follows from Lemmas 4.9 and 4.10.

(2) By (1), $\rgap$ contains either a vertex $v$ of valence 2 or a
boundary vertex of valence at most 4.  In the first case the result
follows from Lemma 4.3.  In the second case by Lemma 2.7(2) $v$ is
incident to at most $2\nb + 4$ positive edges, hence at least $2\nb-4
> \nb$ negative edges, so by Lemma 2.8 $\ga$ has a Scharlemann cycle.
\qed

\section {The case $n_1, n_2 > 4$   }

In this section we will complete the proof that the generic case $n_1,
n_2 > 4$ cannot happen.  We assume throughout the rest of the section
that $n_1, n_2 > 4$.  Let $G$ be a component of $\rgap$.  By Corollary
4.11 $G$ is one of the graphs in Figure 4.2.  We need to rule out all
these possibilities.  Recall that a component of $\rgap$ is of type
(k) if it is the graph in Figure 4.2(k).

Here is a sketch of the proof.  We first show (Lemma 5.4) that $\rgap$
cannot have two boundary vertices of valence 2, hence no component of
$\rgap$ is of type (5)--(8).  Types (3) and (11) will be ruled out in
Lemmas 5.6 and 5.7, so we are left with types (1), (2), (4), (9) and
(10).  Lemma 5.8 will show that each vertex of a type (10) component
is a label of Scharlemann cycle, which implies that all vertices of
$\rgap$ are labels of Scharlemann cycles, except the valence 4 vertex
in a type (4) component.  Since $\rgap$ has at most two Scharlemann
labels of each sign, we see that each $\rgap$ is a union of two type
(4) components.  This will be ruled out in Lemma 5.10, completing the
proof of the theorem.

Each vertex $u_i$ in $\ga$ has $\Delta$ edge endpoints labeled $j$.
Define $\sigma(u_i, v_j)$ to be the number of those on positive edges
minus the number of those on negative edges.  In other words, it is
the sum of the signs of the edges with an endpoint labeled $j$ at
$u_i$.

Define a vertex $u$ of $\rgap$ to be {\it small\/} if it is either of
valence 2 or is a boundary vertex of valence 3.  Note that a component
of type (1) or (3) in Figure 4.2 has one small vertex, a component of
type (10) has no small vertex, and all others have two small vertices.

\begin{lemma} (1) $\sigma(u_i, v_j) = - \sigma(v_j, u_i)$.

(2) If $v_j$ is a small vertex in $\rgbp$ then $\sigma(u_i, v_j) \geq
    0$ for all $i$.
    
(3) If $\rgap$ has a boundary vertex $u_i$ of valence 2, then
$\sigma(u_i, v_j) < 0$ for all but at most two $j$, at most one
for each sign.

(4) If $\rgap$ has a boundary vertex of valence 2, then $\rgbp$ has at
most one small vertex of each sign.
\end{lemma}

\proof
(1) This follows from the parity rule Lemma 2.2(1).  

(2) If $v_j$ has valence 2 in $\rgbp$ then each label $i$ appears at
most twice among the positive edge endpoints.  If $v_j$ is a boundary
vertex of valence 3 in $\rgbp$ then by Lemma 2.7(3) it is incident to
at most $2n_a$ adjacent positive edges in $\gb$, hence again each $i$
appears at most twice among the positive edge endpoints.  Since $\Delta
\geq 4$, the result follows.

(3) If $u_i$ is a boundary vertex of valence 2 then by Lemma 2.7(2)
there are at most $n_b + 2$ adjacent positive edges, so at most two
labels appear more than once among the positive edge endpoints, and if
there are two then they are adjacent, so there is only one for each
sign.

(4) This follows immediately from (2) and (3).
\qed

\begin{lemma} Suppose $\rgap$ has a boundary vertex $u_i$ of
valence 2.  Then all components of $\rgbp$ are of type (1), (3) or
(10).  Moreover, for each sign there is at most one component with
vertices of that sign which is of type (1) or (3).
\end{lemma}

\proof This follows immediately from Lemma 5.1(4) and the fact that
a component of type (1) or (3) has one small vertex, a component of
type (10) has no small vertex, and all others have two small vertices.
\qed

\begin{lemma} Let $v_j$ be a vertex of a type (10) component
$G$ of $\rgbp$.

(1) $v_j$ is incident to at most $2n_a+2$ positive edges in $\gb$.

(2) $\sigma(u_i, v_j) \geq 0$ for all but at most two $u_i$, one for
    each sign.
\end{lemma}

\proof (1) By Lemma 2.7(2) the four families of adjacent parallel positive
edges incident to $v_j$ contain $m \leq 2(n_a+2)$ edges.  If $m >
2n_a+2$, then in particular one of the families contains more than
$n_a/2$ edges, so it contains a Scharlemann bigon.  By Lemma 2.2(4) and
Lemma 2.2(1) the labels at the endpoints of a loop at $v_j$ must have
different parity, which rules out the possibility $m=2n_a+3$.
Hence $m=2n_a+4$.  Note that in this case there are at least 4
parallel loops $\{e_1, ..., e_4\}$, where $e_1$ is the outermost edge
on the annulus containing $G$.  By looking at the labels at the
endpoints of these loops, we see that $e_2, e_3$ form a Scharlemann
bigon, which contradicts Lemma 2.2(6) because $\{e_1, e_4\}$ is then
an extended Scharlemann cycle.  

(2) Since $v_j$ is a boundary vertex of $G$, the positive edges
incident to $v_j$ are adjacent.  Therefore (1) implies that
$\sigma(v_j, u_i) \leq 0$ for all but at most two $i$, hence by Lemma
5.1(1) we have $\sigma(u_i, v_j) \geq 0$ for all but at most two
$u_i$, and if there are two such $u_i$ then they are of opposite sign.
\qed

\begin{lemma} $\rgap$ cannot have two parallel boundary
vertices of valence 2; in particular, no component $G$ of $\rgap$ is
of type (5), (6), (7) or (8).  \end{lemma}

\proof Suppose to the contrary that $\rgap$ has two boundary vertices
$u_{i_1}, u_{i_2}$ of valence 2, and of the same sign.  By Lemma 5.2,
each component $G'$ of $\rgbp$ is of type (1), (3) or (10).  If $G'$
is of type (3) then it has a boundary vertex of valence 2, so
applying Lemma 5.2 to this vertex (with $\rgap$ and $\rgbp$ switched),
we see that $G$ must be of type (1), (3) or (10), which is a
contradiction.  Therefore $G'$ must be of type (1) or (10).

By Lemma 5.1(3), $\sigma(u_{i_1}, v_k) < 0$ for all but at most two
$v_k$.  Similarly for $\sigma(u_{i_2}, v_k)$.  Since $n_b > 4$, there
is a vertex $v'$ such that $\sigma(u_r, v') < 0$ for both $r=i_1,
i_2$.  On the other hand, if $v'$ is on a component $G'$ and if $G'$
is of type (1) then by Lemma 5.1(2) we have $\sigma(u_r, v') = -
\sigma(v', u_r) \geq 0$ for all $u_r$, while if $G'$ is of type (10)
then Lemma 5.3(2) says $\sigma(u_r, v') \geq 0$ for either $r=i_1$ or
$i_2$ because $u_{i_1}$ and $u_{i_2}$ are of the same sign.  This is a
contradiction.  \qed

Note that a vertex $u$ on a component $G$ of $\rgap$ is a boundary
vertex if it lies on one outermost essential cycle $C_1$ of $G$ but
not the other one.  In this case there is a unique component $G'$ of
$\rgap$ and a unique outermost essential cycle $C_2$ on $G'$ such that
$C_1 \cup C_2$ bounds an annulus on $\hat F_a$ whose interior contains
no vertex of $\ga$.  We say that $G'$ and $C_2$ are {\it adjacent\/}
to $u$.

\begin{lemma} Let $u_i$ be a vertex on a type (10) component
$G$ of $\rgap$.  If $u_i$ is not a label of a Scharlemann cycle in
$\gb$, then 

(i) the component $G'$ of $\rgap$ adjacent to $u_i$ is of
type (1), (3) or (10);

(ii) $u_i$ is incident to exactly $2n_{b}-2$ negative edges; and

(iii) $\rgbp$ has only two components, each of type (4) or (11).
\end{lemma}

\proof We assume that $u_i$ is not a label of a Scharlemann cycle.
Let $G'$ and $C$ be the component and outermost cycle adjacent to
$u_i$.  If $C$ has a boundary vertex $u_j$ of valence at most 3,
then by Lemma 4.3 $u_j$ is a label of a Scharlemann cycle.  Since
$u_j$ is a boundary vertex and there is no vertex between $C$ and the
outermost cycle on $G$ containing $u_i$, the edges of the above
Scharlemann cycle must connect $u_j$ to $u_i$, hence $u_i$ is also a
label of the Scharlemann cycle, which is a contradiction.  Also, if
$G'$ is of type (2) then by Lemma 4.1 each of its vertices is a label
of a Scharlemann cycle.  Recall that the edges of a Scharlemann cycle
in $\gb$ cannot lie in a disk on $\hat F_{a}$, hence the edges of one
of the Scharlemann cycles must connect a vertex on $C$ to $u_i$, which
again is a contradiction.  Therefore $C$ does not have a boundary
vertex of valence at most 3, and it is not on a type (2)
component.  Examining the graphs in Figure 4.2, we see that $G'$ must
be of type (1), (3) or (10).  Moreover, if it is of type (3) then $C$
is the loop there.  In any case, $C$ contains only one vertex.

Let $t$ be the number of negative edges incident to $u_i$.  Since $C$
has only one vertex $u_j$, $u_i$ is incident to at most two families
of negative edges $\hat e_1, \hat e_2$, all connecting $u_i$ to $u_j$,
so by Lemma 2.7(1) $t \leq 2\nb$.  On the other hand, by Lemma 5.3
$u_i$ is incident to at most $2\nb + 2$ positive edges, so $t \geq
2\nb -2$.  Therefore we have $2\nb \geq t \geq 2\nb -2$.

First assume $t=2\nb$.  Then each of $\hat e_1$ and $\hat e_2$
contains exactly $\nb$ edges.  Since $u_i$ is not a label of
Scharlemann cycle, by Lemma 2.4 these $2\nb$ edges are mutually
non-parallel on $\gb$, hence $\rgbp$ has at least $2\nb$ edges.  On
the other hand, by Lemma 4.1 it cannot have more than $2\nb$ such
edges, hence $\rgbp$ has exactly $2\nb$ edges, each containing exactly
one edge in $\hat e_1 \cup \hat e_2$.  Counting the number of edges on
each graph in Figure 4.2, we see that each component of $\rgbp$ must
be of type (10) or (11).  Also, a component of type (11) has a vertex
$v_k$ of valence 5 in $\rgbp$, so the above implies that the label $k$
appears 5 times among the endpoints of edges in $\hat e_1 \cup \hat
e_2$, which is absurd.  This rules out the possibility for a component
to be of type (11).  Now notice that these two families of $\nb$
parallel edges have the same transition function, hence if some edge
has the same labels on its two endpoints, then they all do.  It
follows that no component can be of type (10) because it has both loop
and non-loop edges.  This completes the proof for the case
$t=2\nb$.

If $t = 2\nb-1$ then one of $\hat e_1, \hat e_2$ contains $n_{b}$
edges and the other contains $\nb-1$ edges.  Examining the labels at
the endpoints of these edges we see that if an edge in $\hat e_1$ has
labels of the same parity at its two endpoints then an edges in $\hat
e_2$ would have labels of different parities at its endpoints, and
vice versa.  This contradicts the parity rule (Lemma 2.2(1)).

We can now assume $t = 2\nb -2$.  Without loss of generality we may
assume that the labels of the endpoints of $\hat e_1 \cup \hat e_2$
appear as $1,2,...,\nb,1,..., \nb-2$ on $\bdd u_i$ when traveling
clockwise, and we assume that the first $\nb$ are endpoints of $\hat
e_1$.  (The other cases are similar.)  Let $e^k_p$ ($k=1,2$) be the
edge in $\hat e_k$ with label $p$ at $u_i$, and assume that the label
of $e^1_1$ on $u_j$ is $1+r$ for some $r$.  Then one can check that
the label of $e^1_p$ on $u_j$ is $p+r$, and the label of $e^2_p$ on
$u_j$ is $p+r+2$.  (All labels are integers mod $\nb$.)  Hence for any
$p$ between $3$ and $\nb$, the edges $e^1_p$ and $e^2_{p-2}$ have the
same label $p+r$ at $u_j$.  On $\gb$ this implies that there are two
positive edges, connecting $v_p$ to $v_{p+r}$ and $v_{p+r}$ to
$v_{p-2}$, so $v_p$ are $v_{p-2}$ are in the same component of
$\rgbp$.  Since this is true for all $p$ between $3$ and $\nb$, it
follows that $\rgbp$ has only two components.

By Lemmas 2.8 and 2.2(4) $\rgbp$ has the same number of positive
vertices and negative vertices, hence each component $G$ has at least
three vertices.  This rules out the possibility for $G$ to be of type
(1), (2), (3), (9) or (10).  Combined with Lemme 5.4 we see that each
component of $\rgbp$ is of type (4) or (11).  \qed

\begin{lemma} No component of $\rgap$ is of type (3).
\end{lemma}

\proof By Lemma 5.2 if $\rgap$ has a component of type (3) then each
component of $\rgbp$ is of type (1), (3) or (10), and there is at most
one component of type (1) or (3) for each sign.  Since $n_b > 4$ and a
component of type (1) or (3) has at most 2 vertices, there is at least
one component $G$ of $\rgbp$ of type (10) and at least one other
component $G'$ of the same sign.  On the other hand, by Lemma 5.5 each
vertex of $G$ is a label of a Scharlemann cycle, and by Lemmas 4.1 and
4.3 at least one vertex of $G'$ is a label of a Scharlemann cycle, so
there are at least three labels of Scharlemann cycles of the same
sign, contradicting Lemma 2.3(4).  \qed

\begin{lemma}  No component of $\rgap$ is of type (11).
\end{lemma}

\proof An outermost cycle on a component $G$ of type (11) contains two
parallel vertices $u_i$ and $u_j$, where $u_i$ is of valence 3 and
hence the label of a Scharlemann bigon (Lemma 4.3), and $u_j$ has
valence 5.  If $\{e_1, e_2\}$ is a Scharlemann bigon on $\gb$ with
label pair $\{i, i+1\}$, say, then on $\ga$ these edges form an
essential curve containing the vertices $u_i$ and $u_{i+1}$, which
separates $u_j$ from all other vertices of opposite sign, hence all
negative edges incident to $u_j$ have their other endpoints on
$u_{i+1}$, and they are all parallel.  Thus $u_j$ has at most $\nb$
negative edges, and hence at least $3\nb$ adjacent positive edges.  In
particular, each label appears at least three times among endpoints of
positive edges at $u_j$.  Dually, each vertex $v_k$ in $\gb$ is
incident to at least three negative edges labeled $j$ at $v_k$.  If
$v_k$ is a boundary vertex, then this implies that it is incident to
at least $2\na+1$ negative edges, so by Lemma 4.1 it is a label of a
Scharlemann cycle.

By Lemmas 5.4 and 5.6 a component of $\rgbp$ is of type (1), (2), (4),
(9), (10) or (11).  By the above and Lemma 4.1 all vertices of $\gb$
except those with valence 4 in type (4) components are labels of
Scharlemann cycles.  Since $\nb>4$ and there are at most two
Scharlemann labels for each sign, we see that $\rgbp$ has only two
components, each of type (4), so $\nb = 6$, and $\rgbp$ has 10
positive edges.  By Lemma 2.5 $\rgb$ has at most $3\nb - 10 = 8$
negative edges.  On the other hand, we have shown that $u_j$ in $\ga$
is incident to at least $3\nb = 18$ positive edges; since no three of
them are parallel in $\gb$, $\rgb$ has at least $18/2 = 9$ negative
edges, which is a contradiction.  \qed

\begin{lemma}  Each vertex of a type (10) component of $\rgap$
is a label of a Scharlemann bigon.
\end{lemma}

\proof Suppose that a vertex $u_i$ of a type (10) component of $\rgap$
is not a label of a Scharlemann bigon.  By Lemmas 5.5 and 5.7 $\rgbp$
is a union of two type (4) components, so $n_b=6$, $\rgb$ has 10
positive edges, and no more than $3\nb - 10 = 8$ negative edges.

By Lemma 5.5(ii) $u_i$ is incident to $(\Delta-2)n_b+2 = 6\Delta - 10$
positive edges (loops counted twice).  By Lemma 2.7(1) no three of
these are parallel in $\gb$, hence they represent at least $3\Delta -
5$ negative edges in $\rgb$.  Therefore $\Delta = 4$, and we have at
least $7$ negative edges in $\rgb$.  We need to find two more to get a
contradiction.

By Lemma 5.5(i) and Lemma 5.6 the component $G$ of $\rgap$ adjacent to
$u_i$ is of type (1) or (10), so the outermost cycle of $G$ adjacent
to $u_i$ has a single vertex $u_j$ and a single edge $E_0$.  We claim
that $E_0$ contains at least two edges of $\gb$.  

If $G$ is of type (10), then $u_j$ is incident to four families of
positive edges in $\ga$, with a total of $2n_b + 2 = 14$ edges, where
loops are counted twice.  By Lemma 2.3(3) each family contains no more
than 4 edges, so the loop edge $E_0$ contains at least $(14 - 2 \times
4)/2 = 3$ edges of $\gb$.  If $G$ is of type (1) then since no three
negative edges incident to $u_j$ are parallel in $\gb$, and since
$\rgbp$ has only 10 edges, we see that $u_j$ is incident to at most 20
negative edges, hence $E_0$ contains at least $(24-20)/2=2$ edges.
This completes the proof of the above claim.

Let $e'_1, e'_2$ be the two edges in $E_0$ closest to $u_i$.  By Lemma
2.2(2) they are not parallel on $\gb$.  We claim that on $\gb$ neither
of them is parallel to any edge incident to $u_i$, hence $\rgb$
contains at least $7+2=9$ negative edges.  This will be a
contradiction as we have shown above that $\rgb$ has at most 8
negative edges.  

By Lemma 5.5(ii) there are exactly $2n_b - 2 = 10$ negative edges
$e_1, ..., e_{10}$ connecting $u_i$ to $u_j$.  Without loss of
generality we may assume that the sequence of labels of the endpoints
of these edges at $u_i$ is $1,...,6,1,...,4$, counting clockwise, and
the labels of their endpoints at $u_j$ are $r+2, r+3, ..., r-1$,
counting counterclockwise.  Thus $\{e'_1, e'_2\}$ is a Scharlemann
bigon with label pair $\{r, r+1\}$.

Since $e'_i$ is a loop, by Lemma 2.3(5) if it is parallel in $\gb$ to
an edge $e$ incident to $u_i$ then $e$ is also a loop.  Note that
$e'_i$ and $e$ must have the same label pair.  Let $E_3$ be the loop
of $\rga$ based at $u_i$.  It has at most four edges $e''_1, e''_2,
e''_3, e''_4$, with label pairs $\{5,6\}, \{6,5\}, \{1, 4\}, \{2,3\}$,
respectively.  By Lemma 2.3(2) we have $\{r, r+1\} \neq \{2, 3\}$,
hence if $e'_i$ is parallel to some $e''_j$ then $\{r, r+1\} =
\{5,6\}$, so $r = 5$, and hence the label sequence of the above
negative edges at $u_j$ is also $1, ..., 6, 1, ..., 4$.

The 10 edges $e_1, ..., e_{10}$ are divided into two families $E_1,
E_2$.  Since $|E_i| \leq 6$, we have $|E_1| = 4$, 5, or 6.  If
$|E_1|=5$ then the edge $e_1$ would have label $1$ at $u_i$ and label
$6$ at $u_j$.  Since $v_1$ and $v_6$ on $\gb$ are antiparallel, this
is impossible by the parity rule.  If $|E_1| = 4$ then $e_1$ has the
same label $1$ at its two endpoints, which contradicts the fact that
$\gb$ has no loop.  Similarly if $|E_1|=6$ then $e_{7}$ has the same
label $1$ at its two endpoints, which is again a contradiction.  This
completes the proof of the Lemma.  \qed

\begin{lemma} Each $\rgap$ is a union of two type (4)
components.  \end{lemma}

\proof By Lemmas 5.4, 5.6 and 5.7, each component $G$ of $\rgap$ is of
type (1), (2), (4), (9) or (10).  By Lemmas 4.1, 4.3 and 5.8, we see
that all vertices $u_i$ of $G$ are labels of Scharlemann bigons,
unless $G$ is of type (4) and $u_i$ is the vertex of valence 4 in $G$.
Since $\na>4$ and $\rgap$ has at most two vertices which are labels of
Scharlemann bigons for each sign, we see that $\rgap$ consists of exactly two
components, each of type (4).  \qed

\bigskip
\leavevmode

\centerline{\epsfbox{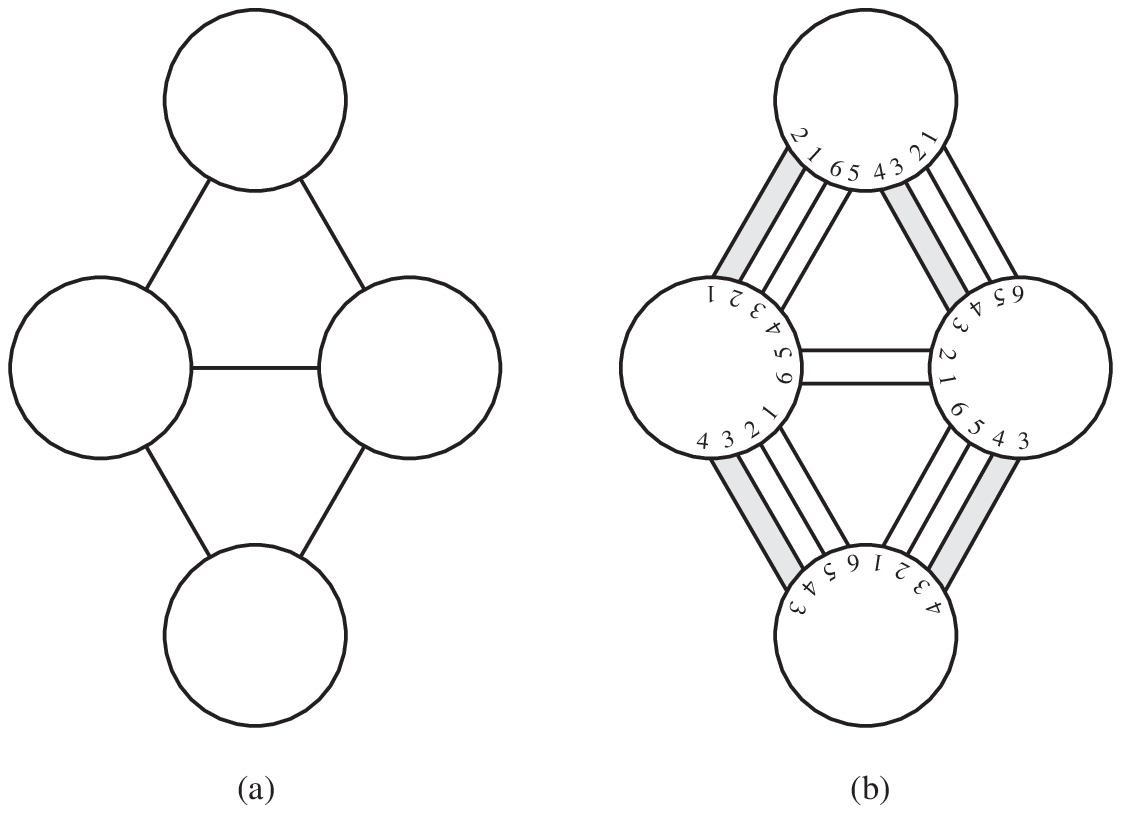}}
\bigskip
\centerline{Figure 5.1}
\bigskip

\begin{lemma} One of the $\rgap$ is not a union of two type (4)
components.  \end{lemma}

\proof Assume that each $\rgap$ is a union of two type (4) components.
Each vertex of $\ga$ has valence $\Delta n_b \geq 24$, hence $\ga$ has
at least 72 edges.  Since a positive edge in $\Gamma_a$ is a negative
edge in $\Gamma_b$, we may assume that $\Gamma_1$ has no more negative
edges than positive edges, so $\Gamma^+_1$ has at least 36 positive
edges.  Thus one component $G$ of $\Gamma^+_1$ has at least 18 edges.
Denote by $\hat G$ the reduced graph of $G$.  It is of type (4), so it
is obtained from the graph in Figure 5.1(a) by identifying the top and
bottom vertices.

Let $E_1, ..., E_5$ be the edges of $\hat G$.  Denote by $|E_i|$ the
number of edges of $G$ in $E_i$, and call it the {\it weight\/} of
$E_i$.  By Lemma 2.3(3), each $|E_i| \leq 4$.  Since $G$ has at least
18 edges, up to relabeling the weights of the edges are at least
$(4,4,4,4,2)$ or $(4,4,4,3,3)$.

Let $D$ be a triangle face of $\hat G$, and let $E_1$, $E_2, E_3$
be the edges of $D$.  We will also use $D$ to denote the
corresponding triangle face in $G$.  If $|E_i| = 4$ then by Lemma 2.4
$E_i$ contains a Scharlemann bigon, which must be at one end of the
family of parallel edges in $E_i$.  We say that the Scharlemann bigon
in $E_i$ is {\it adjacent to $D$\/} if one of its edges is on the
boundary of $D$.

\bigskip
\leavevmode

\centerline{\epsfbox{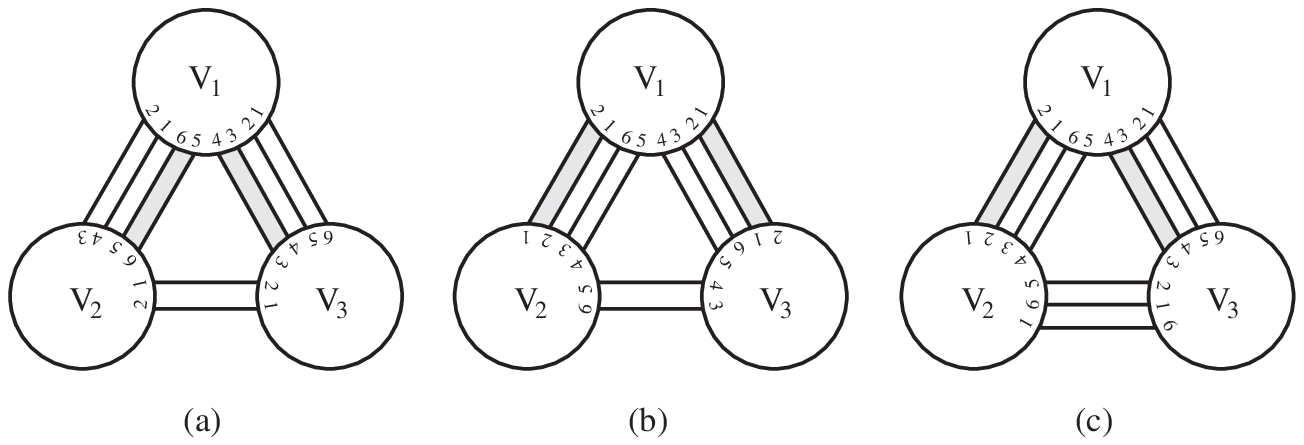}}
\bigskip
\centerline{Figure 5.2}
\bigskip

\noindent
{\bf Sublemma} {\it If $|E_1| = |E_2| = 4$, then (i) $|E_3|=2$, and
(ii) exactly one of $E_1$ and $E_2$ has its Scharlemann bigon adjacent
to $D$. }
\bigskip

\proof Let $V_1$ be the fat vertex incident to both $E_1$ and $E_2$.
Without loss of generality we may assume that the labels on $\bdd V_1$
are as shown in Figure 5.2(a), where $E_1$ is the upper right family of
edges.  Note that the positions of the Scharlemann bigons in $E_1,
E_2$ determine the labels on $\bdd V_2$ and $\bdd V_3$.

If both $E_1$ and $E_2$ have their Scharlemann bigons adjacent to $D$,
then the labels are as shown in Figure 5.2(a), in which case we have
three Scharlemann bigons with disjoint label pairs, contradicting
Lemma 2.3(4).  If both Scharlemann bigons of $E_1, E_2$ are
non-adjacent to $D$, then the labels are as shown in Figure 5.2(b), in
which case the edges adjacent to those of $D$ form an extended
Scharlemann cycle, which contradicts Lemma 2.2(6).  This proves (ii).

We may now assume without loss of generality that the Scharlemann
bigon of $E_1$ is adjacent to $D$ while that of $E_2$ is not
adjacent to $D$.  See Figure 5.2(c).  In this case the label pair
of the Scharlemann bigon in $E_1$ is $\{3,4\}$.  If $|E_3| \geq 3$
then $E_3$ contains a Scharlemann bigon with label pair $\{6, 1\}$.
This contradicts Lemma 2.3(2), completing the proof of the sublemma.
\qed

If the weights of the $E_i$ are $(4,4,4,3,3)$, or if the weights are
$(4,4,4,4,2)$ and the horizontal edge in Figure 5.1(a) has weight 4,
then the boundary edges of one of the triangles in Figure 5.1(a) have
weights $(4,4,3)$ or $(4,4,4)$, which contradicts the sublemma.
Therefore the edges of $G$ are exactly as shown in Figure 5.1(b).  As
in the proof of the sublemma, we may assume that the labels at the
three vertices in the upper triangle of $G$ are as shown in Figure
5.1(b).  The Scharlemann bigons in the upper triangle have label pairs
$\{3,4\}$ and $\{1,2\}$, hence by Lemma 2.3(4) $G$ cannot have a
Scharlemann bigon on label pair $\{5,6\}$.  Therefore the labels of the
endpoints of the lower-right edges must be as shown in Figure 5.1(b).
This determines the labels at the lower vertex.  But then neither
Scharlemann bigon in the lower triangle is adjacent to the triangle,
contradicting the sublemma.  \qed

\begin{prop} The case that both $n_1, n_2 > 4$ is
impossible.  \end{prop}

\proof
This follows from the contradiction between Lemma 5.9 and Lemma 5.10.
\qed

\section {Kleinian graphs }

In Sections 6 -- 11 we will improve Proposition 5.11 to show that $n_i
\leq 2$ for $i=1$ or $2$.  For the most part we will assume that $n_a
= 4$.  In this section we prove some useful lemmas.  In particular,
Lemmas 6.2 -- 6.5 study kleinian graphs.  Lemma 6.2 gives basic
properties of kleinian graphs, which will also be used later in
studying the case $n_a = 2$.

\begin{defn} The graph $\ga$ is said to be {\it
kleinian\/} if $\hat F_a$ bounds a twisted $I$-bundle over the Klein
bottle $N(K)$ such that each component of $N(K) \cap V_a$ is a $D^2
\times I$, and each component of $N(K) \cap F_b$ is a bigon.
\end{defn}

By Lemma 2.12, if $M(r_a)$ contains a Klein bottle $K$ intersecting
$K_a$ at $n_a/2$ points then $\bdd N(K)$ is an essential torus
intersecting $K_a$ at $n_a$ points, hence in this case we may assume
that $\hat F_a = \bdd N(K)$, where $N(K)$ is a small regular
neighborhood of $K$; in particular, $\ga$ is kleinian.  In this case
$N(K)$ is called the {\it black region}, and all faces of $\gb$ lying
in this region are called {\it black faces}, and the others {\it white
faces}.  We assume that the vertices of $\ga$ have been labeled so
that $u_{2i-1} \cup u_{2i}$ lie on the same component of $V_a \cap
N(K)$.  The following lemma lists the main properties of kleinian
graphs.

\begin{lemma} Suppose $\ga$ is kleinian.  Then

(1) each black face of $\gb$ is a bigon;

(2) each family of parallel edges in $\gb$ contains an even number of edges;

(3) $\gb$ has no white Scharlemann disk, hence any Scharlemann cycle
of $\gb$ has label pair $\{k, k+1\}$ with $k$ odd;

(4) there is a free involution of $\hat F_a$, which preserves $\ga$,
sending $u_{2i-1}$ to $u_{2i}$ and preserving the labels of edge
endpoints.  \end{lemma}

\proof (1) follows from the definition.  (2) follows from (1) because
if there is a family containing an odd number of edges then one side
of that family would be adjacent to a black face, which is not a
bigon.

(3) Each edge of a white face is adjacent to a black bigon, so if
there is a white Scharlemann disk then the edges of the Scharlemann
cycle and the adjacent edges would form an extended Scharlemann cycle,
which would be a contradiction to Lemma 2.2(6).

(4) We may assume that the Dehn filling solid torus $V_a$ and the
surface $F_b$ intersect $N(K)$ in $I$-fibers.  Thus the involution of
$\hat F_a$ obtained by mapping each point to the other end of the
$I$-fiber gives rise to the required involution of $\ga$.  \qed

\begin{lemma} Suppose $n_a = 4$.  Then $\ga$ is kleinian if
each vertex of $\ga$ is a label of a Scharlemann bigon in $\gb$.
\end{lemma}

\proof  Without loss of generality we may assume that $\gb$ has a
$(12)$ Scharlemann bigon.  By assumption there is a Scharlemann bigon
with $3$ as a label.  If there is no $(34)$ Scharlemann bigon
then this Scharlemann bigon must have label pair $(23)$.  Similarly
the Scharlemann bigon with $4$ as a label must have label pair
$(14)$.  We may therefore relabel the vertices of $\ga$ so that the
label pairs of the above Scharlemann bigons are $(12)$ and $(34)$
respectively.  

Shrinking the Dehn filling solid torus to its core, the Scharlemann
bigons become M\"obius bands $B_{12}$ and $B_{34}$ in $M(r_a)$.  The
union of these M\"obius bands, together with an annulus on $\hat F_a$,
becomes a Klein bottle which can be perturbed to intersect the core of
the Dehn filling solid torus at $2=n_a/2$ points.  By the convention
after Definition 6.1, $\hat F_a$ should have been chosen so that $\ga$ is
kleinian.  
\qed

\begin{lemma} Suppose $n_a=4$.  Then $\ga$ is kleinian if one
of the following holds.

(1) $\gb$ has a family of 4 parallel positive edges.

(2) $\gb$ is positive.  

(3) $\rgbp$ has a full vertex $v_j$ of valence at most 7.

(4) $\rgbp$ contains 4 adjacent families of positive edges with a
total of at least 12 edges.
\end{lemma}

\proof (1) Each label appears exactly twice among the edge endpoints
of a family of four parallel positive edges, hence by Lemma 2.4 it is a
label of a Scharlemann bigon.

(2) If $\gb$ is positive then every vertex $u_i$ of $\ga$ is incident
to at least $4n_b$ negative edges, two of which must be parallel in
$\gb$ because by Lemma 2.5 $\rgb$ contains at most $3n_b$ edges.
Hence by Lemma 2.4 these two edges form a Scharlemann bigon with $i$ as
a label.  Since this is true for all $i$, $\ga$ is kleinian by Lemma
6.3.

(3) Consider the subgraph $G$ of $\rga$ consisting of negative edges.
Then the signs of the vertices around the boundary of a face of $G$
alternate, hence each face has an even number of edges. Using an Euler
characteristic argument one can show that $G$ contains at most $2n_a =
8$ edges.  By (2) we may assume $\gb$ is not positive, so by Lemma
2.3(1) no 3 $j$-edges are parallel on $\ga$, hence $G$ has exactly 8
negative edges, each containing exactly 2 $j$-edges, with one $j$
label at each ending vertex.  Since each vertex $u_i$ has 4
$j$-labels, we see that $u_i$ is incident to exactly 8 $j$-edges, two
of which must be parallel in $\gb$ because $val(v_j, \rgbp) \leq 7$.  By
Lemma 2.4 they form a Scharlemann bigon with $i$ as one of its labels.

(4) By (1) we may assume that each family contains exactly 3 edges, so
the labels at the endpoints of the middle edge in each family are the
labels of a Scharlemann bigon.  It is easy to see that the 4 endpoints
of the middle edges at the vertex are mutually distinct, hence include
all labels.  
\qed

\bigskip
\leavevmode

\centerline{\epsfbox{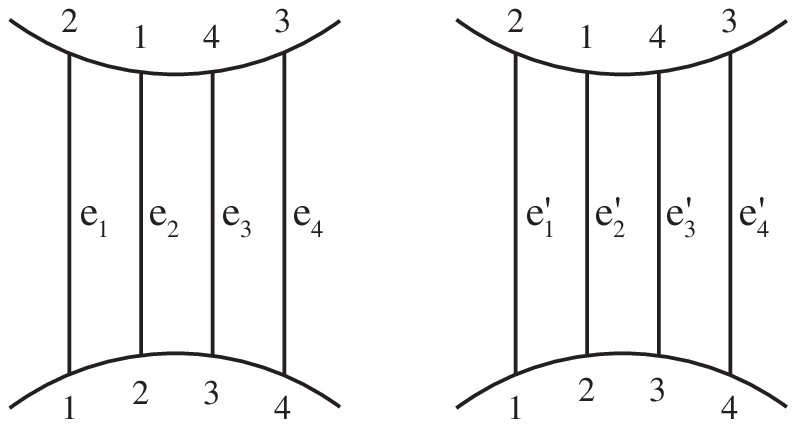}}
\bigskip
\centerline{Figure 6.1}
\bigskip

\begin{lemma} Suppose $n_a = 4$.  Let $e_1\cup e_2 \cup e_3
\cup e_4$ and $e'_1\cup e'_2 \cup e'_3 \cup e'_4$ be two families of
parallel edges in $\gb$ as shown in Figure 6.1.  Then $e_i$ is
parallel to $e'_i$ on $\ga$ for all $i$.  \end{lemma}

\proof Since $e_1\cup e_2$ and $e_3\cup e_4$ form two disjoint
essential cycles on $\hat F_a$ by Lemma 2.2(5), any $(12)$-edge must
be parallel to $e_1$ or $e_2$ and any $(34)$-edge parallel to $e_3$ or
$e_4$ on $\ga$.  Note also that if $e_1$ is parallel to $e'_2$ on
$\ga$ then $e_2$ must be parallel to $e'_1$ (instead of $e'_2$) on
$\ga$ as otherwise $e_1, e_2$ would be parallel on both graphs.
Therefore if the result is not true then either $e_1$ is parallel to
$e'_2$ or $e_4$ is parallel to $e'_3$, so there is a subset $e'_r \cup
... \cup e'_s$ of the second family containing less than $4$ edges,
such that $e'_r \cup e_1$ and $e'_s \cup e_4$ are parallel pairs on
$\ga$.  This contradicts Lemma 2.19
\qed

\begin{lemma} Suppose $n_a=4$ and $\gb$ is non-positive.

(1) No vertex $v_j$ of $\gb$ can have two families of 4 positive edges
with the same label sequence on $\bdd v_j$.  In particular, $v_j$
cannot have two adjacent families of 4 positive edges.

(2) If $\ga$ is kleinian, then two adjacent families of positive edges
of $\gb$ contain at most 6 edges, three contain at most 10, and four
contain at most 12.  

(3) A full vertex of $\rgbp$ has valence at least 6.
\end{lemma}

\proof (1) If there are two families of 4 positive edges with the same
label sequence on $\bdd v_j$ then by Lemma 6.5 the two starting edges
$e_1, e'_1$ of these families will be parallel in $\ga$.  If $e_1,
e'_1$ have label $i$ at $v_j$ then on $\ga$ they have the same label
$j$ at $u_i$, so there are $n_b +1$ parallel negative edges at $u_i$,
and hence by Lemma 2.3(1) $\gb$ would be positive, a contradiction.

(2) By Lemma 6.2(2) the number of edges in each family of positive
edges is either 2 or 4, so by (1) two adjacent families contain a
total of at most 6 edges.  The other two cases follow from this.

(3) Otherwise by Lemma 6.4(3) $\gb$ is kleinian, so the weight of each
positive family of $\gb$ is either 2 or 4.  If some full vertex $v_i$
has valence $5$ or less in $\rgbp$ then it has two adjacent edge of
weight 4, contradicting (1).
\qed

A bigon is called a {\it non-Scharlemann bigon\/} if it is not a
Scharlemann bigon.

\begin{lemma} Suppose $n_a = 4$ and $\ga$ is kleinian.  

(1) Exactly one edge on the boundary of a triangle face of $\rgbp$
represents a non-Scharlemann bigon.  Each of the other two represents
either a Scharlemann bigon or a union of two Scharlemann bigons.

(2) If some vertex $v_i$ is incident to two edges of weight 4 in
$\rgbp$ then any other edge of $\rgbp$ incident to $v_i$ represents a
non-Scharlemann bigon.
\end{lemma}

\proof (1) Let $\hat e_1, \hat e_2, \hat e_3$ be the edges of a
triangle face $\delta$ of $\rgbp$.  By Lemma 6.2(2) each edge of
$\rgbp$ represents 2 or 4 edges.  From the labeling of the edges
around $\delta$ one can see that there are exactly one or three $\hat
e_i$ which are neither a Scharlemann bigon nor a union of two
Scharlemann bigons.  If there are three then they form an extended
Scharlemann cycle, which is impossible by Lemma 2.2(6).  Hence there
must be exactly one such $\hat e_i$.

(2) Otherwise $v_i$ would be incident to 5 Scharlemann bigons, three
of which have the same label pair, say $\{1,2\}$.  Then on $\ga$ there
are six $i$-edges connecting $u_1$ to $u_2$, which form at most two
families because there is a Scharlemann cocycle containing $u_3,
u_4$.  It follows there there are three $i$ labels at the endpoints of
a family, so it contains more than $n_b$ edges, contradicting Lemma
2.3(1).  
\qed

Suppose $\Delta = 4$.  Then a label $j$ is a {\it jumping label\/} at
$u_i$ if the signs of the four $j$-edges incident to $u_i$ alternate.

\begin{lemma} Suppose $\Delta = 4$.  Then a label $i$ is a
jumping label at $v_j$ if and only if $j$ is a jumping label at $u_i$.
In particular, if $v_j$ is a boundary vertex of $\rgbp$ then $j$ is
not a jumping label at any $u_i$.  \end{lemma}

\proof This follows from the Jumping Lemma 2.18.  Let $x_1, ...,
x_4$ be the four points of $u_i \cap v_j$.  Since $\Delta = 4$, the
jumping number must be $\pm 1$.  Therefore they appear in this order
on both $\bdd u_i$ and $\bdd v_j$, appropriately oriented.  If $j$
is a jumping label at $u_i$ then we may assume $x_1, x_3$ are positive
edge endpoints and $x_2, x_4$ are negative edge endpoints on $\bdd
u_i$, which by the parity rule implies that $x_1, x_3$ are negative
edge endpoints and $x_2, x_4$ are positive edge endpoints on $\bdd
v_j$, hence $i$ is a jumping label at $v_j$.  \qed

\begin{lemma}  Suppose $n_a=4$, $n_b\geq 4$, and $\ga$ is
non-positive.  Then $\hat F_a$ is separating.  In particular, $u_1$ is
parallel to $u_3$ and antiparallel to $u_2$ and $u_4$.
\end{lemma}

\proof The result follows from Lemmas 2.8 and 2.2(4) if some vertex
$u_i$ is incident to more than $n_b$ negative edges.  In particular,
since each family of positive edges contain no more than $n_b$ edges,
the result is true if $val(u_i, \rgap) \leq 2$ for some $i$.  Hence we
may assume that $val(u_i, \rgap) > 2$ for all $i$.  One can check that
in this case each component of $\rgap$ must be as shown in Figure
4.2(9), (10) or (11).  In each case $\rgap$ has a boundary vertex
$u_i$ of valence at most 4, so if $n_b > 4$ then by Lemma 2.7(b) the 4
families of positive edges contain at most $2(n_b+2) < 3n_b$ edges,
hence $u_i$ is incident to more than $n_b$ negative edges and the
result follows.  Similarly if 
$\rgap$ has a boundary vertex of valence at most $3$ then by Lemma
2.7(c) it is incident to at most $2n_b$ positive edges in $\ga$ and the
result follows.    Therefore we may assume $n_b = 4$.

A vertex $u_i$ on a component in Figure 4.2 (9) or (10) is a boundary
vertex of valence 3 or 4 in $\rgap$, so by Lemma 2.7(b)--(c) it is
incident to less than $3n_b$ positive edges, and hence more than $n_b$
negative edges, unless $n_b = 4$ and $val(u_i, \rgap) = 4$.  In
particular $\rgap$ must be of type (10) in Figure 4.2.  In this last
case by Lemma 6.4(4) $\gb$ is kleinian, so by Lemma 6.2(2) each family
of positive edges of $\ga$ contains either 2 or 4 edges.  Since there
is a total of at least 12 edges and by Lemma 6.6(2) two adjacent
families contain at most 6 edges, the weights of the four edges of
$\rgap$ incident to $u_i$ must be $(4,2,4,2)$ successively.  However
since the first and the last belong to a loop in $\rgap$, their
weights must be the same, which is a contradiction.  \qed

\section {If $n_a=4$, $n_b \geq 4$ and $\rgap$ has a small
component then $\ga$ is kleinian.
}

A component of $\rgap$ is {\it small \/} if it has at most two edges;
otherwise it is {\it large\/}.  In this section we will show that if
$n_a=4$, $n_b \geq 4$ and $\rgap$ has a small component then $\ga$ is
kleinian.  It is easy to see that the assumption implies that either
$val(u_1, \rgap) \leq 1$, or $val(u_1,\rgap)=val(u_3, \rgap) = 2$ up
to relabeling.  (See the proof of Proposition 7.6.)  The two cases are
handled in Lemmas 7.3 and 7.5, respectively.

\begin{lemma} Suppose $\ga$ contains a loop edge at $u_3$.
Then $\gb$ cannot contain both $(12)$- and $(14)$-Scharlemann bigons.
\end{lemma}

\proof The loop $e$ at $u_3$ must be essential, otherwise it would
bound some disk containing some vertex and hence one of the
Scharlemann cocycles in its interior, which contradicts Lemma 2.2(5).
Now the $(12)$- and $(14)$-Scharlemann bigons in $\gb$ form two
essential cycles in $\ga$ disjoint from $e$, so they must be isotopic
on $\hat F_a$, bounding a disk face containing no vertices of $\ga$ in
its interior.  This is a contradiction to Lemma 2.13.  \qed

\begin{lemma} Suppose $n_a = 4$ and $n_b \geq 4$.  If $val(u_1,
\rgap) \leq 1$ and $\rgbp$ has a boundary vertex $v_j$ of valence at
most 3, then $\ga$ is kleinian.  \end{lemma}

\proof By Lemma 6.9 $u_1$ is parallel to $u_3$ and antiparallel to
$u_2, u_4$.  Since $u_1$ is incident to at most 1 family of positive
edges, it is incident to at least three negative $j$-edges at $u_1$,
so $v_j$ has at least three positive edge endpoints labeled 1.  Hence
$v_j$ being a boundary vertex implies that it has at least 9 positive
edges.  If $v_j$ is incident to 10 or more positive edges of $\gb$
then it has a family of 4 parallel positive edges and hence $\ga$ is
kleinian.  Therefore we may assume that it has exactly 9 positive
edges, divided into three families of parallel edges, each family
containing exactly three edges.  See Figure 7.1.

\bigskip
\leavevmode

\centerline{\epsfbox{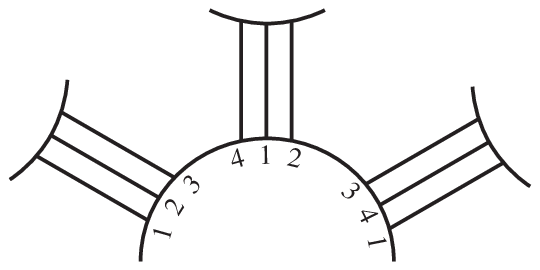}}
\bigskip
\centerline{Figure 7.1}
\bigskip

Since $n_a = 4$, each of these families contains a Scharlemann bigon,
so the labels of the middle edge in the family are labels of a
Scharlemann bigon.  It follows that $1,2,4$ are labels of Scharlemann
bigons.  Thus if the result is not true then $3$ is not a label of
Scharlemann bigon, and $\gb$ contains both $(12)$- and
$(14)$-Scharlemann bigons.  

There are 7 adjacent negative edges at $v_j$, so three of
then have labels $1$ or $3$ at $v_j$.  These cannot all be parallel in
$\ga$ as otherwise there would be three $j$-edges in a family and
hence the family would contain more than $n_b$ edges, contradicting
Lemma 2.3(1) and the fact that $\gb$ is not positive.  On $\rgap$ this
implies that there are at least two edges with endpoints on $\{u_1,
u_3\}$, hence $val(u_1, \rgap) \leq 1$ implies that there is a loop
$\hat e$ based at $u_3$.  Since $\gb$ contains both $(12)$- and
$(14)$-Scharlemann bigons, this is a contradiction to Lemma 7.1. 
\qed

\begin{lemma} Suppose $n_a = 4$ and $n_b \geq 4$.  If $val(u_i,
\rgap) \leq 1$ for some $i$ then $\ga$ is kleinian.  \end{lemma}

\proof If $\gb$ is positive then $\ga$ is kleinian by Lemma 6.4(2).
Therefore we may assume that $\gb$ is non-positive.  By Lemmas 2.3(3)
and 2.7(1) each family of parallel edges in $\ga$ contains at most
$n_b$ edges.  Also, notice that since $u_i$ is incident to more than
$n_b$ negative edges, by Lemmas 2.8 and 2.2(4) the surface $\hat F_a$
is separating, hence $u_i$ is parallel to $u_j$ if and only if $i$ and
$j$ have the same parity.

Without loss of generality we may assume that $val(u_1, \rgap)\leq 1$.
Assume $\ga$ is not kleinian.  Then by Lemma 7.2 $\rgbp$ has no
boundary vertex of valence $3$, and by Lemma 6.4(3) it has no interior
vertex of valence at most 7.  Also, each vertex $v_j$ of $\rgbp$ has
valence at least 3 because it is incident to at least three positive
edges with label $1$ at $v_j$, which by Lemma 2.3(3) must be mutually
non-parallel.  Therefore by Lemma 2.11 all vertices of $\rgbp$ are
boundary vertices of valence 4.

If $n_b > 4$ then by Lemma 2.3(3) the family of positive edges at
$u_1$ contains at most $n_b/2 + 2 < n_b$ edges, so some $v_j$ is
incident to 4 positive edges with label $1$ at $v_j$, which implies
that $v_j$ has at least 13 positive edges in four families, so one of
the families contains 4 edges and hence $\ga$ is kleinian by Lemma
6.4(1).  Similarly if $\Delta > 4$ then $\ga$ is kleinian.

\bigskip
\leavevmode

\centerline{\epsfbox{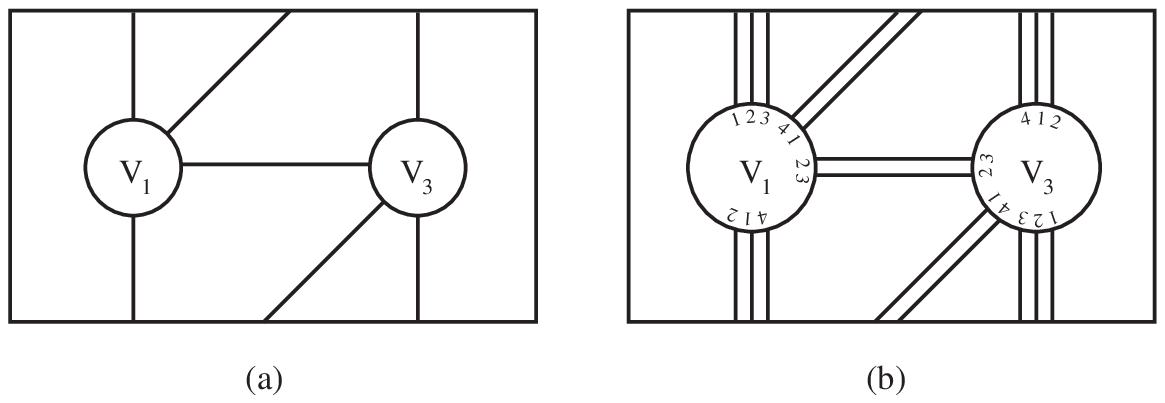}}
\bigskip
\centerline{Figure 7.2}
\bigskip

Now suppose $\Delta = n_b = 4$.  Then $val(v_j, \rgbp)=4$ for all $j$
implies that each component of $\rgbp$ has two loops and two non-loop
edges, as shown in Figure 7.2(a).  By the parity rule a loop based at
$v_j$ has labels of different parity on its two endpoints, hence one
sees that the number of positive edge endpoints of $\gb$ at each $v_j$
is even.  By Lemma 6.4(4) we may assume that $v_j$ has less than 12
positive edges, hence the above implies that each $v_j$ is incident to
exactly 10 positive edges.  If some $v_j$ is incident to only one loop
in $\gb$ then each of the non-loop family incident to $v_j$ contains 4
edges and we are done.  If some $v_j$ is incident to two parallel
loops in $\gb$ then they form a Scharlemann bigon with label pair
$\{1,2\}$, say.  Each of the two non-loop families contains three edges,
hence the middle edge endpoint is a label of a Scharlemann bigon.
Examining the labeling we see that all labels are Scharlemann bigon
labels.

We now assume that each $v_j$ is incident to three parallel loop
edges.  See Figure 7.2(b).  The two outermost loops form a Scharlemann
bigon with $1$ as one of its labels.  There are 6 adjacent negative
edges at $v_1$, so three of then have labels $1$ or $3$ at $v_1$.  By
the same argument as in the last paragraph of the proof of Lemma 7.2
we may assume that the two Scharlemann bigons at $v_1$ and $v_3$ have
the same label pair $(12)$.  The labeling of edge endpoints around
$v_1$ and $v_3$ in a component of $\gbp$ is now as shown in Figure
7.2(b).

Because of the parity rule, the 4 non-loop edges cannot be divided
into a family of 1 and another family of 3 edges, so they must form
two pairs of parallel edges.  From the labeling in Figure 7.2(b) one
can see that they form two Scharlemann bigons with label pairs
$\{2,3\}$ and $\{4,1\}$, respectively.  The result now follows from
Lemma 6.3.
\qed

We now assume that $val(u_1, \rgap) = val(u_3, \rgap) = 2$.  Then
$\rgap$ contains either a cycle $C$ containing both $u_1, u_3$, or
it has two cycle components $C, C'$ containing $u_1, u_3$,
respectively.

\begin{lemma}  If $val(u_i, \rgap) = 2$ then $i$ is a label of
a Scharlemann bigon in $\gb$. 
\end{lemma}

\proof Let $k$ be the number of interior vertices in $\rgbp$.  Let $m$
be the number of edges in $\rgbp$.  We claim that $m_1 \leq 2n_b+k$.
Formally adding edges to $\rgbp$ if necessary we may assume that any face
$A$ between two adjacent components of $\rgbp$ is an annulus.  It is
easy to see that if $\bdd A$ contains $p$ vertices then we can add $p$
edges to make each face on $A$ a triangle.  Therefore we can add at
least $n_b - k$ edges to $\rgbp$ to create a graph $G$ on the torus
$\hat F_b$ whose faces are all triangles.  By an Euler characteristic
argument we see that $G$ has $3n_b$ edges, hence $\rgbp$ has at most
$3n_b - (n_b - k) = 2n_b + k$ edges, and the claim follows.

Now let $m'$ be the number of negative edges of $\ga$ incident to
$u_i$.  Note that if $m < m'$ then two negative edges at $u_i$ are
parallel in $\gb$ and we are done.  By Lemma 2.3(3) each positive
family $\hat e$ in $\ga$ contains at most $(n/2) + 2$ edges.
Moreover, if $k>0$ then some label does not appear on endpoints of
edges in $\hat e$, so $\hat e$ has at most $n_b/2$ edges.  Since $u_i$
is incident to two families of positive edges, we have $m' \geq 4n_b -
2(n_b/2) = 3n_b$ if $k>0$, and $m' \geq 4n_b - 2(n_b/2 + 2) = 3n_b -
4$ if $k=0$.  Since $m\leq 2n_b + k$, we have $m < m'$ (and hence $i$
is a label of a Scharlemann bigon in $\gb$), unless $k=0$, $n_b = 4$
and $m = m' = 8$.  

In this last case ($k=0$, $n_b = 4$ and $m = m' = 8$), all vertices of
$\rgbp$ are boundary vertices, hence by Lemma 6.8 there is no jumping
label at $u_i$.  On the other hand, since $m' = 8$, each positive
family at $u_i$ have 4 edges, so the two positive families cannot be
adjacent by Lemma 6.6(2); hence there is a label $j$ such that the
two negative edges labeled $j$ at $u_i$ are separated by the two
positive edges labeled $j$ at $u_i$, so $j$ is a jumping label at
$u_i$, which is a contradiction.  \qed

\begin{lemma} If $val(u_1, \rgap) = val(u_3, \rgap)=2$ then
$\ga$ is kleinian.  \end{lemma}

\proof By Lemma 7.4 $u_1, u_3$ are labels of Scharlemann bigons.  If
some vertex, say $u_4$, is not a label of Scharlemann bigon then there
must be $(12)$- and $(23)$-Scharlemann bigons in $\gb$.  By Lemma 7.4
we have $val(u_4, \rgap) > 2$, so there is a loop edge $e$ of
$\rgap$ based at $u_4$.  This is a contradiction to Lemma 7.1 (with
labels permuted).
\qed

\begin{prop} If $\rgap$ has a small component then (1)
$\ga$ is kleinian, and (2) $\rgap$ has at most 4 edges.  \end{prop}

\proof Let $G$ be a small component of $\rgap$.  If $G$ contains only
one vertex $u_1$ and two edges then it cuts the torus into a disk
containing the other three vertices.  It is easy to see that in this
case there is a vertex of valence at most 2 in $\rgap$, which by Lemma
2.3(3) is incident to at most $2n_b$ edges, hence at least $2n_b$
negative edges.  By Lemma 2.8  $\gb$ has a Scharlemann cycle, so the
surface $\hat F_a$ is separating.  Therefore $u_3$ is parallel to
$u_1$ and is antiparallel to $u_2$ and $u_4$.  It follows that $u_3$
is incident to no positive edges, so by Lemma 7.3 $\ga$ is kleinian.
If $G$ is not as above then either it contains a vertex of valence at
most 1, or it is a cycle, in which case (1) follows from Lemmas 7.3
and 7.5.

Since $\ga$ is kleinian, by Lemma 6.2(4) there is a free involution of
$\ga$ sending $u_i$ to $u_{i+1}$, hence the number of edges ending at
$\{u_2, u_4\}$ is the same as the number of edges ending at $\{u_1,
u_3\}$, which is at most two in all cases discussed above.  Hence (2)
follows.  \qed

\section {If $n_a=4$, $n_b \geq 4$ and $\gb$ is non-positive then
$\rgap$ has no small component 
}

Denote by $X$ the union of $\rgbp$ and all its disk faces.

\begin{lemma} Suppose $n_a=4$, $n_b \geq 4$, $\gb$ is
non-positive, and $\rgap$ has a small component.  Then

(1) each vertex of $\gb$ is incident to at most 8 negative edges;

(2) if $v_j$ is incident to more than 4 negative edges then $j$ is a
label of a Scharlemann bigon;

(3) if $v_j$ is a boundary vertex of valence 3 in $\rgbp$ then it is
incident to either 6 or 8 negative edges, and $j$ is a label of a
Scharlemann bigon;

(4) $val(v_j, \rgbp) \geq 3$ if $v_j$ is a boundary vertex, and $\geq
2$ otherwise;
    
(5) each component of $X$ is either (a) a cyclic union of disks
    and (possibly) arcs, or (b) a cycle, or (c) an annulus.
\end{lemma}

\proof
Since $\rgap$ has a small component, by Proposition 7.6 $\ga$ is
kleinian, and $\rgap$ has at most 4 edges.

(1) If $v_j$ is incident to 9 negative edges then three of them are
parallel on $\ga$ because $\rgap$ has at most four edges, which
contradicts Lemma 2.3(1).

(2) If $v_j$ is incident to 5 negative edges then two of them form a
Scharlemann bigon in $\ga$ because $\rgap$ has only four edges by
Proposition 7.6.

(3) Since $\ga$ is kleinian, by Lemma 6.2(2) $v_j$ is incident to an
even number of negative edges.  Each family of positive edges contains
at most four edges, and by Lemma 6.6 two adjacent families contain at
most 6 edges, hence the three positive families at $v_j$ contain at
most 10 edges.  The result now follows from (1) and (2).

(4) By (1) $v_j$ is incident to at least 8 positive edges, which are
divided into at least two families, and if two then they cannot be
adjacent by Lemma 6.6(1).

(5) If a component of $X$ is contained in a disk then by Lemma 2.9 it
would have either a boundary vertex of valence at most 2, which is
impossible by (4), or six boundary vertices of valence 3, which is a
contradiction because by (3) each such vertex is a label of
Scharlemann bigon while by Lemma 2.3(4) $\rgbp$ has at most two labels
of Scharlemann bigons for each sign.  Therefore no component of $X$ is
contained in a disk on the torus $\hat F_b$.  Since $\gb$ is not
positive, this implies that each component of $X$ is contained in an
annulus but not a disk on $\hat F_b$.

If there is a sub-disk $D$ of $X$ such that $D \cap \overline{X-D}$ is a
single point $v$ then either $\rgbp \cap D$ contains a boundary vertex
of valence 2 other than $v$, or 3 boundary vertices of valence 3 other
than $v$, which again leads to a contradiction as above.  \qed

Let $X_1$ be a component of $X$, and let $v_1$ be a boundary vertex on
the left cycle $C_l$ of $X_1$, as defined in Section 4.  Then there
is another component $X_2$ of $X$ such that the annulus $A$ between
$C_l$ and the right cycle $C'_r$ of $X_2$ has interior disjoint from
$\rgbp$.  Denote by $m_j$ the number of negative edges incident to
$v_j$, and by $m'= m'_1$ the number of negative edges on $A$ which are
not incident to $v_1$.

\begin{lemma} Suppose $n_a=4$, $n_b \geq 4$, $\gb$ is
non-positive, and $\rgap$ has a small component.  Let $v_1$ be a
boundary vertex of $X_1$ with $m_1 > 4$.  Then

(1) $v_1$ is a label of a Scharlemann bigon;

(2) $m' = 0$ if $m_1 = 8$; 

(3) $m' \leq 2$ if $m_1 > 4$; 

(4) $C_l$ contains no other boundary vertices of valence at most 4.
\end{lemma}

\proof (1) Since $\rgap$ has only four edges, two of the negative
edges at $v_1$ form a Scharlemann bigon on $\ga$.

(2) If $m_1 = 8$ then since $\rgap$ has only 4 edges, the 8 negative
edges at $v_1$ form 4 Scharlemann cocycles, which must all go to the
same vertex $v_2$ on $C'_r$ because $v_1$ is a boundary vertex and the
cocycles are essential loops.  These cocycles separate $C_l$ from
$C'_r$, hence all negative edges in $A$ incident to a vertex of $C_l -
v_1$ must have the other endpoint on $v_2$.  On the other hand, by
Lemma 8.1(1) $v_2$ is incident to at most 8 negative edges, and by the
above all of them must connect $v_2$ to $v_1$.  Hence $m' = 0$.

(3) By Proposition 7.6 $\ga$ is kleinian, so by Lemma 6.2(2) $m_1$ is
even; hence by (2) we may assume that $m_1 = 6$.  Since $\rgap$ has
only 4 edges, the 6 negative edges incident to $v_1$ contain at least
2 Scharlemann cocycles, which connect $v_1$ to some $v_2$ on $C'_r$.
If $v_2$ is incident to $8$ negative edges in $A$ then as in (2) these
edges form 4 Scharlemann cocycles, which must all connect to the same
vertex $v_1$ and hence $m' = 0$.  By Lemma 6.2(2) each family of
parallel edges in $\gb$ has an even number of edges, so $v_2$ cannot
be incident to 7 negative edges in $A$.  If $v_2$ is incident to 6 or
less negative edges in $A$ then by the above 4 of them connect to
$v_1$, so there are at most 2 connecting to $C_l - v_1$, hence $m'
\leq 2$.

(4) A boundary vertex on $C_l - v_1$ of valence at most $4$ in $\rgbp$
is incident to at most 12 positive edges by Lemma 6.6(2), and hence at
least 4 negative edges, which must lie in $A$ because it is a boundary
vertex on $C_l$.  This is contradicts (3).  \qed

\begin{lemma} Suppose $n_a=4$, $n_b \geq 4$, $\gb$ is
non-positive, and $\rgap$ has a small component.  If a component $X_1$
of $X$ contains a boundary vertex $v_1$ of valence 3, then $X_1$ is an
annulus containing exactly two vertices, both of which are of valence
3 and are labels of Scharlemann bigons.  \end{lemma}

\proof By Lemma 6.6(2) $v_1$ is incident to at least 6 negative edges.
Consider the three possible types of $X_1$ in Lemma 8.1(5).  It cannot
be a cycle because it has a boundary vertex $v_1$.  If $X_1$ is an
annulus or a cyclic union of disks and arcs then by Lemma 8.2(4),
$C_l-v_1$ has no boundary vertex of valence at most 4 in $\rgbp$,
which implies that there is a boundary vertex $v_3$ of valence 3 on
the right circle $C_r$ of $X_1$, hence for the same reason $C_r - v_3$
contains no boundary vertex of valence at most 4.  By Lemma 8.2(3)
there are at most 4 negative edges incident to $C_l \cup C_r - \{v_1,
v_3\}$, so there is no (non-boundary) vertex of valence 2 on $X_1$.
Thus either 

(i) $X_1$ is an annulus containing only the two vertices $v_1$ and
$v_3$; or

(ii) $X_1$ is an annulus containing exactly four vertices and the
other two are boundary vertices of valence 5; or

(iii) $X_1$ is as in Figure 8.1 (a) or (b).  

Case (i) gives the conclusion of the lemma because, as in the proof of
Lemma 8.2, a boundary vertex of valence 3 must be a label of a
Scharlemann bigon.  We need to show that (ii) and (iii) are
impossible.

\bigskip
\leavevmode

\centerline{\epsfbox{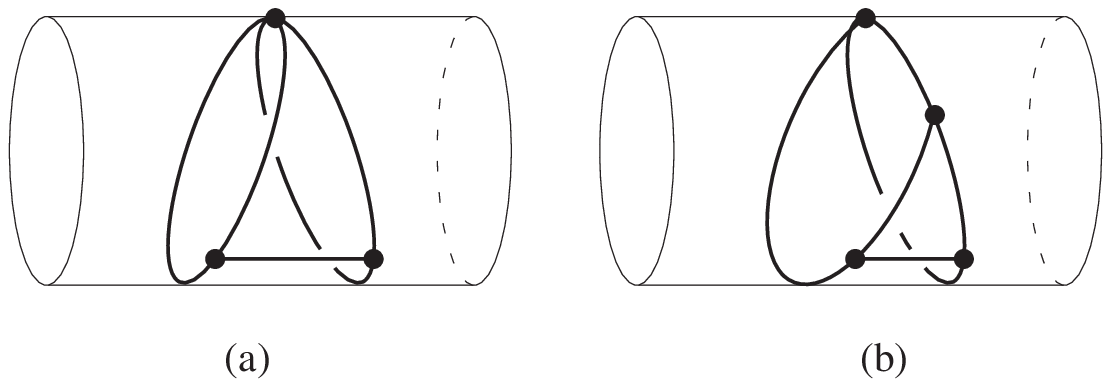}}
\bigskip
\centerline{Figure 8.1}
\bigskip

In case (ii), let $v_5, v_7$ be the boundary vertices of valance 5 on
$X_1$, with $v_5 \subset C_l$.  Note that all faces of $X_1$ are
triangles.  If $v_1$ had 8 negative edges then by Lemma 8.2(2) $v_5$
would have no negative edge, which is impossible by Lemma 6.4(3).
Therefore $v_1$ has exactly 6 negative edges and 10 positive edges, so
by Lemma 6.6(1) the weights of the edges of $X_1$ incident to $v_1$
must be $(4,2,4)$.  Now the middle edge of $X_1$ at $v_5$ has weight 2
and the two boundary edges have weight $4$, so again by Lemma 6.6(1)
the weights around $v_5$ must be $(4,2,2,2,4)$.  By Lemma 6.7(2) the
edges of weight 2 must be non-Scharlemann.  This is a contradiction
because a triangle with a corner at $v_5$ bounded by two weight 2
edges has two non-Scharlemann bigons on its boundary while by Lemma
6.7(1) it has only one.

In case (iii), we assume $X_1$ is as in Figure 8.1(a).  The other case
is similar.  Let $v_1, v_3$ be the vertices of valence $3$ in the
figure, and let $v_5$ be the other vertex.  If $m_1 = 8$ then by Lemma
8.2(2)$v_5$ is incident to no negative edges on the side of $C_l$, and at
most 2 negative edges on the side of $C_r$.  Therefore the four
positive edges of $\rgbp$ at $v_5$ are adjacent to each other,
representing a total of $14$ edges.  It follows that there are two
adjacent families of positive edges, each containing 4 edges, which is
a contradiction to Lemma 6.6(1).  If $m_1 = m_3 = 6$ then each of
$v_1$ and $v_3$ is incident to 10 positive edges, so by Lemma 6.6(1)
the weights of the edges of $\rgbp$ at $v_1$ and $v_3$ are $(4,2,4)$,
in which case $v_5$ again has two adjacent families of 4 positive
edges each, contradicting Lemma 6.6(1).  \qed

\begin{lemma} Suppose $n_a=4$, $n_b \geq 4$, $\gb$ is
non-positive, and $\rgap$ has a small component.  If a component $X_1$
of $X$ does not contain a boundary vertex of valence 3, then 

(1) $X_1$ is either a cycle or an annulus containing exactly two
vertices; and

(2) all vertices of $X_1$ are labels of Scharlemann bigons.
\end{lemma}

\proof By Lemma 6.6(3) all interior vertices of $X$ have valence at
least 6.  Since $X_1$ has no boundary vertex of valence 3, by Lemma
2.11 it is either a cycle, or an annulus with all interior vertices of
valence 6, all boundary vertices of valence 4, and all faces
triangles.  The result is true when $X_1$ is a cycle because any
vertex of valence 2 has more than 4 negative edges and hence is a
label of Scharlemann cycle.  Therefore we assume that $X_1$ is an
annulus.  If $X_1$ has an interior vertex $v_1$ then by Lemmas 6.6(4)
the weights of the edges of $X_1$ around $v_1$ must be $(4,2,4,2,2,2)$
and any edge of weight 2 represents a non-Scharlemann bigon.  Thus the
triangle with a corner at $v_1$ bounded by two weight 2 edges has the
property that it has at least two edges representing non-Scharlemann
bigons, which is a contradiction to Lemma 6.7(1).  Therefore $X_1$ has
no interior vertex.

First assume that some vertex $v_1$ on the left cycle $C_l$ of $X_1$
is incident to less than 12 positive edges, and hence more than 4
negative edges.  Since all vertices of $X_1$ are of valence 4 in
$\rgbp$, by Lemma 8.2(3) $C_l$ has no other vertices.  In this case
$C_r$ also contains only a single vertex $v_3$.  By Lemma 6.6(2) $v_3$
cannot have more than 12 positive edges, and if 12 then the weights
around it are $(4,2,4,2)$.  However this cannot happen because the
first and the last numbers are the weights of the loop edge at $v_3$
and hence must be the same.  Therefore $v_3$ has less than 12 positive
edges.  Now by Lemma 8.1(2) $v_3$ and $v_1$ are labels of Scharlemann
bigons, and the result follows.

Now assume that each vertex of $X_1$ is incident to 12 positive edges.
This implies that the weights of the edges of $X_1$ around each vertex
are either $(4,2,4,2)$ or $(2,4,2,4)$, so any pair of adjacent edges
have different weights.  However, this is impossible because two of
the three edges of a triangle face must have the same weight.
\qed

\begin{lemma} Suppose $n_a=4$, $n_b \geq 4$, $\gb$ is
non-positive, and $\rgap$ has a small component.  Then

(1) $n_b = 4$.

(2) $val(u_i, \rgap) \geq 2$ for all $i$.

(3) Each component of $\rgap$ is a loop.
\end{lemma}

\proof (1) By Lemmas 8.3 and 8.4 each vertex of $\gb$ is a label of a
Scharlemann bigon, and by Lemma 2.3(4) there are at most 4 such
labels.  Hence $n_b = 4$.

(2) Suppose $u_1$ is incident to at most one edge of $\rgap$.  By
Lemmas 8.3 and 8.4 each component $G$ of $\rgbp$ consists of either
(i) a cycle, or (ii) two loops and one non-loop edge, or (iii) two
loops and two non-loop edges.  Since there are at least 3 negative
$j$-labels at $u_1$, there are at least 3 positive $1$-labels at each
$v_j$, hence (i) cannot happen.  Moreover, since each $v_j$ is a
boundary vertex containing at least three positive $1$-labels, it has
more than 8 positive edges.

Suppose $G$ is of type (ii).  Then the label $1$ appears three times
on positive edge endpoints around each of the two vertices of $G$,
hence it appears a total of $6$ times among the three families of
positive edges in $G$, so by Lemmas 2.4 and 6.2(3) there is a
$(12)$-Scharlemann bigon among each of these families.  Since a loop
and a non-loop edge cannot be parallel in $\ga$ (Lemma 2.3(5)), these
represents at least four edges of $\rga$ connecting $u_1$ to $u_2$,
which cut the torus $\hat F_a$ into two disks.  On the other hand, by
the parity rule a loop at a vertex $v$ of $G$ must have labels of
different parity on its two endpoints, so the total number of positive
edges at $v$ is at least 10, divided into three families, hence one of
the families has four parallel edges, which contains a
$(34)$-Scharlemann bigon, giving a pair of edges on $\ga$ lying in the
interior of the disks above which must therefore be parallel.  This is
a contradiction to the fact that a Scharlemann cocycle is essential
(Lemma 2.2(5)).

Now suppose $G$ is of type (iii).  If some $v_j$ is incident to $12$
positive edges then by Lemma 6.6(1) the weights of the positive edges
around $v_j$ are $(4,2,4,2)$, which is impossible because the first
and the last weights are for a loop and hence must be the same.  Since
$v_j$ is incident to more than $8$ positive edges, it is incident to
exactly 10 positive edges, so the weights are $(2,4,2,2)$ or
$(2,2,4,2)$ around each vertex.  The two loops must be
$(12)$-Scharlemann bigons in order for each vertex to have 3 edge
endpoints labeled $1$.  This completely determines the labeling of the
edge endpoints up to symmetry.  Examining the labeling one can see
that the family of 4 parallel edges form an extended Scharlemann
cycle, which is a contradiction to Lemma 2.2(6).

(3) By Proposition 7.6 $\ga$ is kleinian, hence the torus $\hat F_a$
is separating, so each edge of $\rgap$ has endpoints on vertices whose
subscripts have the same parity.  By (2) a small component of $\rgap$
must be a loop $C$.  Let $u_1$ be a vertex of $C$.  If $C$ does not
contain $u_3$ then the component of $\rgap$ containing $u_3$ must also
be a loop because it contains no other vertices, and it cannot contain
more than one edge as otherwise some component of $\rgap$ would lie in
a disk and hence would have a vertex of valence at most 1,
contradicting (2).  Thus the graph $G$ consisting of $u_1, u_3$ and
all edges with endpoints on them is either one loop or two disjoint
loops.  By Lemma 6.2(4) there is a involution of $\rgap$ mapping $u_1$
and $u_3$ to $u_2$ and $u_4$ respectively, hence it maps $G$ to $\rgap
- G$.  Therefore the components in $\rgap - G$ are also loops.
\qed

\begin{prop} Suppose $n_a=4$, $n_b \geq 4$, and $\gb$ is
non-positive.  Then $\rgap$ has no small component.  \end{prop}

\proof Suppose to the contrary that $\rgap$ has a small component.  By
Lemma 8.5 we have $n_b = 4$ and each component of $\rgap$ is a loop.
Thus each $u_i$ is incident to at most $8$ positive edges, so $\ga$
has no more positive edges than negative edges.

By Lemmas 8.3 and 8.4, each component $X_1$ of $X$ is either a circle
or an annulus containing two vertices of $\rgbp$.  First assume that
$X_1$ is a circle, so it is a small component of $\rgbp$.  Applying
Proposition 7.6 and Lemma 8.5 with $\ga$ and $\gb$ reversed, we see
that $\gb$ is kleinian, and all components of $\rgbp$ are also
circles, hence $\gb$ also has the property that it has no more
positive edges than negative edges.  Applying the parity rule we see
that both graphs have the same number of positive edges and negative
edges.  In particular, each family of positive edges contains exactly
4 edges, which by Lemma 6.2(3) must consist of a (12)-Scharlemann bigon
and a (34)-Scharlemann bigon.  Dually it implies that all negative
edges connect $u_1$ to $u_2$ or $u_3$ to $u_4$, so there are 4
families of negative edges, each containing exactly 4 edges.  Now the
two positive families at $u_1$ contain 4 edges each, and, whether
separated by the two negative families or not, their endpoints at
$u_1$ have the same label sequence.  This contradicts Lemma 6.6(1).

We may now assume that $X$ consists of two annular components $X_1,
X_2$, each containing two vertices.  Assume $v_1, v_3 \in X_1$.  As in
the last paragraph of the proof of Lemma 8.5(2), in this case each
vertex of $\gb$ is incident to 8 or 10 positive edges, therefore by
Lemma 8.1(2) each vertex is a label of Scharlemann bigon, hence $\gb$
is also kleinian.

If $val(v_1, \gbp) = 10$ then $v_1$ is incident to 6 negative edges,
which are divided into two families of parallel edges on the annulus
bounded by the loops at $v_1$ and $v_2$.  Since $\gb$ is kleinian,
each family contains an even number of edges, hence the number of
edges in these two families are $4$ and $2$, respectively.  Examining
the labels at the endpoints of these edges, we see that two edges with
the same label $i$ at $v_1$ have different labels at $v_2$.  On
$\ga$ this means that there are both loop and non-loop positive edges
incident to $u_i$, which is a contradiction to the fact that $\rgap$
consists of cycles only.

We have shown that $val(v_j, \gbp) = 8$ for all $j$.  Thus there are
16 positive edges on $\ga$, so each of the two positive families
incident to $u_1$ contains 4 edges.  Since all vertices of $\gb$ are
boundary vertices, there is no jumping label on any vertex of $\gb$,
hence by Lemma 6.8 there is no jumping label at $u_1$, so the two
families of positive edges must be adjacent.  This is a contradiction
to Lemma 6.6(2).  \qed

\section {If $\gb$ is non-positive and $n_a=4$ then $n_b \leq 4$
}

Note that if $\ga$ is positive then each vertex of $\gb$ is a label of
a Scharlemann bigon and hence $n_b \leq 4$ by Lemma 2.3(4).  By
Proposition 8.6 the statement in the title is true if $\rgap$ has a small
component.  Therefore we may assume that $\rgap$ consists of two large
components $G_1$ and $G_2$, each of which must be one of the graphs of
type (3), (9) or (10) in Figure 4.2.  As before, denote by $X_a$ the
union of $\rgap$ and all its disk faces, and by $X_i$ the components
of $X_a$ containing $G_i$, $i=1,2$.  Denote $n = n_b$.

\begin{lemma} Suppose that $\gb$ is non-positive, $n_a=4$ and
$n > 4$.  Then $\rgbp$ contains no interior vertex.
\end{lemma}

\proof Otherwise $\gb$ has a vertex $v_i$ which is incident to
positive edges only.  By Lemma 2.3(1) no three of these edges are
parallel on $\ga$, so $\rga$ contains at least $\Delta n_a / 2 \geq 8$
negative edges, and hence at most $3n_a - 8 = 4$ positive edges by
Lemma 2.5, so $\rgap$ has a small component, contradicting our
assumption.  \qed

\begin{lemma} Suppose that $\gb$ is non-positive, $n_a=4$ and
$n > 4$.  Suppose $X_a$ is a disjoint union of two annuli.  Let $G$
be the subgraph of $\gb$ consisting of positive $1$-edges and all
vertices.  Then $G$ cannot have two triangle faces $D_1, D_2$ with an
edge in common.  \end{lemma}

\proof Since $X_a$ is a disjoint union of two annuli, all negative
edges of $\ga$ incident to $u_1$ must have the other endpoint on the
same vertex, say $u_2$, and vice versa.  On $G$ this means that every
edge has label pair $(12)$, and all positive edges with an endpoint
labeled $1$ or $2$ are in $G$.  Thus no edge in the interior of $D_i$
has label $1$ or $2$ at any of its endpoints.  Up to symmetry the
labels on the boundary of the two triangles must be as shown in Figure
9.1.  Since the labels 3 and 4 must appear between two label 1 at a
vertex, one of the triangles, say $D_1$, must contain some $(34)$
edges.  Since all the vertices are parallel, one can see that the
labels $3$ and $4$ appear at each corner of $D_1$, hence there are
three edges inside of $D_1$.  Since there is no trivial loop, they
must form an extended Scharlemann cycle, contradicting Lemma 2.2(6).
\qed

\bigskip
\leavevmode

\centerline{\epsfbox{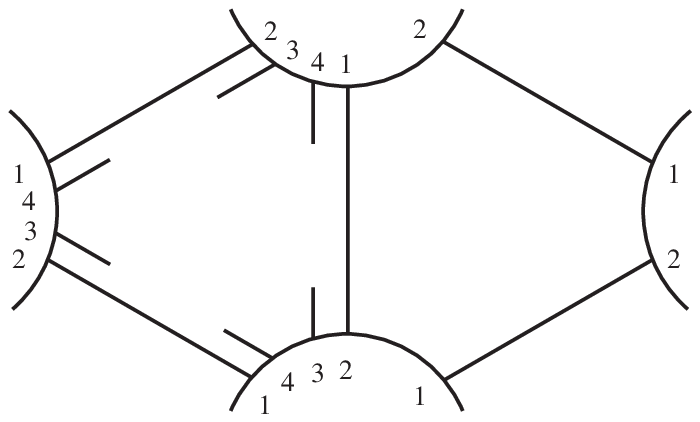}}
\bigskip
\centerline{Figure 9.1}
\bigskip

\begin{lemma} Suppose that $\gb$ is non-positive, $n_a=4$ and
$n > 4$.  If $X_a$ is a disjoint union of two annuli then $\ga$ is
kleinian.  \end{lemma}

\proof By Lemma 6.3 it suffices to show that each vertex $u_i$ of
$\ga$ is a label of a Scharlemann bigon.  Let $t_i$ be the number of
negative edges at $u_i$.  Since there are at most two families of
negative edges incident to $u_i$, we have $t_i \leq 2n$ by Lemma
2.3(1).  On the other hand, since $u_i$ is a boundary vertex of
valence at most 4 in $\rgap$, by Lemma 2.7 it is incident to at most
$2(n+2) = 2n + 4$ positive edges.  If there are $2n+4$ then $u_i$ is
incident to at least four loops and the four outermost loops form an
extended Scharlemann cycle, which is impossible by Lemma 2.2(6).  Also
by the parity rule the loops have labels of different parity at its
two endpoints, hence the number of positive edges at $u_i$ must be
even.  It follows that $u_i$ is incident to either $2n$ or $2n + 2$
positive edges, hence $t_i = 2n$ or $2n-2$.

Assume that $u_1$ is not a label of a Scharlemann bigon.  Then the
negative edges at $u_1$ are mutually non-parallel positive edges in
$\gb$, hence $\rgbp$ has at least $2n - 2$ edges.  Denote by $G$ the
subgraph of $\gb$ consisting of positive $1$-edges.  Let $Y$ be the
union of $G$ and its disk faces, and let $k$ be the number of boundary
vertices of $G$.

First assume $t_1 = 2n$.  By Lemma 9.1 $G$ has no interior vertex, and
clearly it has no isolated vertex, so we can apply Lemma 2.10(1) to
conclude that $k \geq t_1 - n = n$.  Since $G$ only has $n$ vertices,
we must have $k=n$, so all vertices of $G$ are boundary vertices, and
hence there is no cut vertex.  In this case $Y$ contains exactly $n$
boundary edges, so it has at least one (actually $n$) interior edge
$e$.  Since equality holds for the above inequality, by Lemma 2.10(1)
all faces of $Y$ are triangles.  Therefore $e$ is the common edge of
two adjacent triangle faces, which is a contradiction to Lemma 9.2.
Therefore this case is impossible.

Now assume $t_1 = 2n-2$.  In this case the two outermost loops at
$u_1$ form a Scharlemann bigon, so by Lemma 2.2(4) $\hat F_b$ is
separating, hence two vertices of $\gb$ are parallel if and only if
they have the same parity.  Therefore we can define $G_1$ (resp.\
$G_2$) to be the union of the components of $G$ containing $v_i$ with
odd (resp.\ even) $i$.  Similarly for $Y_1$ and $Y_2$.  Then $G_1$
contains all the negative edges at $u_1$ with odd labels, and $G_2$
those with even labels.  Therefore each $G_i$ contains exactly $n - 1$
edges.

The $2n + 2$ positive edges at $u_1$ form at least $n+1$ negative
edges in $\rgb$ because any family of $\gb$ contains at most 2 such
edges.  Hence $\rgb$ contains at least $(n+1) + (2n - 2) = 3n - 1$
edges.  Since a reduced graph on a torus contains at most $3n$ edges
(Lemma 2.5), we may add at most one edge to make the faces of the
graph all triangles.  Hence $\rgb$ has at most one 4-gon and all other
faces are triangles.  In particular, one of the $G_i$, say $G_1$, has
the property that all its faces are triangles.

Let $V$ and $E$ be the number of vertices and edges of $G_1$, and let
$E_b, V_b$ be the number of non-interior edges and boundary vertices,
respectively.  Note that $V-V_b$ is the number of cut vertices, and
$E-E_b$ is the number of interior edges.  We have shown that $V=n/2$
and $E=n-1$.

By Lemma 2.10(1) we have $V_b \geq E-V = (n-1) - n/2$, hence
$G$ has $V - V_b \leq V -(E-V) = (n/2) - (n-1-n/2) = 1$ cut vertex.
If there is no cut vertex then the number of non-interior edges is the
same as the number of vertices $V$, i.e.\ $E_b = V$.  If it has a cut
vertex $v$ then the equality $V_b = E - V$ holds, so by Lemma 2.10(1),
$v$ has exactly two corners not on disk faces, which implies $v$ is
incident to at most 4 non-interior edges, while every other vertex is
incident to exactly two non-interior edges, hence $E_b \leq V+1$.  In
either case we have $E - E_b \geq E - (V+1) = (n-1)-(n/2+1) \geq 1$,
so $G$ has at least one interior edge $e$.  Since all faces of $Y_1$
are triangles, $e$ is incident to two triangle faces of $G_1$, which
is a contradiction to Lemma 9.2.  \qed

\begin{lemma} Suppose that $\gb$ is non-positive, $n_a=4$ and
$n > 4$.  Then $X_a$ is not a disjoint union of two annuli.
\end{lemma}

\proof
Assume to the contrary that $X_a$ is a union of two annuli.  Let $t_i$
be the number of negative edges incident to $u_i$.  As in the proof of
Lemma 9.3, we have $t_i = 2n-2$ or $2n$.

First assume that $t_1 = 2n-2$.  Let $\hat e_1, \hat e_2$ be the two
families of edges in $\ga$ connecting $u_1$ to $u_2$.  Note that a
$(12)$-Scharlemann bigon in $\gb$ must have one edge in each of $\hat
e_1$ and $\hat e_2$.  By Lemma 9.3 $\ga$ is kleinian, so all
$(12)$-edges belong to Scharlemann bigons in $\gb$, hence each edge
$e_i$ in $\hat e_1$ is parallel in $\gb$ to an edge $e_i'$ in $\hat
e_2$, and the label of $e_i$ at $u_1$ is the same as that of $e_i'$ at
$u_2$.  In particular, $\hat e_1$ and $\hat e_2$ have the same number
of edges, hence each contains exactly $n-1$ edges.  Without loss of
generality we may assume that the label $n$ does not appear at the
endpoints at $u_1$ of edges in $\hat e_1$.  By the above, $n$ does not
appear at the endpoints at $u_2$ of edges of $\hat e_2$, hence the
labels must be as in Figure 9.2.  However, in this case the edge
labeled $1$ at $u_1$ has its other endpoint labeled $n$, which is
a contradiction to the parity rule.

\bigskip
\leavevmode

\centerline{\epsfbox{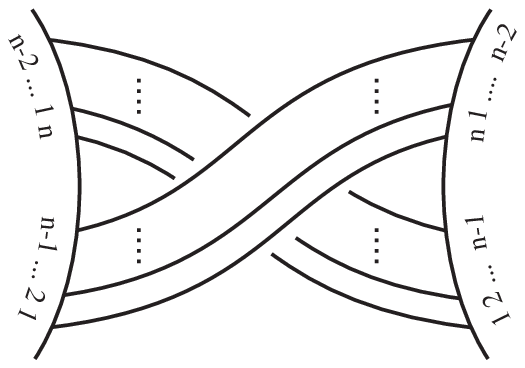}}
\bigskip
\centerline{Figure 9.2}
\bigskip

We now assume $t_1 = 2n$.  Then the two families of negative edges
from $u_1$ to $u_2$ have the same transition function $\varphi$.
Since there is a $(12)$-Scharlemann bigon, $\varphi^2 = id$, so the
length of each $\varphi$-cycle is $1$ or $2$.  Since $n>4$, it follows
that the edges of $\hat e_1$ form at least 3 cycles on $\gb$, which is
a contradiction to Lemma 2.14(2).
\qed

\begin{lemma} Suppose that $\gb$ is non-positive, $n_a=4$ and
$n > 4$.  If $G_1$ is of type (3), then $G_2$ is also of type (3),
and the two loops of $\rgap$ do not separate the two vertices $u_3,
u_4$ which are not on the loops.  \end{lemma}

\proof Let $u_3$ be the vertex of $G_1$ which has valence 2.  If $G_2$
is not of type (3), or if the two loops $\hat e_1 \cup \hat e_2$ of
$\rgap$ separates $u_3$ and $u_4$ then all negative edges incident to
$u_3$ have their other endpoint on the same vertex $u_2$ of $G_2$, and
there are only two such families.  Hence $u_3$ is incident to only
four families of parallel edges, so by Lemmas 2.3(1) and 2.3(3) we
must have $n = 4$ and each family contains exactly four edges.
Since the two families of positive edges are adjacent, this
contradicts Lemma 6.6(1).  \qed

\begin{lemma} Suppose that $\gb$ is non-positive, $n_a=4$ and
$n > 4$.  Then $\rgap$ cannot be a union of two type (3) components.
\end{lemma}

\proof Suppose that $\rgap$ is a union of two type (3) components, so
$\rga$ has 6 positive edges $\hat e_1, ..., \hat e_6$.  By Lemma 9.5
the two loops do not separate the two vertices which are not on the
loops, so the edges appear as in Figure 9.3.  By Lemmas 2.3(1), 2.3(3)
and 6.6(1) each vertex has valence at least 5, so there is one edge
$\hat e_{7}$ from $u_2$ to $u_3$, two edges $\hat e_8, \hat e_9$ from
$u_1$ to $u_2$, and one edge $\hat e_{10}$ from $u_1$ to $u_4$.  There
are one or two edges $\hat e_{11}$ and $\hat e_{12}$ connecting $u_3$
to $u_4$.  See Figure 9.3.

\bigskip
\leavevmode

\centerline{\epsfbox{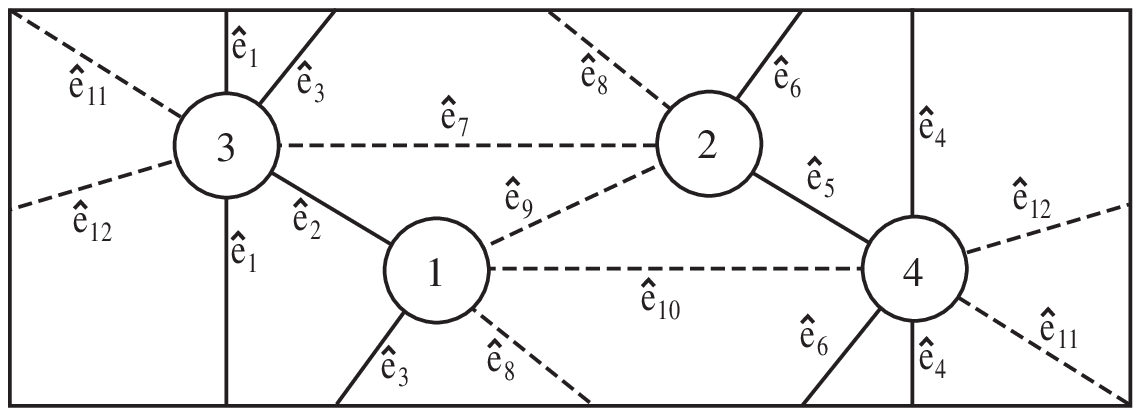}}
\bigskip
\centerline{Figure 9.3}
\bigskip

Denote by $w_i$ the number of edges in $\hat e_i$.  By Lemma 2.7(2) the
two positive edges $\hat e_2, \hat e_3$ contain at most $n + 2$
edges, and the three negative families at $u_1$ contains at most $3n$
edges, hence $\Delta = 4$.  Since each family contains at most $n$
edges, we have the following inequalities.
\begin{eqnarray}
& n \leq w_2 + w_3 \leq n+2  \\
& 3n-2 \leq w_8 + w_9 + w_{10} \leq 3n  \\
& n-2 \leq w_i \leq n  \qquad \text{for $i=7,8,9,10$} 
\end{eqnarray}
Since $u_1$ has at least $n$ adjacent positive edge endpoints and $n$
adjacent negative edge endpoints, each vertex of $\gb$ is incident to
a positive edge and a negative edge, hence $\rgbp$ has no isolated or
interior vertex.

\medskip

Claim 1.  {\it $\rgbp$ has at least $2n -2$ edges, hence $\ga$
has at most two jumping labels.} 

Let $e \in \hat e_8$, and $e' \in \hat e_{10}$.  Then $e$ and $e'$ are
not parallel in $\gb$ because if they were then they would form a
Scharlemann bigon and hence have the same label pair on their
endpoints, which is not the case because their label pairs are $\{1,
2\}$ and $\{1,4\}$ respectively.  Therefore the number of edges in
$\rgbp$ is at least 
$$w_8 + w_{10} \geq 3n-2 - w_9 \geq 2n - 2$$

Since $\rgbp$ has no interior or isolated vertex, by Lemma 2.10(1) it
has at most $2n - (2n - 2) = 2$ non-boundary vertices.  Since
these are the only vertices containing jumping labels, by Lemma 2.18
they are the only possible jumping labels of $\ga$.

\medskip

Claim 2.  {\it $w_{11} + w_{12} \geq 2n - 2$.}

Note that since $w_9, w_{10} \leq n$, we have
\begin{eqnarray*}
8n & = val(u_1,\ga) + val(u_2,\ga) \\
& = (w_2+w_3+ w_{10}) + (w_5+w_6+w_7) + 2 (w_8 + w_9)\\
& \leq (w_2+w_3+w_7) + (w_5+w_6+w_{10}) + 4n 
\end{eqnarray*}
Thus either $w_2 + w_3 + w_7 \geq 2n$ or $w_5 + w_6 + w_{10} \geq
2n$.  Because of symmetry we may assume 
\begin{eqnarray}
w_2 + w_3 + w_7 \geq 2n 
\end{eqnarray}

Divide the edge endpoints on $\bdd u_3$ into $P_1, P_2, P_3, P_4$, as
shown in Figure 9.4.  Denote by $k_i$ the number of edge endpoints in
$P_i$.  A label that appears twice in one of the $P_i$ will be called
a {\it repeated label.\/} Note that if $P_i$ contains a repeated label
then $k_i > n$. Note also that a non-jumping label is a repeated
label.  Thus by Claim 1 there are at least $n-2$ repeated labels among
all the $P_i$.  Since $k_1 = w_7 \leq n$, there is no repeated label
in $P_1$.

\bigskip
\leavevmode

\centerline{\epsfbox{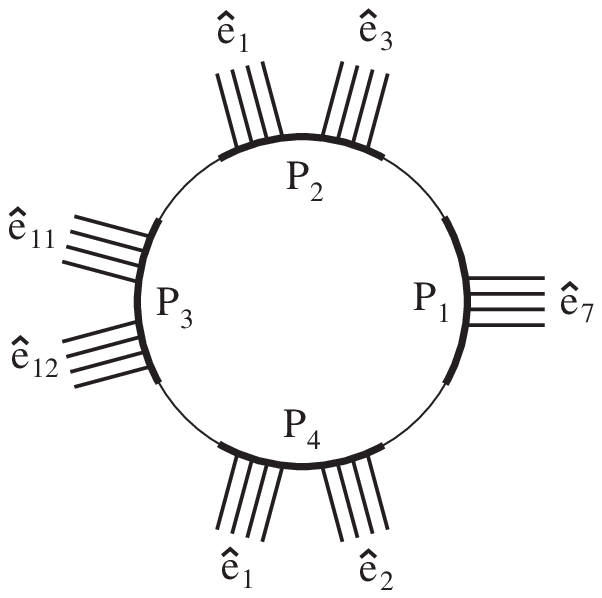}}
\bigskip
\centerline{Figure 9.4}
\bigskip

First assume that $k_3 > n$.  If both $k_2, k_4 \leq n$ then all
repeated labels are in $P_3$, so $k_3 \geq 2n -2$ and we are done.
Because of symmetry we may now assume to the contrary that $k_2 > n$.
Thus one of $\hat e_1, \hat e_3$ contains more than $n/2$ edges,
so by Lemmas 2.4 and 2.2(4) $n$ is even, hence $n\geq 6$.  Note that
in this case $k_4 \leq n$, as otherwise we would have $4n = \sum k_i
\geq (n-2) + 3(n+1) > 4n$, a contradiction.  Thus $P_2, P_3$ contain
all the repeated labels and there are at least $n-2$ of them, so $k_2
+ k_3 \geq 3n-2$.  Since $k_2 = w_1 + w_3 \geq n+1$ and $w_3 \leq
(n/2)+2$, we have $w_1 \geq (n/2)-1$, so $w_3-w_1 \leq 3$.  By
equation (4) above, we have 
\begin{eqnarray*}
4n & = val(u_3, \ga)  = w_7 + k_2 + k_3 + (w_1 + w_2)  \\
& = (k_2 + k_3) + (w_7 + w_3 + w_2) - (w_3 - w_1)  \\
& \geq (3n-2) + 2n - 3 = 5n - 5 > 4n \end{eqnarray*}
This is a contradiction, which completes the proof for the case
$k_3>n$. 

We now assume $k_3 \leq n$.  Then each non-jumping label is a repeated
label in either $P_2$ or $P_4$.  By Lemma 2.7(2) we have $k_2, k_4
\leq n+2$, so each of $P_2, P_4$ contains at most $2$ repeated labels.
Since there are at least $n-2$ repeated labels, we have $n-2 \leq 4$,
i.e., $n\leq 6$.  As above, $P_2$ or $P_4$ containing a repeated label
implies that $n$ is even, so $n = 6$.  In this case $k_2 = k_4 = 8$,
so by Lemma 2.3(3) each of $\hat e_1, \hat e_2, \hat e_3$ contains
exactly 4 edges.

By equation (3) above we have $ n-2 \leq k_1 = w_7 \leq n$.  If $k_1 =
n-1$ or $n-2$ then the labels that do not appear in $P_1$ are repeated
labels of both $P_2$ and $P_4$, so there are at most 3 distinct
repeated labels, and hence at least $n-3 = 3$ jumping labels,
contradicting Claim 1.  If $k_1 = n$ then $k_1 + k_2 + k_4 = 3n+4$, so
the endpoints of the 4 edges of $\hat e_1$ at $P_2$ have the same
label sequence as those of $\hat e_1$ at $P_4$.  Hence these edges
form an extended Scharlemann cycle, which is impossible by Lemma
2.2(6).

\medskip

Claim 3.  {\it $\ga$ is kleinian.}

Clearly each of $u_1, u_2$ is incident to more than $2n$ edges.  Since
$w_7 + w_8 + w_9 = 4n - (w_5 + w_6) \geq 3n-2$, we have $w_7 \geq
n-2$.  By Claim 2 $w_{11}+w_{12} \geq 2n-2$, hence $w_{11}+w_{12}+w_7
\geq 3n-4 > 2n$, so $u_3$ is also incident to more than $n$ negative
edges.  Similarly for $u_4$.  By Lemmas 9.1 and 2.10(2) $\rgbp$ has at
most $2n$ edges, so two of the negative edges at $u_i$ are parallel on
$\gb$, hence $u_i$ is a label of Scharlemann bigon.  Since this is
true for all $i$, $\ga$ is kleinian.

\medskip
Claim 4.  $w_{11} + w_{12} = 2n - 2$.

Since $\ga$ is kleinian, by Lemma 6.2(4) there is a free involution
$\eta$ of $\ga$ mapping $u_1$ to $u_2$ and $u_3$ to $u_4$.  Thus
$\eta$ must map $\hat e_{10}$ to $\hat e_7$ and hence $w_7 =
w_{10}$.  We now have $w_2 + w_3 + w_7 = w_2 + w_3 + w_{10} \geq 4n -
w_8 - w_9 \geq 2n$.  Since $w_1 \geq 1$, we see that $w_{11} + w_{12}
= 4n - (2w_1 + w_2 + w_3 + w_7) \leq 2n - 2$.  The result now
follows from Claim 2.

\medskip

The involution $\eta$ maps $\hat e_{11}$ to $\hat e_{12}$.  Hence
Claim 4 implies that $w_{11} = w_{12} = n - 1$.  Since $\eta$ is
label preserving, the label sequence of $\hat e_{11} \cup \hat e_{12}$
on $\bdd u_3$ is the same as that on $\bdd u_4$, so we may assume
without loss of generality that the label sequences are as shown in
Figure 9.2.  One can see that in this case the transition function of
$\hat e_{11}$ defined in Section 2 is transitive, which implies that
all vertices of $\gb$ are parallel, contradicting the assumption.
This completes the proof of the lemma.  \qed

\begin{prop} Suppose $n_a = 4$ and $\gb$ is
non-positive.  Then $n_b \leq 4$.  \end{prop}

\proof Consider $\rgap$.  If $\rgap$ has a small component then by
Proposition 8.6 we have $n_b \leq 4$.  If $\rgap$ has no small
component then the component $G$ containing $u_1$ must also contain
$u_3$, and it is either of type (3), (9) or (10).  The result follows
from Lemma 9.4 if both components are of type (9) or (10), and from
Lemmas 9.5 and 9.6 if at least one component is of type (3).
\qed

\section {The case $n_1 = n_2 = 4$ and $\Gamma_1, \Gamma_2$
non-positive }

In this section we assume that $n_1 = n_2 = 4$ and $\ga$ is
non-positive for $a=1,2$.  We will show that this case cannot happen.
Denote by $X_a$ the union of $\rgap$ and all its disk faces.  By
Theorem 8.6 $\rgap$ has no small component, so each component of
$\rgap$ is of type (3), (9) or (10) in Figure 4.2.

\begin{lemma} Suppose $n_1 = n_2 = 4$, and both $\Gamma_1$ and
$\Gamma_2$ are non-positive.  Then at least one component of $\hat
\Gamma^+_1$ or $\hat \Gamma^+_2$ is of type (3).  \end{lemma}

\proof Suppose to the contrary that all components of $\hat
\Gamma^+_1$ and $\hat \Gamma^+_2$ are of type (9) or (10).  Then each
component of $X_a$ is an annulus, hence any vertex $u_i$ of $\ga$ is
incident to at most 2 families of negative edges.  By Lemma 2.3(1)
each negative family contains at most 4 edges, so $u_i$ is incident to
at most 8 negative edges, and hence the number of negative edges is no
more than the number of positive edges in $\ga$.  Since this is true
for $a=1,2$ and since a positive edge in one graph is a negative edge
in the other, the numbers of positive and negative edges of $\ga$ must
be the same, hence each vertex must be incident to exactly 8 positive
and 8 negative edges, so each negative family contains exactly 4
edges.  Since $\gb$ contains a loop, one of the negative families in
$\ga$ contains a co-loop and hence is a set of 4 parallel co-loops,
which is a contradiction to the 3-Cycle Lemma 2.14(2).  
\qed

\begin{lemma} Suppose $n_1 = n_2 = 4$, and both $\Gamma_1$ and
$\Gamma_1$ are non-positive.  If both components of $\rgap$ are of
type (3), and no component of $\rgbp$ is of type (3), then $\ga$ is
kleinian.  \end{lemma}

\proof Note that in this case all vertices of $\gb$ are boundary
vertices.  Let $u_1, u_2$ be the vertices of valence $2$ in $\rgap$.
By Lemma 6.6(1), $u_1$ is incident to fewer than $8$ positive edges,
hence there are three negative edges in $\ga$ incident to $u_1$ having
the same label $j$ at $u_1$.  On $\gb$ this implies that the vertex
$v_j$ is incident to at least three positive edges with label $1$ at
$v_j$; since $v_j$ is a boundary vertex, it is incident to at least 9
positive edges.  Since a loop at $v_j$ must have labels of different
parity on its two endpoints, we see that $val(v_j, \gbp)$ is even.  If
$val(v_j, \gbp) = 12$ then by Lemma 6.4(4) $\ga$ is kleinian.  Hence
we may assume that $val(v_j, \gbp) = 10$.  By Lemma 6.4(1) we may
assume that each family of positive edges at $v_j$ has at most 3
edges.  This implies that $v_j$ is incident to 2 or 3 loops.
Examining the labels we see that the two outermost loops form a
Scharlemann bigon, with $1$ as a label.  For the same reason, $2$ is a
label of a Scharlemann bigon in $\gb$.  If there is no
$(12)$-Scharlemann bigon then there must be $(14)$- and
$(23)$-Scharlemann bigon, so $\ga$ is kleinian and we are done.
Therefore we may assume that $\gb$ contains a $(12)$-Scharlemann
bigon.

The $(12)$-Scharlemann bigon and $\rgap$ cuts $\hat F_a$ into faces.
There is now only one edge class in these faces which connects $u_1$ to
$u_4$, hence by Lemma 2.2(5) $\gb$ contains no $(14)$-Scharlemann
bigon.  Similarly there is no $(23)$-Scharlemann bigon.  It follows
that all Scharlemann bigons of $\gb$ have label pair $(12)$.  In
particular, the two outermost loops at $v_j$ must form a
$(12)$-Scharlemann bigon.

We have shown above that the vertex $v_j$ has 2 or 3 loops.  If it has
2 loops then the weights of the positive families at $v_j$ are
$(2,3,3,2)$, and the middle label of a family of $3$ is a label of
Scharlemann bigon, which implies that both $3$ and $4$ are labels of
Scharlemann bigons, which is a contradiction.  Hence $v_j$ has exactly
3 loops $e_1, e_2, e_3$, and 4 non-loop edges divided to 2
non-Scharlemann bigons $e_4\cup e_5$ and $e_6 \cup e_7$.

As shown above, $e_1 \cup e_2$ is a (12)-Scharlemann bigon, so up to
symmetry we may assume that the edges $e_4, e_5$ have labels $4$ and
$1$ at $v_j$.  Since these edges do not form a Scharlemann bigon, the
labels at their other endpoints must be $3$ and $2$ respectively, so
$e_5$ is a $(12)$-edge, which must be parallel on $\ga$ to one of the
two $(12)$-loops $e_1, e_2$ because $\ga$ has only two families
connecting $u_1$ to $u_2$.  This is a contradiction because by Lemma
2.3(5) a loop and a non-loop edge cannot be parallel on $\ga$.
\qed

\begin{lemma} Suppose $n_1 = n_2 = 4$, and both $\Gamma_1$ and
$\Gamma_2$ are non-positive.  Then all components of $\hat \Gamma^+_1$
and $\hat \Gamma^+_2$ are of type (3), and $\hat \Gamma_1$ and $\hat
\Gamma_2$ are subgraphs of the graph shown in Figure 10.1.
\end{lemma}

\proof By Lemma 10.1 we may assume that $\rgap$ has a component $C$ of
type (3).  Let $u_1$ be the valence $2$ vertex in $C$.  If $u_1$ is
incident to at most $4$ families of parallel edges then each family
contains exactly 4 edges, but since the two positive families are
adjacent, this would be a contradiction to Lemma 6.6(1).  Therefore
$u_1$ is incident to at least 5 families of edges.  Note that if the
other component $C'$ of $\rgap$ is not of type (3), or if it is of
type (3) but the loop of $C'$ separates $u_1$ from the valence 2
vertex $u_2$ of $C'$ then $u_1$ would have only two families of
negative edges, which is a contradiction.  It follows that $C'$ is
also of type (3), so the graph $\rga$ is a subgraph of that in Figure
10.1.

For the same reason if some component of $\rgbp$ is of type (3) then
$\rgb$ is a subgraph of that in Figure 10.1 and we are done.
Therefore we may assume that no component of $\rgbp$ is of type (3).
By Lemma 10.2 $\ga$ is kleinian.

Consider the edges $\hat e_1, ..., \hat e_6$ of $\rga$ incident to
$u_3$.  See Figure 10.1.  Let $p_i$ be the weight of $\hat e_i$.  By
Lemmas 6.2(2) and 2.3(1), $p_i$ is even and at most 4.  Note that
$p_1, p_2, p_3$ are non-zero, so by Lemma
6.6(1) we know that $p_1+p_2$ and $p_1 + p_3$ are between 4 and 6.  If
$p_4=0$ then $u_2$ would have valence 4 in $\rga$, which would lead to
a contradiction as above.  Hence $p_4 > 0$.  If $p_5+p_6 = 0$ then
$val(u_3, \rgap)=5$, so either $p_1=4$ or $p_2=p_3=p_4 = 4$; either case
contradicts Lemma 6.6(1).  Therefore $p_5 + p_6 \geq 2$.

One can now check, from the labeling around the boundary of $u_3$,
that if $p_4 = 2$ then both labels at the edge endpoints of $\hat e_4$
at $u_3$ are jumping labels at $u_3$, and if $p_4 = 4$ then all labels
at the endpoints of $\hat e_5 \cup \hat e_6$ at $u_3$ are jumping
labels.  This is a contradiction because all vertices of $\rgbp$ are
boundary vertices and hence by Lemma 6.8 $u_3$ should have no jumping
label.  
\qed

\bigskip
\leavevmode

\centerline{\epsfbox{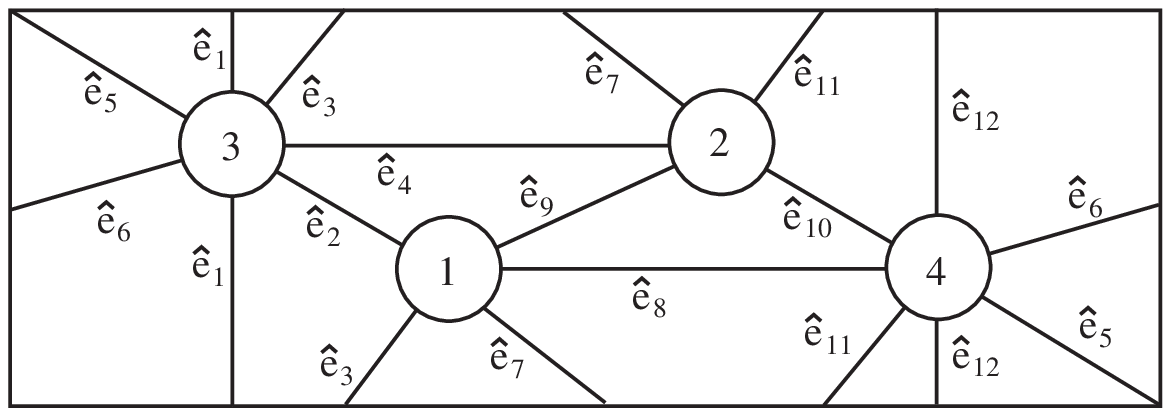}}
\bigskip
\centerline{Figure 10.1}
\bigskip

By Lemma 10.3 we may assume that both $\hat \Gamma_1, \hat
\Gamma_2$ are subgraphs of that in Figure 10.1.  Label the edges of
$\rga$ as in Figure 10.1, and let $p_i$ be the weight of $\hat e_i$.
Label $\rgb$ similarly using $\hat e'_i$.

\begin{lemma}  $\Delta = 4$, $p_2 + p_3 = 6$, $p_8=p_4=2$,
and $p_7=p_9=4$.  Moreover, edges in $\hat e_4 \cup \hat e_8$ are
co-loops, while those in $\hat e_7 \cup \hat e_9$ are not.
\end{lemma}

\proof If $p_2+p_3 > 6$ then one of the $\hat e_2, \hat e_3$ contains
4 edges, so by Lemma 6.4(1) $\ga$ is kleinian.  By Lemma 6.2(2) the
$p_i$ are even, so $p_2+p_3>6$ implies that $p_2=p_3 = 4$, which is a
contradiction to Lemma 6.6(2).  Therefore we have $p_2 + p_3 \leq 6$.
Since each $p_i \leq 4$, by counting edges at $u_1$ we have
$\Delta=4$, and $p_7 + p_8 + p_9 \geq 10$.

Recall that if a pair of negative edges of $\ga$ incident to $u_i$ are
parallel in $\gb$ then they form a Scharlemann bigon in $\gb$, hence
must have the same label pair in $\gb$, so they are incident to the
same pair of vertices in $\ga$.  Therefore no edge in $\hat e_7$ or
$\hat e_9$ is parallel in $\gb$ to an edge in $\hat e_8$.  Since
$\rgbp$ has at most 6 edges, we have $\max \{p_7, p_9\} + p_8 \leq 6$.
Since $p_i \leq 4$, this gives $p_7 + p_8 + p_9 \leq 10$, hence $p_2 +
p_3 = \Delta n_b - (p_7 + p_8 + p_9) \geq 6$.  Together with the
inequalities above, we have $p_2 + p_3 = 6$, $p_7 + p_8 + p_9 = 10$.
Since $\max \{p_7, p_9\} + p_8 \leq 6$ and $p_i \leq 4$, this holds
only if $p_8=2$ and $p_7 = p_9 = 4$.  Similarly $p_4 = 2$.

We have shown that edges in $\hat e_7 \cup \hat e_8$ belong to
distinct families in $\rgbp$.  Since $\rgbp$ has at most 6 edges,
$\hat e_7 \cup \hat e_8$ represent all edges in $\rgbp$.  If some edge
in $\hat e_i$ is a co-loop then all of them are.  Therefore the edges
in $\hat e_7$ cannot be co-loops because $\rgbp$ has only two loops.
It follows that the edges in $\hat e_8$ must be co-loops.  Similarly,
edges of $\hat e_4$ are co-loops, and those in $\hat e_9$ are not.
\qed

\begin{prop} Suppose both $\Gamma_1$ and $\Gamma_2$ are
non-positive, and $n_a = 4$.  Then $n_b < 4$. \end{prop}

\proof By Proposition 9.7 we have $n_b \leq 4$.  Assume to the
contrary that $n_b = 4$.  Since the two edges in $\hat e_8$ are
co-loops, they have labels $3,4$ at $u_1$.  Consider the three
negative edges $e_7, e_8, e_9$ such that $e_i \in \hat e_i$, and they
all have label $3$ at $u_1$.  In $\gb$ these are $1$-edges at $v_3$.
Since $e_8$ is a loop, it belongs to $\hat e'_1$.  The other two edges
are non-loop positive edges on $\gb$, so they belong to $\hat e'_2
\cup \hat e'_3$.  Applying Lemma 10.4 to $\gb$, we see that $\hat e'_2
\cup \hat e'_3 \cup \hat e'_4$ contains $8$ edges, so the two
edges $e_7, e_9$ are adjacent among the four edges labeled $1$ at
$v_3$ in $\gb$.  Since they are not adjacent among the four edges
labeled $3$ at $u_1$ in $\ga$, this is a contradiction to the Jumping
Lemma 2.18.  \qed

\section {The case $n_a = 4$, and $\gb$ positive 
}

In this section we assume that $n_a = 4$ and $\gb$ is positive.  We
will determine all the possible graphs for this case.  Recall from
Lemma 6.4(2) that in this case $\ga$ is kleinian, so the weights of
edges of $\rgb$ are all even.

\begin{lemma} Suppose $n_a = 4$ and $\gb$ is positive.  

(1) Two families of 4 parallel edges with the same label sequence at a
given vertex $v_j$ of $\gb$ connect $v_j$ to the same vertex $v_k$.

(2) There are at most three families of 4 parallel edges with the same
label sequence at any vertex $v_j$, and if $n_b > 2$ then there are at
most two such.

(3) if $\Delta = 4$ then $val(v_j, \rgb) \geq 5$ for all $j$; 

(4) if $\Delta = 5$ then $val(v_j, \rgb) \geq 6$ for all $j$;

(5) two weight 4 edges $\hat e_1, \hat e_2$ of $\rgb$ adjacent at a
vertex $v_j$ form an essential loop on $\hat F_b$.  \end{lemma}

\proof 
(1) If there are two families of 4 parallel edges with the same label
sequence $1,2,3,4$ at $v_j$ then by Lemma 6.5 the initial edges
$e_1, e'_1$ of the two families are parallel in $\ga$, with the same
label $j$ at the vertex $u_1$, hence the other endpoints of $e_1$
and $e'_1$ must also have the same label $k$, which implies that in
$\gb$ the two families have the same endpoints.

(2) If there were four then the leading edges of the $(12)$
Scharlemann bigons in these families are parallel in $\ga$, so they
belong to a family of at least $3n_b+1$ parallel edges connecting $u_1$
to $u_2$.  Since $\rgb$ has at most $3n_b$ edges by Lemma 2.5, two of
these edges would be parallel on both graphs, which is a contradiction
to Lemma 2.2(2).

If $n_b > 2$ and there are three families of 4 parallel edges with the
same label sequence at $v_j$ then as above there would be a family of
$2n_b+1$ parallel edges in $\ga$, which contradicts Lemma 2.22(3).

(3) and (4) follow immediately from (2).

(5) By (1) these two edges have their other endpoints at the same
vertex $v_k$, hence form a loop $C=\hat e_1 \cup \hat e_2$ on $\hat
F_b$.  If $C$ is not essential then we can choose $C$ to be an
innermost such cycle.  $C$ bounds a disk $D$ on $\hat F_b$, which must
contain some vertex because $\rgb$ is reduced.  If some vertex in the
interior of $D$ has valence 5 then it is incident to two adjacent
weight 4 edges, which would form another inessential loop,
contradicting the choice of $C$.  Hence all vertices in the interior
of $D$ have valence at least 6.  By Lemma 2.9 in this case there
should be at least three vertices on $\bdd D$, which is a
contradiction.  \qed

$\hat F_a$ separates $M(r_a) = M \cup V_a$ into the {\it black\/} and
{\it white\/} sides $X_B$ and $X_W$.  Since $\gb$ is kleinian (Lemma
6.4(2)), the black side $X_B$ is a twisted $I$-bundle over the Klein
bottle.  A face of $\gb$ is white if it lies in the white region
$X_W$, otherwise it is black.  In the next two lemmas we assume that
$\gb$ contains a white bigon and a white 3-gon as in Figure 11.1.

\bigskip
\leavevmode

\centerline{\epsfbox{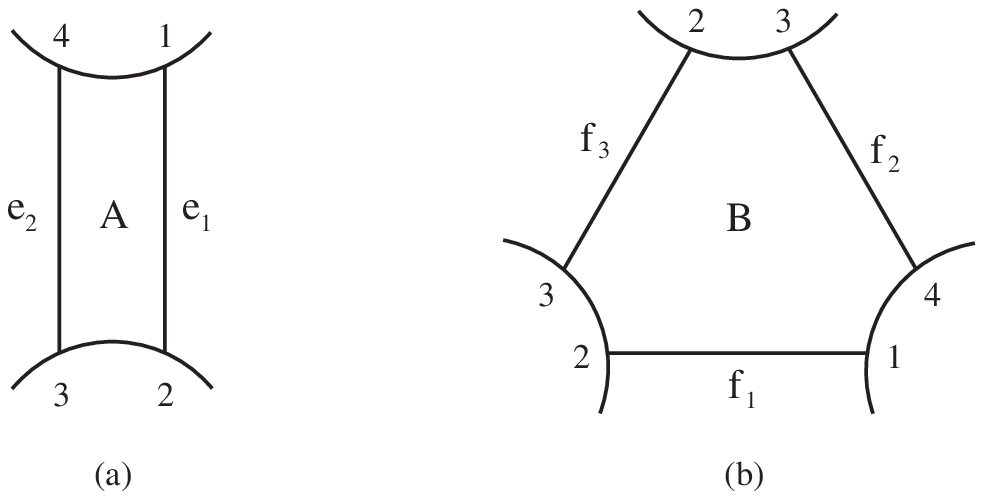}}
\bigskip
\centerline{Figure 11.1}
\bigskip

\begin{lemma} Suppose $n_a = 4$ and $\gb$ is positive.
Suppose $\gb$ contains a white bigon $A$ and a white 3-gon $B$ as in
Figure 11.1.  Then up to homeomorphism of $\hat F_a$, the edges of $A$
and $B$ appear on $\hat F_a$ as shown in Figure 11.2.
\end{lemma}

\proof Let $V_{23}$ and $V_{41}$ be the components of $V_a \cap X_W$
that run between $u_2, u_3$ and $u_4,u_1$, respectively.  Let $Y$ be a
regular neighborhood of $\hat F_a \cup V_{23} \cup V_{41} \cup A \cup
B$.  Then $\bdd Y = \hat F_a \cup T$, where $T$ is a torus in $M$, and
hence either $X_W = Y$, or $X_W$ is the union of $Y$ and a solid torus
along $T$.

Take a regular neighborhood $D$ of $e_1 \cup e_2 \cup f_3$ as
``base point'' for $\pi_1(\hat F_a) \cong \Bbb Z \times \Bbb Z$ and for
$\pi_1(X_W)$.  (See Figure 11.2).  The cores of the 1-handles $V_{23}$
and $V_{41}$ represents elements $x, y$ respectively of $\pi_1(X_W)$,
and $\pi_1(X_W)$ is generated by $\Pi = \pi_1(\hat F_a)$ together with
$x$ and $y$.  Note that since $\hat F_a$ is essential in $M(r_a)$,
$\Pi$ is a proper subgroup of $\pi_1(X_W)$.

\bigskip
\leavevmode

\centerline{\epsfbox{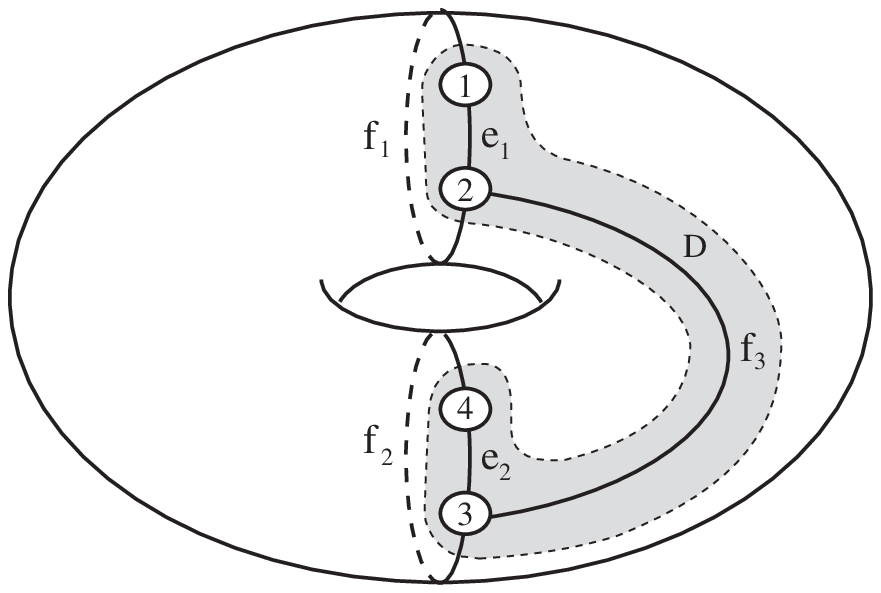}}
\bigskip
\centerline{Figure 11.2}
\bigskip

The bigon $A$ gives the relation $xy = 1$, so the $23$-corners and the
$41$-corner of $\bdd B$ represent $x$ and $x^{-1}$ respectively.  The
$32$-edge on $\bdd B$ lies in $D$ and hence represents $1$.  If one of
the other two edges $f_1, f_2$ of $B$ represents $1 \in \pi_1(\hat
F_a)$, then $B$ gives a relation of the form $x = \gamma$ for some
$\gamma \in \pi_1(\hat F_a)$, hence $\pi_1(X_W) = \Pi$, a
contradiction.  Therefore the edge $f_1$ is as shown in Figure 11.2,
and represents a nontrivial element $\alpha \in \pi_1(\hat F_a)$, when
oriented from $u_1$ to $u_2$.  Similarly, the edge $f_2$ represents a
non-trivial element $\gamma$, say, of $\pi_1(\hat F_a)$, when oriented
from $u_3$ to $u_4$.  Since $e_1 \cup f_1$ and $e_2\cup f_2$ are
disjoint, we must have $\gamma = \alpha$ or $\alpha^{-1}$.  The union
$V_{23} \cup V_{41} \cup N(A)$ is a single $1$-handle attached to
$\hat F_a$, and $\bdd B$ is a simple closed curve on $\hat F_a \cup
\bdd H$.  One can see that $\alpha x^2 \alpha^{-1} x^{-1}$ cannot be
realized by a simple closed curve.  It follows that $\gamma = \alpha$,
so $f_2$ appears on $\hat F_a$ as in Figure 11.2.  \qed

\bigskip
\leavevmode

\centerline{\epsfbox{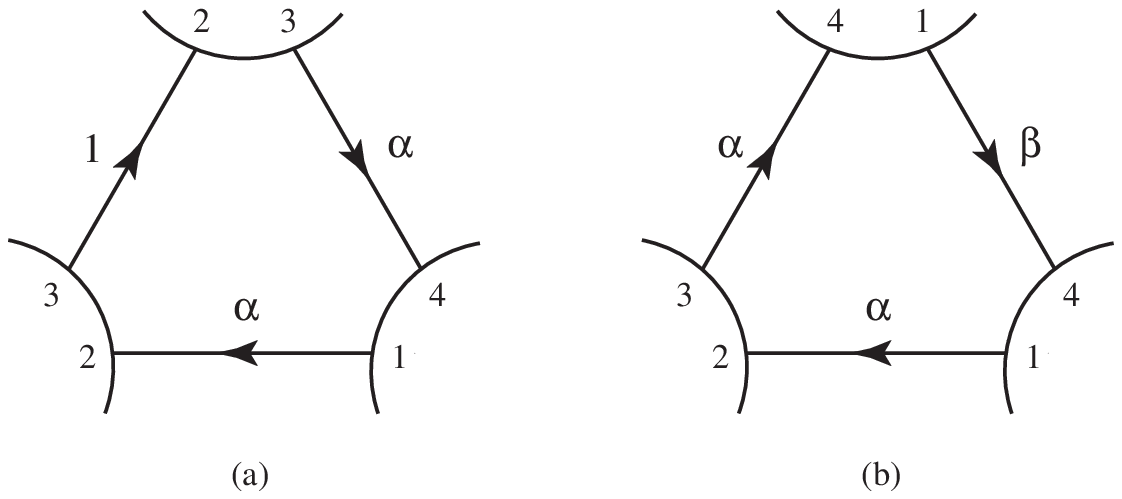}}
\bigskip
\centerline{Figure 11.3}
\bigskip

If we orient an edge $e$ of $\gb$ then the corresponding oriented edge
of $\ga$ represents an element $\gamma$ of $\pi_1(\hat F_a)$, and we
will label the edge $e$ with $\gamma$.  Thus the edge labels of the
bigon $A$ are both $1$, and the edge labels of the 3-gon $B$ are as in
Figure 11.3(a).  Note that any $12$-edge, oriented from $1$ to $2$, or
any $34$-edge, oriented from $3$ to $4$, has label $1$ or $\alpha$.
Also $\pi_1 (\hat F_a)$ has a basis $\{\alpha, \beta\}$, where $\beta$
is represented by an arc joining $u_1$ and $u_4$, disjoint from $\bdd
A$ and $\bdd B$, as shown in Figure 11.4.
\bigskip
\leavevmode

\centerline{\epsfbox{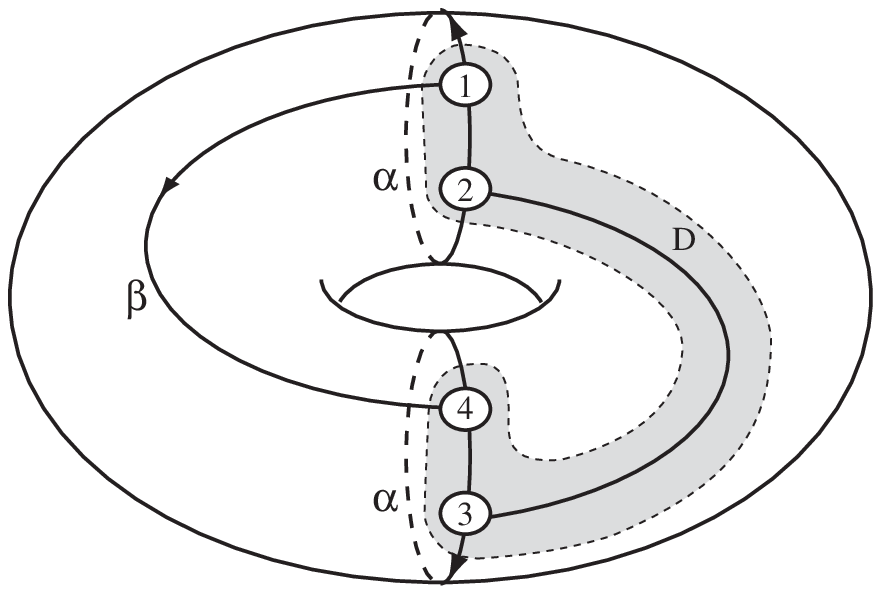}}
\bigskip
\centerline{Figure 11.4}
\bigskip

\begin{lemma} Suppose $n_a = 4$ and $\gb$ is positive.
Suppose $\gb$ contains a white bigon $A$ and a white 3-gon $B$ as in
Figure 11.1.  

(1) All edges of white bigons have label $1$.

(2) $\beta$ can be chosen so that any 3-gon must have edge labels as
shown in Figure 11.3(a) or (b).  \end{lemma}

\proof (1) Since the edges $e_1 \cup f_1$ and $e_2 \cup f_2$ on $\bdd
A$ and $\bdd B$ form two parallel essential circles on $\hat F_a$, any
$12$-edge of $\gb$ must be parallel to either $e_1$ or $f_1$, and any
$34$-edge of $\gb$ must be parallel to either $e_2$ or $f_2$.

Let $A'$ be another white bigon.  Applying Lemma 11.2 to $A'$ and $B$
gives that the $12$-edge of $A'$ and the $12$-edge $f_1$ of $B$ are
not parallel, and similarly the $34$-edge of $A'$ and the $34$-edge
$f_2$ of $B$ are not parallel.  Hence both edges of $A'$ are labeled
$1$ since they must be parallel to $e_1$ and $e_2$, respectively.

(2) Since $\gb$ has no extended Scharlemann cycle, each triangle face
$B'$ has either one or two $23$-corners.  If $B'$ has two $23$-corners
then applying Lemma 11.2 to $A$ and $B'$ shows that the $12$-edge of
$B'$ is not parallel to $e_1$, so it must be parallel to $f_1$ and
hence is labeled $\alpha$, as in Figure 11.3(a).  Similarly the
$34$-edge of $B'$ is also labeled $\alpha$.  If the $32$-edge of $B'$
is labeled $\gamma$ then $B$ and $B'$ together give the relation
$\gamma = 1$, so $\gamma = 1$, as in Figure 11.3(a).

If $B'$ has only one $23$-corner, let $f'_1, f'_2, f'_3$ be the $12$-,
$34$- and $14$-edges of $B'$, respectively.  Applying Lemma 11.2 to
$A$ and $B'$ gives that $f'_1$ is not parallel to $e_1$, so by the
above it must be parallel to $f_1$ and hence is labeled $\alpha$.
Similarly $f'_2$ is parallel to $f_2$ and is also labeled $\alpha$.
The two loops $e_1 \cup f_1$ and $e_2 \cup f_2$ cut $\hat F_a$ into
two annuli $A_L$ and $A_R$, where $A_R$ contains $f_3$; see Figure
11.4.  If the $14$-edge $f'_3$ of $B'$ lies in $A_R$ then it is
labeled $\alpha$, and $B'$ gives the relation $x^{-1} \alpha x \alpha
x^{-1} \alpha = 1$.  It is easy to see that, together with the
relation $x^2 \alpha x^{-1} \alpha = 1$ coming from $B$, this implies
$a = x^2$, so $x^5=1$, and hence $\alpha^5 = 1$, which is a
contradiction to the fact that $\alpha^5$ is a nontrivial element in
$\pi_1(\hat F_a)$ and hence is nontrivial in $\pi_1(X_W)$.  Therefore
$f'_3$ lies in $A_L$, as shown in Figure 11.4.  Let $\beta$ be the
corresponding element of $\pi_1(\hat F_a)$; then the edge labels of
$B'$ are as in Figure 11.3(b).  If $B''$ is any other $3$-gon with one
$23$-corner then the argument above in the case of $3$-gons with two
$23$-corners, using $A$ and the present $B'$, shows that the edge
labels of $B''$ are also as shown in Figure 11.3(b).  \qed

\bigskip
\leavevmode

\centerline{\epsfbox{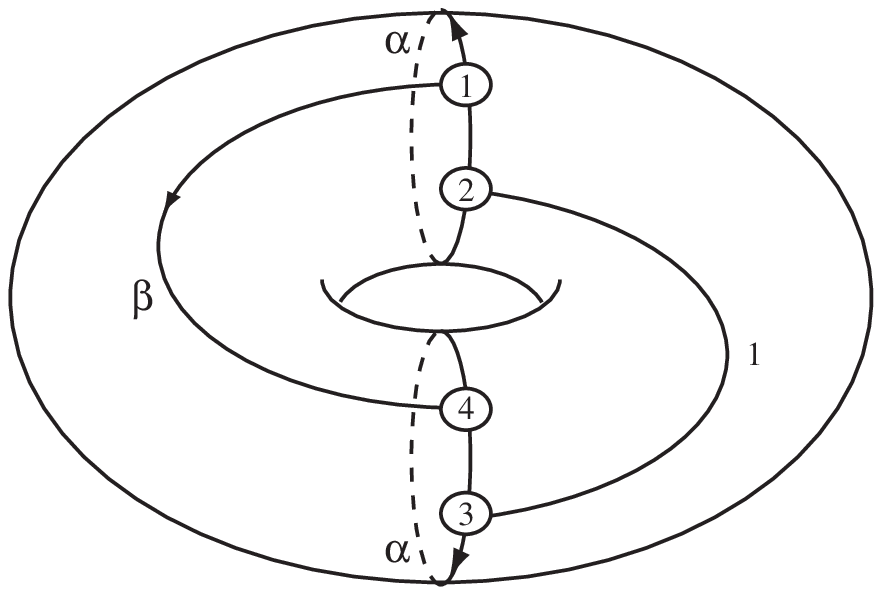}}
\bigskip
\centerline{Figure 11.5}
\bigskip

\begin{cor} Suppose $n_a = 4$ and $\gb$ is positive.

(1) Let $G$ be the subgraph of $\ga$ consisting of edges of white
bigons and white 3-gons on $\gb$.  Then the reduced graph $\hat G$ is
a subgraph of that in Figure 11.5.

(2) If $\gb$ contains a white bigon then it cannot contain a (black)
Scharlemann bigon which is flanked on each side by a (white) 3-gon.

(3) If $\gb$ contains a white bigon then it cannot contain three
3-gons occurring as consecutive white faces at a vertex.
\end{cor}

\proof  
(1) This follows immediately from Lemma 11.3.

(2) The edges of a (black) Scharlemann bigon are either $(12)$- or
$(34)$-edges, so by Lemma 11.3(2) both edges of the Scharlemann bigon
are labeled $\alpha$ and hence are parallel on $\ga$, contradicting
Lemma 2.2(2).

(3) By (2) the two black bigons between the white 3-gons are
$(12,34)$-bigons, as shown in Figure 11.6.  But then the middle 
3-gon would be a white Scharlemann cycle, contradicting Lemma 6.2(3).  
\qed

\bigskip
\leavevmode

\centerline{\epsfbox{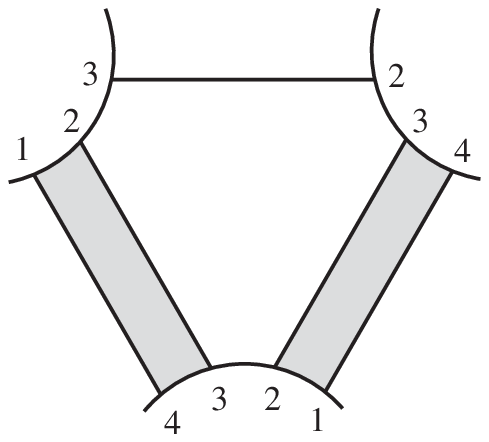}}
\bigskip
\centerline{Figure 11.6}
\bigskip

Let $G$ be a reduced graph on a torus $T$ with disk faces.  One can
endow $T$ with a singular Euclidean structure by letting each edge
have length 1 and each $n$-gon face a regular Euclidean $n$-gon.  The
cone angle $\theta(v)$ at a vertex $v$ of $G$ is the sum of the angles
of the corners incident to $v$.  Such a structure is {\it
hyperbolic\/} if $\theta(v) \geq 2\pi$ for all $v$, and $\theta(v) >
2\pi$ for some $v$.  The following lemma says that no singular
Euclidean structure on $T$ is hyperbolic.

\begin{lemma} Let $G$ and $\theta(v)$ be defined as above.
Then either $\theta(v) < 2\pi$ for some $v$, or $\theta(v) = 2\pi$ for
all $v$.  \end{lemma}

\proof Denote by $V, E, F$ the numbers of vertices, edges and faces of
$G$, respectively.  Let $\theta(c)$ be the angle at a corner $c$ of
the graph.  If $\sigma$ is a face of $G$, denote by $|\sigma|$ the
number of edges of $\sigma$.  Since $\sigma$ is a regular
$|\sigma|$-gon, for each corner $c \in \sigma$ we have $\theta(c) =
\pi(1 - 2/|\sigma|)$.  In the following, the first sum is over all
vertices $v$ of $G$, and the second is over all corners $c$.  Grouping
corners by faces $\sigma$, we get
\begin{eqnarray*}
\sum_v \theta(v) & = \sum_c \theta(c) = \sum_{\sigma} \sum_{c \in \sigma} \theta(c) 
= \sum_{\sigma} \sum_{c \in \sigma} \pi ( 1 - \frac{2}{|\sigma|})  \\
& = \pi (\sum_{\sigma} |\sigma| - \sum_{\sigma} \sum_{c \in \sigma} \frac 2{|\sigma|}) 
= \pi (2E - 2F)  = 2\pi(E - F) = 2 \pi V.
\end{eqnarray*}
Therefore $ \sum_v (2\pi - \theta(v)) = 0$, and the result follows.
\qed

\begin{lemma} Let $G$ be a reduced graph on a torus $T$ such
that $val(v) \geq 5$ for all $v$.  Then either

(1) there exists a vertex of valence 5 with at least four 3-gons
    incident; or

(2) there exists a vertex of valence 6 and all vertices of valence 6
    have all incident faces 3 gons; or 

(3) all faces of $G$ are 3-gons or 4-gons, and every vertex has
    valence 5 and has exactly three 3-gons incident.
\end{lemma}

\proof  We have $\theta(v) > 2\pi$ if $val(v)>6$, $\theta(v)
\geq 6 \times \pi/3 = 2\pi$ if $val(v) = 6$.  
Assuming (1) is not true, then we also have
$$\theta(v) \geq 3\times \frac{\pi}3 + 2 \times \frac{\pi}{4} = 2\pi$$
if $val(v)=5$.  Thus there is no vertex with cone angle $\theta(v) <
2\pi$, so by Lemma 11.5 we see that $\theta(v) = 2\pi$ for all $v \in
G$, hence there is no vertex of valence more than 6, all faces
incident to vertices of valence 6 are 3-gons, and exactly 3 faces
incident to a vertex of valence 5 are 3-gons and the other two are
4-gons.  Therefore either (2) or (3) holds.
\qed

\bigskip
\leavevmode

\centerline{\epsfbox{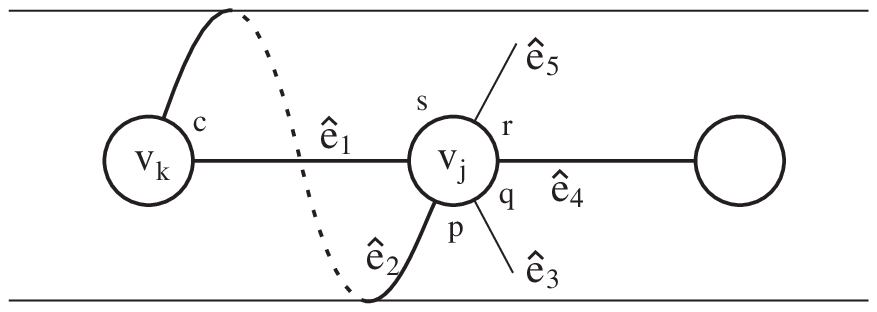}}
\bigskip
\centerline{Figure 11.7}
\bigskip

\begin{lemma} Suppose $n_a = 4$, $n_b > 2$, and $\gb$ is
positive.  Then no vertex $v_j$ of $\rgb$ with $val(v_j)=5$ has four
corners belonging to 3-gons.  \end{lemma}

\proof Let $v_j$ be a vertex of $\rgb$ with $val(v_j)=5$.  By Lemma
11.1(4) we must have $\Delta = 4$.  The weights of the edges at $v_j$
are $4,4,4,2,2$.  By Lemma 11.1(2) the three weight 4 edges are not
consecutive, hence the order around $v_j$ is $(4,4,2,4,2)$.  Label
these edges by $\hat e_1, ..., \hat e_5$, respectively, where $\hat
e_3, \hat e_5$ have weight 2.  By Lemma 11.1(5) the two edges $\hat
e_1, \hat e_2$ form an essential cycle on the torus $\hat F_b$, hence
the graph looks like that in Figure 11.7.

Let $c$ be the corner at $v_k$ between these two edges $\hat e_1, \hat
e_2$, as shown in Figure 11.7.  We claim that $c$ contains no other
edge endpoint.  Let $e$ and $e'$ be the edges in $\hat e_1$ and $\hat
e_2$ with label $1$ at $v_j$.  Let $P,Q$ be the endpoints of $e,e'$ at
$v_j$, and let $R,S$ be the endpoints of $e,e'$ at $v_k$,
respectively.  By Lemma 6.5 these edges are parallel on $\ga$, so they
connect the same pair of vertices $u_1, u_r$ for some $r$.  On $\gb$
this means that $P,Q$ have the same label $1$, and $R,S$ have the same
label $r$.  

Since $e,e'$ are parallel negative edges on $\ga$, we have
$d_{u_1}(P,Q) = d_{u_r}(R,S)$, therefore the four points $P,Q,R,S$
satisfy the assumptions of Lemma 2.16(1), hence by the lemma we have
$d_{v_j}(P,Q) = d_{v_k}(R,S)$.  Without loss of generality assume that
the orientations on $\bdd v_j, \bdd v_k$ are counterclockwise on
Figure 11.7.  Then one can see that $d_{v_j}(P,Q) = 4$, hence
$d_{v_k}(R,S) = 4$, which implies that there are only 3 edge endpoints
from the endpoint of $e$ to that of $e'$ on $\bdd v_k$, so there is no
edge endpoint in the corner $c$ in Figure 11.7.  This proves the
claim.

Label the corners at $v_j$ as shown in Figure 11.7.  The above implies
that the corner $p$ and $s$ belong to the same face $\sigma$, so if
$v_j$ is incident to at least four 3-gons then $\sigma$ must be a
3-gon, hence $\hat e_3 = \hat e_5$ is a loop.  Now the corners $q$ and
$r$ belong to the same face $\sigma'$, which cannot be a 3-gon, hence
the result follows.
\qed

\begin{lemma}  Suppose $n_a = 4$, $n_b > 2$, and $\gb$ is
positive.  Suppose $\Delta = 4$.  

(1) All faces of $\rgb$ are 3-gons or 4-gons, every vertex has valence
5 and has exactly three 3-gons incident, and the weight sequence of
the edges incident to the vertex is $(4,4,2,4,2)$.  In particular,
$n_b$ is even.

(2) The two weight 2 edges at any vertex form a loop, which is
incident to a 3-gon whose other two edges are of weight 4.

(3) Each edge in a weight 2 family of $\gb$ has label pair $(23)$ or
$(14)$.  \end{lemma}

\proof (1) By Lemma 11.6 and Lemma 11.1(3) and (4), $\rgb$ is one of
the three types stated there.  Lemma 11.7 shows that $\rgb$ cannot be
of type (1).  If $\rgb$ has a vertex of valence 6 then the weights are
$4,4,2,2,2,2$, hence there are two consecutive edges of weight 2.  By
Corollary 11.4(3) the three faces incident to these two edges cannot
all be 3-gons, hence $\rgb$ cannot be of type (2) in Lemma 11.6.  It
follows that $\rgb$ is of type (3) in Lemma 11.6, so the weights
of the edges at every vertex of $\rgb$ are $4,4,4,2,2$.  Thus the
number of weight 4 edge endpoints in $\rgb$ is $3n_b$, which must be
even, hence $n_b$ is even.  By Lemma 11.1(2) the three weight 4 edges
cannot all have the same label sequence, hence the weight sequence is
$(4,4,2,4,2)$ at each vertex.

(2) By Lemma 11.1(5) the two adjacent weight 4 edges at $v_j$ connects
$v_j$ to a vertex $v_k$ and form an essential loop on $\hat F_b$.  The
other weight 4 edge at $v_j$ connect to some vertex $v_r$, whose two
other weight 4 edges connect to another vertex $v_s$ and form an
essential loop.  These five weight 4 edges cut off a 6-gon containing
the four weight 2 edges at $v_j$ and $v_r$.  The 6-gon cannot contain
any vertex in its interior because each vertex is incident to two
weight 4 edges forming an essential cycle on $\hat F_b$.  Therefore
the four weight 2 edges at $v_j$ and $v_r$ form two loops.

(3) The loop $\hat e$ of $\rgb$ at $v_j$ cuts off a 3-gon in the 6-gon
above.  Let $e_1$ be the edge of $\hat e$ which is on the boundary of
a 3-gon face $\sigma$ of $\gb$.  Then the other two edges of $\sigma$
belong to families of 4 edges and hence must have label pair $(12)$
and $(34)$ respectively, so the labels on $\bdd \sigma$ are as shown
in Figure 11.3, and $e_1$ has label pair $(23)$ or $(14)$.  Since the
other loop edge of $\hat e$ is parallel to $e_1$, it has label pair
$(14)$ or $(23)$, respectively.  \qed

\bigskip
\leavevmode

\centerline{\epsfbox{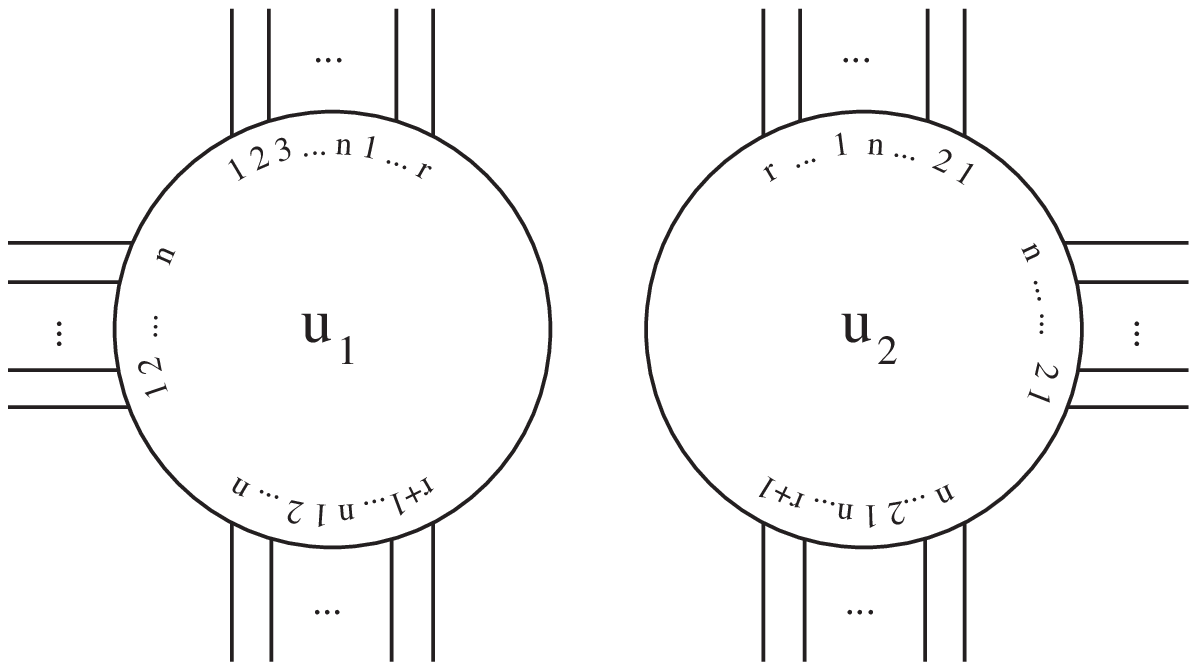}}
\bigskip
\centerline{Figure 11.8}
\bigskip

\begin{prop} Suppose $n_a = 4$ and $\gb$ is positive.  

(1) If $\Delta = 4$ then $n_b = 2$, and the graphs are as shown in
Figure 11.9.

(2) If $\Delta = 5$ then $n_b = 1$, and the graphs are as
shown in Figure 11.10. 
\end{prop}

\proof (1) Put $n=n_b$.  If $n = 1$ then the assumption $\Delta=4$
implies that the weights of the edges of $\rgb$ are either $(4,4,0)$
or $(4,2,2)$.  In either case a family of four edges form an extended
Scharlemann cycle, which is impossible.

Now assume $n > 2$.  Let $C = \hat e_1 \cup \hat e_2$ be the two edges
in $\rga$ connecting $u_1$ to $u_2$.  Since $\hat F_a$ contains both
$(12)$- and $(34)$-Scharlemann cocycles, any edge of $\gb$ with label
pair $(12)$ is parallel to an edge of $C$ on $\ga$.  By Lemma 11.8(1)
each vertex of $\rgb$ has valence $5$ and hence is incident to 3
weight 4 edges, so $\rgb$ has $3n/2$ weight 4 edges.  Each weight 4
edge contributes 2 edges to $C$, one for each $\hat e_i$, and by Lemma
11.8(3) no weight 2 edge of $\rgb$ contributes to $C$.  Thus each
$\hat e_i$ represents exactly $3n/2$ edges, and each label $j$ appears
exactly three times among the edge endpoints of $\hat e_1 \cup \hat
e_2$ at $u_1$.  Hence if the edge endpoints of $\hat e_1$ at $u_1$ are
labeled $1, ..., n, 1, ..., r$, then the labels of those in $\hat e_2$
must be $r+1, ..., n, 1, ..., n$, where $r=n/2$.  It follows that the
$n$ edge endpoints at $u_1$ that do not belong to $C$ must be on one
side of $C$, so up to relabeling we may assume that the labels at
$u_1$ are as shown in Figure 11.8(a).  Let $\varphi$ be the involution
on $\hat F_a$ given by Lemma 6.2(4).  Then $\varphi$ maps $\hat e_1$
to $\hat e_2$ and is label preserving, so the labels at $u_2$ must be
as shown in Figure 11.8(b).  Now the transition function of $\hat e_1$
maps $1$ to $r+1$, which has period $2$.  Since $n > 2$, this function
is not transitive, contradicting Lemma 2.3(1).

We have shown that $n=2$.  The graph $\rgb$ is now a subgraph of that
in Figure 13.1, with vertices labeled $v_1, v_2$ instead of $u_1,
u_2$.  Let $w_i$ be the weight of $\hat e_i$.  By Lemma 6.4(2) $\ga$
is kleinian, and by Lemma 6.2(2) the $w_i$ are all even.  By Lemma
11.1(3) we have $w_5 > 0$.  If $w_5 = 4$ then $\hat e_5$ containing no
extended Scharlemann bigon implies that either $w_1 + w_2 = 6$ or $w_3
+ w_4 = 6$, so $w_i =4$ for some $i\leq 4$, in which case $\hat e_i$
contains a $(12)$-Scharlemann bigon, whose edges, by the above, must
be parallel in $\ga$ to the $(12)$-edges of $\hat e_5$, which is a
contradiction to Lemma 2.3(5).  Therefore we must have $w_5 = 2$.  For
the same reason, the two loops in $\hat e_5$ cannot be a Scharlemann
bigon, so we must have $w_1 = w_2 = 4$ and $w_3+w_4=4$ up to symmetry.
Let $e_1\cup ... \cup e_4$ and $e'_1 \cup ... \cup e'_4$ be the edges
in $\hat e_1, \hat e_2$ respectively, such that $e_i, e'_i$ have label
$i$ at $v_1$.  By Lemma 6.5 $e_i, e'_i$ are parallel on $\ga$, with
the same label $1$ at $u_i$.  Therefore there is another edge between
them, which must belong to $\hat e_3 \cup \hat e_4$.  If $w_3=w_4 =2$
then one can check that these edges would have label pairs $(14)$ and
$(23)$, which is a contradiction.  Therefore we may assume $w_3=4$ and
$w_4=0$.  The graph $\gb$ is now as shown in Figure 11.9(b).

As shown above, there are 12 edges on $\gb$ with label pair $(12)$ or
$(34)$, divided into 4 families of 3 edges each on $\ga$.  Label the
edges of $\gb$ as in the figure.  Up to symmetry we may assume the
edge $A$ is as shown in Figure 11.9(a).  Since $\Delta=4$, we may
assume that the jumping number $J=1$.  The 1-edges around $v_1$ appear
in the order $A,E,G,P$.  By Lemma 6.5 $G,P$ are parallel on $\ga$.
This determines the position of these edges as well as the orientation
of $u_1$.  The $2$-edges at $v_2$ appear in the order $E,R,G,P$,
and $E$ on $\ga$ has already been determined above, hence the position
of $R,G,P$ must appear around $\bdd u_2$ are shown in the figure.
Other edges on $\ga$ can be determined similarly.  

\bigskip
\leavevmode

\centerline{\epsfbox{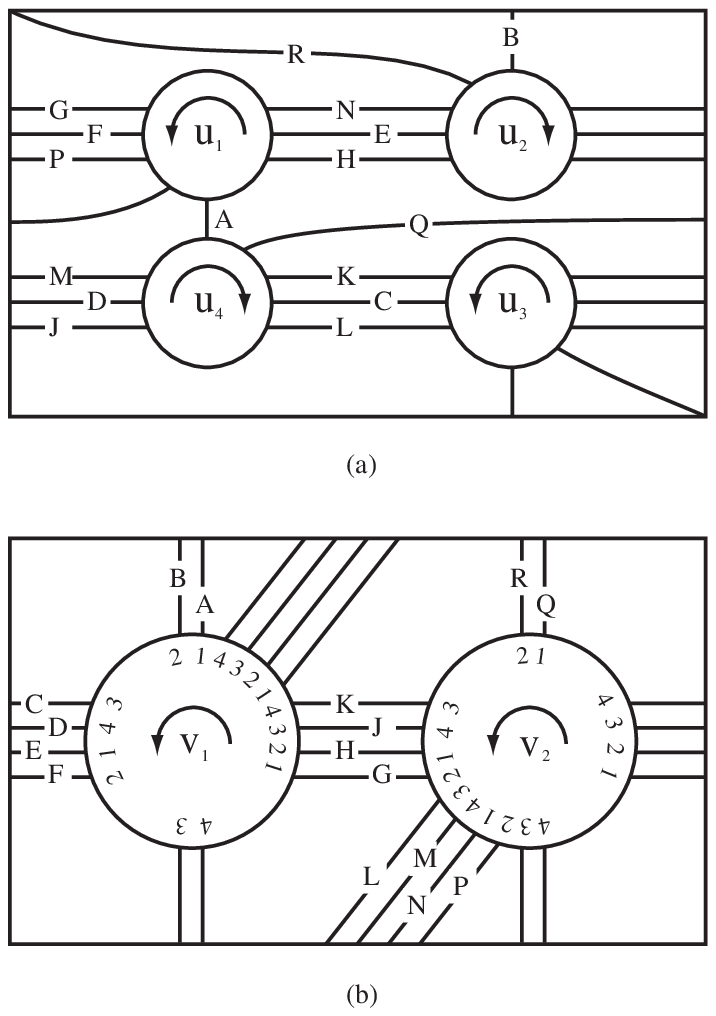}}
\bigskip
\centerline{Figure 11.9}
\bigskip

(2) The proof for $\Delta = 5$ is similar but simpler.  In this case
by Lemma 11.1(4) each vertex $v_j$ of $\rgb$ has valence 6, and the
edges at $v_j$ have weights $4,4,4,4,2,2$.  Thus all white faces are
bigons and 3-gons, so by Corollary 11.4(1) $\rga$ is a subgraph of
that in Figure 11.5.  By Lemma 11.1(2) the two weight 2 edges at any
vertex are non-adjacent, thus any edge $e$ in a weight 2 family is on
the boundary of a 3-gon whose other two edges have label pairs $(12)$
and $(34)$, so the label pair of $e$ must be $(14)$ or $(23)$ and
hence $e$ is not a vertical edge in Figure 11.5.  On the other hand,
each weight 4 edge of $\rgb$ contributes one edge to each vertical
family in Figure 11.5, hence each vertical edge has weight exactly
$2n_b$.  As above, one can show that the transition function defined
by a family of vertical edges is the identity function, which by Lemma
2.3(1) implies that $n_b = 1$, so the graph $\gb$ must be as shown in
Figure 11.10(b).  By the above discussion, $\ga$ is as shown in Figure
11.10(a).  Label edges as in Figure 11.10(b).  By Lemma 6.5 the edges
$A,H$ are parallel on $\gb$, hence we may assume the jumping number
$J=1$.  One can now easily determine the labels of the edges of $\ga$.
\qed

\bigskip
\leavevmode

\centerline{\epsfbox{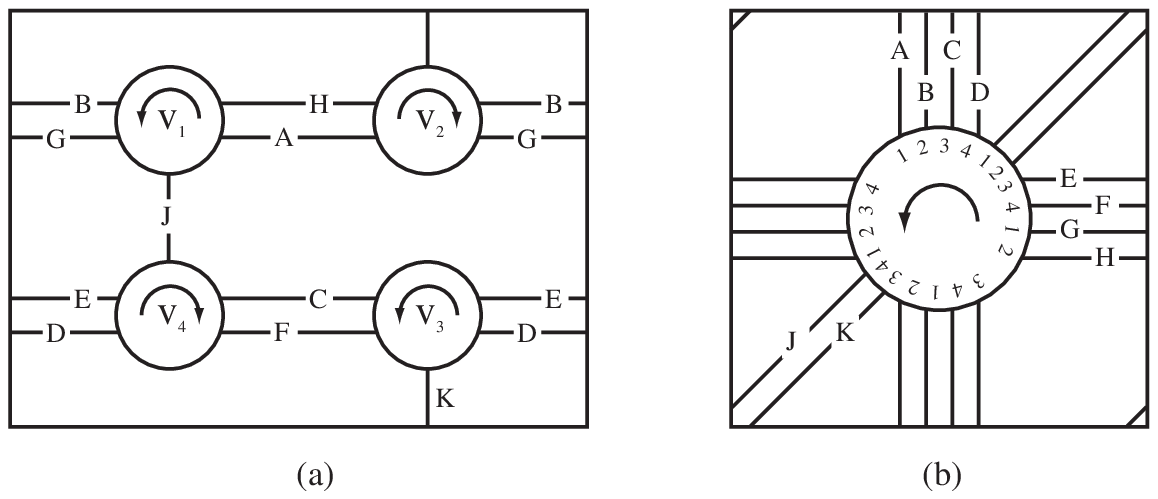}}
\bigskip
\centerline{Figure 11.10}
\bigskip

\begin{prop}  Suppose $n_a \leq n_b$.  Then $n_a \leq 2$.
\end{prop}

\proof By Proposition 5.11 we have $n_a \leq 4$.  Assume $n_a=4 \leq
n_b$. Then by Proposition 11.9, $\gb$ is non-positive. Therefore, by
Proposition 9.7, $n_a=n_b=4$, and, by Proposition 11.9 again, $\ga$ is
also non-positive, contradicting Proposition 10.5.

Suppose $n_a = 3$.  Then $\hat F_a$ is non-separating.  If $\ga$ is
non-positive then some vertex $u_1$, say, has different sign to the
other two vertices $u_2, u_3$.  One of the vertices has valence at
most two in $\rgap$, so it is incident to at most $2n_b$ positive
edges, and hence at least $2n_b$ negative edges.  By Lemma 2.8 this
implies that $\gb$ has a Scharlemann cycle, so by Lemma 2.2(4) $\hat
F_a$ is separating, which is a contradiction.  If $\ga$ is positive
then by Lemma 3.1 we have $n_b \leq 4$.  By Lemma 2.23 $n_1, n_2$
cannot both be odd, hence $n_a=3$ implies that $n_b$ is even, so we
must have $n_b = 4$.  Now applying Proposition 11.9 with $n_a, n_b$
reversed, we get $n_a \leq 2$, which is a contradiction.  \qed

\section {The case $n_a=2$, $n_b \geq 3$, and $\gb$ positive
}

The next few sections deal with the case that $n_a=2$ and $n_b \geq
3$.  The main result of this part is Proposition 16.8, which shows
that there are only a few possibilities for the graphs $\ga, \gb$.

Throughout this section we will assume that $n_a=2$, $n_b \geq 3$ and
$\gb$ is positive.  We will show that this case is impossible.  To
simplify notation, denote $n_b$ by $n$.  Note that $\rga$ has at most
four edges $\hat e_1, ..., \hat e_4$, all connecting $u_1$ to $u_2$.
We will always assume that the first edge of $\hat e_1$ has label $1$
at $u_1$.  Write $\ga = (a_1, a_2, a_3, a_4)$, where $a_i$ is the
weight of $\hat e_i$.  Let $s_i = s(\hat e_i)$ be the transition
number of $\hat e_i$ from $u_1$ to $u_2$.  In the following lemma the
subscripts are mod $4$ integers.

\begin{lemma}  
(1) $s_{i+1} \equiv s_i - a_i - a_{i+1}$ mod $n$;

(2) $s_{i+2} \equiv s_i - a_{i+1} + a_{i+3}$ mod $n$.  In particular,
$s_i \equiv s_{i+2}$ if and only if $a_{i+1} \equiv a_{i+3}$ mod $n$.
\end{lemma}

\proof (1) Orient edges from $u_1$ to $u_2$, and denote by $e(h),
e(t)$ the head and tail of an edge $e$.  Let $e, e'$ be the first edge
of $\hat e_i, \hat e_{i+1}$, respectively.  Let $x$ be the label of
$e(t)$, and $y$ the label of $e'(h)$.  Then traveling from $e(t)$ to
$e'(t)$ on $\bdd u_1$ then to $e'(h)$ through $e'$ gives $y \equiv x+
a_i + s_{i+1}$ mod $n$, while traveling through $e$ to $e(h)$ then
along $\bdd u_2$ to $e'(h)$ gives $y \equiv x + s_i - a_{i+1}$ mod
$n$.  Hence $s_{i+1} \equiv s_i - a_i - a_{i+1}$ mod $n$.

(2) Applying (1) twice gives $s_{i+2} \equiv s_i - a_i - a_{i+1} -
a_{i+1} - a_{i+2} \equiv s_i - a_{i+1} + a_{i+3}$ mod $n$.  \qed

\begin{lemma}  Let $e, e'$ be edges of $\gb$ joining a pair of
distinct vertices, such that $e \cup e'$ is null-homotopic in $\hat
F_b$.  If $e$ belongs to a family of at least $n$ parallel edges in
$\ga$ then $e$ and $e'$ are parallel on $\gb$.
\end{lemma}

\proof
Let $D$ be the disk in $\hat F_b$ bounded by $e \cup e'$.  The family
of $n$ parallel edges of $\ga$ containing $e$ gives a set of essential
loops on $\hat F_b$, corresponding to the orbits of the associated
permutation.  It follows that $D$ contains no vertices in its
interior, and hence $e$ and $e'$ are parallel on $\gb$.
\qed

\begin{lemma} Suppose $\ga$ contains bigons $e_1 \cup e_2$
and $e'_1 \cup e'_2$, such that $e_1,e'_1$ have label $i$ at $u_1$ and
$j$ at $u_2$, and $e_2,e'_2$ have label $i+1$ at $u_1$ and $j+1$ at
$u_2$.  Suppose either

(i) $j \neq i \pm 1$, or

(ii) $\ga$ contains a pair of edges $e_3, e'_3$ with the same
label pair $(r, s)$, such that $r,s \notin \{i, i+1, j, j+1\}$ and
$C_3 = e_3 \cup e'_3$ form an essential loop on $\hat F_b$.

Then $C_1 = e_1 \cup e'_1$ and $C_2 = e_2 \cup e'_2$ are inessential
on $\hat F_b$.
\end{lemma}

\proof
(i) In this case $C_1 \cap C_2 = \emptyset$, so they cannot be
essential and yet non-homotopic on $\hat F_b$, hence by Lemma
2.20 they must be inessential.

(ii) In this case $C_1, C_2$ lie in the interior of the annulus
obtained from $\hat F_b$ by cutting along $C_3$, which again implies
that $C_1, C_2$ cannot be essential and yet non-homotopic on $\hat
F_b$.
\qed

An edge is a {\it border edge\/} if it is the first or last edge in a
family of parallel edges.

\begin{lemma}  Suppose $s_k \neq \pm 1$.  Then

(1) $a_k \leq n+1$, and 

(2) if $a_k = n+1$ and $e'$ is an edge of $\hat e_j$ which has the
same label pair as that of a border edge $e_1$ of $\hat e_k$, then
$e'$ is a border edge.
\end{lemma}

\proof
(1) Assume $a_k \geq n+2$.  Label the first $n+2$ edges of $\hat e_k$
successively as $e_1, e_2, ..., e_n, e_{n+1}, e_{n+2}$.  Let $e'_i =
e_{n+i}$.  Since $s_k
\neq \pm 1$, $e_1, e_2, e'_1, e'_2$ satisfy Condition (i) in
Lemma 12.3, so $e_1 \cup e'_1$ is an inessential loop on $\hat
F_b$.  By Lemma 12.2 this implies that $e_1$ and $e'_1$ are parallel
on $\gb$, and hence parallel on both graphs, which is a contradiction.

(2) If $e'$ is not a border edge then the bigon $e_1 \cup e_2$ and
one of the two bigons containing $e'$ satisfy the assumption of
Lemma 12.3(i), hence $e'$ is parallel to $e_1$ on $\gb$.
Similarly, using the bigon between $e_n, e_{n+1}$ and the other bigon
containing $e'$ one can show that $e'$ is also parallel to $e_{n+1}$
on $\gb$, hence $e_1, e_{n+1}$ are parallel on both $\ga$ and $\gb$,
which is again a contradiction.
\qed

\begin{lemma}  Let $\hat e, \hat e'$ be families of at least
$n$ parallel edges in $\ga$, and let $i,j,k,l \in \Bbb Z_n$ be
distinct.  Then $\ga$ cannot contain both

(i) $ij$-edges $e_1, e_2, e_3$ with $e_1, e_2 \in \hat e$ and $e_3$
non-equidistant with $e_1, e_2$; and 

(ii) $kl$-edges $e'_1, e'_2 \in \hat e'$.
\end{lemma}

\proof
The edges $e_1, e_2, e_3$ are pairwise non-parallel in $\gb$.  Since
$e_1, e_2 \in \hat e$, no pair of $e_1, e_2, e_3$ cobounds a disk in
$\hat F_b$ by Lemma 12.2.  Hence $e'_1, e'_2$ cobound a disk in $\hat
F_b$.  Since $e'_1 \in E'$, $e'_1, e'_2$ are parallel on $\gb$ by
Lemma 12.2.  This contradicts Lemma 2.2(2).
\qed

Two families $\hat e, \hat e'$ of $\ga$ are {\it A-conjugate\/} it
there are $e \in \hat e$ and $e' \in \hat e'$ such that they are
anti-parallel on $\gb$ when oriented on $\ga$ from $u_1$ to $u_2$.
They are {\it P-conjugate\/} if the $e, e'$ above are parallel on
$\gb$ as oriented edges.  They are {\it conjugate\/} if they are
either A-conjugate or P-conjugate.

\begin{lemma} (1) There exist $\hat e_i, \hat e_j$ on $\rga$
which are A-conjugate.  Moreover, if $a_i<(\Delta-3)n$ or
$a_j<(\Delta-3)n$ then there is another such pair.  (The two pairs may
have one family in common.)

(2) If $\hat e_i, \hat e_j$ are A-conjugate then $s_i \equiv - s_j$
mod $n$; moreover, if $a_i \geq n$ or $a_j \geq n$ then $s_i \not
\equiv s_j$ mod $n$.  \end{lemma}

\proof (1) Since there are $\Delta n$ edges while $\rgb$ has at most
$3n$ edges, $\gb$ has at least $(\Delta - 3)n$ bigons.  The two edges
of a bigon in $\gb$ belong to a pair of A-conjugate families $\hat
e_i, \hat e_j$ on $\ga$.  If $a_i<(\Delta-3)n$ or $a_j<(\Delta-3)n$
then these families cannot contain all the bigons on $\gb$, hence
there must be another A-conjugate pair.

(2) If $\hat e_i, \hat e_j$ are A-conjugate then by definition
there exist $e \in \hat e_i$ and $e' \in \hat e_j$ which are
anti-parallel on $\gb$, hence the label of $e'$ at $u_2$ is the same
as that of $e$ at $u_1$, and vice versa.  Therefore $s_i \equiv -s_j$
mod $n$.

If we also have $s_i \equiv s_j$ mod $n$ then $s_i=0$ or $n/2$, which
is a contradiction to the 2-Cycle Lemma 2.14(3). 
\qed

\begin{lemma} Let $\hat e = e_1 \cup ... \cup e_p$ and $\hat
e' = e'_1 \cup ... \cup e'_q$ be two families of $\ga$, where $p\leq
q$.

(1) If $\hat e$ and $\hat e'$ are conjugate then $p\equiv q$ mod $2$.

(2) If $\hat e$ and $\hat e'$ are conjugate and $q\geq p$ then each
edge $e_r$ is parallel to the edge $e'_{r+c}$, where $c=(q-p)/2$;
hence the set of edges in $\hat e'$ which are parallel to those in
$\hat e$ lie exactly in the middle of $\hat e'$.

(3) If $p+q \equiv 0$ mod $2n$ and $\hat e, \hat e'$ are adjacent on
$\rga$ then they are not A-conjugate.  

(4) If $p+q \equiv 0$ mod $2n$, $\hat e, \hat e'$ are adjacent on
$\rga$, and $J\neq \pm 1$ then they are not conjugate.  
\end{lemma}

\proof (1) By definition there are edges $e_i, e'_j$ which are
parallel in $\gb$.

First consider the case that $\hat e, \hat e'$ are adjacent.
Denote by $e(k)$ the endpoints of $e$ at $u_k$.  We may assume that
the first edge $e'_1$ of $\hat e'$ is adjacent to the last edge $e_p$
of $\hat e$ on $\bdd u_1$.  Then the distance from $e_i(1)$ to
$e'_j(1)$ is 
\begin{eqnarray}
d_{u_1}(e_i, e'_j) = j + p - i 
\end{eqnarray}
On $\bdd u_2$
$e'_q(2)$ is adjacent to $e_1(2)$, so we have 
\begin{eqnarray}
d_{u_2}(e'_j, e_i) = i + q - j 
\end{eqnarray}
Since $e_i, e'_j$ are parallel positive edges on $\gb$, they
are equidistant, hence by Lemma 2.17 we have $d_{u_1}(e_i, e'_j) =
d_{u_2}(e'_j, e_i)$, which gives 
\begin{eqnarray}
2(j-i) = q-p 
\end{eqnarray}
and 
\begin{eqnarray}
2d = p+q 
\end{eqnarray}
where $d = d_{u_1}(e_i, e'_j) = d_{u_2}(e'_j, e_i)$.  Equation (C)
gives $q-p \equiv 0$ mod $2$.

Now suppose $\hat e, \hat e'$ are not adjacent.  Let $\hat e''$ be the
family whose endpoints on $\bdd u_1$ are between $e_p(1)$ and
$e'_1(1)$.  Then on $\bdd u_2$ the endpoints of $\hat e''$ are also
exactly the ones between $e'_q(2)$ and $e_1(2)$.  Thus if $\hat e''$
has $k$ edges then the equations (5) and (6) above become $d = j + k
+ p - i$ and $d = i + k + q - j$.  Therefore again we have $2(j-i) =
q-p$, and the result follows.

(2) From equation (7) we have $j = i + (q-p)/2$.  If $i>1$ then
the above and the condition $q\geq p$ imply that $j>1$.  Since $e_i$
is parallel to $e'_j$ on $\gb$, by Lemma 2.20 applied to the bigons
$e_{i-1}\cup e_i$ and $e'_{j-1} \cup e'_j$, the loop $e_{i-1} \cup
e'_{j-1}$ is null-homotopic on $\hat F_b$, hence by Lemma 12.2
$e_{i-1}$ is parallel to $e'_{j-1}$ on $\gb$.  Similarly, if $i<p$
then the edge $e_{i+1}$ is parallel to $e'_{j+1}$ on $\gb$.  By
induction it follows that every edge $e_k$ in $\hat e$ is parallel to
the edge $e'_{k+(q-p)/2}$.

(3) Assume without loss of generality that $\hat e = \hat e_1$ and
$\hat e' = \hat e_2$.  If they are A-conjugate then the label of
$e'_j$ at $u_1$ is the same as the label of $e_i$ at $u_2$, so $d =
d_{u_1} (e_i, e'_j) \equiv s_1$ mod $n$.  Hence equation (8) and the 
assumption $p+q\equiv 0$ mod $2n$ gives $s_1 \equiv d \equiv
0$ mod $n$.  Since $n\geq 3$, this is a contradiction to Lemma 2.14(2).

(4) By (3) $\hat e, \hat e'$ are not A-conjugate.  Assume they are
P-conjugate.  

If $p+q=4n$ then by Lemma 2.22(3) we have $p=q=2n$, and by (2) each
$e_i$ of $\hat e_1$ is parallel to the corresponding edge $e'_i$ of
$\hat e_2$ on $\gb$ for $i=1,...,2n$.  Since $\hat e, \hat e'$ are not
A-conjugate, $e_i, e'_i$ are parallel as oriented edges, with
orientation from $u_1$ to $u_2$.  Hence there is another edge $e''_i$
between them, which cannot belong to $\hat e \cup \hat e'$ as
otherwise there would be two edges parallel on both graphs,
contradicting Lemma 2.2(2).  This gives at least $6n$ edges on $\ga$,
which is a contradiction.

Now assume $p+q=2n$.  Let $e, e'$ be the edges of $\hat e, \hat e'$
which are parallel on $\gb$ as oriented edges, so they have the same
label $k$ at $u_1$ for some $k$.  The condition $p+q=2n$ implies that
$e, e'$ are adjacent among edges labeled $k$ at $u_1$.  Since $J \neq
\pm 1$, they are non-adjacent on $\gb$ among edges labeled $1$ at
$v_k$, hence they belong to a family of at least 5 parallel edges,
which is a contradiction to Lemma 2.2(2) because $\rga$ has at most 4
edges.  \qed

\begin{lemma} If the jumping number $J=\pm 1$ (in particular
if $\Delta = 4$), then $\ga$ has at most $n+1$ parallel edges. 
\end{lemma}

\proof Suppose for contradiction that $\hat e_1$, say, contains edges
$e_1, ..., e_{n+2}$.  By Lemma 12.4(1) we may assume that $s_1 = 1$,
so the label sequences of these edges are $(1,2,..., n, 1,2)$ at
$u_1$, and $(2,3,...,n,1,2,3)$ at $u_2$.  By Lemma 2.22(1) we may
assume that the subgraph of $\gb$ consisting of these edges is as
shown in Figure 2.3.  Up to symmetry we may assume that the
orientation of $\bdd v_i$ is counterclockwise on Figure 2.3.

Orient edges from $u_1$ to $u_2$.  Denote by $h_i, t_i$ the head and
tail of $e_i$, respectively.  For $i>1$, $h_{i-1}$ and $t_i$ both have
label $i$ on $\ga$, so they are on $\bdd v_i$.  Define $d_i =
d_{v_i}(t_i, h_{i-1})$, where $i = 2, ..., n+2$.  Note that $d_i = 1$
implies that the corner from $t_i$ to $h_{i-1}$ on $\bdd v_i$ contains
no edge endpoint.

\medskip
CLAIM 1. {\it $d_i = d_j$ for $2 \leq i,j \leq n+2$.}
\medskip

\proof  Isotoping on $T_0$ along the positive direction of $\bdd u_i$
moves $h_1$ to $h_2$ and $t_2$ to $t_3$, so the distance on $\bdd v_2$
from $h_1$ to $t_2$ should be the same as that on $\bdd v_3$ from
$h_2$ to $t_3$, i.e., $d_2 = d_3$.  (Alternatively one may apply Lemma
2.16 to obtain the result.)  Similarly we have $d_i = d_{i+1}$
for $2\leq i \leq n+1$.
\qed

\medskip
CLAIM 2. {\it $d_i = 1$ for $2 \leq i \leq n+2$.}
\medskip

\proof By assumption we have $J = \pm 1$, so either $d_{v_2}(h_1,
h_{n+1}) = 2$ (when $J=1$), or $d_{v_1}(t_{n+1}, t_1) = 2$ (when
$J=-1$).  In the first case, from Figure 2.3 we see that the tail
of $e_{n+2}$ is the only edge endpoint at the corner from $h_1$ to
$h_{n+1}$, hence $d_{n+2} = 1$.  Similarly in the second case the head
of $e_n$ is the only edge endpoint on $\bdd v_1$ from $t_{n+1}$ to
$t_1$, hence $d_{n+1} = d_{v_1}(t_{n+1}, h_n) = 1$.  In either case by
Claim 1 we have $d_i = 1$ for all $i$ between $2$ and $n+2$.
\qed

Let $D$ be the disk face indicated in Figure 2.3.  Claim 2 shows that all
corners of $D$ except $c_1, c_2, c_3$ shown in the figure contain no
edge endpoints.

When $J=1$, we have $d_{v_1}(t_1, t_{n+1}) = 2$, so there is one edge
endpoint in $c_1$.  Similarly there is one edge endpoint in $c_3$.
Since $d_{v_2}(t_{n+2}, t_2) = 2\Delta - 2 \geq 6$, there are at least
4 edge endpoints in $c_2$.  Thus there would be some trivial loops
based at $v_2$, contradicting the assumption that $\gb$ has no trivial
loops.

When $J = -1$, $d_{v_1}(t_1, t_{n+1}) = d_{v_3}(h_2, h_{n+2}) =
2\Delta - 2 \geq 6$, and $d_{v_2}(t_{n+2}, t_2) = 2$, so there are at
least 5 edge endpoints in each of $c_1$ and $c_3$, and no edge
endpoints in $c_2$.  It follows that $D$ contains at least 5 interior
edges, all parallel to each other, two of which would then be parallel
on both graphs, contradicting Lemma 2.2(2).  \qed

\begin{lemma}  $\ga$ has at most $n+2$ parallel edges.  
\end{lemma}

\proof Assume to the contrary that $\hat e_1 \supset e_1 \cup ... \cup
e_{n+3}$.  By Lemma 12.4(1) we may assume without loss of generality
that $s_1 = 1$.  By Lemma 2.22(1) the first $n+2$ edges appear in
$\gb$ as shown in Figure 2.3

First assume $n\geq 4$. Orient edges of $\ga$ from $u_1$ to $u_2$, and
denote by $e(h), e(t)$ the head and tail of an edge $e$, respectively.
From Figure 2.3 we see that $(e_2(h), e_{n+2}(h), e_3(t))$ is a
positive triple, hence by Lemma 2.21(2) the triple $(e_3(h), e_{n+3}(h),
e_4(t))$ is also positive, so the head of $e_{n+3}$ lies on
the inner circle in Figure 2.3.  Note that $e_{n+2}$ shields this edge
endpoint from the outside circle of the annulus in Figure 2.3, hence
the tail of $e_{n+3}$ also lies in the inner circle in the figure,
therefore $e_{n+3}$ is parallel to $e_3$ on $\gb$, which is a
contradiction as they cannot be parallel on both graphs.

Now consider the case $n=3$.  By Lemma 2.22(3) we may assume $a_i \leq
6$, and $a_1 = n+3 = 6$.  By Lemma 12.8 we may assume that $\Delta=5$
and the jumping number $J \neq \pm 1$.  Also, $a_j \neq 5$, otherwise
by Lemma 12.7(1) the 11 edges in $\hat e_1 \cup \hat e_j$ would be
mutually non-parallel on $\gb$, contradicting the fact that $\rgb$ has
at most $3n = 9$ edges (Lemma 2.5).  One can now check that the
following are the only possible values of $(a_1, a_2, a_3, a_4)$ up to
symmetry, where $*$ indicates any possible value.  Let $s=s_1$.  Then
the other $s_i$ can be calculated using Lemma 12.1.  The second
quadruple in the following list indicates the values of $(s_1, s_2,
s_3, s_4)$. 
$$
\begin{array}{lll}
(1) & \qquad (6,6,*,*) & \qquad  (s,\; s,\; *,\; *) \\
(2) & \qquad (6,1,6,2) & \qquad  (s,\; s-1,\; s+1,\; s-1) \\
(3) & \qquad (6,1,4,4) & \qquad  (s,\; s-1,\; s,  \; s+1) \\
(4) & \qquad (6,4,3,2) & \qquad  (s,\; s-1,\; s+1,\; s-1) \\
(5) & \qquad (6,4,1,4) & \qquad  (s,\; s-1,\; s,\; s+1) \\
(6) & \qquad (6,3,2,4) & \qquad  (s,\; s,\; s+1,\; s+1 ) \\
(7) & \qquad (6,3,4,2) & \qquad  (s,\; s,\; s-1,\; s-1) \\
(8) & \qquad (6,3,3,3) & \qquad (s,\; s,\; s,\; s) \end{array}
$$

In case (1) by Lemma 12.7(4) the 12 edges in $\hat e_1 \cup \hat e_2$
are mutually non-parallel on $\gb$, which is impossible because $\rgb$
contains at most $3n = 9$ edges.  Case (8) is impossible by Lemma
12.6.  Also, Lemma 2.14(3) implies that $s_i \not \equiv 0$ mod $3$
if $a_i \geq 2$, which can be applied to exclude cases (2), (4) and
(5).

In case (6) and (7), by Lemma 12.7(1) $\hat e_2$ is not conjugate to
$\hat e_1, \hat e_3$ or $\hat e_4$, so $\hat e_1 \cup \hat e_2$
represents all $9$ edges in $\rgb$, hence each of $\hat e_3$ and $\hat
e_4$ must be conjugate to $\hat e_1$.  By Lemma 12.7(2) the two middle
edges of $\hat e_1$ are parallel to the middle edges in each of $\hat
e_3$ and $\hat e_4$, so $\hat e_3, \hat e_4$ are conjugate.  Since
$a_3+a_4 = 6 = 2n$ and $J\neq \pm 1$, this is a contradiction to Lemma
12.7(4).

In case (3), by Lemma 2.14(3) we have $s=1$.  Since $a_1 + a_3 = a_1 +
a_4 = 10$ while $\rgb$ has at most 9 edges, each of $\hat e_3, \hat
e_4$ must have an edge parallel to some edge of $\hat e_1$ on $\gb$.
By Lemma 12.7(2) this implies that each edge of $\hat e_3 \cup \hat
e_4$ is parallel to one of the 4 middle edges in $\hat e_1$.  Note
that the edge $e'$ in $\hat e_2$ is a loop based at $v_1$ in $\gb$,
which cannot be parallel to any other edge on $\gb$.  Therefore $\gb$
has exactly 7 families.  Moreover, if we let $e_1, e_6$ be the first
and last edges in $\hat e_1$ then each of $e_1, e_6, e'$ forms a
single family.

Now consider the graph in Figure 2.4.  Clearly there is only one
possible position for $e'$, which has exactly one endpoint on the
corner from the tail of $e_1$ to the head of $e_6$.  By the above
there are no other edge endpoints in this corner, which is a
contradiction because the label of the tail of $e_1$ is $1$ while the
label of the head of $e_6$ is $2$, so the number of edge endpoints
between them must be even.  
\qed

\begin{lemma} Suppose $n\geq 4$.  Then $\Delta = 4$ and
$\ga$ has at most $n+1$ parallel edges.  
\end{lemma}

\proof We need only show that $\Delta=4$.  The second statement will
then follow from Lemma 12.8.

Suppose to the contrary that $\Delta = 5$.  First assume that
$a_i < n+2$ for all $i$.  Then $\Delta n = 5n \leq 4(n+1)$, so $n=4$,
and $\ga = (5,5,5,5)$.  By Lemma 2.3(1) $s_1$ is coprime with $n=4$,
so we may assume without loss of generality that $s_1 = 1$.  Thus the
label sequences of $\hat e_1$ are $(1,2,3,4,1)$ at $u_1$ and
$(2,3,4,1,2)$ at $u_2$.  One can check that the label sequences of
$\hat e_3$ are $(3,4,1,2,3)$ at $u_1$, and $(4,1,2,3,4)$ at $u_2$.
This contradicts Lemma 12.5 with $e_1, e_2$ the two $12$-edges in
$\hat e_1$, $e_3$ the $12$-edge in $\hat e_3$, and $e'_1, e'_2$ the
two $34$-edges in $\hat e_3$.

We may now assume without loss of generality that $a_1 > n+1$.  By
Lemma 12.9 we must have $a_1 = n+2$.  By Lemma 12.4(1) we have $s_1 =
\pm 1$.

\medskip
CLAIM 1. {\it $a_2, a_4 \leq n+1$}.
\medskip

\proof  Suppose $a_2 = n+2$.  Then $s_2 = \pm 1$ by Lemma 12.4(1).  Also
by Lemma 12.1 we have 
$$ s_1 - s_2 \equiv a_1 + a_2 \equiv 4 \qquad \text{mod $n$}
$$
Hence either $n=4$ and $s_1 = s_2$ ($=1$ say), or $n=6$, $s_1 = -1$
and $s_2 = 1$.  In either case one can check that there is a pair of
parallel $12$-edges $e_1, e_2$ in $\hat e_1$, a $12$-edge $e_3$ in
$\hat e_2$ which is not equidistant to $e_1, e_2$, and a pair of
parallel $34$-edges $e'_1, e'_2$ in $\hat e_2$.  This is a
contradiction to Lemma 12.5.  

Hence $a_2 \leq n+1$.  A symmetric argument shows that $a_4 \leq n+1$.
\qed

We now have $5n \leq 2(n+2) + 2(n+1)$, giving $n\leq 6$.  Also if
$n=6$ then $\ga = (8,7,8,7)$.

\medskip
CLAIM 2. {$n = 5$.}
\medskip

\proof Otherwise we have $n=4$ or $6$.  If $a_2 = n+1$ then Lemma 12.1
gives $s_2 = s_1 - a_1 - a_2 = \pm 1 - (n+1) - (n+2) \equiv 0$ mod
$2$, which is a contradiction to the fact that the transition function
of a family of more than $n$ edges must be transitive (Lemma 2.3(1)).
Therefore we have $a_2 \leq n$.  Similarly for $a_4$.  This rules out
the case $n=6$.

When $n=4$ we must have $\ga = (6,4,6,4)$.  Assume without loss of
generality that $s_1 = 1$.  We now apply Lemma 12.5 with $e_1, e_2$
the $12$-edges in $\hat e_1$, $e_3$ the $12$-edge in $\hat e_3$, and
$e'_1, e'_2$ the $34$-edges in $\hat e_3$.
\qed

\medskip
CLAIM 3. {If $n = 5$ then $a_3 \neq 7$.}
\medskip

\proof Otherwise by Claim 1 we have $(a_2, a_4) = (6,5)$ or $(5,6)$,
so $a_4 - a_2 = \pm 1$.  We may assume that $s_1 = 1$.  By Lemma 12.1
we have
$$ s_3 \equiv s_1 - a_2 + a_4 = 1 \mp 1 = \text{ $0$ or $2$ mod $5$}
$$
which contradicts the fact that $s_3 = \pm 1$ mod $n$ (Lemma 12.4(1)).
\qed

The only possibility left is that $n = 5$ and $\ga = (7,6,6,6)$.
We may assume $s_1 = 1$.  Then this can be ruled out by applying Lemma
12.5 with $e_1, e_2$ the $12$-edges in $\hat e_1$, $e_3$ the $12$-edge
in $\hat e_3$, and $e'_1, e'_2$ the $45$-edges in $\hat e_3$.  
\qed

\begin{lemma}
(a) $\ga$ has at most $n+1$ parallel edges.  

(b) $\Delta = 4$.
\end{lemma}

\proof (a) This follows from Lemmas 12.8 and 12.10 if either $J=\pm
1$, or $\Delta=4$, or $n\geq 4$.  Hence we may assume that $\Delta=5$,
$J\neq \pm 1$, and $n = 3$.  By Lemma 12.9 we have $a_i \leq n+2 = 5$.
Thus the possible values of $(a_1, a_2, a_3, a_4)$ are given below.
The second quadruple gives $(s_1, s_2, s_3, s_4)$, calculated as
functions of $s = s_1$, using Lemma 12.1.  
$$
\begin{array}{lll}
(1) & \qquad (5,5,5,0) & \qquad  (s,\; s-1,\; s+1,\; - ) \\
(2) & \qquad (5,5,4,1) & \qquad  (s,\; s-1,\; s-1,\; s  ) \\
(3) & \qquad (5,3,2,5) & \qquad  (s,\; s+1,\; s-1,\; s+1) \\
(4) & \qquad (5,4,5,1) & \qquad  (s,\; s,  \; s,  \; s  ) \\
(5) & \qquad (5,4,4,2) & \qquad  (s,\; s,  \; s+1,\; s+1) \\
(6) & \qquad (5,4,3,3) & \qquad  (s,\; s,  \; s-1,\; s-1) \\
(7) & \qquad (5,4,2,4) & \qquad  (s,\; s,  \; s  ,\; s  ) \\
(8) & \qquad (5,3,5,2) & \qquad  (s,\; s+1,\; s-1,\; s+1) \\
(9) & \qquad (5,3,4,3) & \qquad  (s,\; s+1,\; s , \; s-1) 
\end{array}
$$

Cases (1), (3), (8) and (9) are impossible because there is an $i$ such
that $a_i\geq 2$ and $s_i = 0$, contradicting Lemma 2.14(3).  Cases (4)
and (7) contradict Lemma 12.6.

In case (2), by Lemma 12.7(1) the edges in $\hat e_3$ are not parallel
to those in $\hat e_1 \cup \hat e_4$ on $\gb$, and by Lemma 12.7(4)
the edge in $\hat e_4$ is not parallel to those in $\hat e_1$.  Thus
the 10 edges in $\hat e_1 \cup \hat e_3 \cup \hat e_4$ are mutually
non-parallel on $\gb$, contradicting the fact that $\rgb$ has at most
$3n$ edges (Lemma 2.5).  Similarly, in case (5) the edges in $\hat e_1
\cup \hat e_3 \cup \hat e_4$ are mutually non-parallel on $\gb$, and
in case (6) the edges in $\hat e_2 \cup \hat e_3 \cup \hat e_4$ are
mutually non-parallel on $\gb$, which lead to the same contradiction.

(b) Assume $\Delta=5$.  By Lemma 12.10 we have $n=3$, and by (a) we
have $a_i \leq 4$, hence the weights of $\hat e_i$ must be $(4,4,4,3)$
up to symmetry, and the transition numbers are $(s, s+1, s-1, s+1)$.
This is a contradiction to Lemma 2.14(2) because one of the families
has $s_i = 0$ and hence is a set of co-loops.  \qed

\begin{lemma} Suppose $\Delta = 4$.  Let $e, e'$ be edges of
$\ga$ with label $i$ at vertex $u_1$ and $j$ at $u_2$, $i\neq j$,
where the $i$-labels of $e, e'$ at $u_1$ are not adjacent among all
$i$-labels at $u_1$.  Suppose also that $e$ belongs to a family of at
least $n$ parallel edges of $\ga$.  Then $e\cup e'$ forms an essential
loop on the torus $\hat F_b$.  \end{lemma}

\proof Note that in this case the jumping number $J=\pm 1$, so the
assumption that the $i$-labels of $e, e'$ at $u_1$ are not adjacent
implies that the $1$-labels of $e_1, e'_1$ at the vertex $v_i$ in
$\gb$ are not adjacent among all $1$-labels.  By assumption $e$
belongs to a family of at least $n$ parallel edges of $\ga$, so if
$e\cup e'$ is inessential on $\hat F_b$ then by Lemma 12.2 $e_1$ and
$e'_1$ are parallel on $\gb$, which gives rise to at least $5$
parallel edges in $\gb$, contradicting Lemma 2.2(2) because $\rga$ has
at most $4$ edges.  \qed

Up to symmetry we may assume that $a_1 \geq a_3$, $a_2 \geq a_4$, and
$a_1 + a_3 \geq a_2 + a_4$.  Since $a_i \leq n+1$, the possibilities
for $\ga$ are listed below.  The second quadruple indicates the values
of $s_i$, calculated in terms of $s=s_1$ using Lemma 12.1.
$$
\begin{array}{lll}
(1) & \quad (n+1,\; n+1,\; n+1,\; n-3) & \qquad  (s,\;s-2,\;s-4,\;s-2) \\
(2) & \quad (n+1,\; n+1,\; n,\; n-2)   & \qquad  (s,\; s-2,\; s-3,\; s-1) \\
(3) & \quad (n+1,\; n+1,\; n-1,\; n-1) & \qquad  (s,\; s-2,\; s-2,\; s) \\
(4) & \quad (n+1,\; n,\; n+1,\; n-2)   & \qquad  (s,\; s-1,\; s-2,\; s-1) \\
(5) & \quad (n+1,\; n,\; n,\; n-1)     & \qquad  (s,\; s-1,\; s-1,\; s) \\
(6) & \quad (n+1,\; n,\; n-1,\; n)     & \qquad  (s,\; s-1,\; s,\; s+1) \\
(7) & \quad (n+1,\; n-1,\; n+1,\; n-1) & \qquad  (s,\;s,\;s,\;s) \\
(8) & \quad (n,\; n,\; n,\; n) & \qquad (s,\;s,\;s,\;s) 
\end{array}
$$

\begin{lemma}  Cases (4), (5), (6), (7), (8) are impossible.
\end{lemma}

\proof In case (4) $\hat e_1, \hat e_3$ are not A-conjugate to $\hat
e_2, \hat e_4$ by Lemma 12.7(1).  Since $a_4 < n$, by Lemma 12.6(1)
$\hat e_2, \hat e_4$ cannot be the only A-conjugate pair, hence $\hat
e_1$ must be A-conjugate to $\hat e_3$.  Since they have the same
number of edges, by Lemma 17.2 the first edge $e$ of $\hat e_1$ is
parallel to the first edge $e'$ of $\hat e_3$.  Since $e, e'$ have
labels $1,2$ at $u_1$, respectively, the label of $e$ at $u_2$ is $2$,
hence $s=1$.  Now $\hat e_2$ is a family of at least 3 co-loops,
contradicting Lemma 2.14(2).

In case (5), by Lemma 12.7(1) $\hat e_1$ can only be conjugate to
$\hat e_4$ and $\hat e_2$ to $\hat e_3$, but since $a_4<n$, by Lemma
12.6(1) $\hat e_2$ must be A-conjugate to $\hat e_3$.  Since $a_2+a_3
= 2n$, this is a contradiction to Lemma 12.7(3).

For the same reason, in case (6) $\hat e_2$ must be A-conjugate to
$\hat e_4$.  By Lemma 12.7(2) the first edge $e_1$ of $\hat e_2$ must
be parallel to the first edge $e'_1$ of $\hat e_4$. Examining the
labeling we see that they have labels $2$ and $1$ at $u_1$,
respectively, so the label of $e_1$ at $u_2$ is $1$, hence $s_2 =
-1$.  It follows that $s_1=s=0$, which is a contradiction to Lemma
2.3(1).  

Cases (7) and (8) are impossible by Lemma 12.6.
\qed

\begin{lemma}  Case (1) is impossible.
\end{lemma}

\proof  
Since $a_4<n$, by Lemma 12.6(1) two of the first three families are
$A$-conjugate.  Up to symmetry we may assume that $\hat e_1$ is
$A$-conjugate to $\hat e_2$ or $\hat e_3$.  By Lemma 12.7(2) the first
edges of the above conjugate pair must have the same label pair.
Examining the labels of these edges on $u_1$ we see that $s=1$ if
$\hat e_1$ is $A$-conjugate to $\hat e_2$, and $s=2$ if $\hat e_1$ is
$A$-conjugate to $\hat e_3$.  The second case cannot happen because
then $\hat e_2$ would be a set of at least 3 co-loops, contradicting
Lemma 2.14(2).

The graph $\ga$ is now shown in Figure 12.1.  If $n\geq 5$ then there
are bigons $e_1 \cup e_2$ in $\hat e_2$ and $e'_1, e'_2$ in $\hat e_4$
with labels $4,5$ at $u_1$ and $5,6$ at $u_2$ ($6=1$ when $n=5$).
Note also that there is a pair of parallel $23$-edges $e_3\cup e'_3$
in $\hat e_2$.  By Lemma 12.3(ii), these conditions imply that $e_1
\cup e'_1$ is inessential on $\hat F_b$, which contradicts Lemma
12.12.

When $n=4$, there is a pair of $14$-edges $e_1, e_2$ in $\hat e_1$,
a $14$-edge $e_3$ in $\hat e_4$, and a pair of $23$-edges in $\hat
e_2$.  Note that $e_3$ is not equidistant to $e_1, e_2$.  This leads
to a contradiction to Lemma 12.5.

When $n=3$, $s_i=0$ for some $i=1,2,3$, so one of the first three
families contains 4 co-loop edges, which is a contradiction
to the 3-Cycle Lemma.  
\qed

\bigskip
\leavevmode

\centerline{\epsfbox{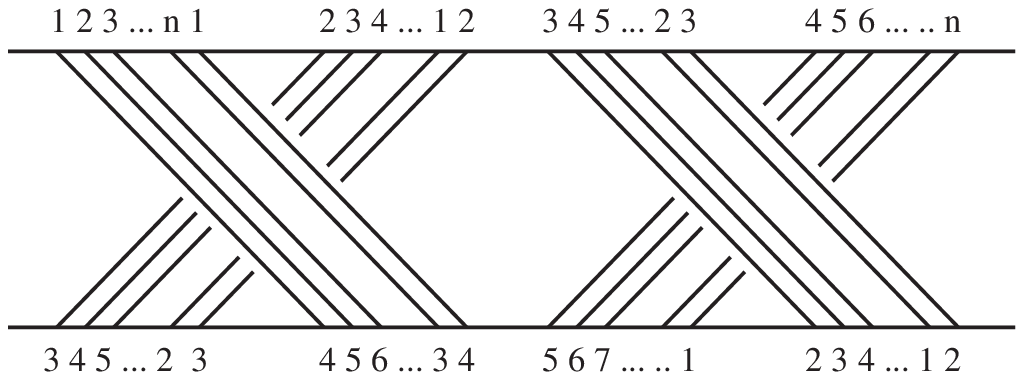}}
\bigskip
\centerline{Figure 12.1}
\bigskip

\begin{lemma}  Case (2) is impossible.
\end{lemma}

\proof In this case $a_1 \equiv a_2 \not \equiv a_3 \equiv a_4$ mod
$2$, so by Lemma 12.7(1) no edge in $\hat e_1 \cup \hat e_2$ is 
parallel to an edge in $\hat e_3 \cup \hat e_4$.  Since $\gb$ contains
at least $n$ bigons while $\hat e_3 \cup \hat e_4$ contributes
at most $a_4 = n-2$ bigons on $\gb$, it follows that some edge in
$\hat e_2$ is parallel to an edge in $\hat e_1$ on $\gb$.  Since $a_1
= a_2$, by Lemma 12.7(2) this implies that the first edge $e_1$ of
$\hat e_1$ is parallel to the first edge $e'_1$ of $\hat e_2$ on
$\gb$.  In particular, they must have the same label pair.  Since
$e_1$ has label $1$ at $u_1$ and $e'_1$ has label $2$ at $u_1$, we see
that $e_1$ has label $2$ at $u_2$, hence $s=1$.  Since $s_4 = s-1 =
0$, this is a contradiction to Lemma 2.14(3) unless $a_4 = n-2 <
2$, i.e.\ $n \leq 3$.  

Now suppose $n=3$.  Let $e_1, e_2$ be the two $12$-edges in $\hat
e_1$.  Note that there is a $12$-edge $e_3$ in $\hat e_3$, which by
the above is not parallel to any edge in $\hat e_1$, hence $e_1, e_2,
e_3$ cut $\hat F_b$ into a disk.  Now $\hat e_4$ is a loop based at
$v_3$ in $\gb$, so it must be a trivial loop.  This is a contradiction
because $\gb$ contains no trivial loop.
\qed

\begin{lemma}  Case (3) is impossible.
\end{lemma}

\proof We claim that $s=1$.  By Lemma 12.6(1) one of $\hat e_1, \hat
e_2$ is A-conjugate to some other $\hat e_j$.  Because of symmetry we
may assume that $\hat e_1$ is conjugate to some $\hat e_j$.  If $j=2$
then by Lemma 12.7(2) the first edge $e_1$ of $\hat e_1$ is parallel
on $\gb$ to the first edge $e'_1$ of $\hat e_2$, which has label $2$
at $u_1$, hence $e_1$ has label $2$ at $u_2$, so $s = s_1 = 1$.
Similarly, if $j=3$ then by Lemma 12.7(2), $a_1=n+1$ and $a_3=n-1$
implies that the second edge $e_2$ of $\hat e_1$ is parallel on $\gb$
to the first edge of $\hat e_3$, which has label $3$ at $u_1$, hence
$e_2$ has label pair $(23)$, which again implies that $s=1$.  The case
$j=4$ is impossible by Lemma 12.6(2).  The graph $\ga$ is now shown in
Figure 12.2.

\bigskip
\leavevmode

\centerline{\epsfbox{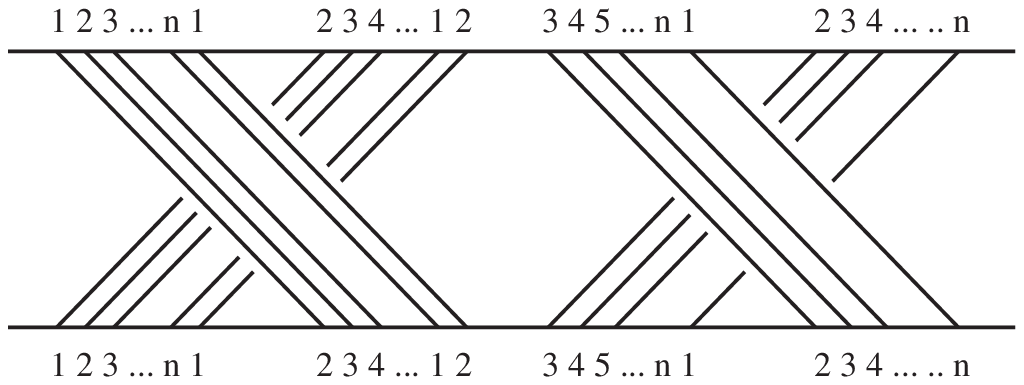}}
\bigskip
\centerline{Figure 12.2}
\bigskip

There are four edges $e'_1, ..., e'_4$ with label pair $(2,3)$, where
$e'_i \in \hat e_i$.  One can check on Figure 12.2 that they are all
equidistant to each other.  We claim that they are all parallel in
$\gb$.

The first $n$ edges of $\hat e_1$ form a loop $C$ on $\hat F_b$.  Let
$a_1, a_2, a_3$ be the first three edges of $\hat e_1$, oriented from
$u_1$ to $u_2$, and let $a_i(t), a_i(h)$ be the tail and head of
$a_i$, respectively.  Then as in the proof of Lemma 2.22(1), one can
show that $d_{v_2}(a_1(h), a_2(t)) = d_{v_3}(a_2(h), a_3(t))$.  In
other words, the corners at $v_2, v_3$ on one side of the above loop
contain the same number of edge endpoints.  Since $e'_2$ is
equidistant to $e'_1 = a_2$, we have $d_{v_2}(a_2(t), e'_2(h)) =
d_{v_3}(e'_2(t), a_2(h))$, hence the two endpoints of $e'_2$ lie on
the same side of the loop $C$.  It follows that $e'_2$ is parallel to
$e'_1$.  Similarly, $e'_3, e'_4$ are also parallel to $e'_1$.  This
proves the claim above.

Among the four parallel edges $e'_1,..., e'_4$, at least one of $e'_3,
e'_4$ is adjacent to $e'_1$ or $e'_2$ on $\gb$.  Because of symmetry
we may assume without loss of generality that $e'_3$ is adjacent to
$e'_1$ or $e'_2$.  Relabel it as $e_3$.

Note that $e_3$ is a border edge.  There is a face $D$ of $\ga$ with
$\bdd D = e_1 \cup e_2 \cup e_3 \cup e_4$, see Figure 12.3.  Let
$\alpha$ be the arc in $D$ connecting the middle points of $e_2, e_4$.
Since $e_3 = e'_3$ is parallel and adjacent to $e'_1$ or $e'_2$ and
$e'_1, e'_2$ are non-border edges in $\ga$, the face $D$ has a bigon
as a coupling face.  (See Section 2 for definition.)  It follows from
Lemma 2.15 that the surface $F_a$ can be isotoped rel $\bdd$ so that
the new intersection graph $\ga'$ is obtained from $\ga$ by deleting
$e_2, e_4$ and replacing them with two edges parallel to $e_1, e_3$
respectively.  The first family of $\ga'$ has $n+2$ edges, which is a
contradiction to Lemma 12.8.  \qed

\bigskip
\leavevmode

\centerline{\epsfbox{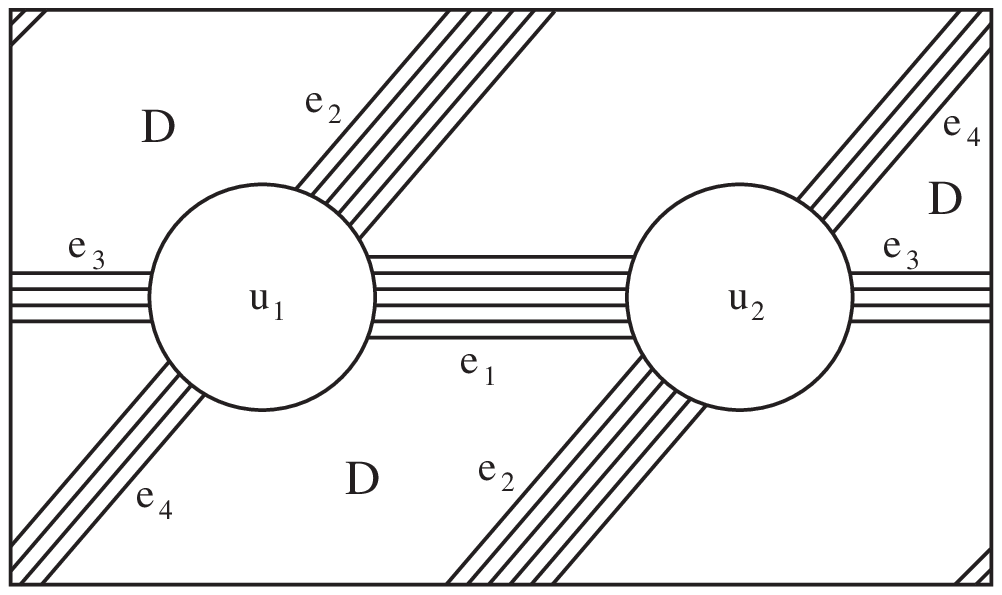}}
\bigskip
\centerline{Figure 12.3}
\bigskip

\begin{prop} The case that $n_a=2$, $n=n_b \geq 3$ and 
$\gb$ positive, is impossible.
\end{prop}

\proof 
We have shown that $\ga$ has 8 possibilities.  These have 
been ruled out in Lemmas 12.13 -- 12.16.
\qed

\section {The case $n_a = 2$, $n_b > 4$, $\Gamma_1, \Gamma_2$
non-positive, and $\text{max}(w_1 + w_2,\,\, w_3 + w_4) = 2n_b-2$
}

Suppose $n_a \leq 2$ and $n=n_b \geq 4$.  In Section 12 it has been
shown that $\gb$ cannot be positive.  In sections 13--16 we will
discuss the case that $\gb$ is non-positive.  The result will be given
in Propositions 14.7 and 16.8.

As before, we will use $n$ to denote $n_b$.

\begin{lemma} Suppose $n_a = 2$, $n \geq 4$, and $\Gamma_1,
\Gamma_2$ are non-positive.

(1) The reduced graph $\rga$ is a subgraph of the graph shown in
Figure 13.1.

(2) Let $w_i$ be the weight of $\hat e_i$.  Then up to relabeling we
may assume $w_3 + w_4 \leq w_1 + w_2$, and $w_1 + w_2 = 2n -2$ or
$2n$.  \end{lemma}

\bigskip
\leavevmode

\centerline{\epsfbox{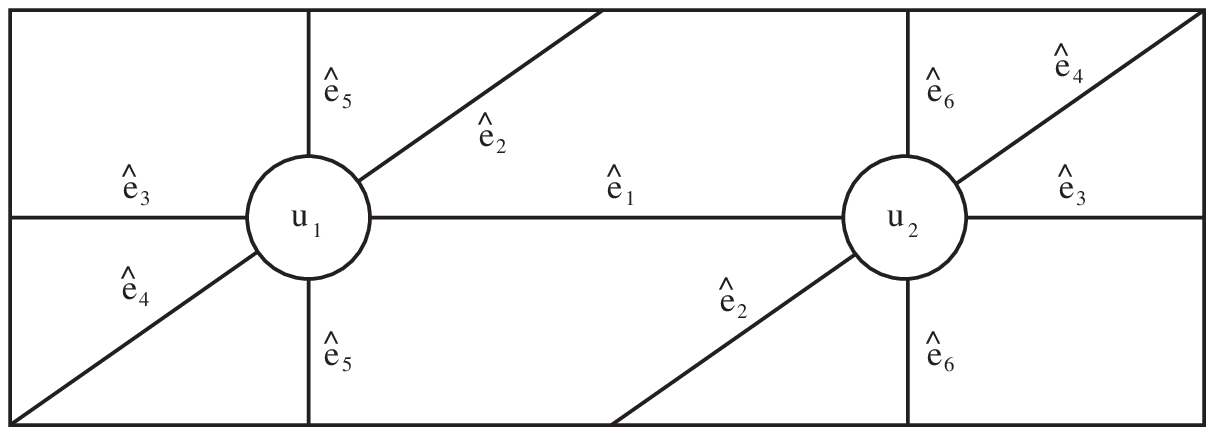}}
\bigskip
\centerline{Figure 13.1}
\bigskip

\proof (1) First note that the number of loops in $\ga$ at $u_1$ is the
same as that at $u_2$, because they have the same valence and the same
number of non-loop edges.  If $\rga$ has two loops based at $u_1$ then
they cut the torus into a disk, so there is no loop at $u_2$, which
would be a contradiction to the above.  Therefore $\rga$ has at most
one loop edge at each vertex.  If there is no loop at $u_i$ then
$\rga$ has at most four edges connecting $u_1$ to $u_2$.  If there is
one loop of $\rga$ at each $u_i$ then these cut the torus into two
annuli, each containing at most two edges of $\rga$.  In either case
$\rga$ is a subgraph of that in Figure 13.1.

(2) Up to relabeling we may assume that $w_1 + w_2 \geq w_3 + w_4$.
Since $\gb$ is non-positive, by Lemma 2.3(1) we have $w_i \leq n$,
hence $w_1 + w_2 \leq 2n$.

Assume $w_1 + w_2 = 2n - k$ and $k>0$.  Then $w_5 \geq (4n - (w_1 +
... + w_4))/2 \geq k$.  Note that the set of $k$ edges $e_1 \cup ...
\cup e_k$ of $\hat e_5$ adjacent to $\hat e_1 \cup \hat e_2$ has the
same set of labels at each of its two ends.  Hence by Lemma 2.4 $\hat
e_5$ contains a Scharlemann bigon, so by Lemma 2.2(4) and the parity
rule $k$ must be even.  If $k\geq 4$ then $e_1 \cup ... \cup e_k$
contains an extended Scharlemann bigon, which is a contradiction to
Lemma 2.2(6).  Hence $k=2$.  \qed

If $\gb$ has a Scharlemann cycle then by Lemma 2.2(4) the surface $\hat
F_a$ is separating, cutting $M(r_a)$ into a black region and a white
region.  Two Scharlemann cycles of $\gb$ have the same color if the
disks they bound lie in the same region.

\begin{lemma}  Suppose $n_a=2$ and $n \geq 1$.  

(1) If $e_1 \cup e_2$ and $e'_1 \cup e'_2$ are two Scharlemann bigons
of $\gb$ of the same color, then either (i) up to relabeling $e_i$ is
parallel to $e'_i$ on $\ga$ for $i=1,2$, or (ii) $\ga$ is kleinian,
and the four edges $e_1, e_2, e'_1, e'_2$ are mutually non-parallel on
$\ga$.

(2)  If $\gb$ has four parallel positive edges then $\ga$ is kleinian.
\end{lemma}

\proof (1) If the four edges are in two families of $\ga$ then (i)
holds.  If they are in three families, i.e., $e_1$ is parallel to
$e'_1$ but $e_2$ is not parallel to $e'_2$, then the nontrivial loop
$e_2 \cup e'_2$ on $\hat F_a$ is homotopic in $M(r_a)$ to the trivial
loop $e_1 \cup e'_1$, which contradicts the incompressibility of $\hat
F_a$.

Now assume that they are mutually non-parallel.  Let $G$ be the
subgraph of $\ga$ consisting of these four edges and the two vertices
of $\ga$.  Let $B$ be the side of $\hat F_a$ which contains the two
Scharlemann bigons.  Shrinking the Dehn filling solid torus $V_a$ to
its core $K_a$ and cutting $B$ along $\hat F_a$ and then along the two
Scharlemann bigons, we obtain a manifold whose boundary consists of
the two disk faces of $G$ and two copies of the two Scharlemann
bigons, which is a sphere, so by the irreducibility of $M(r_a)$ it
bounds a 3-ball.  It follows that $B$ is a twisted $I$-bundle over a
Klein bottle $K$, and $K$ intersects $K_a$ at a single point.
Therefore by Lemma 2.12 $\ga$ is kleinian.

(2) If $\gb$ has four parallel positive edges then they form two
Scharlemann bigons of the same color.  By Lemma 2.2(2) no two of these
edges are parallel on $\ga$, hence by (1) $\ga$ is kleinian.
\qed

In the remainder of this section we assume that $\ga, \gb$ are
non-positive, $n_a = 2$, $n = n_b > 4$, and $\text{max}(w_1 +
w_2,\,\, w_3 + w_4) = 2n-2$.  We may assume without loss of generality
that $w_3 + w_4 \leq w_1 + w_2 = 2n-2$.  Since $w_1 + ...  + w_4 +
2w_5 = \Delta n \geq 4n$, we have $w_5 = w_6 \geq 2$.  Let $\alpha_1
\cup \alpha_2$ (resp.\ $\beta_1 \cup \beta_2$) be the two edges of
$\hat e_5$ (resp.\ $\hat e_6$) adjacent to $\hat e_1 \cup \hat e_2$.
Note that these are Scharlemann bigons, hence $\hat F_b$ is
separating, and $n$ is even.  Without loss of generality we may assume
that $\alpha_1 \cup \alpha_2$ is a $(12)$-Scharlemann bigon.  Assume
that $\beta_1 \cup \beta_2$ is a $(k,k+1)$-Scharlemann bigon.

\begin{lemma}  (1) $k$ is even if and only if $w_2 = n-1$.

(2) $\{1, 2\} \cap \{k, k+1\} = \emptyset$.
\end{lemma}

\proof (1) This follows from the parity rule.  Orient $u_1$
counterclockwise and $u_2$ clockwise in Figure 13.1.  If $w_2 = n-1$
then the first edge of $\hat e_2$ has label $2$ at $u_1$ and $k+2$ at
$u_2$, so by the parity rule $k$ must be even.  Similarly if $w_2 = n$
or $n-2$ then $k$ is odd.

(2) If $k = 1$ then by (1) we have $w_2 = n$ or $n-2$.  In the first
case the edges of $\hat e_1$ would be co-loops, while in the second
case the edges of $\hat e_2$ would be co-loops.  If $k = 2$ then
$w_1=w_2 = n-1$ and the edges of $\hat e_1$ are co-loops.  Similarly
if $k=n$ then the edges of $\hat e_2$ are co-loops.  Since $n-2>2$,
all cases contradict Lemma 2.14(2) because the above would imply that
there are at least three parallel co-loop edges.  \qed

\begin{lemma} Suppose $w_1 = w_2 = n-1$.  Then for $i=1,2$, the
edges of $\hat e_i$ on $\gb$ form a cycle $C_i$ and a chain $C'_i$
disjoint from $C_i$.  Moreover, the vertices of $C_1$ ($C'_2$) are the
set of $v_j$ with $j$ odd, while the vertices of $C_2$ ($C'_1$) are the
set of $v_j$ with $j$ even.  The cycles $C_1, C_2$ are essential on
$\rgb$.  
\end{lemma}

\proof Let $\varphi_i$ be the transition function of $\hat e_i$.  Let
$h$ be the number of orbits of $\varphi_i$.  Since $\hat e_i$ has
$n-1$ edges, all but one component of the subgraph of $\gb$ consisting
of edges of $\hat e_i$ are cycles.  Therefore $h-1 \leq 2$ by Lemma
2.14(2).  Note also that each orbit contains the same number ($n/h$)
of vertices.  Since $\ga$ has a Scharlemann bigon, $\hat F_b$ is
separating and the number of positive vertices of $\gb$ is the same as
that of negative vertices, hence the number of orbits $h$ is even, so
we must have $h=2$.  Hence $\hat e_i$ forms exactly one cycle
component $C_i$ and one non-cycle component $C'_i$ on $\gb$.  Since
each odd number appears twice at the endpoints of $\hat e_1$, $C_1$
contains $v_j$ with $j$ odd, and $C'_1$ contains those with $j$ even.
For the same reason the edges of $\hat e_2$ form a cycle $C_2$ and a
chain $C'_2$.  Since $n/2$ edges of $\hat e_2$ have even labels at
$u_1$, $C_2$ must contain $v_j$ with $j$ even, while $C'_2$ contains
$v_j$ with $j$ odd.  It follows that $C_1 \cap C_2 = \emptyset$.
\qed

When $w_1 = n-2$ and $w_2 = n$, the edges of $\hat e_2$ form exactly
two cycles $C_1$ and $C_2$ on $\gb$, essential on $\hat F_b$, where
the vertices of $C_1$ ($C_2$) are the $v_j$ with $j$ odd (even).  This
is because by Lemma 2.14(2) they cannot form more than two cycles,
while $\gb$ being non-positive implies that $\hat e_2$ cannot form
only one cycle.  When $w_1 = w_2 = n-1$, let $C_1, C_2$ be the cycles
given in Lemma 13.4.  In either case, let $A_1, A_2$ be the annuli
obtained by cutting $\hat F_b$ along $C_1 \cup C_2$.  Consider the
cycles $\alpha = \alpha_1 \cup \alpha_2$ and $\beta = \beta_1 \cup
\beta_2$ on $\gb$.  Note that either $\alpha$ and $\beta$ are in
different $A_i$, or each of them has exactly one edge in each $A_i$.
We say that $\alpha, \beta$ are {\it transverse to\/} $C_i$ in the
second case.

\begin{lemma} The cycles $\alpha, \beta$ are disjoint, and
transverse to $C_i$.  \end{lemma}

\proof The first statement follows from Lemma 13.3(2), so we need only
show that $\alpha, \beta$ are transverse to $C_i$.

First assume $\Delta = 5$.  Then $w_5 = \frac 12 (\Delta n -
(w_1 + ... + w_4)) \geq \frac n2 + 2$.  By Lemma 2.3(3) we also have
$w_5 \leq \frac n2 + 2$, hence $w_5 = \frac n2 + 2$, in which case the
two outermost bigons of the family $\hat e_5$ are Scharlemann bigons,
with label pair $(12)$ and $(r+1, r+2)$, respectively, where $r=n/2$.
By Lemma 2.3(4) the label pair of $\beta_1 \cup \beta_2$ must be
either $(1,2)$ or $(r+1, r+2)$, and by Lemma 13.3 it cannot be the
former.  Therefore it must be $(r+1, r+2)$.

If $\alpha$ is not transverse to $C_i$, then it is an essential cycle
in one of the annuli, say $A_1$, obtained by cutting $\hat F_b$ along
$C_1 \cup C_2$, so $\beta$ must be an essential cycle in the other
annulus $A_2$.  The two cycles $\alpha$ and $\beta$ separate the
vertices of $C_1$ from $C_2$, except $v_1, v_2, v_{r+1}$ and $v_{r+2}$
which lie on $\alpha \cup \beta$.  On the other hand, the edge $e$ in
$\hat e_5$ adjacent to $\alpha_2$ has label pair $(3,n)$, so there is
an edge on $\gb$ connecting $v_3$ to $v_n$.  Since $n$ is even, the
vertices $v_3, v_n$ belong to different $C_i$, but since $n>4$,
neither $3$ nor $n$ belongs to the set $\{1,2,r+1, r+2\}$, which is a
contradiction.

Now assume $\Delta = 4$.  In this case the jumping number $J(r_a, r_b)
= \pm 1$.  Consider the two negative edges $e', e''$ of $\ga$ with
label $2$ at $u_1$.  Note that their endpoints at $u_1$ are separated
by the label 2 endpoints of $\alpha_1, \alpha_2$, hence by the Jumping
Lemma, on $\gb$ the endpoints of $e', e''$ at $v_2$ are separated by
those of $\alpha$; in other words, $e', e''$ are on different sides of
the cycle $\alpha$.  Assume that $v_2 \in C_2$ is positive.  If
$\alpha$ is not transverse to $C_2$ then all positive edges at $v_2$
must be on one side of $\alpha$ because the other side is shielded by
the cycle $C_1$, which contains only negative vertices.  This is a
contradiction.  Therefore $\alpha$, and hence $\beta$, must be
transverse to $C_i$.  \qed

\begin{lemma} Each edge of $\hat e_1 \cup \hat e_2$ is either
on $C_1 \cup C_2$ or parallel to an edge of $C_1 \cup C_2$ on $\gb$.
\end{lemma}

\proof Let $C_1, C_2$ and $\alpha, \beta$ be as above.  By definition
$C_2$ consists of the edges in $\hat e_2$ with even labels.  Let
$C'_1$ be the edges of $\hat e_1$ with even labels.  Because of
symmetry it suffices to show that each edge of $C'_1$ is parallel to
an edge in $C_2$.

Note that $\alpha \cap C_2 = v_2$.  Let $v_t = \beta \cap C_2$.  (Thus
$t$ is the even label of the Scharlemann bigon $\beta_1 \cup \beta_2$
in $\ga$.)  Since $w_1+w_2 = 2n-2$ and the edges adjacent to $\hat e_1
\cup \hat e_2$ on $\ga$ are the $(12)$-Scharlemann bigon $\alpha_1
\cup \alpha_2$ and the $(t,t+1)$- or $(t-1,t)$-Scharlemann bigon
$\beta_1 \cup \beta_2$, we see that $2, t$ are the only even labels
appear three times among the endpoints of edges in $\hat e_1 \cup \hat
e_2$, hence on $\gb$ the edges of $C'_1$ form a chain with endpoints
at $v_2, v_t$, and possibly some cycle components.  Therefore $C'_1 -
v_2 \cup v_t$ is disjoint from $\alpha \cup \beta \cup C_1$, hence
lies in the interior of the two disks obtained by cutting $\hat F_b$
along $\alpha \cup \beta \cup C_1$.  By Lemma 2.14(1) this implies
that $C'_1$ has no cycle component, and hence is a chain.  Since
$C'_1$ contains all vertices of $C_2$, this also implies that one
component of $C_2 - v_2 \cup v_t$ contains no vertices of $\gb$; in
other words, the two vertices $v_2, v_t$ are adjacent on $C_2$.

Let $q, p$ be the transition number of $\hat e_1, \hat e_2$,
respectively.  Since $C_2$ has an edge connecting $v_2$ to $v_t$, we
have $p \equiv \pm (t-2)$ mod $n$.  Since $C'_1$ is a chain of length
$(n/2)-1$ connecting $v_2, v_t$, we have $((n/2) - 1) q \equiv \pm
(t-2)$ mod $n$.  An edge of $C'_1$ has even labels on both endpoints,
so $q$ is even, hence $((n/2) -1)q \equiv -q$ mod $n$.  It follows
that $p \equiv \pm q$ mod $n$, which implies that each edge $e'$ of
$C'_1$ has its endpoints on adjacent vertices of $C_2$.  Let $e$ be
the edge of $C_2$ connecting these two vertices.  Since $e'$ has
interior disjoint from $\alpha \cup \beta \cup C_1$, it must be
parallel to $e$.  \qed

\begin{prop} The case that $\ga, \gb$ are non-positive,
$n_a = 2$, $n = n_b > 4$, and $w_3 + w_4 \leq w_1 + w_2 = 2n-2$, is
impossible.  \end{prop}

\proof First assume that $w_1 = n-2$ and $w_2 = n$.  By Lemma 13.3, the
label pair of $\beta$ is $(k, k+1)$, where $k$ is odd and $k\neq 1$.
If $k=n-1$ then the edges of $\hat e_2$ are co-loops, which
contradicts Lemma 2.14(2).  Therefore $n-3 \geq k \geq 3$.  

Since the label sequence of $\hat e_1$ at $u_2$ is $k+2, ..., n, 1,
..., k-1$, the above implies that there are adjacent edges $e'_1, e'_n
\in \hat e_1$ with labels $1$ and $n$ at $u_2$, respectively.  By
Lemma 13.6 each edge of $\hat e_1$ is parallel to some edge of $\hat
e_2$ on $\gb$, hence the transition function $\psi_1$ of $\hat e_1$ is
either equal to $\psi_2$ of $\hat e_2$, or $\psi_2^{-1}$, but since
the two edges of $\hat e_1 \cup \hat e_2$ with label $n$ at $u_1$ have
labels $k+1$ and $k-1$ respectively at $u_2$, the first case is
impossible, hence $\psi_1 = \psi_2^{-1}$.  Let $\hat e_2 = e_1 \cup
... \cup e_n$, where $e_i$ has label $i$ at $u_1$.  Since $e_1$ is the
only edge of $\hat e_2$ with label $1$ at $u_1$, it must be the one
that is parallel to $e'_1$ on $\gb$.  Similarly, $e_n$ is parallel to
$e'_n$ on $\gb$.  This is a contradiction to Lemma 2.19. 

Now assume that $w_1 = w_2 = n-1$.  As above, let $\hat e_2 = e_2 \cup
... \cup e_n$, where $e_i$ has label $i$ at $u_1$.  The label sequence
of $\hat e_1$ at $u_2$ is $k+1, k+2, ..., n, 1, ..., k-1$.  By Lemma
13.3, $k$ is even, and $\{1,2\} \cap \{k, k+1\} = \emptyset$, so $n-2
\geq k \geq 4$.  It follows that there are three consecutive edges
$e'_n, e'_1, e'_2$ of $\hat e_1$ such that $e'_i$ has label $i$ at
$u_2$.  For the same reason as above, $e'_2$ is parallel to $e_2$ and
$e'_n$ is parallel to $e_n$ on $\gb$.  Since the number of edges
between $e'_n$ and $e'_2$ is 1 while the number of edges between $e_2$
and $e_n$ is $n-3 > 1$ on $\ga$, this is a contradiction to Lemma
2.19.  \qed

\section {The case $n_a = 2$, $n_b > 4$, $\Gamma_1, \Gamma_2$
non-positive, and $w_1 = w_2 = n_b$  }

In this section we consider the case that $n_a = 2$, $n = n_b > 4$,
$\Gamma_1, \Gamma_2$ non-positive, and $w_1 = w_2 = n$.  We will also
assume without loss of generality that $w_3 \geq w_4$.  Let $\hat e_1
= e_1 \cup ... \cup e_{n}$, $\hat e_2 = e'_1 \cup ... \cup e'_{n}$,
and assume that $e_i, e'_i$ have label $i$ at $u_1$.

Let $r$ be such that the label of the endpoint of $e_1$ on $\bdd u_2$
is $r+1$.  One can check that both $e_i, e'_i$ have label $r+i$ at
$\bdd u_2$.

Since $\gb$ is non-positive, the vertices of $\gb$ cannot all be
parallel, so the edges of $\hat e_1$ form at least two cycles on $\gb$.
By Lemma 2.14(2) they form exactly two cycles $C_1 \cup C_2$ on $\gb$.

\begin{lemma} $\ga$ is not kleinian.  In particular, $\gb$
cannot contain four parallel positive edges.
\end{lemma}

\proof  If $\ga$ is kleinian then by Lemma 6.2(4) there is a free
orientation reversing involution $\phi$ of $(\hat F_a, \ga)$, which
maps $u_1$ to $u_2$, and is label preserving.  If there is no loop on
$\ga$ (i.e.\ $w_5 = w_6 = 0$), then $\Delta = 4$ and $w_i = n$ for all
$i$, so the label sequences of $\hat e_i$ at $u_1$ are all the same.
The above implies that the label sequences of $\hat e_i$ at $u_2$ are
also the same as those at $u_1$, so the transition function $\varphi$
defined by $\hat e_i$ is the identity map and hence all edges of
$\ga$ are co-loops, contradicting the 3-Cycle Lemma 2.14(2).

Now assume $w_5 = w_6 > 0$.  Then $\phi$ maps $\hat e_1 \cup\hat e_2$
to either $\hat e_1 \cup \hat e_2$ or $\hat e_3 \cup \hat e_4$.  In
the first case since $\phi$ is label preserving and orientation
reversing on the torus, the label sequence of $\hat e_1$ at $u_2$ is
the same as that of $\hat e_1$ at $u_1$, hence all edges of $\hat e_1$
are co-loops and we have a contradiction to Lemma 2.14(2).  In the
second case $w_3 = w_4 = w_1 = w_2 = n$, so $\Delta = 5$ and $w_5 =
w_6 = n/2$.  We have assumed that $\hat e_1$ has label sequence
$1,2,...,n$ at $u_1$.  so $\hat e_3$ has the same label sequence at
$u_2$.  Since $w_5 = n/2$, the label sequence of $\hat e_1$ at $u_2$
is $(k+1, k+2, ...n, 1, ...,k)$, where $k = n/2$.  Therefore $\phi$ is
of period 2, so it has $n/2 > 2$ orbits, which again contradicts Lemma
2.14(2).

The second statement follows from the above and Lemma 13.2(2).
\qed

\begin{lemma}  The edges $e_i, e'_i$ are parallel on $\gb$.
\end{lemma}

\proof The cycles $C_1 \cup C_2$ defined at the beginning of the
section cut the torus $\hat F_b$ into two annuli $A_1, A_2$.  Each
$e'_i$ lies in one of the $A_j$ and has the same endpoints as $e_i$,
so if it is not parallel to $e_i$ and $e_i \subset C_1$ then it is
parallel to $C_1 - e_i$.  There are at most two such $e'_i$ for $C_1$,
one in each $A_j$.  Since $n > 4$, $C_1$ contains at least three
edges, hence there exists some $e'_j$ parallel to $e_j \subset C_1$.

Assume $e_i$ is not parallel to $e'_i$, and let $e_j, e'_j$ be
parallel on $\hat F_b$, which exist by the above.  Let $D$ (resp.
$D'$) be the disk on $F_a$ that realizes the parallelism between $e_i,
e_j$ (resp.\ $e'_i, e'_j$), and let $D''$ be the disk between $e_j$
and $e'_j$ on $F_b$.  Shrinking $V_b$ to its core $K_b$, $B = D \cup
D' \cup D''$ becomes a disk in $M(r_b)$ with $\bdd B = e_i \cup
e'_i$, which contradicts the fact that $\hat F_b$ is incompressible in
$M(r_b)$.
\qed

\bigskip
\leavevmode

\centerline{\epsfbox{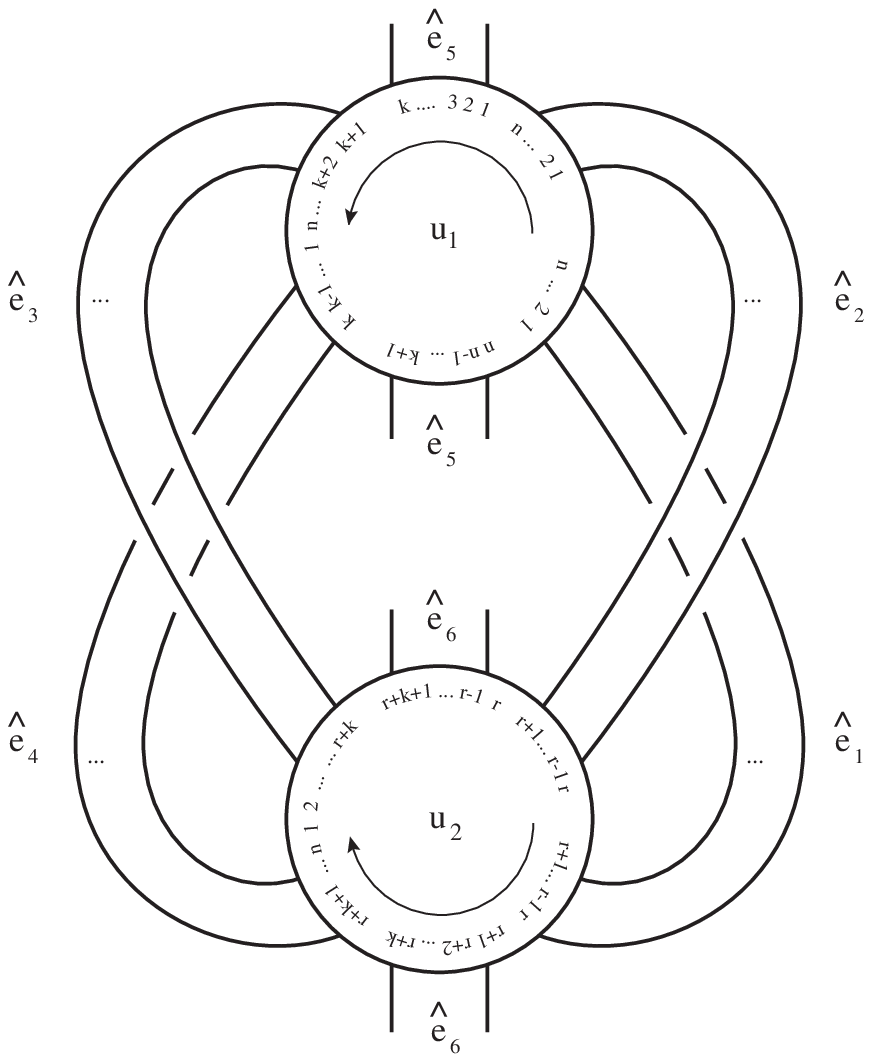}}
\bigskip
\centerline{Figure 14.1}
\bigskip

\begin{lemma} $w_3 + w_4 \neq n$.
\end{lemma}

\proof Assume to the contrary that $w_3 + w_4 = n$.  We have $\Delta =
4$ as otherwise there would be $n$ parallel positive edges in $\hat
e_5$, contradicting Lemma 2.3(3).  Now $w_5 = w_6 = n/2$, so the graph
$\ga$ is as shown in Figure 14.1, where $k=n/2$.  Since $\Delta = 4$, we
may assume that the jumping number is $1$.

Let $i$ be a label such that $1 \leq i \leq k$, so it appears on the
top of the vertex $u_1$ in Figure 14.1.  Consider the vertex $v_i$ of
$\gb$, see Figure 14.2.  By Lemma 14.2 $e_i$ of $\hat e_1$ is parallel
to $e'_i$ of $\hat e_2$.  Since $e_i$ and $e'_i$ have the same label
$1$ at $v_i$, there is an edge of $\hat e_3 \cup \hat e_4$ between
them.  Similarly there are parallel edges $e_j, e'_j$ with label $2$
at $v_i$, and there is another edge between them.  See Figure 14.2.
From the labeling we see that the two negative edges at $v_i$
(corresponding to loops in $\ga$) must be adjacent to each other.

\bigskip
\leavevmode

\centerline{\epsfbox{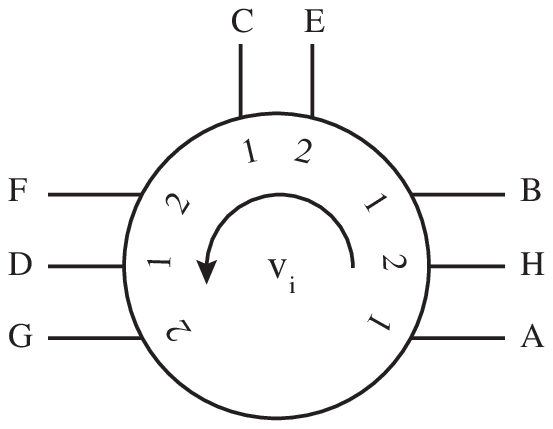}}
\bigskip
\centerline{Figure 14.2}
\bigskip

On $\bdd u_1$ the $i$-labels appear on endpoints of edges in the order
of $A, B, C, D$, where $A \in e_i$, $B \in e'_i$, $C$ is a loop in
$\hat e_5$, and $D \in \hat e_3 \cup \hat e_4$.  Since the jumping
number is $1$, the $1$-labels at $v_i$ also appear in the same order.
In Figure 14.2 this implies that the negative edge $C$ appears on the
top of the vertex.

Now consider the four edges labeled $2$ at $v_i$, denoted by $E, F, G,
H$, where $E$ is the negative edge, which is uniquely determined.
Since $F, G$ are parallel positive edges on $\gb$, they are the $e_j,
e'_j$ given above, belonging to $\hat e_1 \cup \hat e_2$.  On $\bdd
u_2$ this implies that the endpoint of the loop $E$ labeled $i$
appears on the top of $\bdd u_2$ in Figure 14.1.  Since $i$ is any
label between $1$ and $n/2$, it follows that the labels on the top of
$\bdd u_2$ must be $1,2,...,k$, so the integer $r$ in the figure
satisfies $r = k$.  However, in this case the edges of $\hat e_1$
would form cycles of length 2 in $\gb$, which is a contradiction to
Lemma 2.14(2).  \qed

\begin{lemma} $w_3 = n$, and $0 < w_4 < n$.  Moreover, an edge
$e''$ of $\hat e_3$ with label $j$ at $u_2$ is parallel to the edges
$e_j$ and $e'_j$ .
\end{lemma}

\proof We have assumed $w_4 \leq w_3\leq n$.  If $w_4 = n$ then the
argument of Lemma 14.2 applied to $\hat e_3, \hat e_4$ shows that each
edge of $\hat e_3$ is parallel to exactly one edge of $\hat e_4$.  On
the other hand, since the two parallel edges $e_i, e'_i$ in the proof
of Lemma 14.2 have the same label $1$ at the vertex $v_i$, there must
be another edge $e''_i$ in $\hat e_3 \cup \hat e_4$ between $e_i$ and
$e'_i$.  Together with the other edge in $\hat e_3 \cup \hat e_4$
which is parallel to $e''_i$, we get four parallel positive edges in
$\gb$, which contradicts Lemma 14.1.  Therefore $w_4 < n$.

Recall from Lemma 14.2 that the edges $e_i$ and $e'_i$ are parallel in
$\gb$, with the same label $1$ at $u_i$, so there must be another edge
$e''_i \in \hat e_3 \cup \hat e_4$ between them.  Note also that if
$e_i$ has label $i+r$ at $u_2$ then $e''_i$ has the property that it
has label $i$ at $u_2$ and $i+r$ at $u_1$.  This is true for all $i$,
so either $\hat e_3$ and $\hat e_4$ have the same transition function,
or these $e''_i$ all belong to the same family.  The first case
happens only if $w_3 + w_4 \equiv 0$ mod $n$, which is impossible
because by Lemma 14.3 we have $w_3 + w_4 \neq n$, while $w_3 \leq n$
and by the above we have $w_4 < n$.  Therefore all the $e''_i$ belong
to $\hat e_3$.  Since $w_3 \leq n$, this implies that $w_3 = n$.
Again by Lemma 14.3 we have $w_4 \neq 0$.  \qed

\bigskip
\leavevmode

\centerline{\epsfbox{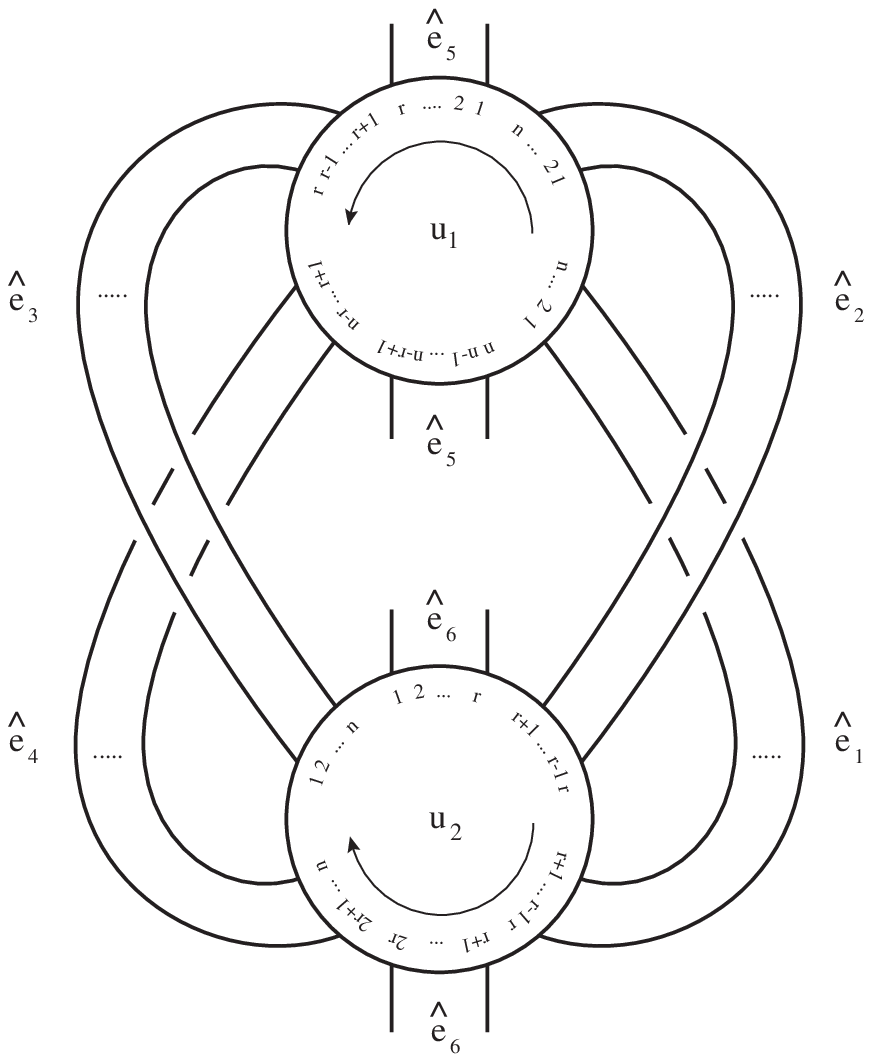}}
\bigskip
\centerline{Figure 14.3}
\bigskip

\begin{lemma}  
The label sequence of $\hat e_3$ is $1,2,...,n$ at $u_2$, and
$1+r, 2+r, ..., n, 1, ..., r$ at $u_1$.  The labels of $\ga$ are as
shown in Figure 14.3.
\end{lemma}

\proof First assume that the label sequence of $\hat e_3$ at $u_2$ is
not $1,2,..., n$.  Then there is a pair of adjacent parallel edges
$e''_{n}, e''_1$ with label $n$ and $1$ at $u_2$, respectively.  By
Lemma 14.4 $e_1 \cup e''_1$ and $e_n \cup e''_n$ are parallel pairs on
$\gb$.  Since $e''_1, e''_n$ are adjacent on $\ga$ while $e_1, e_n$
are not, this is a contradiction to Lemma 2.19.  Therefore the label
sequence of $\hat e_3$ at $u_2$ must be $1,2,..., n$.

Since by Lemma 14.4 the edge $e''_i$ connects $v_i$ to $v_{i+r}$ with
label $2$ at $v_i$ and $1$ at $v_{i+r}$, we see that on $\ga$ it has
label $i$ at $u_2$ and $i+r$ at $u_1$, hence the label sequence of
$\hat e_3$ at $u_1$ is $r+1, ..., n, 1, ..., r$.  The labels of $\hat
e_1, \hat e_2, \hat e_3$ determine those of the loops, and hence those
of $\hat e_4$.  Therefore $\ga$ must be as shown in Figure 14.3.
\qed

\begin{lemma} (1) The jumping number $J = \pm 1$.

(2) Orient the negative edges of $\ga$ from $u_1$ to
$u_2$.  Then on $\gb$ the edges of $\hat e_1$ form two essential
cycles of opposite orientation on $\hat F_b$.  \end{lemma}

\proof (1) Since $\Delta =
4$ or $5$, the jumping number is either $\pm 1$ or $\pm 2$.  Let $e_i,
e'_i$ be the edges of $\hat e_1, \hat e_2$, respectively, with label
$i$ at $u_1$.  If $J = \pm 2$ then these edges are not adjacent among
the $1$-edges at $v_i$ in $\gb$.  Since by Lemma 14.2 they are parallel
in $\gb$, there would be more than $2n_a = 4$ parallel edges in $\gb$,
which contradicts Lemma 2.2(2).  Therefore $J = \pm 1$.  Changing the
orientation of $\hat F_b$ if necessary, we may assume that $J = 1$.  

(2) Now let $C_1, C_2$ be the cycles on $\gb$ consisting of edges of
$\hat e_1$.  We need to show that they are of opposite orientation.

Let $a_1, a_2, \bar a, a_3$ be the edges with label $1$ at $u_1$,
where $a_i \in \hat e_i$ for $i=1,2,3$, and $\bar a \in \hat e_5$.
Note that they appear in this order on $\bdd u_1$.  
Since $J=1$, they also appear in this order on $\bdd v_1$ in $\gb$,
see Figure 14.4.  By the proof of Lemma 14.3 we see that $a_3$ is in
the middle of a pair of parallel positive edges incident to $v_1$,
which is not parallel to $a_1, a_2$, hence the orientation of $C_1$
must be as shown in Figure 14.4, where $C_1$ is represented by the
lower level chain.

\bigskip
\leavevmode

\centerline{\epsfbox{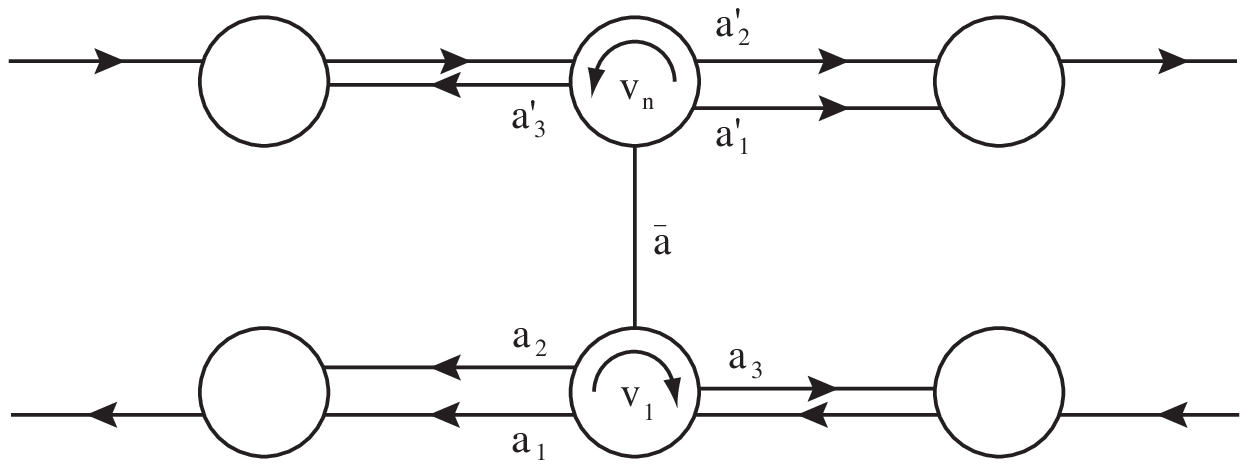}}
\bigskip
\centerline{Figure 14.4}
\bigskip

Now consider the edges labeled $n$ at $u_1$.  There are 5 of them if
$\Delta = 5$, but we only consider $\bar a$ and the edges $a'_1, a'_2,
a'_3$, where $a'_i \in \hat e_i$.  The order of the label $n$
endpoints of these edges on $\bdd u_1$ is $a'_3, \bar a, a'_1, a'_2$,
while the orientation of $v_n$ is opposite to that of $v_1$.
Therefore these edges appear on $\bdd v_n$ as shown in Figure 14.4.
We see that $C_1, C_2$ are of opposite orientation on $\hat F_b$.
\qed

\bigskip
\leavevmode

\centerline{\epsfbox{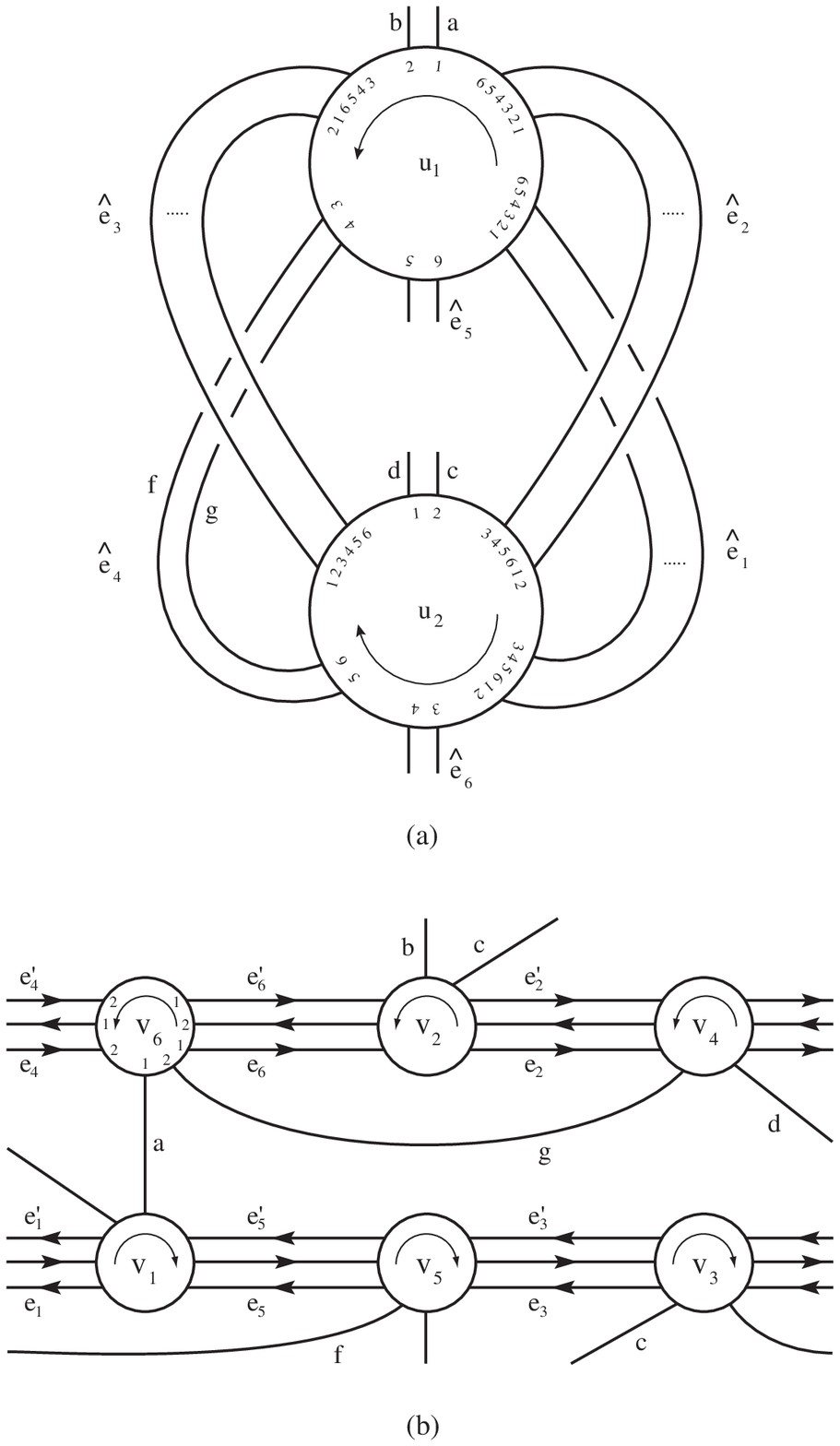}}
\bigskip
\centerline{Figure 14.5}
\bigskip

\begin{prop} Suppose $n_a = 2$, $n > 4$, $\ga, \gb$ are
non-positive, and $w_1 = w_2 = n$.  

(1) On $\gb$ each edge of $\hat e_4$ connects a pair of adjacent
    vertices of some $C_i$, but is not parallel to an edge of $C_i$.

(2) Two adjacent edges of $\hat e_4$ lie in different annuli of $\hat
    F_b - \cup C_i$.

(3) $w_4 = w_5 = w_6 = 2$, $\Delta = 4$, and $n = 6$.  

(4) The graphs $\ga, \gb$ and their edge correspondence are as shown
in Figure 14.5, where $e_i$ (resp.\ $e'_i$) is the edge in $\hat e_1$
(resp.\ $\hat e_2$) with label $i$ at $u_1$, and the edge between
$e_i, e'_i$ is the edge of $\hat e_3$ with label $i$ at $u_2$.
\end{prop}

\proof (1) From Figure 14.3 we see that an edge $e$ of $\hat e_4$ with
label $i$ at $u_1$ has label $i+r$ at $u_2$.  Since the transition
function of $\hat e_1$ also maps $i$ to $i+r$, $v_i$ and $v_{i+r}$ are
connected by the edge $e_i$ of $\hat e_1$, and hence are adjacent on
one of the cycles $C_j$.  This proves the first part of (1).  By
Lemmas 14.2 and 14.4 each edge $e_i$ of $C_j$ is parallel to an edge
$e'_i$ in $\hat e_2$ and an edge $e''_i$ in $\hat e_3$, so if $e$ is
parallel to $e_i$ then there would be four parallel positive edges in
$\gb$, which would contradict Lemma 14.1.

(2) By Lemma 14.6(2) the two cycles $C_1, C_2$ have opposite
orientations.  Without loss of generality we may assume that the
orientations of $C_i$ are as shown in Figure 14.5.  Recall that $C_1$
is the cycle containing the vertex $v_1$.  By Lemma 14.6(1) we may
assume without loss of generality that the jumping number of the
graphs is $1$.  Let $e$ be an edge of $\hat e_4$ with label $k$ at
$u_1$, and let $e_k, e'_k, e''_k$ be the edges of $\hat e_1, \hat e_2,
\hat e_3$ with label $k$ at $u_1$.  Then the endpoints of these edges
appear at $\bdd u_1$ in the order $e_k, e'_k, e''_k, e$, so on
$\bdd v_k$ they appear in the same order.  If $v_k$ is in $C_1$ then
the orientation of $C_1$ points to the left and the orientation of
$v_k$ is clockwise, so $e$ is in the annulus below $C_1$.  If $v_k$ is
in $C_2$ then the orientation of $C_1$ points to the right and the
orientation of $v_k$ is counterclockwise, so again $e$ is in the
annulus below $C_1$.  Since the labels of adjacent edges of $\hat e_4$
belong to different $C_i$ in $\gb$, it follows from the above that
they are in different annuli of $\hat F_b - \cup C_i$.

(3) Since each $C_i$ contains $n/2 > 2$ vertices, there cannot be two
edges on the same side of $C_i$ connecting two different pairs of
adjacent vertices and yet not parallel to an edge of $C_i$.  Hence by
(1) and (2) $\hat e_4$ contains at most two edges.  By Lemma 14.4 $w_4
> 0$, and from the labeling in Figure 14.3 we see that $w_4$ is even.
Therefore $w_4 = 2$.

If $\Delta = 5$ then the loop family of $\ga$ at $u_1$ contains $n-1$
edges.  This contradicts Lemma 2.3(3) for $n>4$.  Hence $\Delta =
4$.

Let $e$ be an edge of $\hat e_4$ with endpoints on $v_i$ and $v_j$ in
$C_2$, lying on the annulus $A$ below $C_2$.  By (1) it is not
parallel to the edge on $C_2$ connecting $v_i, v_j$, so on $A$ it
separates $C_1$ from other vertices of $C_2$, hence there is no edge
in $A$ connecting $C_1$ to vertices of $C_2$ except possibly $v_i$ and
$v_j$.  By Lemmas 14.2 and 14.4 there are three parallel edges for
each edge of $C_i$.  Together with $e$, they contribute $7$ edge
endpoints to each of $v_i$ and $v_j$, therefore $\Delta = 4$ implies
that there are at most two edges in $A$ connecting $C_1$ to $C_2$, one
for each of $v_i, v_j$.  Note that these correspond to loop edges in
$\ga$.  Therefore the two annuli give rise to at most 4 loops in
$\ga$, so $w_5 = w_6 \leq 2$.  Since $n>4$ and $2w_5 + w_4 = (\Delta -
3)n = n$, it follows that $n=6$, and $w_5=w_6 = 2$.

(4) By Lemma 14.5 $\ga$ is the graph in Figure 14.3.  We have $w_4 =
w_5 = w_6 = 2$ and $w_1=w_2=w_3 = n = 6$, hence $\ga$ is as shown in
Figure 14.5(a).

The edges in $\hat e_2, \hat e_3$ are parallel to those in $\hat e_1$,
as shown in Lemmas 14.2 and 14.5, therefore they form families of
three parallel edges, as shown in Figure 14.5(b).  Orientations are
from $u_1$ to $u_2$ on $\ga$, so the tails of these edges are labeled
$1$ and the heads labeled $2$ on $\gb$.  The two edges in $\hat e_4$
connect $v_4, v_6$ and $v_3, v_5$ respectively, and by (1) and (2)
they are not parallel to edges in $C_i$ and lie in different annuli of
$\hat F_b - C_1\cup C_2$, hence we may assume that they look like that
in Figure 14.5(b).  The four edges in $\hat e_5$ and $\hat e_6$ are
now determined by the labeling of the edges and the vertices on
$\gb$.  The labeling of the weight 3 families in $\gb$ are determined
by the single edges and the assumption that the jumping number is $1$.
\qed

\section {$\ga$ with $n_a \leq 2$ }

The next few sections deal with the case that $n_a \leq 2$ and $n_b
\leq 4$.  In this section we set up notation and give some
preliminary results.

We use $G = (b_1, b_2, b_3)$ to denote a graph $G$ on a torus with one
vertex and three families of edges weighted $b_1, b_2, b_3$.
Similarly, denote by $G = (\rho; a_1, ..., a_4)$ a graph $G$ on a
torus which has two vertices, two families of loops of weight $\rho$,
and four families of edges $\he_i$ with weight sequence $a_1, ...,
a_4$ around the vertices.  It is possible that $\rho$ and some of the
$a_i$ may be zero.  When $\rho = 0$ we will simply write $G = (a_1,
..., a_4)$.  Note that the weight sequence is defined up to cyclic
rotation and reversal of order.  When $\rho=0$, any weight $0$ can be
moved around without changing the graph, hence $(2,2,0,0)$ is
equivalent to $(2,0,2,0)$, but $(1;2,0,2,0)$ is different from
$(1;2,2,0,0)$ and $(3,1,3,1)$ is different from $(3,3,1,1)$.  When it
is necessary to indicate whether the vertices of $G$ are parallel or
antiparallel, we write $G =+(\rho; a_1, ..., a_4)$ if the vertices of
$G$ are parallel, and $G=-(\rho; a_1, ..., a_4)$ otherwise.

If $n_a = 2$ then $\ga$ is of the form $(\rho; a_1, ..., a_4)$.  Note
that if $e$ is co-loop then all edges parallel to $e$ are.  Hence we
may define an edge $\he$ in $\rga$ to be co-loop if one (and hence
all) of its edges is co-loop.  Define $\epsilon_i = 0$ if $\he_i$ is a
co-loop, and $\epsilon_i=1$ otherwise.  Note that if $n_b =2$ then
$\epsilon_i$ measures the difference between the labels at the two
endpoints of an edge in $\he_i$, so it is actually the same as the
transition number defined in Section 2.

\begin{lemma} {\bf (The Congruence Lemma.)}  Suppose $n_a =
2$.  Let $\he_i, \he_j$ be edges in $\rga$ with the endpoints on the
same pair of vertices $u_1, u_2$.  Let $a_k$ be the weight of $\he_k$.

(1) If $\rga$ has no loops and $a_i, a_j\neq 0$ then $a_i +
\epsilon_i \equiv a_j + \epsilon_j$ mod $2$.  In other words, $a_i
\equiv a_j$ mod $2$ if and only if $\he_i$ and $\he_j$ are both
co-loop or both non co-loop.

(2) If $\rga$ has loops and $a_i, a_j \neq 0$, then $a_i \equiv
a_j$ mod $2$.

(3) If $\rga$ has loops and the endpoints of $\he_i, \he_j$ at $u_1$
are on the same side of the loop at $u_1$ then $a_i \equiv a_j$
mod $2$.  \end{lemma}

\proof (1) Delete edges of $\rga$ with zero weight.  We need only
prove the statement for adjacent edges $\he_1, \he_2$ of $\rga$ with
non-zero weight.  Let $e_1, ..., e_{a_1}$ and $e'_1, ..., e'_{a_2}$
be the edges in $\he_1$ and $\he_2$, respectively, so that $e'_1$ is
adjacent to $e_{a_1}$ on $\bdd u_1$.  Then $e'_{a_2}$ is adjacent to
$e_1$ on $\bdd u_2$.  Without loss of generality assume that the label
of $e_i$ at $u_1$ is $i$. (Since $n_b=2$, all labels of endpoints of
$e_i, e'_j$ are mod 2 integers.)  Then the label of $e_1$ at $u_2$ is
$1+\epsilon_1$.  On the other hand, the label of $e'_1$ at $u_1$ is
$a_1 + 1$, so the label of $e'_{a_2}$ at $u_1$ is $a_1 + a_2$, and
the label of $e'_{a_2}$ at $u_2$ is $a_1 + a_2 + \epsilon_2$.  See
Figure 15.1.  Since $e'_{a_2}$ is adjacent to $e_1$ on $u_2$, the
label of $e_1$ on $u_2$ is $a_1 + a_2 + \epsilon_2 + 1$.  These two
equations give
$$
1 + \epsilon_1 \equiv a_1 + a_2 + \epsilon_2 + 1 \qquad \text{ mod
} 2$$
It follows that $a_1 \equiv a_2$ if and only if $\epsilon_1
\equiv \epsilon_2$ mod $2$.

\bigskip
\leavevmode

\centerline{\epsfbox{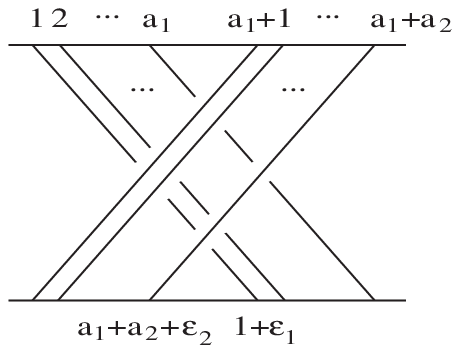}}
\bigskip
\centerline{Figure 15.1}
\bigskip

(2) Again we need only prove the statement for adjacent edges $\he_1,
\he_2$ with non-zero weight.  Since $\ga$ has loops, the two vertices
of $\rgb$ must be antiparallel.  If the two vertices of $\ga$ are
parallel then $\he_1, \he_2$ are both positive on $\ga$ and hence both
negative on $\gb$, so they are both non co-loops.  Similarly if the
two vertices of $\ga$ are antiparallel then $\he_1, \he_2$ are both
co-loops.  Therefore $\epsilon_1 = \epsilon_2$ in either case.  Note
that the endpoints of $\he_1, \he_2$ are on the same side of the loop
at $u_1$ if and only if their other endpoints are on the same side of
the loop at $u_2$.  Since the number of loops at the two vertices are
the same, the distance between the endpoints of $e_{a_1}$ and $e'_1$
on $\bdd u_1$ is the same as that of $e'_{a_2}$ and $e_1$ on $\bdd
u_2$, hence the above argument can be modified to show that if $a_1,
a_2 \neq 0$ then $a_1 \equiv a_2$ mod $2$.  More explicitly, if
there is no loop between the endpoints of $\he_1, \he_2$ then the
above argument follows verbatim, while if there are $k$ loops between
them then the endpoint $e'_{a_2}$ is $a_1 + a_2 + k$ at $u_1$, and
$a_1 + a_2 + k + \epsilon_2$ at $u_2$, and we have 
$a_1 + a_2 + k + \epsilon_2 + k + 1 = 1 + \epsilon_1$, hence the
result follows because $\epsilon_1 = \epsilon_2$ mod $2$.

(3) This follows from (2) if $a_i, a_j \neq 0$.  If $a_i=0$ and
$\he_j$ is on the same side of the loop as $\he_1$ (which is empty),
then since a loop has different labels on its two
endpoints, the number of edges in $\he_2$ must be even, hence $a_2
\equiv a_1 = 0$ mod $2$.  \qed

\begin{lemma} Suppose $\gb$ is positive, and contains a black
bigon $e_1\cup e_2$ and a white bigon $e'_1 \cup e'_2$.  Then on $\ga$
the four edges $e_1, e_2, e'_1, e'_2$ cannot be contained in two
families of parallel edges.  \end{lemma}

\proof Recall that no two edges are parallel on both graphs, so if the
lemma is not true then we may assume that $e_i$ is parallel to $e'_i$
on $\ga$.  Let $B_i$ be the disk on $F_a$ realizing the parallelism,
and let $D, D'$ be the bigon on $F_b$ bounded by $e_1\cup e_2$ and
$e'_1 \cup e'_2$, respectively.  Then $A = D \cup D' \cup B_1 \cup
B_2$ is either a M\"obius band or an annulus.  The first case
contradicts the fact that a hyperbolic manifold $M$ contains no M\"obius
bands.  In the second case $A$ contains a single white bigon and hence
each of its boundary components intersects a curve of slope $r_a$
transversely at a single point.  Since $e_1$ is an essential arc on
both $F_a$ and $A$, $A$ cannot be boundary parallel, and hence is
essential in $M$, which is again a contradiction to the hyperbolicity
of $M$.  \qed

\section {The case $n_a = 2$, $n_b=3$ or $4$, and $\Gamma_1,
\Gamma_2$ non-positive }

Throughout this section we assume that $n_a = 2$, $n_b = 3$ or $4$,
and both $\Gamma_1, \Gamma_2$ are non-positive.  We will show that in
this case there are only three possibilities for the pair $(\ga,
\gb)$, given in Figures 16.6, 16.8 and 16.9.  The following lemma
rules out the possibility that $n_b = 3$.

\begin{lemma} The case $n_a=2$, $n_b=3$ and $\ga, \gb$
non-positive, is impossible.  \end{lemma}

\proof The graph $\ga$ contains at most one loop at each vertex as
otherwise it would contain a Scharlemann bigon, which contradicts
Lemma 2.2(4) because $n_b = 3$ implies that $\hat F_b$ is
non-separating.  There are at most four families of edges on $\ga$
connecting $u_1$ to $u_2$, containing a total of at least $\Delta n_b
- 2 \geq 10$ edges, hence there is a family containing 3 edges $e_1
\cup e_2 \cup e_3$.  These are positive edges in $\gb$, and we may
assume that $e_i$ has label $i$ at $u_1$.  Since one of the vertices
of $\gb$, say $v_1$, is anti-parallel to the other two vertices, the
edge $e_1$ is a loop on $\gb$, so its label on $u_2$ is also $1$.
Since $u_1, u_2$ are antiparallel, we see that the label of $e_i$ at
$u_2$ is $i$ for $i=1,2,3$, hence they are all co-loop edges on $\ga$.
This is a contradiction to the 3-Cycle Lemma 2.14(2).  \qed

We will assume in the remainder of this section that $n_b = 4$.  By
Lemma 13.1 the graph $\rga$ is as shown in Figure 13.1.  Note that
$\hat e_1, \hat e_2$ are on the same side of the loop at each $u_i$.
Denote by $w_i$ the weight of $\hat e_i$, and put $\lambda = w_5 =
w_6$.  Then we can denote $\ga$ by $(\lambda; w_1, w_2, w_3, w_4)$,
and by Lemma 13.1(2) we may assume that $w_3 + w_4 \leq w_1 + w_2 = 6$
or $8$.  By Lemmas 2.3(1) and 2.3(3) we have $\lambda, w_i \leq 4$.
Also, counting the number of edges incident to $u_i$ gives $$
\sum_{i=1}^4 w_i + 2\lambda = 4\Delta $$

\begin{lemma}  (1) If $w_i \geq 3$ then $s_i = 2$, where $s_i$
is the transition number of $\hat e_i$.

(2) $v_1$ is parallel to $v_3$ and antiparallel to $v_2$ and $v_4$.

(3) $(w_1, w_2)$ and $(w_3, w_4)$ cannot be $(3,2)$, $(3,3)$ or
$(3,4)$.  
\end{lemma}

\proof (1) Let $s_i$ be the transition number of $\hat e_i$.  By the
3-Cycle Lemma (2.14(2)) we have $s_i \neq 0$.  If $s_1 = \pm 1$ then
all vertices of $\gb$ would be parallel, which is a contradiction to
the assumption that $\gb$ is non-positive.  Since $n_b = 4$, the only
remaining possibility is that $s_i = 2$.

(2) If $\lambda \geq 3$ then $\ga$ contains a Scharlemann cycle among
the loops, so $\hat F_b$ is separating and the result follows.  If
$\lambda \leq 2$ then the equation $\sum w_i + 2\lambda = 4 \Delta$
gives $w_i \geq 3$ for some $i$.  By (1) and the parity rule, $v_j$ is
parallel to $v_{j+2}$, hence the result follows because $\gb$ is
non-positive.

(3) Assume $w_1 = 3$.  By the equation above, $\lambda > 0$, hence by
Lemma 15.1 $w_2$ is odd.  The transition function of $\hat e_1$ is
given by (1), and it will determine that of $\hat e_2$.  If $w_2 = 3$
then one can check that the transition function of $\hat e_2$ would
map $j$ to $j$, which would be a contradiction to (1).  \qed

\begin{lemma}  $\lambda \geq 2$.
\end{lemma}

\proof First assume $\lambda = 0$.  Then $\ga = (0;4,4,4,4)$.  By
Lemma 16.2(1) all edges of $\ga$ have label pair $(1,3)$ or $(2,4)$,
see Figure 16.1.  Thus $\rgb$ is a union of two cycles, hence all edges
from $v_1$ to $v_3$ in $\gb$ are equidistant.  Since two of these
edges are in $\hat e_1$ and are not equidistant on $\ga$, this is a
contradiction to the Equidistance Lemma 2.17.

\bigskip
\leavevmode

\centerline{\epsfbox{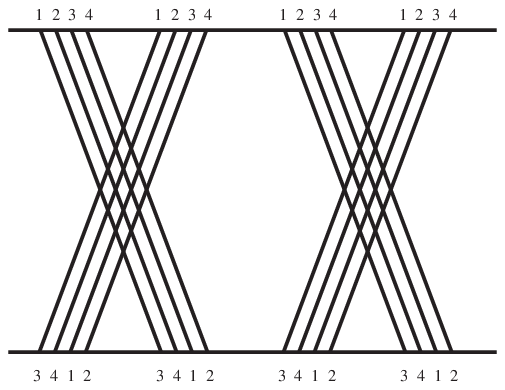}}
\bigskip
\centerline{Figure 16.1}
\bigskip

If $\lambda = 1$, then $w_j \geq 2$ for $j=1,...,4$, so $w_i \neq 3$
by Lemma 16.2(3).  Hence $\ga = (1;4,4,4,2)$.  One can check that the
one of the families of weight 4 would have the same label on the two
endpoints of any of its edges, which is a contradiction to Lemma
16.2(1).  \qed

\begin{lemma} Suppose $w_i = w_j = 4$ and $\hat e_i = e_1 \cup
... \cup e_4$ and $\hat e_j = e'_1 \cup ... \cup e'_4$ satisfy
(i) they have the same label sequence at $u_1$, and (ii)
$e_1$ is equidistant to $e'_1$ on $\ga$.  Then there exist at least 4
non co-loop edges in the other two non-loop families of $\ga$.
\end{lemma}

\proof The graph $\ga$ is as shown in Figure 16.2 for the case $(i,j)
= (1,2)$.  (The proof works in all cases.)  Note that $e_1$ being
equidistant to $e'_1$ implies that $e_k$ is equidistant to $e'_k$ for
$k=1,2,3,4$.  We may assume that the label sequence of $\hat e_i$ and
$\hat e_j$ is $1,2,3,4$ at $u_1$.  By Lemma 16.2(1) the four edges
$e_1 \cup ... \cup e_4$ form two essential cycles on $\gb$, so any
edge on $\gb$ with endpoints $v_1, v_3$ must be parallel to $e_1$ or
$e_3$.  In particular, the edge $e'_1$ has label pair $(1,3)$ and
hence must be parallel to either $e_1$ or $e_3$.  Note that two
parallel positive edges are equidistant.  Since $e'_1$ is equidistant
to $e_1$ and $e_1$ is not equidistant to $e_3$ on $\ga$, it follows
that $e'_1$ is not equidistant to $e_3$ on $\ga$, therefore by the
Equidistance Lemma and the above we see that $e'_1$ must be parallel
to $e_1$ on $\gb$.  Similarly each $e'_k$ is parallel to $e_k$ on
$\gb$.  Since $e'_k$ and $e_k$ have the same label $k$ at $u_1$ on
$\ga$, they have the same label $1$ at $v_k$ in $\gb$, so there must
be another edge $e''_k$ between them.  By the above $e''_k$ cannot be
in $\hat e_i \cup \hat e_j$, hence they belong to the other two
families of non-loop edges in $\ga$, and the result follows.
\qed

\bigskip
\leavevmode

\centerline{\epsfbox{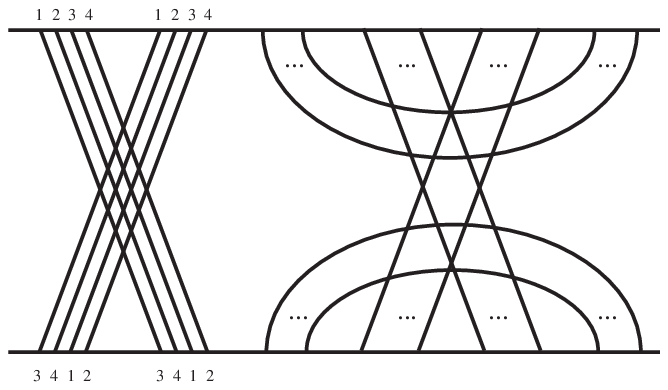}}
\bigskip
\centerline{Figure 16.2}
\bigskip

\begin{lemma}  $\lambda = 3$ is impossible.
\end{lemma}

\proof
Suppose $\lambda = 3$.  Using the Congruence Lemma (Lemma 15.1) and
Lemma 16.2(3) one can show that $\ga$ has the following possibilities.

(1) $\Delta = 5$, $(3;4,4,4,2)$;

(2) $\Delta = 4$, $(3;4,4,2,0)$;

(3) $\Delta = 4$, $(3;4,2,4,0)$;

(4) $\Delta = 4$, $(3;4,2,2,2)$.

In each case, the family of $\hat e_1$ has weight $4$.  We assume that
its label sequence at $u_1$ is $1,2,3,4$.  Then by Lemma 16.2(1) its
label sequence at $u_2$ is $3,4,1,2$, which then completely determines
the labels of $\ga$.  One can check that in case (1) and (3) the
family $\hat e_3$ gives $4$ parallel co-loops, which is a
contradiction to the 3-Cycle Lemma (Lemma 2.14(2)).  Case (2) is
impossible by Lemma 16.4.

It remains to consider case (4).  The graph $\ga$ is shown in Figure
16.3.  The third edge $A$ of $\hat e_1$ and the second edge $B$ of
$\hat e_3$ in the figure both have label pair $(1,3)$.  As in the
proof of Lemma 16.4, this implies that they are parallel on $\gb$.
Since $\rga$ has at most 4 negative edges and at most 2 positive
edges, by Lemma 2.2(2) $\gb$ cannot have more than $2n_a = 4$ parallel
edges, so the endpoints of $A$ and $B$ at $v_3$ are adjacent among the
four edge endpoints labeled $1$ at $v_3$.  Since $\Delta =4$, the
jumping number is $\pm 1$, so the endpoints of $A, B$ at $u_1$ in
$\ga$ are also adjacent among the four edge endpoints labeled $3$ at
$u_1$.  This is a contradiction because this is not the case in Figure
16.3.  \qed

\bigskip
\leavevmode

\centerline{\epsfbox{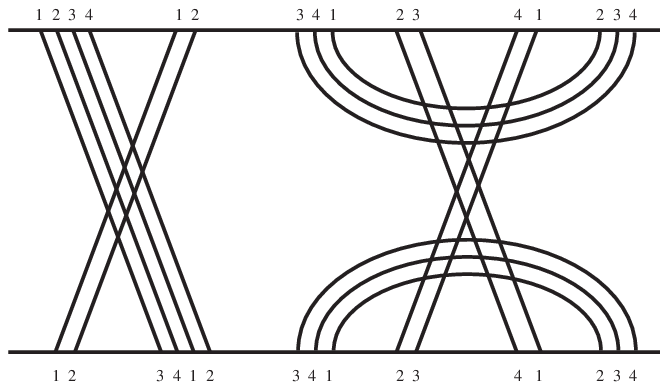}}
\bigskip
\centerline{Figure 16.3}
\bigskip

\begin{lemma} If $\lambda = 4$, then $\ga = (4; 4,2,4,2)$, and
the graphs $\ga, \gb$ and their edge correspondence are as shown in
Figure 16.6.  \end{lemma}

\proof Since $\ga$ does not contain an extended Scharlemann cycle, by
considering the labels at the endpoints of the four loops at $u_1$ we
see that $w_1 + w_2 \equiv w_3 + w_4 \equiv 2$ (mod $4$).  This,
together with Lemmas 15.1 and 16.2(3), give the following
possibilities for $\ga$.

(1) $\Delta = 5$, $\ga = (4; 4,2,4,2)$;

(2) $\Delta = 5$, $\ga = (4; 4,2,2,4)$;

(3) $\Delta = 4$, $\ga = (4; 4,2,2,0)$.

We shall show that (2) and (3) are impossible, and (1) gives the
example in Figure 16.6.

Case (2) can be excluded by Lemma 16.4.  The graph $\ga$ is shown in
Figure 16.4.  Note that the corresponding edges of the two non-loop
families of weight 4 are equidistant in $\ga$, and they have the same
label sequence at $u_1$.  Since the other two non-loop families of
$\ga$ consist of co-loops, this is a contradiction to Lemma 16.4.

\bigskip
\leavevmode

\centerline{\epsfbox{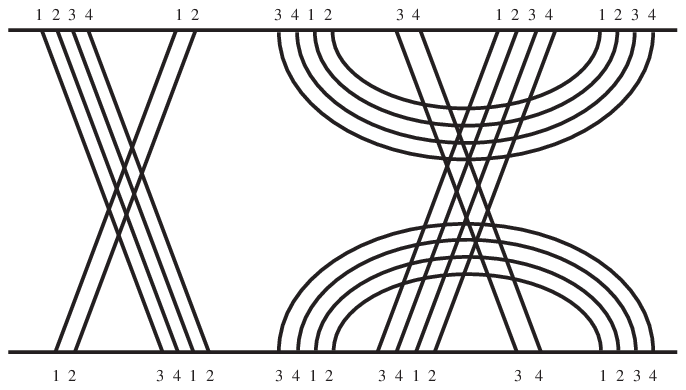}}
\bigskip
\centerline{Figure 16.4}
\bigskip

The graph for case (3) is shown in Figure 16.5.  Note that there is a
loop in $\gb$ based at each vertex $v_i$, so two edges connecting
$v_i$ to different vertices must be on different sides of the loop.
Consider the four edges with label $3$ at $u_1$, indicated by $A, B,
C, D$ in Figure 16.5.  Note that they appear in this order on $\bdd
u_1$.  Since $\Delta = 4$, the jumping number is $\pm 1$, so they must
also appear in such an order on $\bdd v_3$ in $\gb$.

On the other hand, since $A$ connects $v_3$ to $v_1$ while
$B, D$ connect $v_3$ to $v_4$, $A$ must be on a different side of the
loop $C$ at $v_3$ than $B, D$.  Hence when traveling around $\bdd v_3$
in a certain direction the four edges appear in the order $A,C,B,D$
or $A,C,D,B$.  This is a contradiction.  Therefore Case (3) is
impossible.

\bigskip
\leavevmode

\centerline{\epsfbox{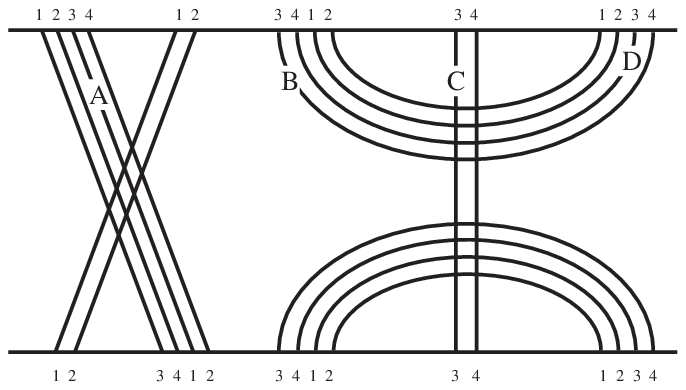}}
\bigskip
\centerline{Figure 16.5}
\bigskip

In case (1), the graph $\ga$ is shown in Figure 16.6(a).  By the same
argument as above, we see that the edges $B \cup E$ and $A \cup C$
must be on different sides of the loop $D$ in $\gb$.  Therefore $B, E$
are adjacent among the 5 edges labeled $1$ at $v_3$.  Since they are
not adjacent among the $3$-edges at $u_1$, the jumping number must be
$\pm 2$.  This completely determines the edges around the vertex $v_3$
up to symmetry, which in turn determine the edges at adjacent vertices
$v_1, v_3$ and then the edges at $v_2$.  The graph $\gb$ is shown in
Figure 16.6(b).  
\qed

\bigskip
\leavevmode

\centerline{\epsfbox{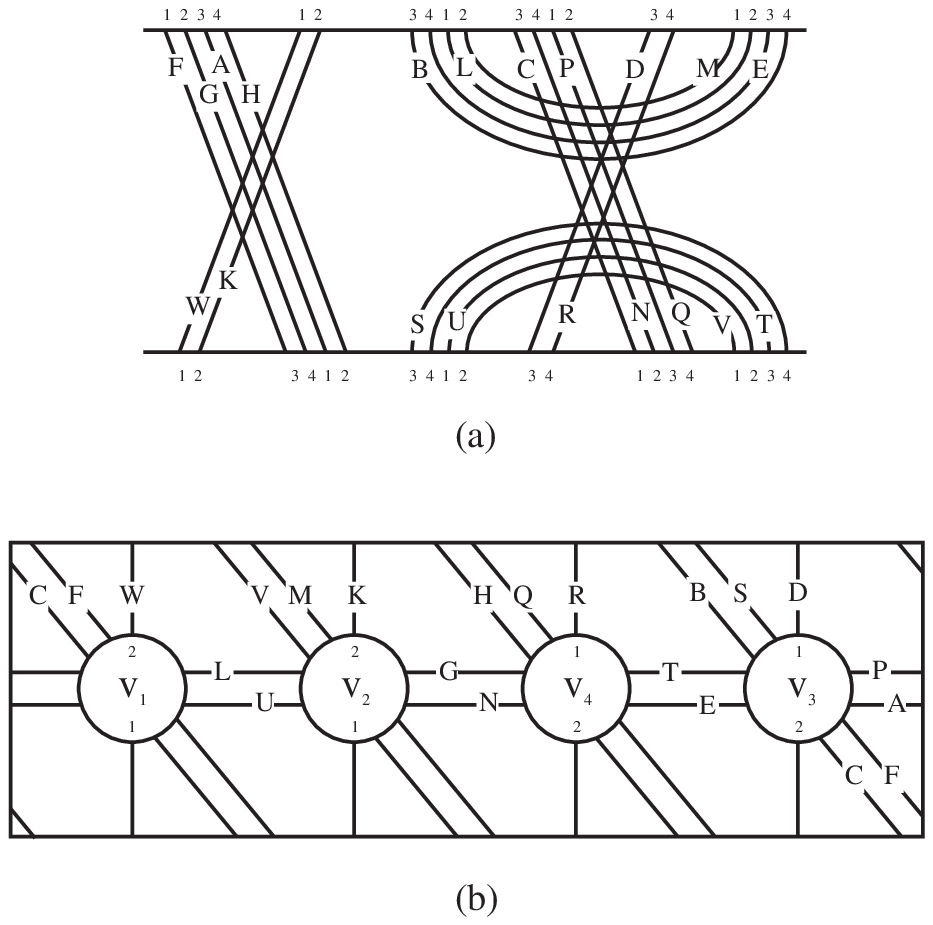}}
\bigskip
\centerline{Figure 16.6}
\bigskip

\begin{lemma}  If $\lambda = 2$, then either $\Delta=5$ and
$\ga = (2;4,4,4,4)$, or $\Delta = 4$ and $\ga = (2;4,4,4,0)$.  The 
graphs $\ga, \gb$ are as shown in Figures 16.8 and 16.9.  
\end{lemma}

\proof
Here the possibilities for $\ga$ are

(1) $\Delta = 5$, $\ga = (2;4,4,4,4)$;

(2) $\Delta = 4$, $\ga = (2;4,4,4,0)$;

(3) $\Delta = 4$, $\ga = (2;4,4,2,2)$;

(4) $\Delta = 4$, $\ga = (2;4,2,4,2)$;

(5) $\Delta = 4$, $\ga = (2;4,2,2,4)$.

The graphs in cases (3) -- (5) are shown in Figure 16.7 (a) -- (c).
In cases (3) and (4) the corresponding edges in the two weight 4
families are equidistant, and the other two non-loop families are
co-loops.  Therefore these cases are impossible by Lemma 16.4.  In
case (5) there are loops at $v_1$ and $v_2$ in $\gb$, and there is a
(34)-Scharlemann bigon in $\ga$ which forms another essential cycle
$C$ in $\gb$.  Consider the two edges of $\ga$ with label $3$ at $u_1$
and label $1$ at $u_2$.  On $\gb$ these edges connect $v_3$ and $v_1$,
and therefore must lie on the same side of $C$. Hence they are
adjacent among the four edges labeled $1$ at $v_3$ because the other
two edges connect $v_3$ to $v_4$.  Since $\Delta = 4$, the jumping
number must be $\pm 1$, so these edges are also adjacent among the
four edges with label $3$ at $u_1$, which is a contradiction because
on Figure 16.7(c) the two edges with label $3$ at $u_1$ and $1$ at
$u_2$ are not adjacent among the four edges labeled $3$ at $u_1$.
Therefore (5) is also impossible.

\bigskip
\leavevmode

\centerline{\epsfbox{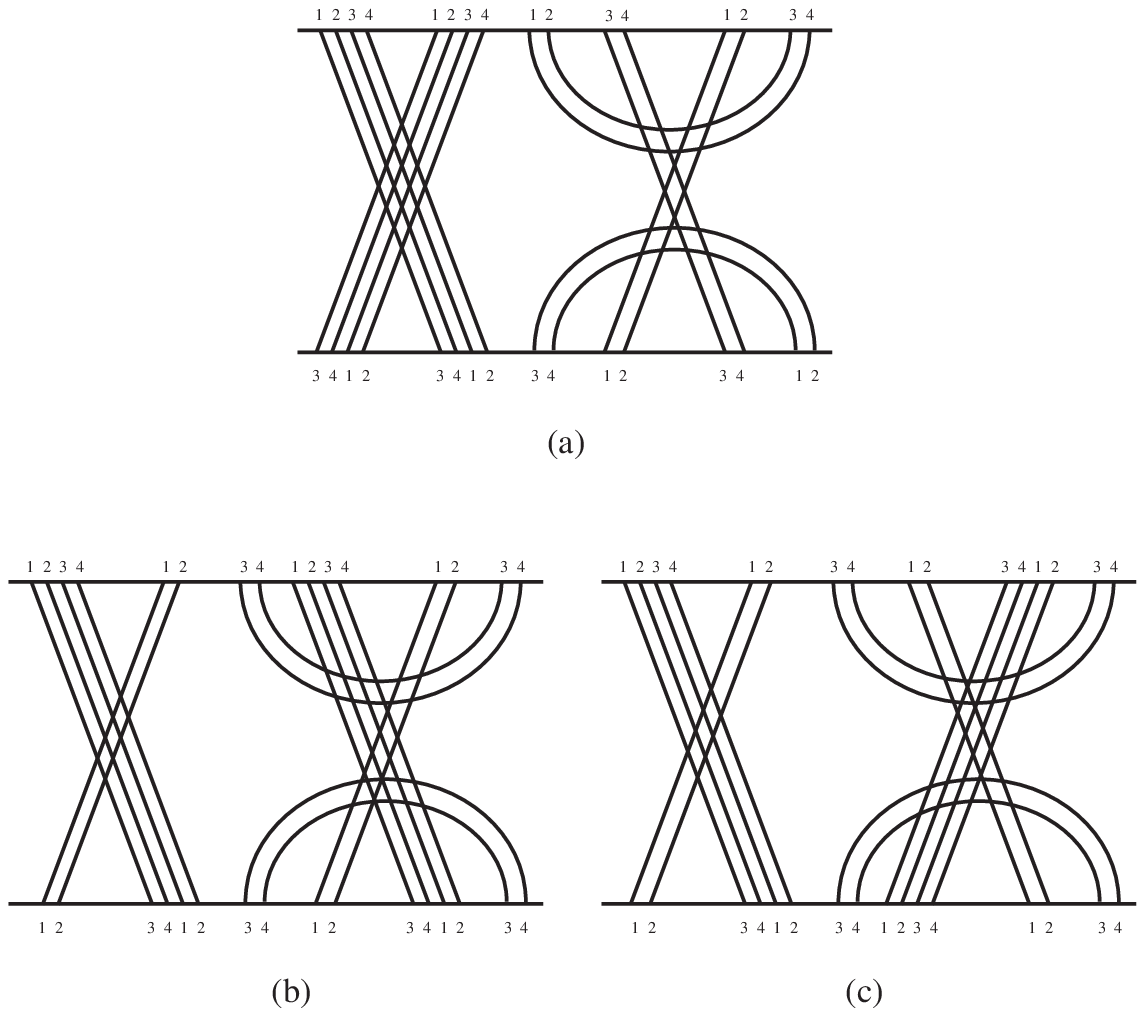}}
\bigskip
\centerline{Figure 16.7}
\bigskip

In case (1) the graph $\ga$ is shown in Figure 16.8(a).  Label the
edges as in the figure, and orient non-loop edges of $\ga$ from $u_1$
to $u_2$.  As in the proof of Lemma 16.4, the $i$-th edge $e_i$ in
$\hat e_1$ must be parallel to the $i$-th edge $e'_i$ in $\hat e_2$ on
$\gb$, and there is an edge of $\hat e_3 \cup \hat e_4$ between them
because $e_i, e'_i$ both have label $1$ at $v_i$.  For the same reason
the $i$-th edge of $\hat e_3$ is parallel to the $i$-th edge of $\hat
e_4$, hence the positive edges of $\gb$ form four families of weight
4.  The two edges $e_i, e'_i$ are adjacent among the five edges
labeled $1$ at $v_i$ in $\gb$, hence the jumping number $J = \pm 1$.
Reversing the orientation of the vertices of $\gb$ if necessary we may
assume $J=1$.  We may also assume that the vertices $v_1, v_3$ are
oriented counterclockwise and $v_2, v_4$ clockwise, otherwise we may
look at $\hat F_b$ from the other side.

Since $\gb$ contains 4 parallel positive edges, by Lemma 13.2(2) $\ga$
is kleinian, so the weight of edges of $\rgb$ are all even.  There are
only two $(14)$-edges $K, W$ in $\ga$, so they must be parallel in
$\gb$.  They may appear in the order $(K,W)$ or $(W,K)$ on $\bdd
v_1$, but there is a homeomorphism of $(\hat F_a, \ga)$ which is label
preserving, interchanging $u_1, u_2$ and mapping $K$ to $W$, hence up
to symmetry we may assume that the order is $(K,W)$.  Thus up to
symmetry we may assume that $K$ and $W$ appear in $\gb$ as shown in
Figure 16.8(b).  This, together with the orientation of the vertices
and the fact that $J=1$, completely determines the edges around $v_1$
and $v_4$, and then the edges around $v_2$ and $v_3$.  See Figure
16.8(b).

\bigskip
\leavevmode

\centerline{\epsfbox{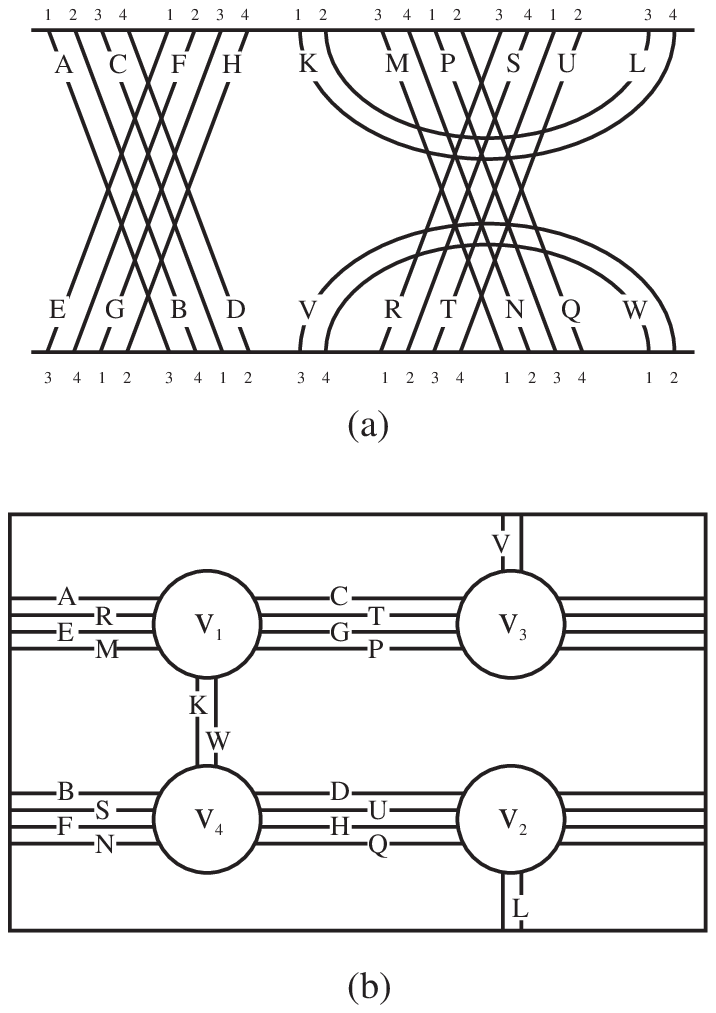}}
\bigskip
\centerline{Figure 16.8}
\bigskip

The graph $\ga$ in case (2) is shown in Figure 16.9(a).  As above, one
can show that each edge $e_i$ in $\hat e_1$ is parallel to $e'_i \in
\hat e_2$ and $e''_i \in \hat e_3$, where $e_i, e'_i$ have label $i$ at
$u_1$ and $e''_i$ has label $i$ at $u_2$.  Orient $v_i$ as above.  Up
to symmetry we may assume $J=1$, and $A,E$ on $\gb$ are as shown in
Figure 16.9(b).  This determines $P$ and the position of $K$ at $\bdd
v_1$, and hence the labels of the $2$-edges at $v_1$.  The $4$-labels
at $u_1$ appear in the order $K,D,H,N$, so on $\gb$ they appear in
this order around $v_4$, clockwise, hence $D,H$ must be to the right
of $v_4$ in the figure.  This also determines the $2$-edges at $v_4$.
In particular, the edges $K$ and $S$ must be non-parallel.  The
remaining two edges $R$ and $L$ can be determined similarly, using
labels at $v_2$ and $v_3$.  See Figure 16.9(b).  \qed

\bigskip
\leavevmode

\centerline{\epsfbox{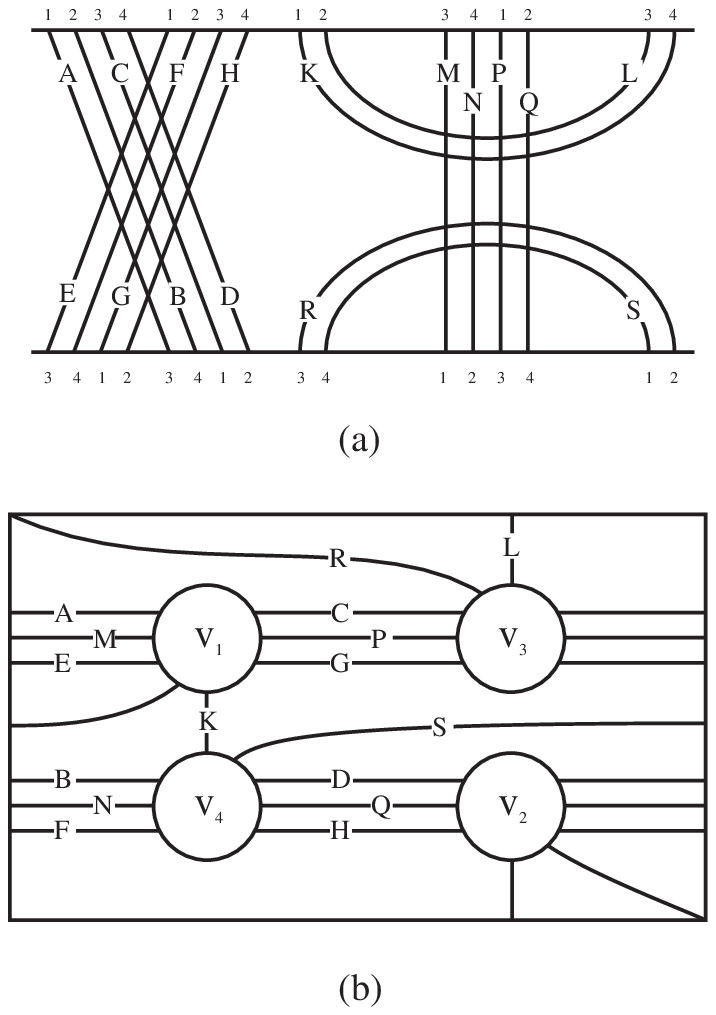}}
\bigskip
\centerline{Figure 16.9}
\bigskip

\begin{prop}  Suppose $n_a \leq 2$ and $n_b \geq 3$.  Then
$\ga, \gb$ and their edge correspondence are given in Figure 11.9,
11.10, 14.5, 16.6, 16.8 or 16.9.
\end{prop}

\proof First assume that $\ga$ is positive.  Then $n_b \leq 4$ by
Lemma 3.2.  By Lemma 2.23 $n_b$ must be even, hence our assumption
implies that $n_b = 4$.  By Proposition 11.10 the graphs are as shown
in Figure 11.9 or 11.10.

Now assume $\ga$ is non-positive.  Then we have $n_a = 2$.  The case
that $\gb$ is positive has been ruled out by Proposition 12.17.  Hence
$\ga, \gb$ are both non-positive.  By Lemma 16.1 $n_b$ cannot be $3$.
By Proposition 14.7 if $n_b > 4$ then $\ga, \gb$ are given in Figure
14.5.  Finally if $n_b = 4$.  Then Lemma 16.3 and 16.5 says that
$\lambda = 4$ or $2$, which are covered by Lemmas 16.6 and 16.7,
respectively, showing that if $\lambda=4$ then the graphs are in
Figure 16.6, and if $\lambda = 2$ then the graphs are the pair in
Figure 16.8 or 16.9.  \qed

\section {Equidistance classes }

The next few sections deal with the case that $n_i \leq 2$ for
$i=1,2$.  In this section we introduce the concept of equidistance
classes.  The main properties are given in Lemmas 17.1 and 17.2, which
will be used extensively in the next few sections.

Define a relation on the set of edges $E_a$ of $\ga$ such that $e_1
\sim e_2$ if and only if (i) they have the same label pair, (ii) they
have the same endpoint vertices, and (iii) they are equidistant.

\begin{lemma}  This is an equivalence relation.
\end{lemma}

\proof We need only show that condition (iii) is transitive, i.e.\ if
$e_1, e_2, e_3$ are edges on a graph $\Gamma$ such that $e_1, e_2$ and
$e_2, e_3$ are equidistant pairs, then $e_1, e_3$ are equidistant.

By definition we have $d_{u_1}(e_1, e_2) = d_{u_2}(e_2, e_1)$, and
$d_{u_1}(e_2, e_3) = d_{u_2}(e_3, e_2)$, hence $d_{u_1}(e_1, e_3) =
d_{u_1}(e_1, e_2) + d_{u_1}(e_2, e_3) = d_{u_2}(e_2, e_1) +
d_{u_2}(e_3, e_2) = d_{u_2}(e_3, e_1)$.  This completes the proof.
\qed

We will call this equivalence relation the {\it ED relation}.  An
equivalence class is then called an {\it ED class}, and the number of
ED classes is called the {\it ED number\/} of $\ga$, denoted by
$\eta_a = \eta(\Gamma_a)$.  We can then define $D_a = D(\ga) = (c_1, ...,
c_{\eta_a})$, where $c_i$ are the number of edges of the equivalence
classes, ordered lexicographically.

\begin{lemma} Let $\ga, \gb$ be intersection graphs.  Then
the edge correspondence between the graphs induces a one to one
correspondence between the ED classes of $\ga$ and $\gb$; in
particular $\eta(\ga) = \eta(\gb)$, and $D(\ga) = D(\gb)$.
\end{lemma}

\proof Note that $e_1, e_2$ satisfy (i) on $\ga$ if and only if they
satisfy (ii) on $\gb$.  The Equidistance Lemma 2.17 now says that a
pair of edges are equivalent on $\ga$ if and only if they are
equivalent on $\gb$.  \qed

\bigskip
\leavevmode

\centerline{\epsfbox{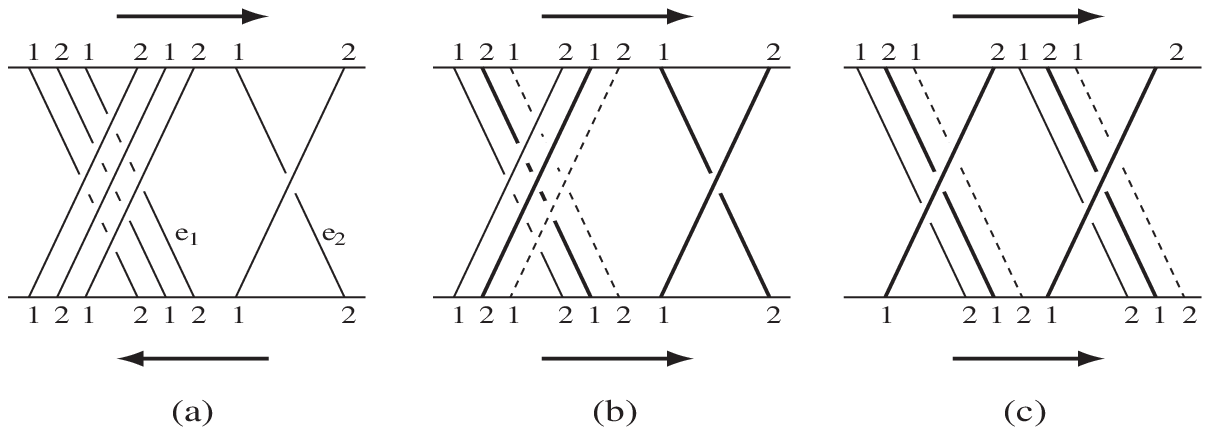}}
\bigskip
\centerline{Figure 17.1}
\bigskip

\noindent
{\bf Example 17.3} (1) Consider a graph $\ga = +(3,3,1,1)$ and
assume $n_b = 2$, see Figure 17.1(a).  In general if $n_b =2$ then all
parallel positive edges are in the same ED class because they have the
same label pairs and they are equidistant.  One can check that
non-parallel edges are not equidistant.  (For example, let $u_1$ be
the top vertex, $u_2$ the bottom vertex, and let $e_1, e_2$ be as
shown in the figure; then $d_{u_1}(e_1, e_2) = 4 \neq 2 = d_{u_2}(e_2,
e_1)$.)  Hence $D(\ga) = (3,3,1,1)$.  Compare this with $\rga =
+(3,1,3,1)$, in which case the two families of 3 edges are
equidistant, and the other two families of weight 1 are equidistant,
hence $D(+(3,1,3,1)) = (6,2)$.

(2) Consider $\ga = -(3,3,1,1)$ or $-(3,1,3,1)$, and suppose that the
edges of $\ga$ are not co-loops (hence conditions (i) and (ii) in the
definition of ED equivalence are satisfied), see Figure 17.1(b) and
(c).  Equidistant edges are indicated in the figure by different kind
of lines.  We can see that $D(-(3,3,1,1)) = D(-(3,1,3,1) = (4,2,2)$.

(3) When $\ga = -(4,2,2,0)$, each of the middle edges of the family of
4 is equidistant to one edge in each of the two weight 2 families, and
the other two edges of the weight 4 families are not equidistant to
any other edges.  Hence $D(-(4,2,2,0)) = (3,3,1,1)$

(4) Suppose $\ga = +(4,2,2,0)$ and all edges have label pair
$(12)$.  Then one can check that each family of parallel edges forms an
ED class, hence $D(\ga) = (4,2,2)$.  

(5) Suppose $\ga = +(2,2,2,2)$ and all edges have label pair $(12)$.
Then one can show that the first family is equidistant to the third
family, but not to the adjacent families.  Hence $D(\ga) = (4,4)$.

(6) Similarly if $\ga = +(4,4,0,0)$ and all edges have the same label
pair then $D(\ga) = (8)$.  

\medskip

\section {The case $n_b = 1$ and $n_a = 2$}

\begin{lemma}  Suppose $n_a = 2$ and $n_b = 1$.  Then 
one of the following holds.

(1) $\ga = -(1,1,1,1)$ and $\gb = (4,0,0)$.

(2) $\ga = -(2,2,0,0)$ and $\gb = (2,2,0)$.

(3) $\ga = -(2,1,1,1)$ and $\gb = (3,1,1)$.  The graphs $\ga, \gb$ and
their edge correspondence are given in Figure 18.2.
\end{lemma}

\proof In this case $\rgb$ has a single vertex, and $\rga$ has two
vertices of opposite orientation and has no loops.  Hence we have $\ga
= (a_1, ..., a_4)$, and $\gb = (b_1, b_2, b_3)$.  We have $b_1 + b_2 +
b_3 = \Delta$.  If $b_i, b_j$ are non-zero and $b_i+b_j$ is odd then
one can check that one of the $\hat e_i, \hat e_j$ is a family of
co-loops, which is a contradiction to the parity rule.  Hence $b_i
\equiv b_j$ mod $2$ for all $b_i, b_j$ non-zero.  Thus if $\Delta= 5$
then up to symmetry we have $\gb = (3,1,1)$, and if $\Delta = 4$ then
$\gb = (4,0,0)$ or $(2,2,0)$.

If $\gb = (4,0,0)$ then the four parallel edges are mutually
non-parallel on $\ga$, hence $\ga = -(1,1,1,1)$.  

If $\gb = (2,2,0)$, one can check that edges in different families are
not equidistant, hence $D(\gb) = (2,2)$.  Since each pair of parallel
edges contributes one edge to each of two families in $\ga$, we have
$\ga = -(2,2,0,0)$, $-(2,1,1,0)$ or $-(1,1,1,1)$.  When
$\ga=-(2,1,1,0)$ the two single edges are equidistant, while each of
the two parallel edges form an ED class, so $D(\ga) = (2,1,1) \neq
D(\gb)$.  Also, when $\ga=-(1,1,1,1)$ we have $D(\ga) = (4)$.
Therefore in this case we have $\ga = -(2,2,0,0)$.

No suppose $\gb = (3,1,1)$.  In this case the three parallel edges are
equidistant, and each of the other two edges is not equidistant to any
other edges.  Hence $D(\gb) = (3,1,1)$.  Since the three parallel
edges in $\gb$ are mutually non-parallel on $\ga$, $\rga$ has at least
three edges.  One can show that $D(-(2,2,1,0)) = (2,2,1) \neq D(\gb)$,
hence $\rga \neq -(2,2,1,0)$.  Therefore $\ga = -(3,1,1,0)$ or
$-(2,1,1,1)$.

In the case that $\ga = -(3,1,1,0)$ and $\gb = (3,1,1)$, the graphs
are as shown in Figure 18.1.  The three parallel edges $B,C,E$ are
equidistant, hence they represent the two weight 1 edges $\hat e_2,
\hat e_3$ and the middle edge of the weight 3 edge $\hat e_1$, so the
other two edges $A, D$ must be as shown in Figure 18.1(a) up to
symmetry.  Since they are non-adjacent at $u_1$ and their label $1$
endpoints are non-adjacent among the five label $1$ edge endpoints at
$v_1$ in $\gb$, the jumping number must be $\pm 1$.  This determines
the edge correspondence between $\ga$ and $\gb$, as shown in Figure
18.1.

The torus $\hat F_a$ cuts $M(r_a)$ into two components.  Let $W$ be
the one containing the bigon $\alpha$ on $F_b$ bounded by $B\cup E$
and the 3-gon $\beta$ bounded by $A\cup C\cup D$.  It can be
constructed by attaching a 1-handle representing part of the Dehn
filling solid torus, then two 2-handles represented by $\alpha,
\beta$, then a 3-cell.  The fundamental group of $W$ is generated by 
the horizontal circle $x$ and the vertical circle $y$ shown in the
figure, and the $1$-handle $z$ from $u_2$ to $u_1$.  On the boundary
of $\alpha, \beta$, $A, B, C, D, E$ represent $1,x,xy,1,1$,
respectively, and each corner represents $z$, hence $\alpha, \beta$
give the relations $zzx = 1$ and $zxyzz = 1$, respectively.  Solving
these in $x$ and $y$ shows that $\pi_1(W) = \Bbb Z$, generated by
$z$.  It follows that $\hat F_a$ is not $\pi_1$-injective in $W$, and
hence is compressible.  This is a contradiction.

\bigskip
\leavevmode

\centerline{\epsfbox{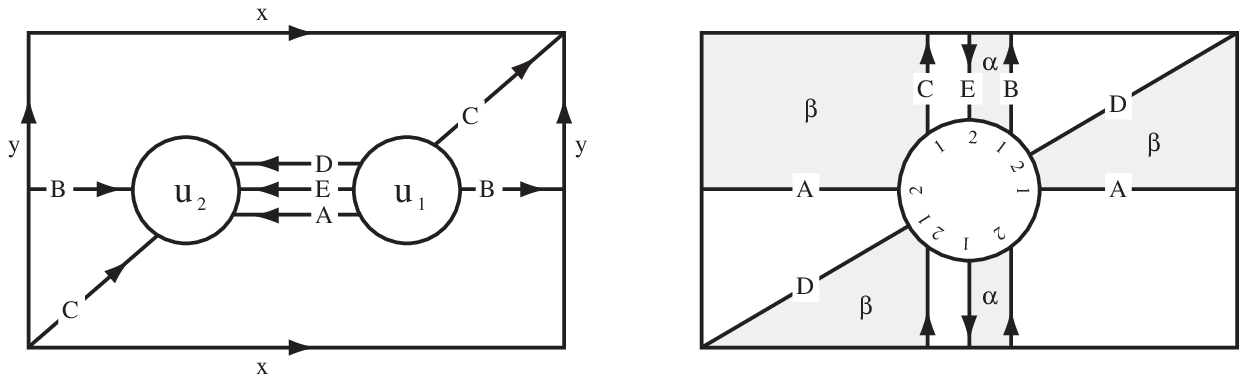}}
\bigskip
\centerline{Figure 18.1}
\bigskip

We now have $\ga = -(2,1,1,1)$ and $\gb = (3,1,1)$.  The three
parallel edges $B,C,E$ are equidistant, hence on $\ga$ they are the
single edges because they are equidistant to each other but not to the
edges in the weight 2 family.  Since the edge endpoints of these are
consecutive on $\bdd v_1$ while the $1$-label endpoint of $E$ at $v_1$
is not adjacent to that of either $B$ or $C$, the jumping number must
be $\pm 2$.  This determines the correspondence of the edges up to
symmetry, see Figure 18.2.  \qed

\bigskip
\leavevmode

\centerline{\epsfbox{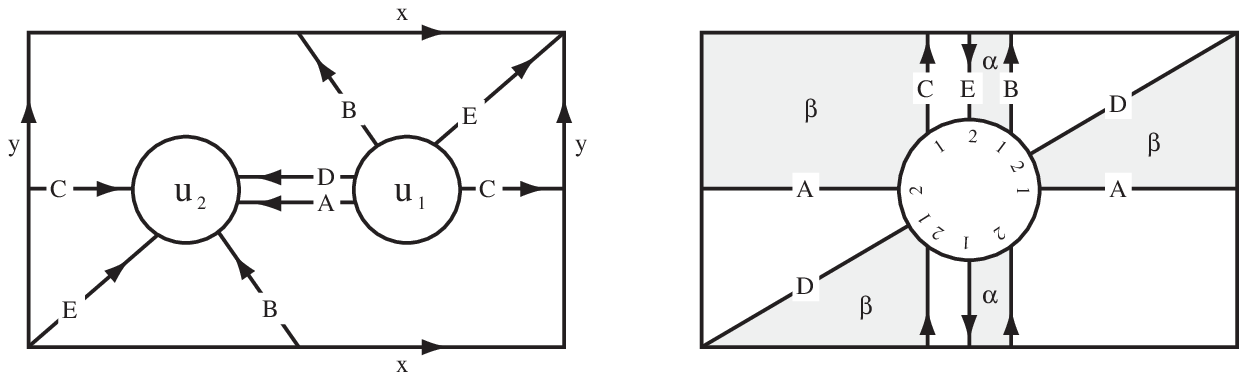}}
\bigskip
\centerline{Figure 18.2}
\bigskip

\section {The case $n_1 = n_2 = 2$ and  $\gb$ positive
}

In this section we assume that $n_1 = n_2 = 2$ and $\gb$ is positive.
Then no edge of $\rga$ is a loop, hence $\ga = -(a_1, ... a_4)$, and
$\gb = +(\rho; b_1, ..., b_4)$.  

When $\rho \neq 0$ we may rearrange the $a_i$ to write $\ga =
-(r_1, ..., r_p \, |\,  s_1, ..., s_q)$, where $r_i$ are the weights of 
the co-loop edges, and $s_j$ are the weights of the non co-loop
edges.

\begin{lemma}  Suppose $n_1 = n_2 = 2$ and $\gb$ is positive.  

(1) All non-zero $b_i$ are of the same parity, all non-zero $r_j$ are
of the same parity, all non-zero $s_k$ are of the same parity, and the
non-zero $r_j$ and $s_k$ are of opposite parity.

(2) $r_i \leq 2$, $s_j \leq 4$, and $\rho + b_k \leq 4$.

(3) $2 \rho + \sum b_i = 2 \Delta$, and $\sum r_i + \sum s_j = 2 \Delta$.
\end{lemma}

\proof
(1) This follows from the Congruence Lemma 15.1.

(2) Since $\rgb$ has at most two loops, each co-loop family of $\ga$
contains at most two edges, hence $r_i \leq 2$.  Similarly, since
$\rgb$ has at most four non-loop edges, $s_j \leq 4$.  On $\gb$,
$\rho$ is the number of edges in a loop family, which is no more than
$p$, the number of co-loop edges in $\rga$.  Similarly, $b_k$ is no
more than $q$, the number of non co-loop edges in $\rga$.  Since
$\rga$ has at most 4 edges, we have $\rho + b_k \leq 4$.
    
(3) This follows from the fact that each vertex of $\ga$ or $\gb$
has valence $2 \Delta$.
\qed

\begin{lemma}  Suppose $n_1 = n_2 = 2$ and $\gb$ is positive.
If $\rho=4$ then $\ga = -(2,2,2,2)$, and $\gb = +(4; 0, 0, 0, 0)$.
\end{lemma}

\proof Since $\rho + b_j \leq 4$ (Lemma 19.1(2)), we have $b_i = 0$
for all $i$, hence from Lemma 19.1(3) we have $\Delta = 4$.  Thus
$\rgb$ is a union of two disjoint loops, each representing a family of
four edges.  Since each family of four parallel edges in $\rgb$
contributes one edge to each family in $\ga$, we have $\ga =
-(2,2,2,2)$.  \qed

\begin{lemma}
Suppose $n_1 = n_2 = 2$ and $\gb$ is positive.  Then $\rho \neq 3$.
\end{lemma}

\proof Suppose $\rho = 3$.  The three loops in a family represent
different classes on $\rga$, so $\rga$ has at least three co-loop
edges.  Since $\gb$ has some non-loop edges, $\rga$ has at least one
non co-loop edge.  It follows that $\rga$ has exactly three co-loop
edges, so $\ga = -(2,2,2\, |\, s_1)$.  Since $\sum r_i + \sum s_j =
2\Delta$ is even, $s_1$ is even, which contradicts Lemma 19.1(1).
\qed

\begin{lemma}
Suppose $n_1 = n_2 = 2$ and $\gb$ is positive.  Then $\rho \neq 2$.
\end{lemma}

\proof On $\rga$ there are non-co-loop edges, so there are at most
three co-loop edges, but since $r_i \leq 2$ and $\sum r_i = 4$ and the
$r_i$'s are of the same parity, there must be exactly two co-loop
edges.  Hence $\ga = -(2,2 \, |\, s_1, s_2)$.  By the Congruence
Lemma, the $s_i$ are odd, hence either $\Delta = 5$ and $\ga = -(2,2\,
|\, 3,3)$, or $\Delta = 4$ and $\ga = -(2,2\, |\, 3,1)$.

If $\Delta=4$ and $\ga = -(2,2\, |\, 3,1)$, then from Lemma 19.1 we
have $b_i \leq 2$, $\sum b_i = 4$, and $b_i \equiv b_j$ mod $2$ if
$b_i, b_j \neq 0$.  These conditions give $\gb = +(2; 1,1,1,1)$,
$+(2; 2,2,0,0)$ or $+(2;2,0,2,0)$.  One can check that in the first
two cases the four non-loop edges of $\gb$ form two equidistance
classes of 2 edges each, so $D_b = (2,2,2,2)$, and in the third case
the four non-loop edges are all equidistant to each other, so $D_b =
(4,2,2)$.  On the other hand, the three parallel edges of $\ga$ belong
to distinct classes, and there are at least two co-loop classes, hence
$\eta(\ga) \geq 5$.  This is a contradiction to Lemma 17.2.

If $\Delta=5$ and $\ga = -(2,2\, |\, 3,3)$, then from Lemma 19.1 we
have $b_i \leq 4-\rho = 2$, $\sum b_i = 6$, and $b_i \equiv b_j$
mod $2$ if $b_i, b_j \neq 0$, so we must have $\gb = +(2; 2, 2, 2,
0)$.  Depending on the weight sequence of the edges of $\rga$, we have
$\ga = -(3,2,3,2)$ or $-(3,3,2,2)$.  If $\ga = -(3,2,3,2)$ then from
the labeling one can see that the two edges with both endpoints
labeled $1$ are not equidistant.  Since these are parallel loops at
$v_1$ of $\gb$, they are equidistant on $\gb$, which is a
contradiction to the Equidistance Lemma 2.17.  Therefore $\ga \neq
(-3,2,3,2)$.

Now suppose $\ga = -(3,3,2,2)$, and $\gb = +(2;2,2,2,0)$.  Then the
graphs are as shown in Figure 19.1.  Consider the edges $A,B,C,D,E$ with
label $1$ at $u_1$ of $\ga$.  These correspond to the 5 edges with
label $1$ at $v_1$ of $\gb$.  Note that on $\ga$ $D,E$ are co-loop
edges, hence on $\gb$ they are the two loops at $u_1$.  Since their
endpoints with label $1$ are not adjacent among the $1$-edges at $v_1$
in $\gb$, the jumping number must be $\pm 2$, so among these edges in
$\gb$, the edge $B$ is the one in $\gb$ which is adjacent to both $D$
and $E$ at $v_1$, as shown in the Figure.  Now consider the five edges
labeled $2$ at $u_2$.  Note that they appear in the order $CABFG$.
Using the same argument as above we see that the edge $A$ is the one
adjacent to both $F$ and $G$, so we would have $A=B$, which is a
contradiction.  
\qed

\bigskip
\leavevmode

\centerline{\epsfbox{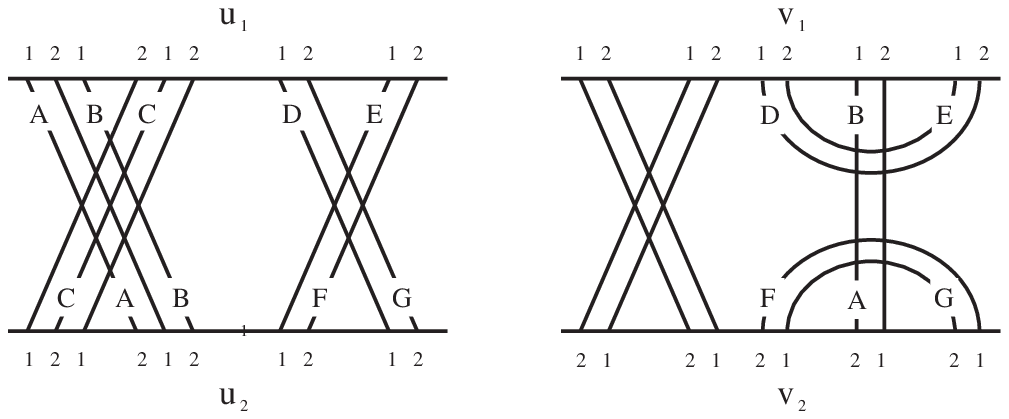}}
\bigskip
\centerline{Figure 19.1}
\bigskip

\begin{lemma}
Suppose $n_1 = n_2 = 2$ and $\gb$ is positive.  Then $\rho \neq 1$.
\end{lemma}

\proof Suppose $\rho = 1$.  Then $\ga$ has two co-loop edges, hence
$\ga = -(1,1\, |\, s_1, s_2)$ or $-(2 \, |\, s_1, s_2, s_3)$.

If $\Delta = 5$ then the second case does not happen since by the
Congruence Lemma $s_i$ would be $0$ or odd and $\sum s_i = 8$, which
would give $s_i > 4$ for some $i$, contradicting Lemma 19.1(2).
Therefore $\ga = -(1,1 \, |\, 4,4)$.  Since each weight 4 family
contributes one edge to each non-loop family of $\gb$, we have $\gb =
+(1; 2,2,2,2)$.  Now the graph $\rgb$ contains both black and white
bigons, whose edges all belong to the two weight 4 families on $\ga$.
This is a contradiction to Lemma 15.2.

If $\Delta = 4$ then $\ga = -(1,1\, |\, 4,2)$ or $-(2\, |\, 3,3,0)$.
The first case cannot happen because the $\he_i$ of weight 4
contributes one edge to each family in $\rgb$, while the edge of
weight 2 contributes one edge to each of two families, so $\gb = +(1;
2,2,1,1)$, which contradicts the Congruence Lemma.  In the second case
for the same reason above we must have $\gb = +(1; 2,2,2,0)$.  Again
there are black and white bigons, which contradicts Lemma 15.2 because
on $\ga$ the edges of these bigons all belong to the two weight 3
families.  \qed

\begin{lemma} 
Suppose $n_1 = n_2 = 2$ and $\gb$ is positive.  If $\rho = 0$ then 
$\Delta = 5$.
\end{lemma}

\proof In this case there is no loop on either graph, hence by Lemma
19.1(1) all non-zero $a_i$ have the same parity, and all non-zero
$b_j$ have the same parity.  Any two edges connect the same pair of
vertices and have the same pair of labels on their two endpoints,
hence by definition they are ED equivalent if and only if they are
equidistant.

Assume $\Delta = 4$.  By the Congruence Lemma each of $\rga$ and $\gb$
is of type $(4,4,0,0)$, $(2,2,2,2)$, $(4,2,2,0)$, $(3,1,3,1)$, or
$(3,3,1,1)$.  Let $e_1 \cup e_2$ be a bigon on $\gb$.  Then $e_1$ and
$e_2$ are equidistant on $\gb$, so by Lemma 2.17 they form an
equidistant pair on $\ga$.  Note that since $e_1$ and $e_2$ are not
loops on $\gb$, on $\ga$ they have different labels on $u_1$.  On the
other hand, one can check that if $\ga = -(4,4,0,0)$ or $-(2,2,2,2)$
then an equidistant pair $e_1, e_2$ on $\ga$ must have the same label
on $u_1$, which is a contradiction.  Therefore $\ga =-(4,2,2,0)$,
$-(3,3,1,1)$ or $-(3,1,3,1)$.  (Note that the above argument does not
apply to $\rgb$ since a pair of parallel edges on $\ga$ is not an
equidistant pair.)  We will rule these out one by one.

\medskip
CLAIM 1.  {\it The case $\ga = - (4,2,2,0)$ is impossible.}

If $\ga = -(4,2,2,0)$ then $b_i \neq 0$ for all $i$, hence $\gb =
+(2,2,2,2)$ or $+(3,1,3,1)$, or $+(3,3,1,1)$.  In the first case all
black (say) faces of $\gb$ are bigons, so Lemma 13.2(2) implies that
$\gb$ is kleinian because $\rga$ has more than two edges.  Since
$\rga$ has a single edge of weight 4, it will be fixed by the free
involution given in Lemma 6.2(4), which is absurd.  By Example 17.3 we
have $D(-(4,2,2,0)) = (3,3,1,1)$, and $D(+(3,1,3,1)) = (6,2)$, hence
$\gb \neq +(3,1,3,1)$.  It follows that $\gb = +(3,3,1,1)$.  The
graphs are shown in Figure 19.2.  Each of $A, D$ on $\ga$ forms an
equidistance class, hence they are the single edges on $\gb$.  Up to
symmetry we may assume that $A, D$ are as shown in Figure 19.2(b).
This and the jumping number $J$ determines the edge correspondence of
the graphs.  The case that $J=1$ is shown in the figure.  When $J=-1$
the edges $G, E$ would be equidistant on $\ga$ but not on $\gb$, which
is impossible.

\bigskip
\leavevmode

\centerline{\epsfbox{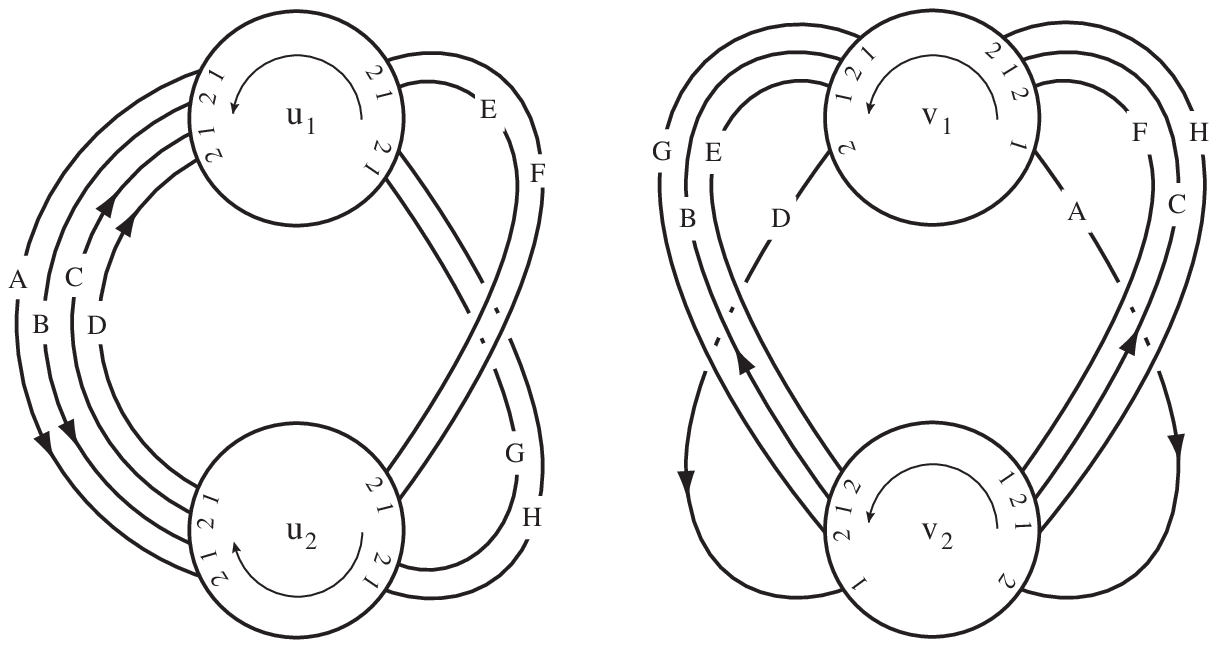}}
\bigskip
\centerline{Figure 19.2}
\bigskip

Let $P_1, P_2$ be the bigon disks on $F_a$ bounded by $A\cup B$
and $C\cup D$, respectively.  Then the union of $P_1, P_2$ and two
disks on $T_0$ form an annulus $Q$.  More explicitly, let $a_1$
(resp.\ $b_1$) be the arc on $\bdd u_1$ (resp.\ $\bdd u_2$) from the
endpoint of $A$ to that of $B$, $a_2$ (resp.\ $b_2$) the arc on $\bdd
u_1$ (resp.\ $\bdd u_2$) from the endpoint of $C$ to that of $D$, $a_3$
(resp.\ $b_3$) on $\bdd v_1$ (resp.\ $\bdd v_2$) from $A$ to $C$, and
$a_4$ (resp.\ $b_4$) on $\bdd v_2$ (resp.\ $\bdd v_1$) from $B$ to
$D$.  Then $a_1 \cup ... \cup a_4$ (resp.\ $b_1 \cup ... \cup b_4$)
bounds a disk $P_3$ (resp.\ $P_4$) on the boundary torus $T_0$.  Now
$Q = P_1 \cup ... \cup P_4$ is an annulus in $M$.  
Note that $\bdd Q$ consists of two simple closed curves $\bdd_1 = A
\cup C \cup a_3 \cup b_3$ and $\bdd_2 = B \cup D \cup a_4 \cup b_4$.

Orient $A, B$ to point from $u_1$ to $u_2$ on $\ga$.  This determines
the orientation of $\bdd_i$.  Note that they are parallel on the
annulus $Q$.  On $\gb$ the orientations of $A, B$ are from label $1$
to label $2$, as shown in the figure.  This determines the orientation
of $C, D$.  It is important to see that $\bdd_1, \bdd_2$ are parallel
as oriented curves on $\hat F_b$.  Let $Q'$ be an annulus on $\hat
F_b$ with $\bdd Q' = \bdd_1 \cup \bdd_2$.  Then $Q \cup Q'$ is a
non-separating torus (not a Klein bottle!) in $M(r_b)$ intersecting
the Dehn filling solid torus at a single meridian disk, which
contradicts the choice of $\hat F_b$.  Therefore this case is
impossible.  

\medskip
CLAIM 2.  {\it The case $\ga = - (3,1,3,1)$ is impossible.  If
$\ga = -(3,3,1,1)$ then $\gb = +(4,2,2,0)$.}

Now suppose $\ga = -(3,3,1,1)$ or $-(3,1,3,1)$ and $\gb = +(4,4,0,0)$,
$+(4,2,2,0)$, $+(2,2,2,2)$, $+(3,3,1,1)$ or $+(3,1,3,1)$.  By Example
17.3 we have $D(-(3,3,1,1)) = -(3,1,3,1) = (4,2,2)$. On the other
hand, by Example 17.3 we also have $D(+(4,4,0,0)) = (8)$,
$D(+(4,2,2,0)) = (4,2,2)$, $D(+(2,2,2,2)) = (4,4)$, $D(+(3,3,1,1)) =
(3,3,1,1)$, and $D(+(3,1,3,1)) = (6,2)$.  Therefore by Lemma 17.2 in
this case we must have $\rgb = +(4,2,2,0)$.  If $\rga = -(3,1,3,1)$
then the four edges in the same ED class all have label $2$ (say) at
$u_1$, which means that on $\gb$ they all have label $1$ at $v_2$, so
they cannot be the four parallel edges in $+(4,2,2,0)$.  Therefore
$\rga \neq -(3,1,3,1)$.

\medskip
CLAIM 3.  {\it The case $\rga = -(3,3,1,1)$ is impossible.}

By Claim 2 we have $\rgb = +(4,2,2,0)$.  The graphs are as shown in
Figure 19.3.  While the graphs are similar to that in Figure 19.2, the
argument is necessarily different because the orientation of the
vertices of $\gb$ here are parallel while that of $\ga$ in Figure 19.2
are antiparallel.  One can check that up to symmetry the edge
correspondence must be as shown in the figure.

We would like to apply Lemma 2.15 to get a contradiction.  To do that,
let $Q$ be the face of $\ga$ bounded by $A\cup B\cup E \cup H$.  The
edge $B$ is parallel to $C$ on $\gb$, and $C$ is a non-border edge on
$\ga$, hence one of the bigons $C\cup H$ or $C\cup F$ is a coupling
face $Q'$ of $Q$ along the edge $B$.  By Lemma 2.15 there is a rel
$\bdd$ isotopy of $F_a$ such that the new intersection graph $\ga'$ is
obtained from $\ga$ by deleting $A$ and $E$ and adding two edges
parallel to $B$ and $F$, respectively.  It follows that $\ga' = -
(4,2,2,0)$.  This is impossible by Claim 1.  Therefore the case $\rga
= -(3,3,1,1)$ is also impossible.
\qed

\bigskip
\leavevmode

\centerline{\epsfbox{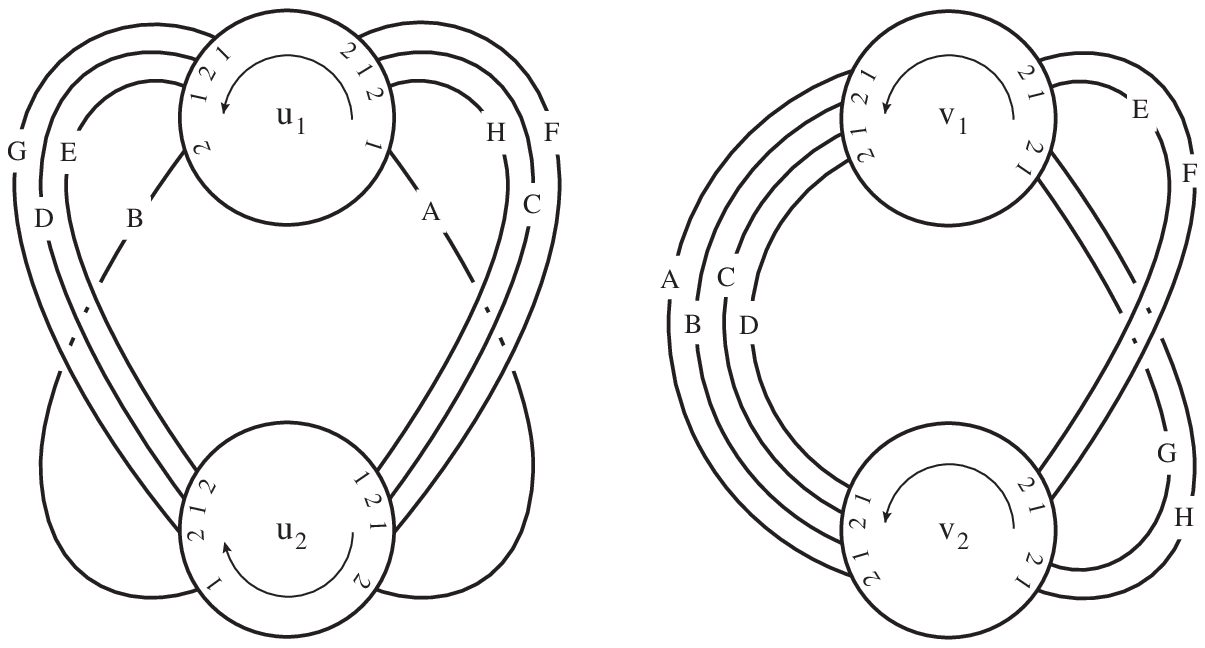}}
\bigskip
\centerline{Figure 19.3}
\bigskip

\begin{lemma}
Suppose $n_1 = n_2 = 2$ and $\gb$ is positive.  If $\rho = 0$ then 
$\ga = -(3,3,3,1)$ and $\gb = +(3,3,3,1)$.  The graphs $\ga, \gb$
and their edge correspondence are shown in Figure 19.4.
\end{lemma}

\proof By Lemma 15.1 all non-zero $a_i$ have the same parity, and all
non-zero $b_j$ have the same parity.  By Lemma 19.6 we have $\Delta =
5$, so each of $\ga$ and $\gb$ is of type $(4,4,2,0)$, $(4,2,2,2)$
or $(3,3,3,1)$.  If some $b_i=4$ then by Lemma 13.2(2) $\ga$ is
kleinian, but since each of the above type has an edge whose weight is
non-zero and different from the others, it must be mapped to itself by
the involution in Lemma 6.2(4), which is a contradiction because it is
supposed to be a free involution on $\hat F_a$.  It follows that $\gb
= +(3,3,3,1)$.  Direct calculation gives $D(+(3,3,3,1)) = (4,3,3)$,
$D(-(4,4,2,0)) = (3,3,2,2)$, $D(-(4,2,2,2)) = (4,4,1,1)$, and
$D(-(3,3,3,1)) = (4,3,3)$.  Hence by Lemma 17.2 we have $\ga =
-(3,3,3,1)$.

\bigskip
\leavevmode

\centerline{\epsfbox{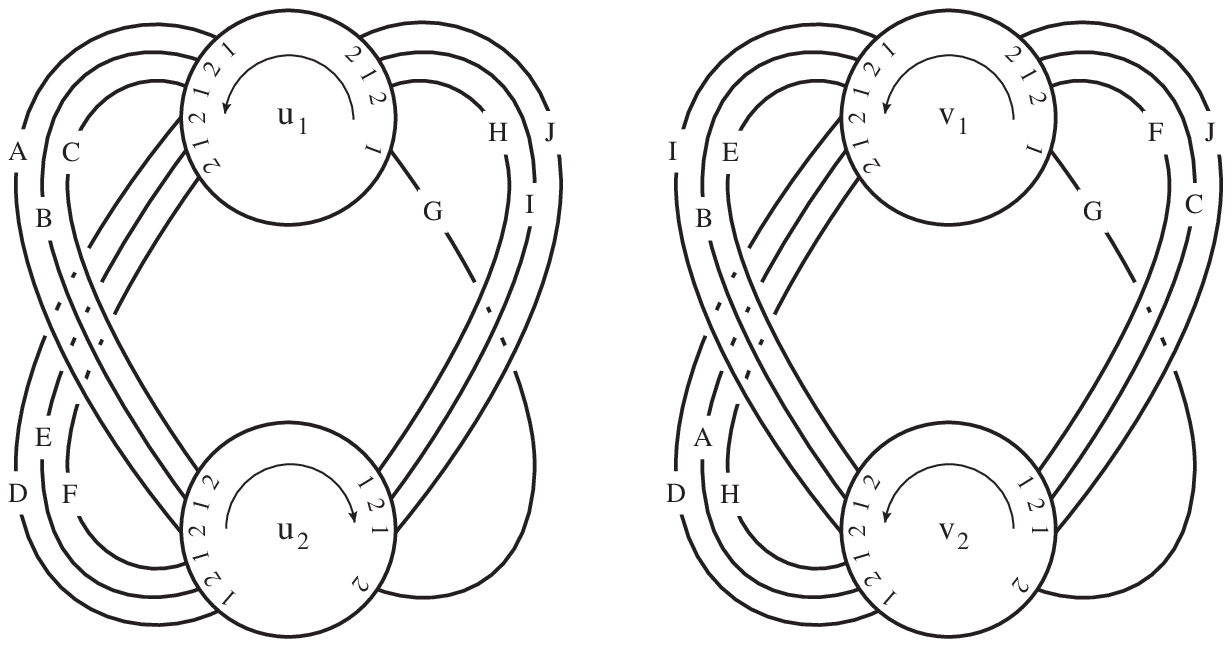}}
\bigskip
\centerline{Figure 19.4}
\bigskip

The graphs $\ga, \gb$ are as shown in Figure 19.4.  Label the edges of
$\ga$ as in the figure.  Relabeling the vertices of $\gb$ if
necessary, we may assume that the labels of edges of $\ga$ at $u_1$
are as shown.  Since $\gb$ has no loops, each edge of $\ga$ has
different labels on its two endpoints, which determines the labels at
$u_2$.  One can check that $\ga$ has three equidistance classes $c_1 =
\{B,E,G,I\}$, $c_2 = \{A,D,H\}$ and $c_3 = \{C,F,J\}$.  Since $\gb$ is
positive, each family belongs to an equidistance class; moreover, one
can check that the single edge is equidistant to the non-adjacent
family of weight 3, which we will denote $\hat e_1$.  Therefore these
must belong to $c_1$.  On $\ga$, $B$ has label $1$ at $u_2$, so on
$\gb$ $B$ has label $2$ at $v_1$.  It follows that $B$ is the middle
edge in $\hat e_1$.  This determines the labels on $v_1$ and $v_2$.
Now the endpoints of $E,G,I \in c_1$ are adjacent on $\bdd u_1$ among
edge endpoints labeled $1$, but they are not all adjacent on $\bdd
v_1$ because the single edge is not adjacent to those in $\hat e_1$
among edges with label $1$ at $v_1$ in $\gb$.  Therefore the jumping
number $J$ cannot be $\pm 1$, so $J=\pm 2$.  Reversing the orientation
of $v_1, v_2$ if necessary we may assume that $J = 2$.  Thus the edges
$E,G,I$ in $\gb$ must be as shown.  The other edges are now determined
by this information.  For example, the edges with label $2$ at $u_1$
appear in the order $B,D,F,H,J$, so on $\gb$ the edges with label
$1$ at $v_2$ appear in the order $B,H,D,J,F$.  \qed

\section {The case $n_1 = n_2 = 2$ and both $\Gamma_1, \Gamma_2$
non-positive }

In this section we assume that $n_1 = n_2 = 2$ and both $\Gamma_1,
\Gamma_2$ are non-positive.  Let $\ga = (\rho_a; a_1, ..., a_4)$, and
$\gb = (\rho_b; b_1, ..., b_4)$.  Without loss of generality we may
assume that $\rho_b \geq \rho_a$.

\begin{lemma} Suppose $n_1 = n_2 = 2$, and $\Gamma_1,
\Gamma_2$ are non-positive.

(1)  $\Delta/2 \leq \rho_b \leq 4$.

(2) $2\rho_a + \sum a_i = 2\Delta$, and $2\rho_b + \sum b_i =
    2\Delta$.

(3) $a_i, b_i \leq 2$.
\end{lemma}

\proof (1) Since a loop in $\ga$ corresponds to a non-loop in $\gb$
and vice versa, we have $\rho_a + \rho_b = \Delta$.  We have assumed
$\rho_b \geq \rho_a$, so $\rho_b \geq \Delta/2$.  Since no two edges
are parallel on both graphs and $\rga$ has at most four non-loop
edges, we also have $\rho_b \leq 4$.
 
(2) This follows from the fact that the valence of a vertex in $\ga$
or $\gb$ is $2 \Delta$.

(3) Since $\ga$ and $\gb$ are non-positive, a non-loop edge in $\ga$
is a loop in $\gb$, hence there are at most two edges in each non-loop
family of $\ga$, i.e.\ $a_i \leq 2$.  Similarly for $b_i$.
\qed

\bigskip
\leavevmode

\centerline{\epsfbox{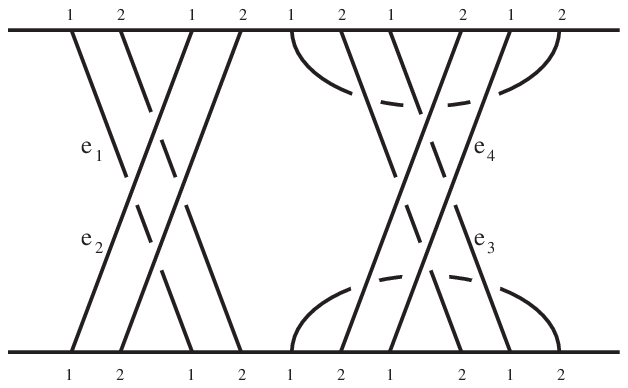}}
\bigskip
\centerline{Figure 20.1}
\bigskip

\begin{lemma} Suppose $n_1 = n_2 = 2$, $\Delta = 5$, and
$\Gamma_1, \Gamma_2$ are non-positive.  Then $\ga = -(2;2,2,2,0)$ and
$\gb = -(3;1,1,1,1)$ or $-(3;2,2,0,0)$.  \end{lemma}

\proof By Lemma 20.1(1) we have $\rho_b = 3$ or $4$.  If $\rho_b = 4$
then each loop family contributes one edge to each non-loop family of
$\ga$, hence $\ga = -(1;2,2,2,2)$.  The four loops $e_1, e_2, e_3,
e_4$ at $v_1$ are equidistant to each other; on the other hand, from
Figure 20.1 one can see that $e_i$ is equidistant to $e_j$ on $\ga$ if
and only if $e_i$ and $e_j$ are on the same side of the loop at $u_1$.
This is a contradiction.  Therefore this case cannot happen.

Now assume $\rho_b = 3$.  Then by the Congruence Lemma we have $\gb =
-(3;2,0,2,0)$, $-(3;2,2,0,0)$ or $-(3;1,1,1,1)$.  Since $\rho_a =
\Delta - \rho_b = 2$, we have $\sum a_i = 10 - 2 \rho_a = 6$.  By
Lemma 20.1 we have $a_i \leq 2$, therefore by the Congruence Lemma we
must have $\ga = -(2;2,2,2,0)$.  The first case for $\gb$ above cannot
happen because the two non-loop $1$-edges are not equidistant in $\gb$
while as parallel loops on $\ga$ they are equidistant on $\ga$.
Therefore $\gb = -(3;2,2,0,0)$ or $-(3;1,1,1,1)$.  
\qed

\bigskip
\leavevmode

\centerline{\epsfbox{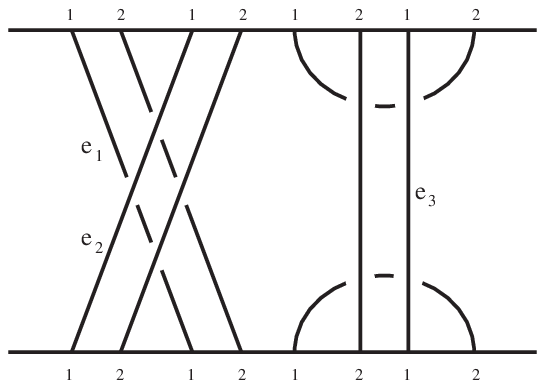}}
\bigskip
\centerline{Figure 20.2}
\bigskip

\begin{lemma} Suppose $n_1 = n_2 = 2$, $\Delta = 4$, and 
$\Gamma_1, \Gamma_2$ are non-positive.  Then one of the following
holds. 

(1) $\ga = -(2,2,2,2)$ and $\gb = -(4;0,0,0,0)$.

(2) Both $\Gamma_1$ and $\Gamma_2$ are of type $-(2;1,1,1,1)$.

(3) Both $\Gamma_1$ and $\Gamma_2$ are of type $-(2;2,0,2,0)$.

(4) Both $\Gamma_1$ and $\Gamma_2$ are of type $-(2;2,2,0,0)$.  
\end{lemma}

\proof If $\rho_b = 4$ then $\gb = -(4;0,0,0,0)$ and each loop family
contributes one edge to each family of $\ga$, hence $\ga =
-(2,2,2,2)$.

If $\rho_b = 3$, then by the Congruence Lemma 15.1 we have $\ga = -(1;
2,2,2,0)$.  Let $e_1, e_2, e_3$ be the three loops at $v_1$.  As
parallel positive edges, they are equidistant on $\gb$.  On $\ga$ they
are as shown in Figure 20.2.  One can check that $e_1$ is equidistant
to $e_2$ but not $e_3$, which is a contradiction to Lemma 2.17.
Therefore $\rho_b \neq 3$.

When $\rho_b=\rho_a = 2$, by the Congruence Lemma each of $\ga$ and
$\gb$ is of type $-(2;1,1,1,1)$ or $-(2;2,2,0,0)$ or $-(2;2,0,2,0)$.
We are done if both $\ga, \gb$ are of the same type.

If $\ga = -(2;2,2,0,0)$ then the two non-loop edges with label $1$ at
both endpoints are adjacent among the four edges labeled $1$ at $u_1$,
hence on $\gb$ the two loops at $u_1$ are adjacent among the four
edges with label $1$ at $v_1$, which implies that $\gb$ cannot be 
$-(2;2,0,2,0)$ or  $-(2;1,1,1,1)$. 

It remains to rule out the possibility that $\ga = -(2;1,1,1,1)$ and
$\gb = -(2;2,0,2,0)$.  In this case the graphs are as shown in Figure
20.3.

Label the edges of $\ga$ as in the figure.  We want to show that this
determines the labels of the edges of $\gb$ up to symmetry.  Since
$\Delta=4$, by changing the orientation of $\hat F_b$ if necessary we
may assume that the jumping number is $1$.  The $1$-edges at $u_1$ are
in the order $A,B,C,D$, so these labels appear in this order at
$v_1$ on $\gb$.  The order of the $1$-edges at $u_2$ is $A,X,C,Y$, so
the $2$-edges at $v_1$ are also in this order, which determines the
edges $X,Y$ on $\gb$.  Finally, the order of the $2$-edges at $u_1$
determines the edges $E,F$ in $\gb$.  Hence the labels of the graphs
are as shown in Figure 20.3.

One way to see that these graphs are not realizable is to consider the
annulus $A$ from $\bdd v_1$ to $\bdd v_2$ along the positive
orientation, draw the segments of $\bdd u_1, \bdd u_2$ on this annulus
and check that these arcs must intersect on $A$, which contradicts the
fact that $\bdd u_1, \bdd u_2$ are parallel curves on the torus $T_0$.
Here is another way.  Consider the endpoints of the edges $D,X$,
labeled $a,b,c,d$ on the two graphs.  We have $$d_{v_1}(a,c) =
d_{v_2}(b,d) = 1$$
so by Lemma 2.16 (applied with $u_i, u_j, v_k, v_l$
replaced by $v_1, v_2, u_1, u_2$ and $P,Q,R,S$ replaced by $a,c,b,d$),
we should have $$d_{u_1}(a,b) = d_{u_2}(c,d)$$
However, on $\ga$ we
have $d_{u_1}(a,b) = 5$ while $d_{u_2}(c,d) = 3$, which is a
contradiction.  \qed

\bigskip
\leavevmode

\centerline{\epsfbox{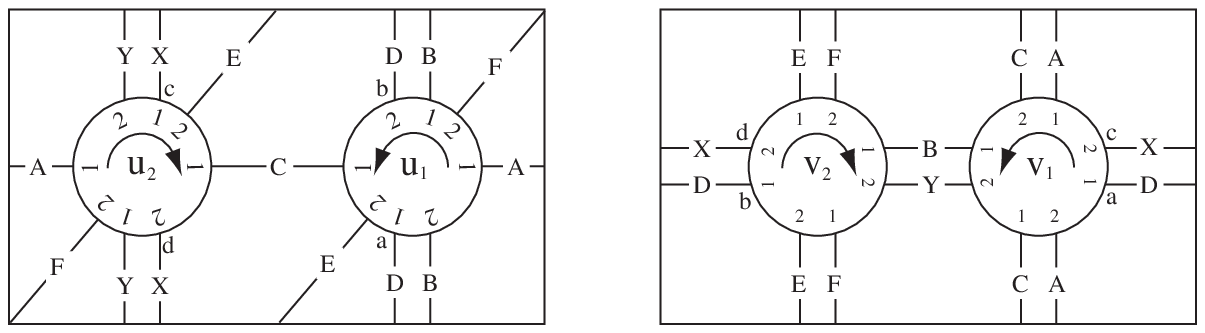}}
\bigskip
\centerline{Figure 20.3}
\bigskip

\begin{prop}  Suppose $n_a, n_b \leq 2$.  Then up to
symmetry $\ga$ and $\gb$ are one of the following pairs.
$$ 
\begin{array}{lll}
(1) \qquad \qquad \qquad  & -(1,1,1,1) \qquad \qquad \qquad   &(4,0,0) \\ 
(2) \qquad \qquad \qquad  & -(2,2,0,0)     & (2,2,0)      \\  
(3) \qquad \qquad \qquad  & -(2,1,1,1)     & (3,1,1)      \\  
(4) \qquad \qquad \qquad  & -(2,2,2,2)    &+(4;0,0,0,0)   \\ 
(5) \qquad \qquad \qquad   & -(3,3,3,1)    &+(3,3,3,1)    \\ 
(6) \qquad \qquad \qquad   & -(2;2,2,2,0)  &-(3;1,1,1,1)  \\ 
(7) \qquad \qquad \qquad   & -(2;2,2,2,0)  &-(3;2,2,0,0)  \\ 
(8) \qquad \qquad \qquad   & -(2,2,2,2)    &-(4;0,0,0,0)  \\ 
(9) \qquad \qquad \qquad   & -(2;1,1,1,1)  &-(2;1,1,1,1)  \\ 
(10) \qquad \qquad \qquad   & -(2;2,0,2,0)  &-(2;2,0,2,0)  \\ 
(11) \qquad \qquad \qquad   & -(2;2,2,0,0)  &-(2;2,2,0,0)    
\end{array}
$$
\end{prop}

\proof This follows from the lemmas in Sections 18--20.  More
precisely, the case $n_b = 1$ is done in Lemma 18.1, which gives
(1)--(2) above; the case $n_a=n_b=2$ and $\gb$ positive is discussed
in Lemmas 19.2--19.7 according to different numbers of loops on $\gb$,
which gives (3)--(5); the case $n_a=n_b=2$ with both graphs
non-positive is discussed in Lemmas 20.2--20.3, with the possibilities
listed in (6)--(11).  \qed

\begin{prop} For each of the cases (3), (5), (6), (9)
and (10) of Proposition 20.4, the correspondence between edges of
$\ga, \gb$ is unique up to symmetry, and is shown in Figures 18.2,
19.4, 20.4, 20.5 and 20.6, respectively.  \end{prop}

\proof
For cases (3) and (5) this follows from Lemmas 18.1 and 19.7.

In case (6) we have $\ga = -(2;2,2,2,0)$ and $\gb = -(3;1,1,1,1)$.
The graphs $\ga, \gb$ are as shown in Figure 20.4.  Label the edges of
$\gb$ as shown in the figure.  By symmetry we may assume that the
labels on the edge endpoints of $\gb$ are as in the figure.  Also up to
symmetry of $\gb$ on the torus $\hat F_b$ we may assume that the
labels on $v_1$ are as in the figure.

The label 1 endpoints of $A,B,C$ are non-adjacent among the $1$-labels
on $\bdd v_1$.  These are non-loops on $\ga$, and one in each family,
hence their endpoints at $u_1$ are also non-adjacent among endpoints
labeled 1.  This forces the jumping number $J$ to be $\pm 1$.  Now on
$\ga$ the edge $A$ must be as shown.  It is easy to see that this
determines the labels on the other edges in $\ga$.

\bigskip
\leavevmode

\centerline{\epsfbox{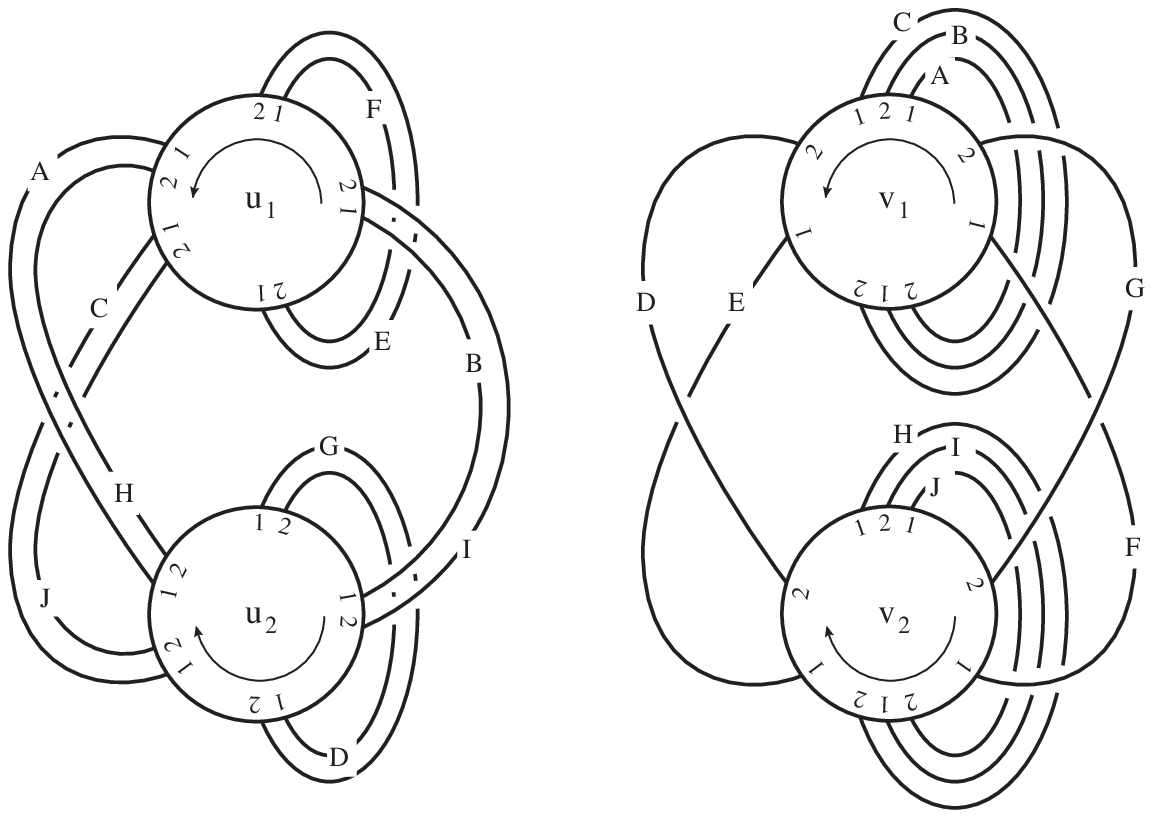}}
\bigskip
\centerline{Figure 20.4}
\bigskip

In case (9) we have $\ga = -(2;1,1,1,1)$, $\gb = -(2;1,1,1,1)$, and
$\Delta = 4$, so we may assume $J=1$.  Label edge endpoints and edges
of $\ga$ as in the figure.  Using symmetry we may assume $A$ to be any
one of the two non-loop edges labeled $1$ at $v_1$.  Then this
determines the labels on the other edges.  See Figure 20.5.

The determination of the edge correspondence for case (10) is similar.
The graphs are shown in Figure 20.6.  \qed

\bigskip
\leavevmode

\centerline{\epsfbox{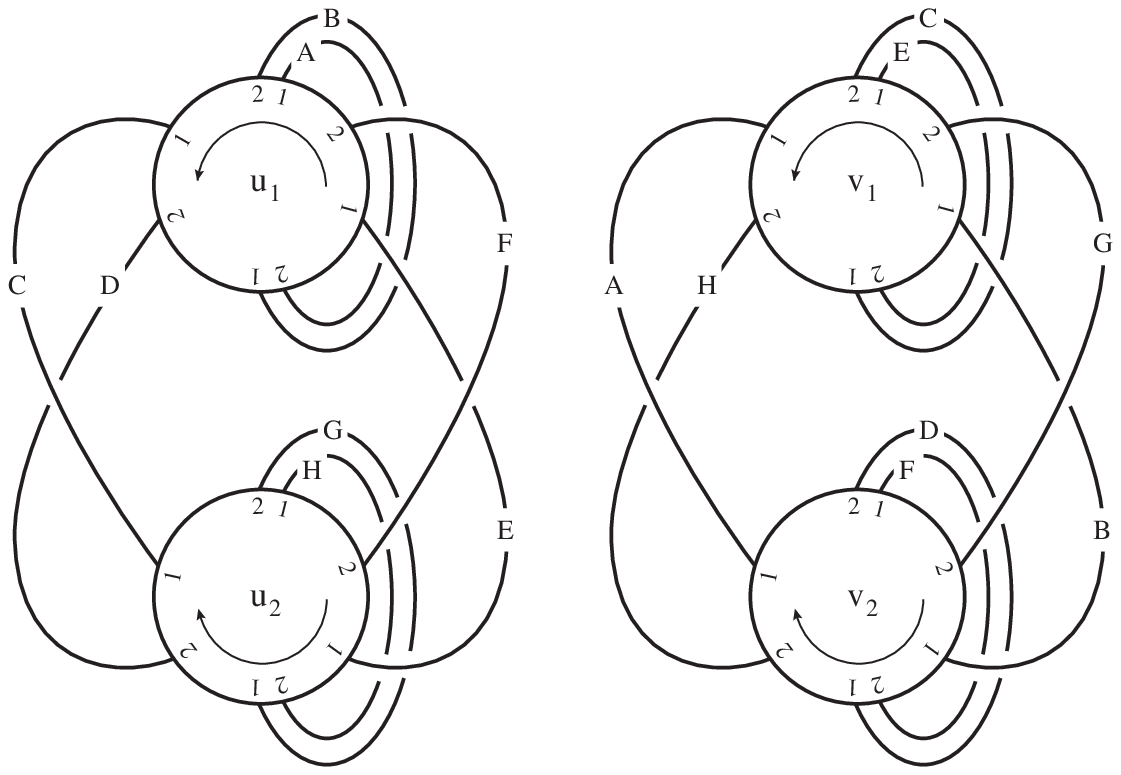}}
\bigskip
\centerline{Figure 20.5}
\bigskip

\bigskip
\leavevmode

\centerline{\epsfbox{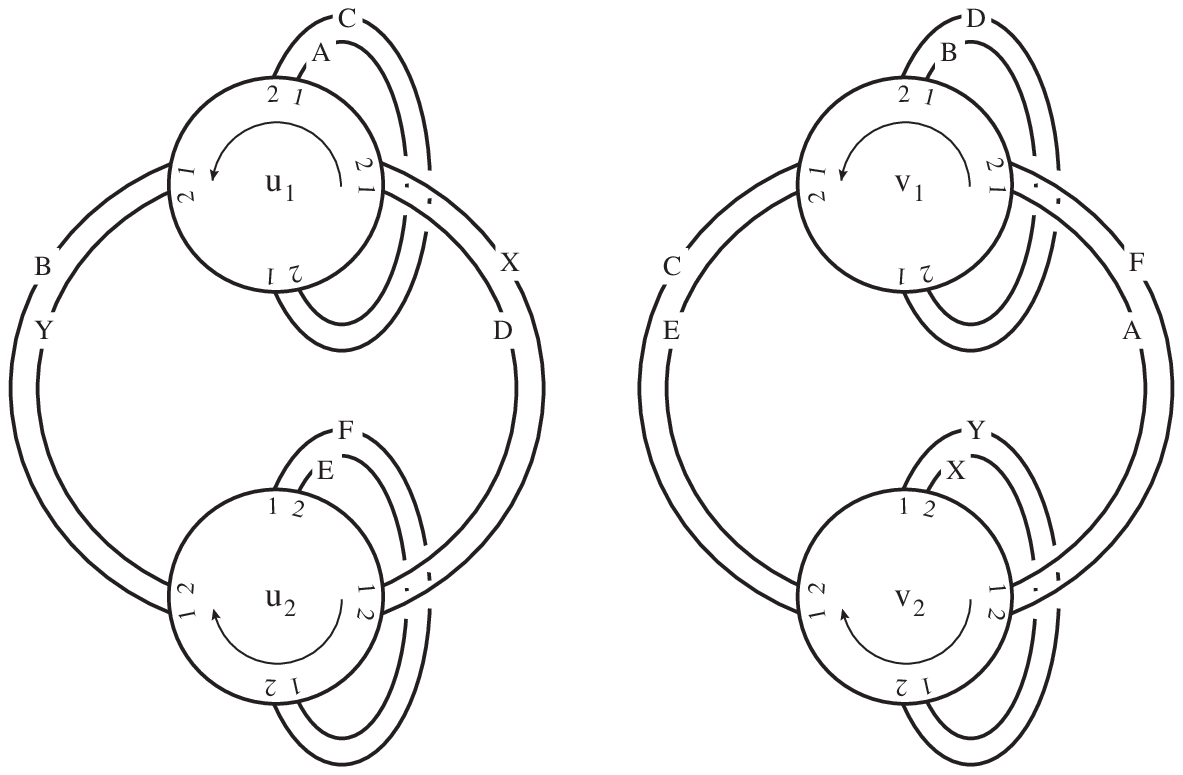}}
\bigskip
\centerline{Figure 20.6}
\bigskip

\section  {The main theorems }

Suppose $M$ is a hyperbolic manifold admitting two toroidal Dehn
fillings $M(r_1), M(r_2)$.  Let $F_a$ be essential punctured tori in
$M$ such that $\bdd F_a$ consists of a minimal number of copies of
$r_a$, and $F_1$ intersects $F_2$ minimally.  Let $X(F_1, F_2)$ be
obtained from $N(F_1 \cup F_2 \cup T_0)$ by capping off its 2-sphere
boundary components with 3-balls.  We will use $X(r_1, r_2)$ to denote
any $X(F_1, F_2)$ above with $\bdd F_a$ of slope $r_a$, and call it a
{\it core\/} of $M$ with respect to the toroidal slopes $r_1, r_2$.
Note that $X(r_1, r_2)$ may not be unique.

\begin{lemma}  Suppose $M$ is a hyperbolic manifold
admitting two toroidal Dehn fillings $M(r_1), M(r_2)$ of distance $4$
or $5$.  Then each of $\bdd X(r_1, r_2)$ and $\bdd M$ is a union of
tori.
\end{lemma}

\proof By the result of the previous sections we see that $\ga, \gb$
are either the graphs in Figures 11.9, 11.10, 14.5, 16.6, 16.8, 16.9,
or one of the pairs given in Proposition 20.4.  

In all figures except 11.9 and 11.10, $\ga$ has two vertices and they
have opposite signs.  Now $X(F_1, F_2)$ can be constructed by adding
thickened faces of $\gb$ to $N(F_a \cup T_0)$, which has two boundary
components of genus 2.  It is easy to check that in all cases $\gb$
has at least one disk face on each side of $F_a$.  The boundary of a
disk face of $\gb$ is always an essential curve on $F_a \cup T_0$.
Adding a 2-handle corresponding to a disk face will change a genus 2
boundary component to one or two tori.  It follows that the boundary
of $X(F_1, F_2)$ is a union of tori.  Since $M$ is irreducible and
atoroidal, each torus boundary component of $X(F_1, F_2)$ either is
boundary parallel, or bounds a solid torus.  Therefore $\bdd M$ is
also a union of tori.

The proof for Figures 11.9 and 11.10 is similar.  In these cases $\ga$
has 4 vertices, so $N(F_a \cup T_0)$ has two boundary components $S_i$
of genus 3.  It suffices to find two faces on each side of $F_a$ whose
boundary curves give rise to non-parallel and non-separating curves on
$S_i$.  For 11.9 one can check that the bigons on $F_b$ bounded by the
edges $K\cup J$ and $H\cup G$ are on the same side of $F_a$ and give
non-parallel boundary curves on $S_1$, say, while the bigon bounded by
the edges $J\cup H$ and the $3$-gon bounded by the edges $K \cup L
\cup R$ give non-parallel non-separating curves on $S_2$.  Hence the
result follows.  For 11.10, use the bigons bounded by $E \cup
F$ and $G \cup H$ on one side, and the bigon $F \cup G$ and the 3-gon
$E \cup K \cup A$ on the other side.
\qed

Consider the three manifolds $M_1, M_2, M_3$ in [GW1, Theorem 1.1].
More explicitly, $M_1$ is the exterior of the Whitehead link, $M_2$ is
the exterior of the 2-bridge link associated to the rational number
$3/10$, and $M_3$ is the exterior of the $(-2, 3,8)$ pretzel link,
also known as the Whitehead sister link.  Each of these manifolds
admits two Dehn fillings $M_i(r_1)$ and $M_i(r_2)$, both toroidal and
annular, with $\Delta = 4$ for $i=1,2$, and $\Delta = 5$ for $i=3$.
Let $T_0$ be the Dehn filling component of $\bdd M_i$, and let $T_1$
be the other component of $\bdd M_i$.  Then for all except a few
slopes $s$ on $T_1$, $M' = M_i(s)$ is a hyperbolic manifold, and it
admits two toroidal Dehn fillings $M'(r_1), M'(r_2)$ of distance $4$
or $5$.  The following lemma shows that several of the cases in
Proposition 20.4 can only be realized by these manifolds.

\begin{lemma} Suppose $\ga$ has a non-disk face.  Then $M =
M_i$ or $M_i(s)$ for some $i=1,2,3$ and slope $s$ on $T_1$, and the
toroidal slopes $r_1, r_2$ are the same as the toroidal/annular slopes
given in [GW1, Theorem 1.1].
\end{lemma}

\proof Let $K$ be a curve on $F_a$ which is essential on $\hat F_a$
and disjoint from $F_a \cap F_b$.  Consider the manifold $X = M-\Int
N(K)$.  If $X$ is hyperbolic then $\hat A = \hat F_a - \Int N(K)$ is
an essential annulus in $X(r_a)$ and $\hat F_b$ is an essential torus
in $X(r_b)$, so by [GW1, Theorem 1.1] $X = M_i$ for some $i=1,2,3$,
and we are done.  Hence we may assume that $X$ is non-hyperbolic.  $X$
is irreducible as otherwise there would be an essential sphere $S$ in
$X$ bounding a 3-ball in $M$ containing $K$, which would be a
contradiction to the fact that $K$ is an essential curve on $\hat
F_a$.  Also $X$ cannot be a Seifert fibered manifold as otherwise $M =
X \cup N(K)$ would be non-hyperbolic.  Since by Lemma 21.1 $\bdd M$ is
a union of tori, the above implies that $X$ must be toroidal.

Since $M$ is atoroidal, an essential torus $T$ in $X$ must be
separating.  Let $T_0$ be the Dehn filling torus component of $M$, and
let $T_1 = \bdd N(K)$.  Recall that $M|T$ denotes the manifold
obtained by cutting $M$ along $T$.  Let $V = V_T$ and $W = W_T$ be the
components of $M|T$, where $W$ is the component containing $T_0$.
Among all essential tori in $X$, choose $T$ so that (a) if there is
some $T$ in $X$ such that $T_1 \subset W_T$, choose $T$ so that $V_T$
contains no essential torus; (b) if every essential torus in $X$
separates $T_0$ from $T_1$, choose $T$ such that $W_T$ contains no
essential torus.

Since $M$ is atoroidal, $T$ is inessential in $M$, hence $V$ is either
(i) a solid torus, or (ii) $T^2 \times I$, or (iii) a 3-ball with a
knotted hole.  Note that in the first two cases $V$ must contain the
curve $K$.  Let $N = V - \Int N(K)$ in the first two cases, and $N =
V$ in the last case.  Let $C = T \cap F_a$.  Using a standard cut and
past argument we may assume that each component of $C$ is essential on
both $T$ and $F_a$.  In case (iii) let $D$ be a compressing disk of
$T$ in $W$.

\medskip 
{\bf Claim 1.}  {\it $C \neq \emptyset$.}
\medskip

\proof If $C = \emptyset$ then $F_a$ lies in $W$, which is
impossible in cases (i) or (ii) because the curve $K$ on $F_a$
lies in $V$.  In case (iii) $D \cap F_a$ is a set of circles and
one can use the incompressibility of $F_a$ in $W$ to isotop $F_a$ so
that it is disjoint from $D$.  But then $D$ is disjoint from $K$, so
$T$ would be compressible in $X$, which is a contradiction.
\qed

\medskip
{\bf Claim 2.}   {\it $C$ is a set of essential curves on $\hat F_a$
parallel to $K$.}
\medskip

\proof Since $C$ is disjoint from $K$, we need only show that each
component $\alpha$ of $C$ is an essential curve on $\hat F_a$.  Assume
to the contrary that $\alpha$ bounds a disk $E$ on $\hat F_a$ and is
innermost on $\hat F_a$.  Then $E$ must contain some boundary
component of $F_a$, hence $E \subset W(r_a)$.  In case (i) $V \cup
N(E)$ is either a 3-ball, or a punctured lens space or $S^1 \times
S^2$, containing the curve $K$, contradicting the fact that $\hat F_a$
is incompressible and $M(r_a)$ irreducible.  In case (ii) $V \cup
N(E)$ is a punctured solid torus, so the irreducibility of $M(r_a)$
implies that $M(r_a)$ is a solid torus, which is absurd because it is
supposed to be toroidal.  In case (iii), for homological reasons $\bdd
E$ and $\bdd D$ must be homotopic on $T$, hence $\bdd E$ is
null-homotopic in $M$, which contradicts the facts that $C$ is
essential on $F_a$ and $F_a$ is incompressible in $M$.  \qed

\medskip {\bf Claim 3.}  {\it Case (iii) cannot happen, i.e.\ 
$V$ is not a 3-ball with a knotted hole.}  
\medskip

\proof We have shown that all components of $C$ are essential curves
on $\hat F_a$ parallel to $K$, and $C \neq \emptyset$.  Let $\alpha$ be
a component of $C$.  Then $K$ is isotopic to $\alpha$ in $M(r_a)$,
but since $\alpha \subset T$ lies in the 3-ball $V \cup N(D)$,
$\alpha$, and hence $K$, is null-homotopic in $M(r_a)$, which
contradicts the fact that $\hat F_a$ is incompressible in $M(r_a)$.
\qed

\medskip {\bf Claim 4.}  {\it $W$ is hyperbolic.}
\medskip

\proof Clearly $W$ is irreducible (since $X$ is) and not a Seifert
fibered space (since $M$ is hyperbolic).  Suppose $W$ contains an
essential torus $T'$.  By Claim 3 we see that $T'$ cannot be of type
(iii), so it must be of type (i) or (ii), which, by our choice of $T$,
implies that every essential torus in $X$ separates $T_0$ from $T_1$.
By the choice of $T$, $W$ must be atoroidal.  \qed

We now continue with the proof of Lemma 21.2.  Let $A$ be a component of
$F_a \cap W$ which contains some boundary components of $F_a$.  By
Claims 1 and 2, the corresponding component $\hat A$ of $\hat F_a \cap
W(r_a)$ is an annulus in $W(r_a)$, which is incompressible because
$\hat F_a$ is incompressible, and not boundary parallel because
otherwise $\hat F_a$ would be isotopic to a torus with fewer
intersections with the Dehn filling solid torus.  Therefore $W(r_a)$
is annular.

Let $P$ be the component of $F_a \cap V$ containing $K$, and let
$\beta$ be a component of $P \cap T$.  Note that $P$ is an annulus.
Since $F_b$ is disjoint from $K$, it can be isotoped to be disjoint
from $P$, hence after isotopy we may assume that $F_b \cap T$ and $F_a
\cap T$ are all parallel to $\beta$ and hence mutually disjoint.  If
$F_b \cap T = \emptyset$ then $\hat F_b$ is an essential torus in
$W(r_b)$, and if $F_b \cap T \neq \emptyset$ then as above, a
component of $\hat F_b \cap W(r_b)$ which intersects the Dehn filling
solid torus is an essential annulus in $W(r_b)$, hence $W(r_b)$ is
either toroidal or annular.  Using Theorem 1.1 of [GW1] in the first
case and Theorem 1.1 of [GW3] in the second case, we see that $W =
M_i$ for $i=1$, $2$, or $3$.

By Claim 3 $V$ is either a solid torus or $T^2 \times I$.  In the
first case $M = M_i(s)$ for some $s$ on $T_1 = \bdd V$, and in the
second case $M = M_i$.
\qed

\begin{defn} 
(1) Define a set of triples $(M_i, r'_i, r''_i)$ as follows.  For
$i=1,2,3$, $(M_i,r'_i,r''_i)$ are the manifolds and the
toroidal/annular slopes given in Theorem 1.1 of [GW1].  $M_4, ...,
M_{14}$ are the manifolds $X(F_1, F_2)$ corresponding to the
intersection graphs given in Figures 11.9, 11.10, 14.5, 16.6, 16.8,
16.9, 18.2, 19.4, 20.4, 20.5 and 20.6, and $r'_i,r''_i$ are the
boundary slopes of the corresponding surfaces $F_1, F_2$.  

(2) Two triples $(M, r', r'')$ and $(N, s', s'')$ are {\it
equivalent}, denoted by $(M,r',r'') \cong (N,s',s'')$, if there is a
homeomorphism from the 3-manifold $M$ to $N$ which sends the boundary
slopes $(r', r'')$ to $(s', s'')$ or $(s'', s')$.
\end{defn}

The following theorem shows that if a hyperbolic manifold $M$ admits
two toroidal Dehn fillings along slopes $r_1, r_2$ of distance 4 or 5
then $(M, r_1, r_2)$ is either one of these triples, or obtained from
such an $M_i$ by Dehn filling on $\bdd M_i - T_0$.

\begin{thm}
  Let $M$ be a hyperbolic 3-manifold admitting two toroidal Dehn
  fillings $M(r_1), M(r_2)$ with $\Delta(r_1, r_2) = 4$ or $5$.  Let
  $n_a$ be the minimal number of intersections between essential tori
  and the Dehn filling solid torus in $M(r_a)$.  Assume $n_a \leq
  n_b$.  Let $(M_i, r'_i, r''_i)$ be the manifolds defined above, and
  let $T_0$ be the boundary component of $M_i$ containing $r'_i,
  r''_i$.  Then

  (1) $n_a \leq 2$, $n_b \leq 6$;

  (2) either $(M,r_1,r_2) \cong (M_i, r'_i, r''_i)$ for some $i=1,...,14$,
  or $(M,r_1,r_2) \cong (M_i(s), r'_i,r''_i)$, where $i\in
  \{1,2,3,14\}$ and $s$ is a slope on $T_1 = \bdd M_i - T_0$; and

  (3) $i \in \{1,2,4,6,9,13,14\}$ if $\Delta = 4$, and $i \in
  \{3,5,7,8,10,11,12\}$ if $\Delta = 5$.
\end{thm}

\proof This is a summary of the results in the previous sections.
Assume $n_a \leq n_b$.  Then by Proposition 11.10 we have $n_a \leq 2$.
By Proposition 16.8 if $n_b \geq 3$ then $X$ is one of those in Figure
11.9, 11.10, 14.5, 16.6, 16.8 or 16.9.

We may now assume $n_a, n_b \leq 2$.  Then by Proposition 20.4 $\ga,
\gb$ is one of the 11 pairs listed there.  One can check that all but
cases (3), (5), (6), (9), (10) have the property that one of $\hat
F_a, \hat F_b$ contains a non-disk face, so by Lemma 21.2 the triple
$(M,r_1,r_2)$ is $(M_i,r'_i, r''_i)$ for some $i=1,2,3$.  Finally, by
Proposition 20.5 the graphs of the above cases are given in Figures
18.2, 19.4, 20.4, 20.5 and 20.6.

(3) follows by counting $\Delta$ for the graph pairs of each of the
manifolds listed in (2).  \qed

\section {The construction of $M_i$ as a double branched cover}

The first three of the 14 manifolds $M_i$ have already been identified
as the exteriors of links in $S^3$.  See [GW1].  The links are shown in
Figure 24.1 Besides $M_4$ and $M_5$, the other nine manifolds $M_6,
..., M_{14}$ have the property that $\ga$ is a graph on $\hat F_a$
with two vertices of opposite signs.  In this section we will
construct, for each $i=6,...,14$, a tangle $Q_i = (W_i, K_i)$, where
$W_i$ is a 3-ball for $i=6,...,13$, and an $S^2 \times I$ for $i=14$,
such that $M_i$ is the double branched cover of $W_i$ with branch set
$K_i$.  It is well known that once we have such a presentation then
the Dehn filling $M_i(r)$ will be the double branched cover of
$Q_i(r)$, where $Q_i(r)$ is obtained by attaching a rational tangle of
slope $r$ to $Q_i$, with coordinates properly chosen.

Here is a sketch of the construction.  Assume $\ga$ is non-positive
and $n_a = 2$, and suppose there is an orientation-preserving
involution $\alpha_1$ on $F_a$ which maps $\bdd u_1$ to $\bdd u_2$ and
preserves $\ga$.  The restriction of $\alpha_1$ on $\bdd F_a$ extends
to an involution $\alpha_2$ on $T_0$ which has four fixed points, and
it preserves the curves $\bdd F_b$ on $T_0$.  Thus $\alpha = \alpha_1
\cup \alpha_2$ is an involution on $F_a \cup T_0$, which has eight
fixed points, four on each of $F_a$ and $T_0$.  Since $\alpha$
preserves $\ga \cup \bdd F_b$, it extends over each disk face of $F_b$
to give an involution on $F_b$.  One can now further extend the
involution $\alpha$ from $F_a \cup F_b \cup T_0$ to a regular
neighborhood $Y$ of $F_a \cup F_b \cup T_0$.  For $i\geq 6$, $M_i$ is
obtained by capping off spherical boundary components of $Y$ by
3-balls, hence $\alpha$ extends to an involution of $M_i$.  Clearly
the quotient of $N(F_a \cup T_0)$ is a twice punctured 3-ball $W_i$.
After attaching 2-handles corresponding to faces of $F_b$ and some
3-balls we see that $W_i$ is a punctured 3-ball.  From the
construction below we will see that $W_i$ is a 3-ball when $i = 6,
..., 13$, and an $S^2 \times I$ when $i=14$.  Denote by $K_i$ the
branch set of $\alpha$ in $W_i$.  Then $Q_i = (W_i, K_i)$ is the
tangle corresponding to the manifold $M_i$, and $M_i$ is the double
branched cover of $Q_i$ in the sense that it is the double branched
cover of $W_i$ with branch set $K_i$.  Attaching a rational tangle of
slope $t$ to $T_i$, we obtain a new tangle $Q_i(t)$ whose double
branched cover is $M(r)$ for some slope $r$ on $T_0 \subset \bdd M$.
This makes it possible to see the essential torus in $M_i(r_a)$ as a
lifting of some surface in $Q_i(t)$.

To illustrate this procedure, we give below a step by step
construction of the tangle $Q_6 = (W_6, K_6)$ for the manifold $M_6$
corresponding to the graphs in Figure 14.5.  The constructions for the
other manifolds are similar.

Denote by $N(C)$ a regular neighborhood of a set $C$ in a 3-manifold,
and by $I$ the interval $[-1, 1]$.

\medskip

STEP 1.  {\it Identify $[N(F_a \cup T_0)/\alpha] - D_2 \times I$ with
$S^2 \times I$.}

\medskip

Recall that $\alpha$ has four branch points on each of $T_0$ and
$F_a$, so $T_0/\alpha = S$ is a 2-sphere, and $F_a/\alpha = D_1$ is a
disk.  Let $D_2$ be a small disk in the interior of $D_1$, disjoint
from $\gb/\alpha$ and the branch points of $\alpha$.  Then $A_1 = D_1
- \Int (D_2)$ is a collar of $\bdd D_1$.  Therefore $N( (F_a \cup
T_0)/\alpha)$ can be written as $$(S \times I) \cup (A_1 \times I)
\cup (D_2 \times I)$$
Note that $A_1 \times I$ is a collar of the
attaching annulus $\bdd D_1 \times I$, hence $X = (S\times I) \cup
(A_1 \times I)$ is homeomorphic to $S^2 \times I$.

One boundary component of $X = S^2 \times I$ is $\bdd_- X = (T_0\times
\{-1\})/\alpha$, and the other boundary component $\bdd_+ X$ can be
written as $D_+ \cup D_- \cup A$, where the two disks $D_+ \cup D_- =
\bdd X \cap (S \times 1)$ lift to two annuli on $T_0 \times 1$, and
$A$ is the annulus $\bdd X \cap (A_1 \times I)$.  We identify $X$ with
$(\Bbb R^2 \cup \{\infty\}) \times I$, so that the disks $D_{\pm}$ are
identified with the squares $I \times [\pm 2, \pm 4]$ on the plane $P
= \Bbb R^2 \times 1$ on $\bdd X$, the annulus $A$ is the closure of
$P- D_+ \cup D_-$, and the core $c_0$ of $A$ is identified with the
closure of the $x$-axis of $P$.  See Figure 22.1.  (Not drawn to
scale.)

\bigskip
\leavevmode

\centerline{\epsfbox{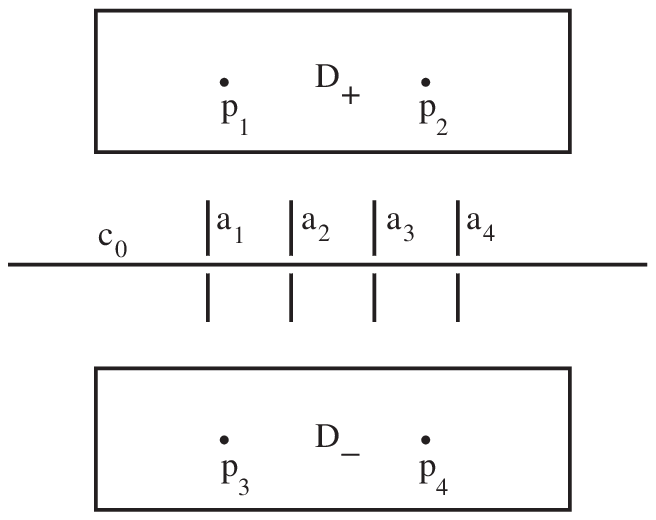}}
\bigskip
\centerline{Figure 22.1}
\bigskip

The branch set of $\alpha$ now consists of eight arcs.  Four of them
come from the fixed points of $\alpha$ on $T_0 \times I$, and are of
type $p_i \times I$ in $X$, where $p_1, p_2 \in D_+$, and $p_3, p_4
\in D_-$.  These will be represented by four dots $p_1, ..., p_4$ on
$P$, two in each of $D_{\pm}$, as shown in Figure 22.1.  The other
four branch arcs of $\alpha$ are of type $a_i = q_i \times I \subset A
\times I$, where $q_1, ..., q_4$ are the branch points of $\alpha$ in
$A = D_1 - \Int (D_2)$.  Note that $a_i$ has both endpoints on $S^2
\times 1$.  We may assume that these project to four vertical arcs on
the annulus $A$ in $P$ above, and we may arrange so that the endpoints
of these arcs have $y$-coordinates $\pm 1$ on the plane $P$.  See
Figure 22.1.  Denote by $a_i(1)$ and $a_i(-1)$ the endpoints of $a_i$
with $y$ coordinates $1$ and $-1$, respectively.

STEP 2.  {\it Draw the arcs $G' = (\ga \times \{\pm 1\})/\alpha$ on
$P$, with edges and edge endpoints labeled.}

\medskip

The graph $\ga \times 1$ on $F_a \times 1$ projects to a set of arcs
$E$ on $D_1 \times 1$.  We may choose the disk $D_2$ above to be
disjoint from $E$.  Then $E$ lies in the annulus $A_+ = A \times 1$.
If a family $\hat e_i$ has $2k$ edges then they project to $k$ edges
on $A_+$ with endpoints on $\bdd D_+$, each circling around the branch
point $a_i(1)$.  If $\hat e_i$ has $2k+1$ edges then the quotient is a
set of $k$ edges as above together with an edge connecting a point on
$\bdd D_+$ to $a_i(1)$.  Up to isotopy we may assume that all edge
endpoints of $E$ on $\bdd D_+$ lie on the horizontal line $y=2$ on
$P$.  Similarly the projection of $\ga \times (-1)$ is a set of arcs
on the annulus $A_- = A \times (-1)$, which is the mirror image of
the arcs $(\ga \times 1)/\alpha$ along the circle $c_0$ on $P$.
Denote by $G'$ the set of arcs above.

For the graph $\ga$ in Figure 14.5(a), the edges in $G'$ are shown in
Figure 22.2.  The edges are labeled by the corresponding edges in
$\ga$.  (We only show a few of the labels in the figure; the others
should be easy to identify.)  Each edge in $G$ is the image of two
edges in $\ga$, hence it has two labels.  (Note that if one of the
families has an odd number of edges then the middle one projects to an
arc in $G'$ with a single label.)  All arcs appear in the region $I
\times [-2,2]$.  The top and bottom lines in the figure represent arcs
on $\bdd D_{\pm}$.  Note that each edge endpoint on $\bdd D_+$
corresponds to one edge endpoint on each of $\bdd u_1$ and $\bdd u_2$.
The labels on the top and bottom lines correspond to the labels on
$\bdd u_1$ in Figure 14.5(a).

\bigskip
\leavevmode

\centerline{\epsfbox{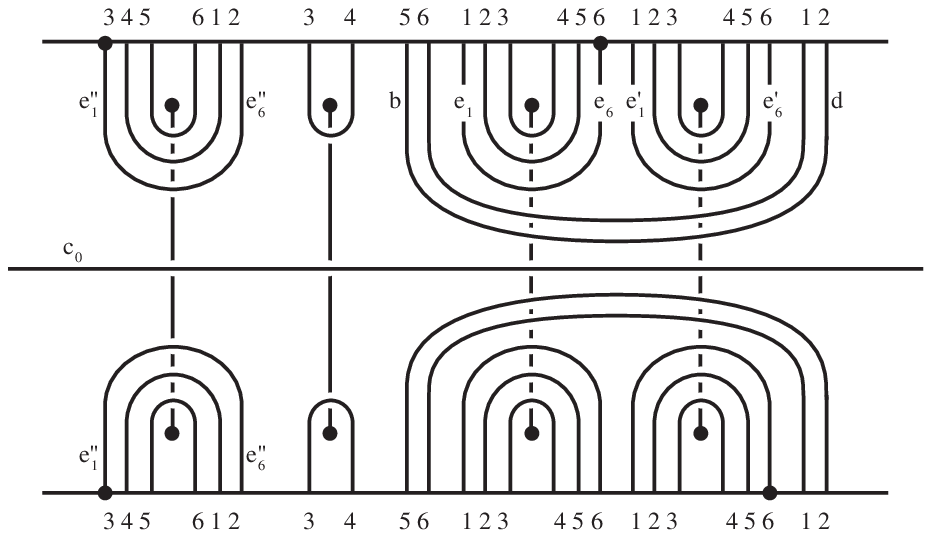}}
\bigskip
\centerline{Figure 22.2}
\bigskip

STEP 3.  {\it Add arcs $G'' = (\bdd F_b \times 1)/\alpha$ on $P$ to
obtain $G = G' \cup G''$.}

\medskip

Recall that the preimage of $D_{\pm}$ are two annuli $A_{\pm}$ on $T_0
\times 1$.  The curves $G'' = (\bdd F_b \times 1 \cap A_{\pm})/\alpha$
is now a set of arcs in $D_{\pm}$, with $\bdd G''$ the union of $\bdd
G'$ and possibly some of the branch points $p_i$ in $D_{\pm}$.  We
need to determine how the endpoints of $G'$ are connected by the edges
of $G''$.

Consider the circle $\bdd v_6$ in Figure 14.5(b).  We may assume that
the segments on $\bdd v_6$ from label 1 to label 2 (in the
counterclockwise direction) project to arcs in $D_+$ while those from
label $2$ to label $1$ project to arcs in $D_-$.  Consider the arc
$\beta$ on $\bdd v_6$ from the tail of $e_6$ to the head of $e''_6$.
(Recall that $e''_6$ is the edge between $e_6$ and $e'_6$ in Figure
14.5(a)).  Note that the tail of $e_6$ projects to the endpoint of
$e_1=e_6$ with label $6$ on $\bdd D_+$ in Figure 22.2.  The other
endpoint $q$ of $\beta$ is the head of $e''_6$, which lies on $\bdd
u_2$.  Since the labels in Figure 22.2 are the ones corresponding to
those on $\bdd u_1$ in Figure 14.5(a), we have to find the
corresponding point on $\bdd u_1$ in order to determine the position
of the edge endpoint $q$ on $\bdd D_+$.  On Figure 14.5(a) the
involution $\alpha$ restricted to $\bdd u_2$ is a vertical
translation, which maps the head of $e''_6$ (i.e.\ the edge in $\hat
e_3$ labeled $6$ at $u_2$) to the tail of $e''_1$, which has label $3$
at $u_1$.  It follows that $q$ is the endpoint of $e''_1=e''_6$ in
Figure 22.2 with label $3$ at $\bdd D_+$.  The two endpoints of
$\beta$ are represented by the two dots on the top line in Figure
22.2.  Similarly, let $\beta'$ be the arc on $\bdd v_6$ from the head
of $e''_6$ to the tail of $e'_6$.  Then it is an arc in $D_-$ with
endpoints on the dots at the bottom line in Figure 22.2.

The arcs $G''$ in $D_{\pm}$ are parallel to each other, and they are
non-trivial in the sense that none of them cuts off a disk in
$D_{\pm}$ that does not contain a branch point of $\alpha$.  Therefore
the above information completely determines the arcs $G''$ as well as
the branch points $p_1, ..., p_4$ of $\alpha$.  (Note that if the
number of edge endpoints between the dots is odd then the middle arc
will have an endpoint on a branch point $p_i$.)  The graph $G = G'
\cup G''$ is now shown in Figure 22.3.

\bigskip
\leavevmode

\centerline{\epsfbox{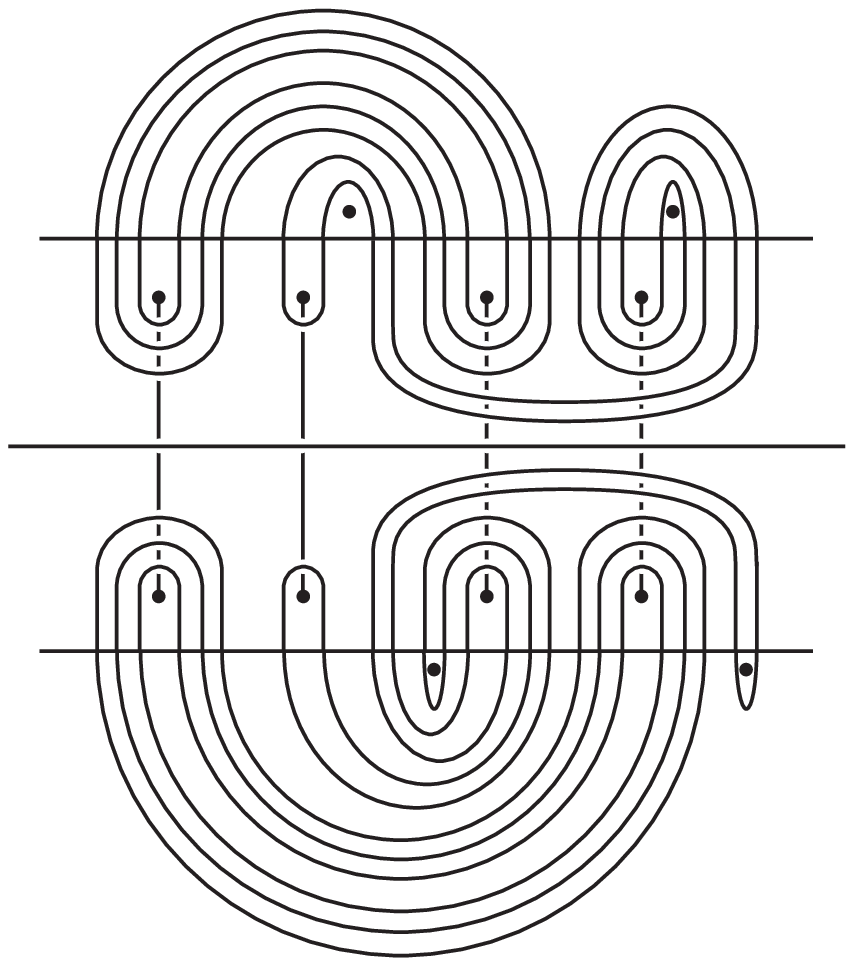}}
\bigskip
\centerline{Figure 22.3}
\bigskip

STEP 4.  {\it Construct the tangle $Q_i$.}

\medskip

Each component $c \neq c_0$ of $G$ lifts to curves on $F_b$ bounding
disk faces $\sigma$ of $\gb$.  The quotient of $\sigma \times I$ is
either a 2-handle attached to $c$ if $c$ is a circle, or a 3-ball
attached to a neighborhood of $c$ if $c$ is an arc.  Examining the
branch set of $\alpha$ in $\sigma \times I$ gives the following
procedure.  We use $X$ to denote the initial manifold at the beginning
of each step below.  In particular, $X = S^2 \times I$ before the
first step.

(1) If $c$ is an arc then $X \cup (\sigma \times I)/\alpha$ is
homeomorphic to $X$.  The new branch set is obtained by adding a
trivial arc in $(\sigma\times I)/\alpha$ joining the two endpoints of
$c$.  Therefore we can simply modify the branch set of $\alpha$ by
pushing $c$ into the interior of $X$.  

(2) If $c \neq c_0$ is a circle component bounding a disk $D_1$ on
$\bdd X$ containing no branch point of $\alpha$, then attaching a
2-handle along $c$ creates a 2-sphere boundary component, which must
bound a 3-ball in $M_i/\sigma$.  Thus after attaching the 2-handle and
the 3-ball the manifold is homeomorphic to $X$, and the homeomorphism
maps the new branch set to the old one.  Therefore in this case we
can simply delete the curve $c$ from $G$.

(3) If $c \neq c_0$ is a circle component bounding a disk $D_1$ on
$\bdd X$ containing one branch point of $\alpha$, then $c$ lifts to a
circle on the boundary of a face $\sigma$ of $\gb$, which necessarily
contains a fixed point of $\alpha$.  Hence the cocore of the
corresponding 2-handle is a branch arc of $\alpha$.  The 2-sphere
boundary component created after attaching the 2-handle contains two
branch points of $\alpha$, hence bounds a 3-ball containing a trivial
arc as branch set of $\alpha$.  Thus after attaching the 2-handle and
the 3-ball the manifold is homeomorphic to $X$, and the branch set of
$\alpha$ has not changed.  As in Case (2), we will simply delete the
curve $c$ from $G$ in this case.

(4) If a circle component $c \neq c_0$ of $G$ bounds a disk $D_1$
containing exactly two branch points of $\alpha$, then after attaching
a 2-handle and a 3-cell, the manifold is homeomorphic to $X$, and the
branch set of $\alpha$ is obtained by adding a trivial arc in the
3-cell joining the two branch points of $\alpha$ in $D_1$.  Therefore
in this case we will add an arc in $D_1$ joining the two branch points
of $\alpha$, push the arc into the interior of $X$ as branch set of
$\alpha$, and then delete the curve $c$.

(5) If $c$ is a circle component of $G$ bounding a disk $D_1$
containing $k > 2$ branch points of $\alpha$, simply attach a
2-handle along $c$.  If $k$ is odd, add an arc in the center of the
2-handle to the branch set of $\alpha$.

(6) Finally, attach a 2-handle along $c_0$, fill each 2-sphere
boundary component containing at most 2 branch points with a 3-ball,
and add a trivial arc in the 3-ball to the branch set if the 2-sphere
contains exactly two branch points.  If the 2-sphere contains four branch
point, shrink it by an isotopy to a small sphere, which projects to a
small disk on the diagram, with four branch arcs attached.  (This
happens only for $Q_{14}$.  See Figure 22.13.)  This completes the
construction of the tangle $Q_i$.

For $M_6$, the above procedure produces the tangle $Q_6 = (W_6, K_6)$
in Figure 22.4(a), where $K_6$ should be considered as a tangle lying
in the half space $Q_6$ (including $\infty$) in front of the
blackboard.  The four boundary points of $K_6$ lie on the blackboard,
which is the boundary of $W_6$.

\bigskip
\leavevmode

\centerline{\epsfbox{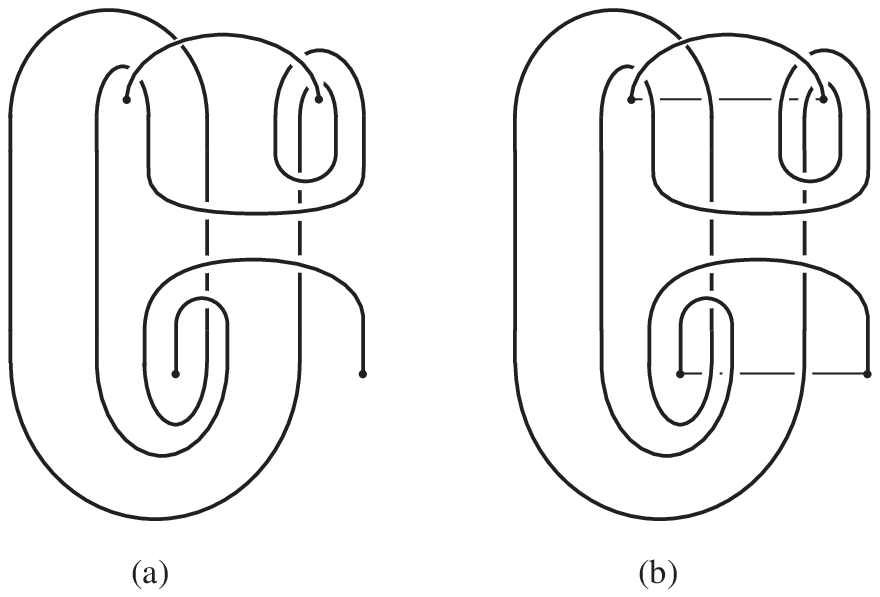}}
\bigskip
\centerline{Figure 22.4}
\bigskip

STEP 5.  {\it Find the tangles $Q_i(t_a) = M_i(r_a)/\alpha$.}

\medskip

There is one branch point $p_i$ in each quadrant of $P' = \Bbb R^2$.
Let $m, l$ be curves on $T_0$ that project to the $y$-axis union
$\infty$ and the $x$-axis union $\infty$ on $P'$, respectively.  This
sets up coordinate systems on $T_0$ and $P'$.  For $t=p/q$ a rational
number or $\infty$, let $M_i(t)$ denote the Dehn filling along slope
$pm + ql$, and $Q_i(t)$ denote the tangle obtained by attaching a
rational tangle of slope $t$ to $P'$.  In other words, $Q_i(t)$ is
obtained by attaching a 3-ball to $Q_i$ on $P'$, and adding two arcs
on $P'$ connecting the branch points of $\alpha$, which lift to curves
of slope $t$ on $T_0$.  Since the attached rational tangle lifts to a
solid torus with meridional slope $t$ on $T_0$, $M_i(t)$ is the double
branched cover of $Q_i(t)$.

By construction $\bdd F_a$ projects to the $x$-axis, hence $M_i(r_1) =
M_i(0)$.  The slope $r_2$ can be obtained by connecting the curves
$G''$ in $D_{\pm}$ by vertical arcs in $A = P' - \cup D_{\pm}$.  For
$M_6$, one can check that the slope $r_2 = 4$.  

\bigskip

Denote by $T(a_1, a_2)$ a Montesinos tangle which is the sum of two
rational tangles of slopes $1/a_1$ and $1/a_2$, respectively, where
$a_1, a_2$ are integers.  Denote by $T(a_1, b_1; a_2, b_2)$ the
collection of pairs $(S^3, L)$ which can be obtained by gluing two
tangles $T(a_i, b_i)$ along their boundary.  Denote by
$X(a_1,a_2)$ the collection of Seifert fiber spaces with orbifold a
disk with two cone points $c_1, c_2$ of index $a_1$ and $a_2$, i.e.\
the cone angle at $c_i$ is $2\pi/a_i$.  Note that the double branched
cover of $T(a_1, a_2)$ is in $X(a_1, a_2)$.  Denote by
$X(a_1,b_1;a_2,b_2)$ the collection of graph manifolds which are the
union of two manifolds $X_1, X_2$ glued along their boundary, where
$X_i \in X(a_i, b_i)$.

Denote by $K_{p/q}$ the two bridge knot or link associated to the
rational number $p/q$.  Denote by $C(p_1, q_1; p_2, q_2)$ the link
obtained by replacing each component $K_i$ of a Hopf link by its
$(p_i, q_i)$ cable $K'_i$, where $q_i$ is the number of times $K'_i$
winds around $K_i$.  Denote by $Y(p_1,q_1;p_2,q_2)$ the double
branched cover of $S^3$ with branch set $C(p_1, q_1; p_2, q_2)$.
Denote by $C(C; p,q)$ the link obtained by replacing one component
$K_1$ of a Hopf link by a Whitehead knot in the solid torus $N(K_1)$,
and the other component $K_2$ by a $(p,q)$ cable of $K_2$.  Let
$Y(C;p,q)$ be the double branched cover of $S^3$ with branch set
$C(C;p,q)$.  Denote by $Z$ the double branched cover of $S^3$ with
branch set the 2-string cable of the trefoil knot shown in Figure
22.12(d).

If $Q_i(r) = (S^3, L)$ then we will sometimes simply write $Q_i(r) =
L$.

\begin{lemma} (1) $Q_6(0) \in T(2,6; 2,3)$, as shown in Figure
22.4(b).  

(2) $Q_6(4) = C(3,1; 2,5)$, as shown in Figure 22.5(b).

(3) $Q_6(\infty) = K_{9/2}$.
\end{lemma}

\proof (1) The tangle $Q_6(0) = (S^3,L)$ is shown in Figure 22.4(b).
A horizontal line at the middle of the diagram corresponds to a
2-sphere $S$ which cuts the link $L$ into two Montesinos tangles
$T(2,6)$ and $T(2,3)$.

(2) The tangle $Q_6(4)$ is shown in Figure 22.5(a), which can be
isotoped to that in Figure 22.5(b).  One can see that it is the link
$C(3,1;2,5)$ in $S^3$.

(3) The tangle $Q_6(\infty)$ is shown in Figure 22.5(c).  One can
check that it is isotopic to the knot $K_{2/9}$ in Figure 22.5(d).
\qed

\bigskip
\leavevmode

\centerline{\epsfbox{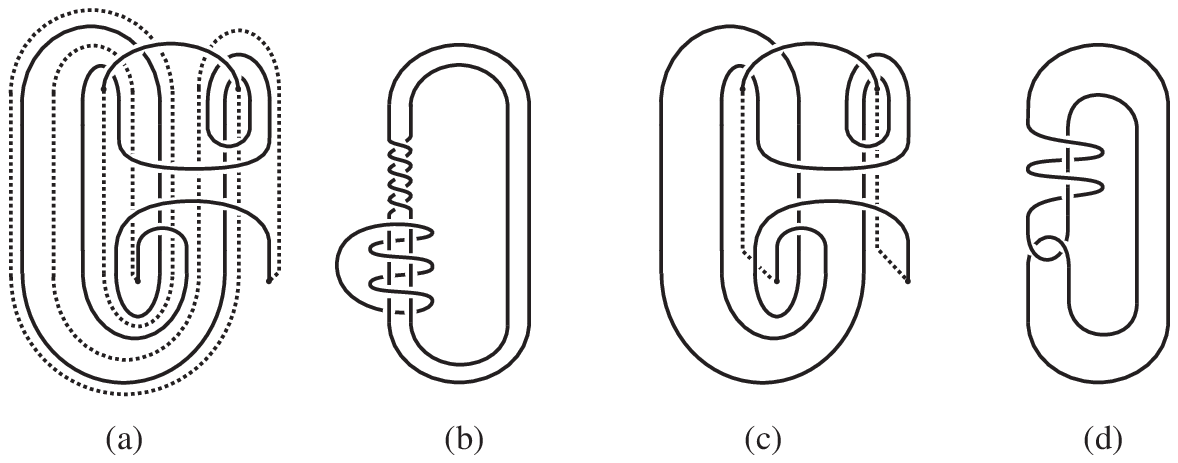}}
\bigskip
\centerline{Figure 22.5}
\bigskip

\begin{lemma} (1) For $i=6,...,14$, each $M_i$ is the double branched
  cover of a tangle $Q_i = (W_i, K_i)$, where $Q_i$ is shown in Figure
  22.4(a) for $i=6$, and in Figure 22.i(b) (with dotted lines removed)
  when $i>6$.

  (2) Each $M_i$ ($i=6, ..., 13$) admits a lens space surgery
  $M_i(r_3)$.  For each $i$, let $r_1, r_2$ be the slopes $r'_i,
  r''_i$ given in Definition 21.3.  Then the manifolds $M_i(r_1)$,
  $M_i(r_2)$ and $M_i(r_3)$ are given in the following table.  
$$
  \begin{array}{lll} M_6(0) \in
    X(2,6;2,3) \quad & M_6(4) = Y(3,1;5,2)
    \quad  & M_6(\infty)  = L(9,2) \\
    M_7(0) \in X(2,3;3,3) \quad & M_7(-5/2) \in X(2,3;2,2)
    \quad  & M_7(\infty)  = L(20,9) \\
    M_8(0) \in X(2,2;2,6) \quad & M_8(-5/4) = Y(3,1;2,5)
    \quad  & M_8(-1)  = L(4,1) \\
    M_9(0) \in X(2,3;2,3) \quad & M_9(-4/3) = Y(3,1;2,4)
    \quad  & M_9(-1)  = L(8,3) \\
    M_{10}(0) \in X(2,3;2,3) \quad & M_{10}(-5/2) = Y(C;2,1)
    \quad  & M_{10}(\infty)  = L(14,3) \\
    M_{11}(0) \in X(2,4;2,4) \quad & M_{11}(-5/2) = Y(C; 2,1)
    \quad  & M_{11}(\infty)  = L(24,5) \\
    M_{12}(0) \in X(2,3;2,3) \quad & M_{12}(5) = Y(3,1;2,3)
    \quad  & M_{12}(\infty)  = L(3,1) \\
    M_{13}(0) \in X(2,3;2,3) \quad & M_{13}(4) = Z \quad &
    M_{13}(\infty) = L(4,1)
\end{array}
$$
\end{lemma}

\proof The result for $M_6$ follows from Lemma 22.1 because $M_6(r)$
is a branched cover of $Q_6(r)$.  The proof for the other cases are
similar.  Each $M_i(r)$ is the double branched cover of $Q_i(r)$ and
the tangle $Q_i(r)$ is a link $L$ in $S^3$.  More explicitly, Figure
22.i(a) shows the curves $G = G' \cup G''$ in Step 3 of the above
construction; Figure 22.i(b) gives the tangle $Q_i$ as well as
$Q_i(r_1)$, which is obtained by attaching a $0$-tangle (the two
horizontal dotted lines) to $Q_i$; Figure 22.i(c) gives $Q_i(r_2)$,
which is simplified to that in Figure 22.i(d); $Q_i(r_3)$ is in Figure
22.i(e), which is simplified to that in Figure 22.i(f) for some $i$.
(The figures are numbered so that Figure 22.i corresponds to the
manifold $M_i$ for $i\geq 7$.  Note that there is no Figure 22.6.)
The manifold $M_{14}(r_3)$ is the double branched cover of
$Q_{14}(r_3) = T(2,2)$ in Figure 22.14(e), and hence is a twisted
$I$-bundle over the Klein bottle.
\qed

\bigskip
\leavevmode

\centerline{\epsfbox{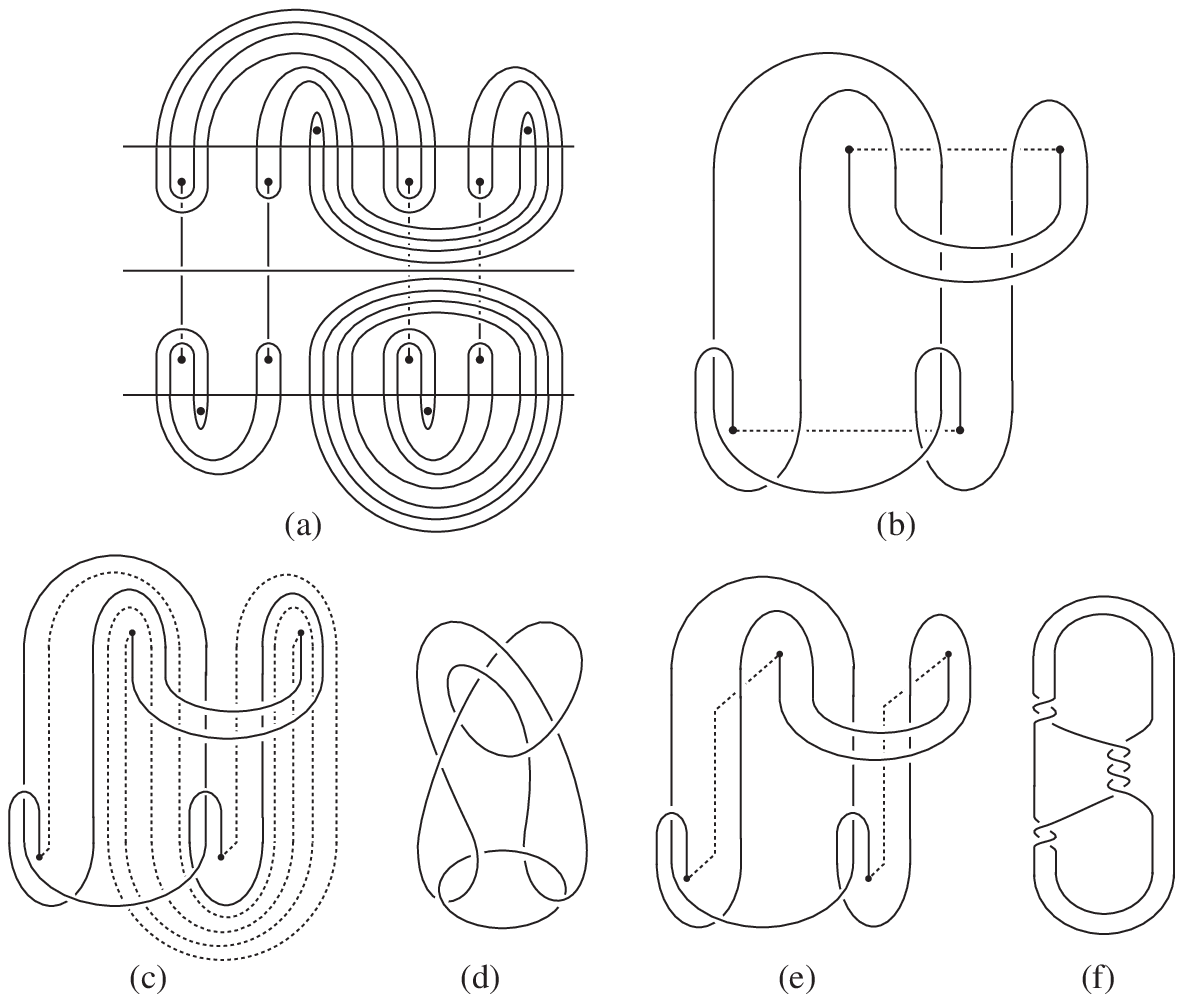}}
\bigskip
\centerline{Figure 22.7}
\bigskip

\bigskip
\leavevmode

\centerline{\epsfbox{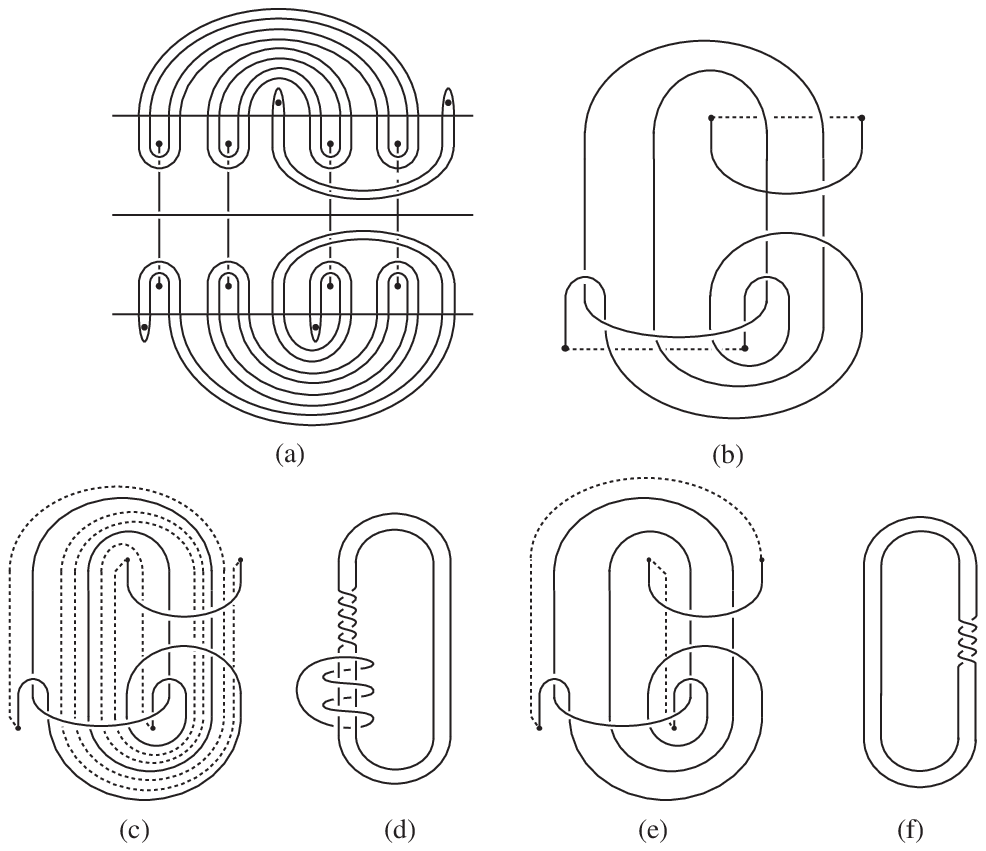}}
\bigskip
\centerline{Figure 22.8}
\bigskip

\bigskip
\leavevmode

\centerline{\epsfbox{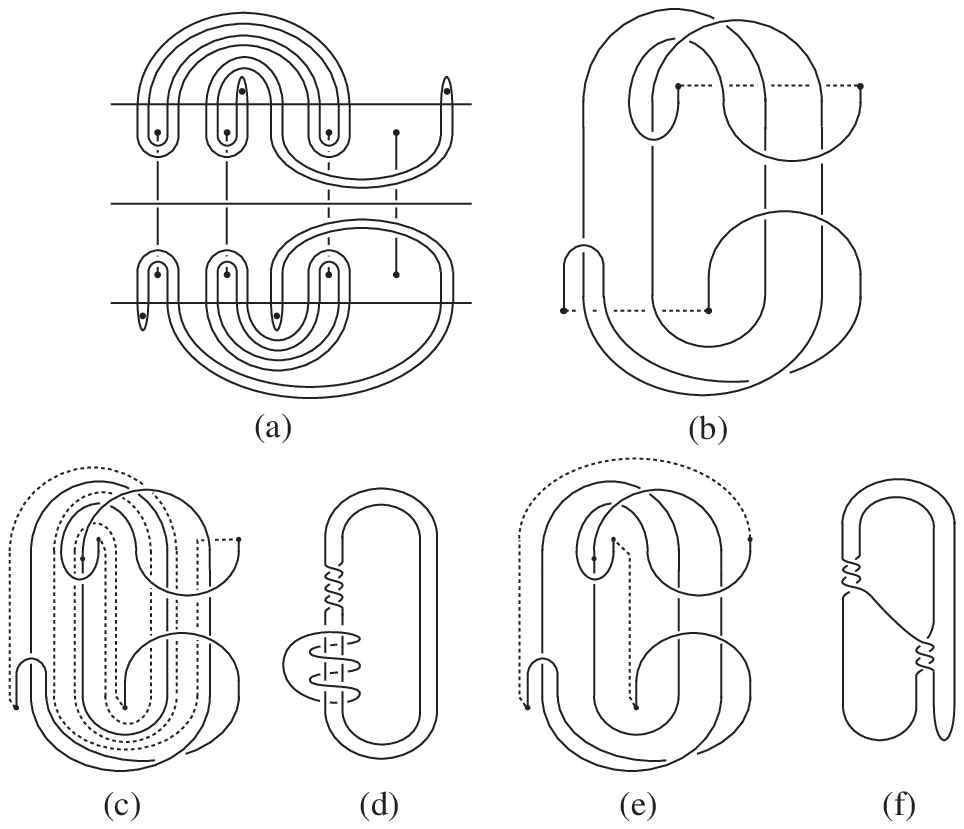}}
\bigskip
\centerline{Figure 22.9}
\bigskip

\bigskip
\leavevmode

\centerline{\epsfbox{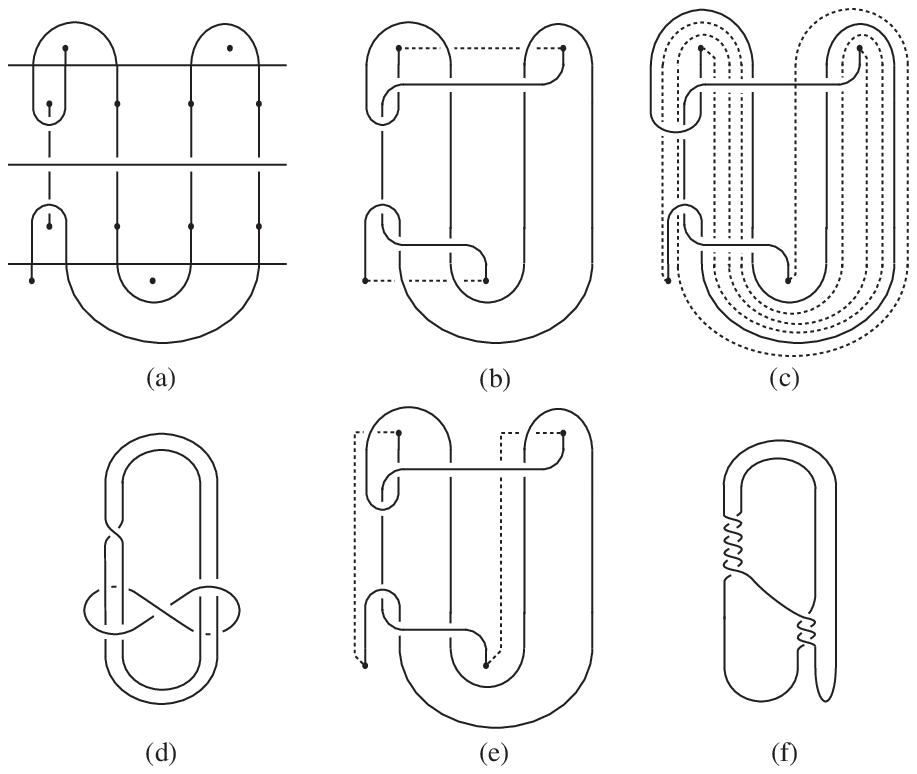}}
\bigskip
\centerline{Figure 22.10}
\bigskip

\bigskip
\leavevmode

\centerline{\epsfbox{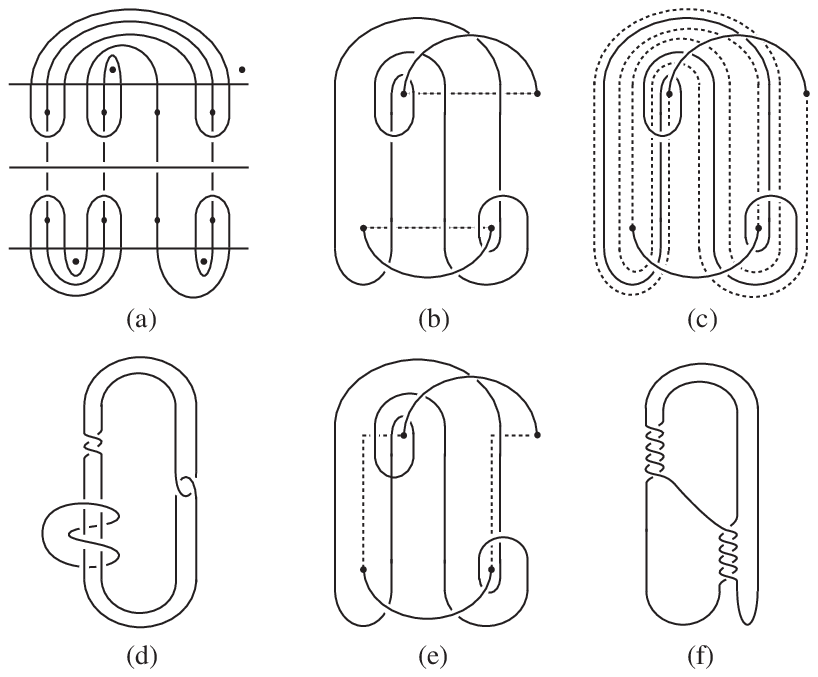}}
\bigskip
\centerline{Figure 22.11}
\bigskip

\bigskip
\leavevmode

\centerline{\epsfbox{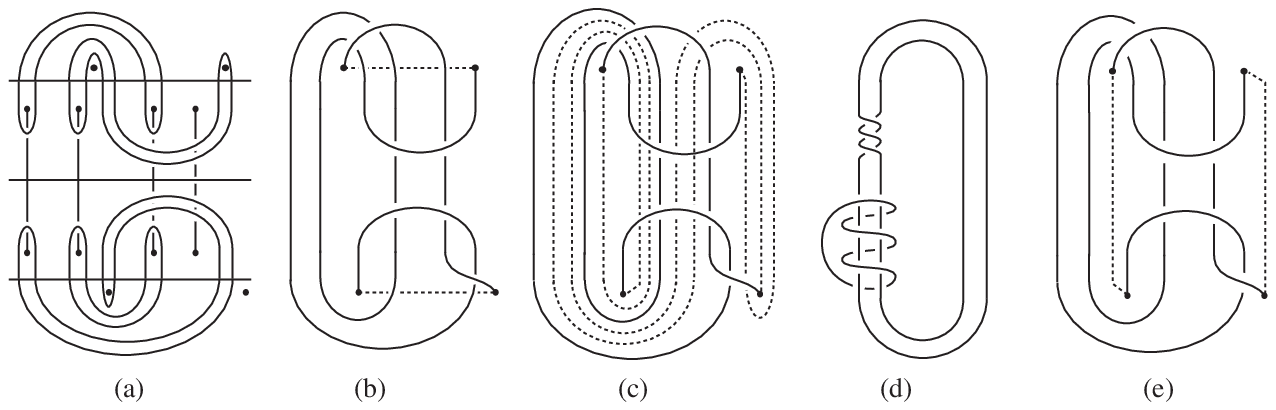}}
\bigskip
\centerline{Figure 22.12}
\bigskip

\bigskip
\leavevmode

\centerline{\epsfbox{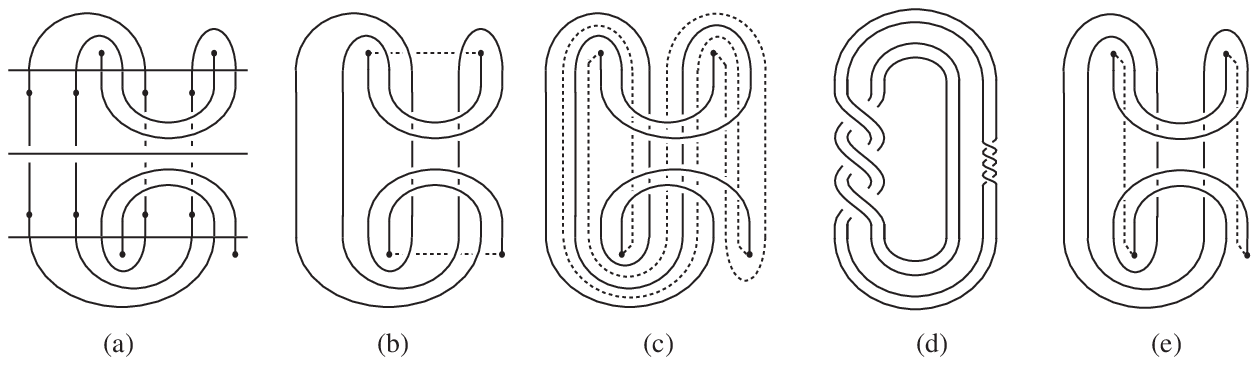}}
\bigskip
\centerline{Figure 22.13}
\bigskip

\bigskip
\leavevmode

\centerline{\epsfbox{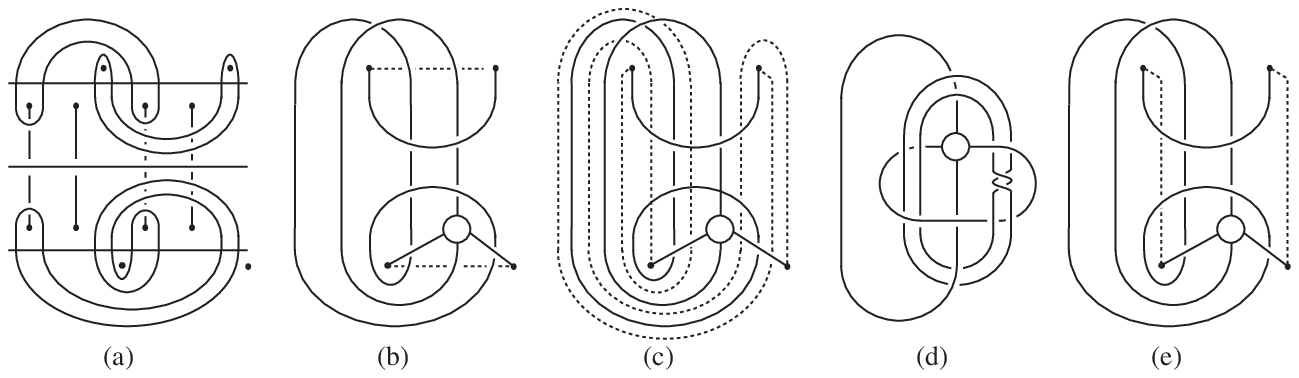}}
\bigskip
\centerline{Figure 22.14}
\bigskip

Recall that a manifold $M$ with a fixed torus $T_0 \subset \bdd M$ is
large if $H_2(M, \bdd M - T_0) \neq 0$.  Teragaito [T2] proved that
there is no large hyperbolic manifold $M$ admitting two toroidal
fillings of distance at least 5.  The following result shows that
there is only one such manifold for $\Delta=4$.

\begin{thm} 
  Suppose $(M,T_0)$ is a large manifold and $M$ is hyperbolic and
  contains two toroidal slopes $r_1, r_2$ on $T_0$ with $\Delta(r_1,
  r_2) \geq 4$.  Then $M$ is the Whitehead link exterior, and
  $\Delta(r_1, r_2)=4$.
\end{thm}

\proof Let $r$ be a slope on $T_0$, $V_r$ the Dehn filling solid torus
in $M(r)$, and $K_r$ the core of $V_r$.  By duality we have $H_2(M,
\bdd M-T_0) \cong H^1(M, T_0)$, which is isomorphic to the free part
of $H_1(M, T_0)$.  Also, $$H_1(M, T_0) \cong H_1(M(r), V_r) \cong
H_1(M(r)) / H_1(K_r).$$
Put $G(M,r) = H_1(M(r)) / H_1(K_r)$.  Then we need only show that
$G(M_i,r)$ is a (possibly trivial) torsion group for $i=2,...,14$ and
$r$ some slope on $T_0$.

For $i=2,3$, $M_i$ is the exterior of a closed braid $K_i$ in a solid
torus $V$.  Let $r$ be the meridian slope of $K_i$.  Then $G(M_i, r) =
\Bbb Z_p$, where $p$ is the winding number of $K_i$ in $V$.

For $i=6,...,13$, by Lemma 22.2 $M_i$ has a lens space filling
$M_i(r_3)$.  Therefore $G(M_i, r_3)$ is a quotient of the finite cyclic
group $H_1(M_i(r_3))$ and hence is a torsion group.  Similarly for the
four manifolds in [Go] with toroidal slopes of distance at least $6$.

For $i=14$, take a regular neighborhood of $u_1 \cup u_2 \cup D$ on
$\hat F_a$ as a base point.  See Figure 20.6(a).  Then
$H_1(M_{14}(r_a))$ is generated by $x,y,s_1,s_2$, where $x$ is the
element of $H_1(\hat F_a)$ represented by the edges $C$ on Figure
20.6(a), oriented from the label $2$ endpoint to the label $1$
endpoint, $y$ is represented by $B$, oriented from $u_1$ to $u_2$, and
$s_i$ by the part of the core of the Dehn filling solid torus running
from $u_i$ to $u_{i+1}$ with respect to the orientation of $\bdd F_b$.
Then the bigons $B\cup D$, $C\cup E$ and the 4-gon bounded by $C\cup
D\cup E \cup Y$ on $F_b$ give relations $2s_1 - y = 0$, $2x = 0$, and
$y + 2x = 0$.  The other faces of $\Gamma_b$ are parallel to these.
To calculate $G(M_{14}, r_a) = H_1(M_{14}(r_a)) / H_1(K_a)$ we further
add the relation $s_1 + s_2 =0$.  One can now check that $G(M_{14},
r_a) = \Bbb Z_2 \oplus \Bbb Z_2$, and the result follows.

For $i=4$, choose a regular neighborhood of $v_1 \cup v_2 \cup J$ in
Figure 11.9(b) as a base point.  Then $H_1(M_4(r_b))$ is generated by
$x,y,s_1,s_2$, where $x,y$ are represented by the edges $L,C$ in
Figure 11.9(b), oriented from $v_1$ to $v_2$, and $s_i$ by the part of
the core of the Dehn filling solid torus from $v_i$ to $v_{i+1}$.  The
faces bounded by $L\cup C$, $C\cup K$ and $Q\cup K \cup M \cup A$ give
the relations $y - s_1 + x + s_2 = 0$, $s_1 - x - s_2 = 0$, and $s_2 -
s_1 + y = 0$.  Together with the relation $s_1+s_2=0$ from $H_1(K_b)=0$,
these give $G(M_4, r_b) = \Bbb Z_2$.

For $i=5$, $H_1(M_5(r_b))$ is generated by $x,y,s$, where $x,y$ are
represented by edges $E$ and $C$ on Figure 11.10(b), oriented from
label $3$ to label $4$, and $s$ is represented by the core of the Dehn
filling solid torus.  Then the bigon $A\cup H$ and the annulus bounded
by $A\cup G \cup C \cup E$ on Figure 11.10(a) containing $J$ give the
relations $x+y=0$ and $2x - 2y = 0$.  Adding the relation
$s=0$ gives $G(M_5, r_b) = \Bbb Z_4$.  \qed

\section {The manifolds $M_i$ are hyperbolic}

The manifolds $M_1, M_2, M_3$ in Definition 21.3 are known to be
hyperbolic, see [GW1, Theorem 1.1].  In this section we will show that
the other 11 manifolds $M_i$ in Definition 21.3 are also hyperbolic.  See
Theorem 23.14 below.

A knot $K$ in a solid torus $V$ is a $(p,q)$ knot if it is isotopic to
a $(p,q)$ curve on $\bdd V$ with respect to some longitude-meridian
pair on $\bdd V$.  In particular, the winding number of $K$ in $V$ is
$p$.

\begin{lemma} (1) If $i\in \{6,...,13\}$ and $j=1$ or $i\in
\{6,7,8,9,12,13\}$ and $j=2$, then $M_i(r_j)$ contains an
essential torus $T$ cutting it into two Seifert fiber spaces
$E_1, E_2$.

(2) For $i=10, 11$, $M_i(r_2)$ contains a non-separating essential
torus cutting $M_i(r_2)$ into a Seifert fiber space whose orbifold is
an annulus with a cone point of index 2.

(3) For $i=6,...,13$ and $j=1,2$, $M_i(r_j)$ is irreducible, and
contains no hyperbolic submanifold bounded by an incompressible torus.
\end{lemma}

\proof (1) By Lemma 22.2, for $i=6,...,13$, $M_i(r_1)$ is of type
$X(a_1, b_1; a_2, b_2)$, which is the union of two Seifert fiber
spaces of types $X(a_1, b_a)$ and $X(a_2, b_2)$, hence the result is
true for $M_i(r_1)$.  Similarly it is true for $M_7(r_2)$.

For $i=13$, $Q_i(r_2) = (S^3, L)$, where $L$ is the link in Figure
22.13(d), which consists of two parallel copies of the trefoil knot.
The two components of $L$ bound an annulus $A$.  Cutting $S^3$ along
$A$ gives the trefoil knot exterior $E$, and $A$ becomes a torus $T$.
The double branched cover of $Q_{13}(r_2)$ is obtained by gluing two
copies of $E$ along $T$.  Hence the result is true because $E$ is a
Seifert fiber space and $T$ is incompressible in $E$.

For $i=6,8,9,12$, $Q_i(r_2) = (S^3, L)$ is of type $C(p_1,q_1;
p_2,q_2)$, so there is a torus $T'$ cutting $S^3$ into two solid tori
$V_1, V_2$, such that each $L_j = L \cap V_j$ is a $(p_j,q_j)$ knot in
$V_j$ for some $p_j>1$.  Note also that in these cases at least one of
the $p_j$ is odd, which implies that $T'$ lifts to a single torus $T$,
cutting $M_i(r_2)$ into two components $W_1, W_2$, such that $W_j$ is
a double branched cover of $(V_j, L_j)$.  The $(p_j,q_j)$ fibration of
$V_j$ now lifts to a Seifert fibration of $W_j$, hence the result
follows.

(2) For $i=10,11$, $Q_i(r_2) = (S^3, L)$, and there is a torus $T'$
cutting $S^3$ into two solid tori $V_1, V_2$, such that $L_1 = V_1
\cap L$ is a $(2,1)$ knot, and $L_2 = V_2 \cap L$ is a Whitehead knot
in the solid torus $V_2$.  Since both winding numbers of $L_j$ are
even, $T'$ lifts to two tori in $M_i(r_2)$.  Let $W_i$ be the lifting
of $V_i$.  A meridian disk of $V_1$ lifts to an annulus in $W_1$,
hence $W_1$ is a $T^2 \times I$ (not a twisted $I$-bundle over the
Klein bottle because $\bdd W_1$ has two components).  Let $T$ be the
core of this $T^2 \times I$.  Then it cuts $M_i(r_2)$ into the
manifold $W_2$.  We need to show that $W_2$ is Seifert fibered.

Let $D$ be a meridian disk of $V_2$ which intersects $L_2$ at two
points.  Then $(V_2, L_2) = (B_1, L'_1) \cup (B_2, L'_2)$, where $B_1
= N(D)$, $B_2$ is the closure of $V_2 - B_1$, and $L'_k = L_2 \cap
B_k$.  Note that each $L'_k$ is a trivial tangle in $B_k$, hence its
double branched cover $V'_k$ is a solid torus.  One can check that
each component of $V'_1 \cap V'_2$ is a longitudinal annulus on $\bdd
V'_1$, and it is an annulus on $\bdd V'_2$ with winding number $2$ in
$V'_2$.  Therefore $W_2 = V'_1 \cup V'_2$ is a Seifert fiber space
whose orbifold is the union of a $D^2$ and a $D^2(2)$ glued along two
boundary arcs.  Since $W_2$ has two torus boundary components, the
orbifold must be an annulus with a single cone point of index $2$.

(3) Let $T$ be the essential torus in $M_i(r_j)$ given in the above
proof.  Then it cuts $M_i(r_j)$ into one or two bounded Seifert fiber
spaces, which are irreducible.  Since $T$ is incompressible,
$M_i(r_j)$ is also irreducible.  The second statement follows from the
fact that the JSJ (Jaco-Shalen-Johannson) decomposition of an
irreducible closed 3-manifold is unique.  \qed

A $(p',q')$ knot $K$ in a solid torus $V$ is also called a $0$-bridge
knot.  In this case there is an essential annulus in $V - \Int N(K)$
with one boundary component in each of $\bdd V$ and $\bdd N(K)$.  This
defines a longitude $l$ for $K$, which is unique if $K$ is not the
core of $V$.  A $(p,q)$ cable of a 0-bridge knot $K$ is a knot on
$\bdd N(K)$ which represents $pl + qm$ in $H_1(\bdd N(K))$, where $m$
is a meridian of $K$.  We refer the readers to [Ga1] for the definition
of a $1$-bridge braid in $V$.

\begin{lemma} Suppose $X$ is an irreducible,
$\bdd$-irreducible, compact, orientable $3$-manifold with $\bdd X =
T_1 \cup T_0$ a pair of tori.  Let $r_1, r_2$ be distinct slopes on
$T_0$ such that $X(r_1)$, $X(r_2)$ are both $\bdd$-reducible.  Let
$K_a$ be the core of the Dehn filling solid torus in $X(r_a)$.  Then
one of the following holds, up to relabeling of $r_i$.

(1) Each $X(r_a)$ is a solid torus, $K_a$ is a $0$- or $1$-bridge
braid in $X(r_a)$, and $\Delta(r_1, r_2) = 1$ if it is not a
$0$-bridge knot.

(2) $X(r_1)$ is a solid torus, and $X(r_2) = (S^1 \times D^2) \#
L(p,q)$ with $p\geq 2$.  $K_1$ is a $(p,q)$ cable of a $(p',q')$ knot
in $X(r_1)$, and $r_2$ is the cabling slope of $K_1$ in $X(r_1)$.
Moreover, if $m_a$ is the slope on $T_1$ bounding a disk in $X(r_a)$,
then $\Delta(m_1, m_2) = pp'$.  \end{lemma}

\proof If both $X(r_1), X(r_2)$ are irreducible then they are solid
tori and (1) holds by [Ga1, Theorem 1.1] and [Ga2, Lemma 3.2].  Now
assume $X(r_2)$ is reducible.  Then by [Sch, Theorem 6.1] $K_1$ is a
$(p,q)$ cable of some knot $K'$ in $X(r_1)$ with respect to some
meridian-longitude pair $(m,l)$ of $K'$, and $r_2$ is the cabling
slope.  In this case $X(r_2)$ is a connected sum $W_1 \cup L(p,q)$,
where $W_1$ is obtained by surgery on the knot $K'$ in $X(r_a)$ along
the cabling slope $r' = pl+qm$.  Denote by $m'$ the meridian slope of
$K'$.  Then $\Delta(m', r') = p > 1$.

Denote by $K'(s)$ the manifold obtained by $s$-surgery on $K'$ in $V =
X(r_1)$.  The assumption on $X$ implies that $T_1$ is incompressible
in $V-K'$, and $V-K'$ is irreducible.  By the above, $T_1$ is
compressible in both $K'(m') = V$ and $K'(r') = W_1$, and
$\Delta(m',r')>1$, hence by [Wu2, Theorem 1.1] and [CGLS, Theorem
2.4.3], either $V-\Int N(K') = T^2 \times I$, or there is an annulus
$A$ in $V-\Int N(K')$ with one boundary component on $T_1$ and another
boundary component a curve of slope $r$ on $T' = \bdd N(K')$,
satisfying $\Delta(r, m') = \Delta(r, r') = 1$.  In either case $K'$
is isotopic to a curve on $T_1$ and hence is a $0$-bridge knot.  Since
$V-K$ is irreducible, this implies that $V$ is also irreducible.
Therefore $V = X(r_1)$ is a solid torus, $K'$ is a $(p',q')$ knot in
$V$ for some $(p',q')$, and $r$ is the cabling slope of $K'$ when
$p'>1$.

If $p'=1$, i.e.\ $V-\Int N(K') = T^2 \times I$, then $K$ is the
$(p,q)$ cable of the core $K'$ of $V$, and $r_2$ is the cabling slope
of $K$.  We have $X(r_2) = V \# L(p,q)$, and the slope $m_2$ on $T_1$
which bounds a disk in $X(r_2)$ is the $(p,q)$ curve on $T_1$, hence
$\Delta(m_1, m_2) = p = pp'$.

Now assume $p' > 1$.  Choose the longitude $l$ of $K'$ to be the
cabling slope $r$ of $K'$ given above.  Since $r' = pr+qm$, the
equation $\Delta(r, r')=1$ above implies that $q=\pm 1$.  Reversing
the orientation of $m'$ if necessary we may assume $q=1$.  Hence $K$
is a $(p,1)$ cable of a $(p',q')$ knot in $V$, and $r_2$ is the
cabling slope.  It is easy to see that the meridian slopes $m_a$ of
$X(r_a)$ satisfy $\Delta(m_1, m_2) = p p'$, and the result
follows.  \qed

\begin{lemma} Let $i=6,...,13$.  Let $\alpha$ be the covering
transformation of the double branched cover $M_i \to Q_i$.  

(1) $M_i$ is irreducible, not Seifert fibered, and contains no
non-separating torus. 

(2) If $M_i$ is not hyperbolic then it contains a separating essential
torus $T$ such that $T$ is $\alpha$-equivariant, and the component $W$
of $M_i|T$ which does not contain $T_0$ is either Seifert fibered or
hyperbolic.  \end{lemma}

\proof (1) If $M_i$ is reducible then the summand which does not
contain $T_0$ is a summand of $M_i(r_j)$ for all $r_j$, but since
$M_i(r_3)$ is a lens space while $M_i(r_1)$ does not have a lens space
summand, this is impossible.  

By Lemma 22.2 $M_i(r_3)$ is a lens space, so if $M_i$ is Seifert
fibered then the orbifold of $M_i$ is a disk with two cone points,
hence $M_i(r_1)$ is either a connected sum of two lens spaces or a
Seifert fibered space with orbifold a sphere with at most three cone
points.  This is impossible because by Lemma 22.2 $M_i(r_1)$ is of
type $X(a_1, b_1; a_2, b_2)$ with some $a_i$ or $b_i$ greater than
$2$, which is irreducible and contains a separating essential torus,
at least one side of which is not an $I$-bundle.

If $M_i$ contains a non-separating torus then the lens space
$M_i(r_3)$ would contain a non-separating surface, which is absurd.

(2) If $M_i$ is non-hyperbolic then by (1) it has a non-trivial JSJ
decomposition.  By [MeS] we may choose the JSJ decomposition surfaces
$F$ to be $\alpha$-equivariant.  If we define a graph $G$ with the
components of $M_i|F$ as vertices and the components of $F$ as edges
connecting adjacent components of $M_i|F$, then the fact that $M_i$
contains no non-separating torus implies that $G$ is a tree.  Let $T$
be a component of $F$ corresponding to an arc incident to a vertex $v$
of valence $1$ in $G$.  Then $T$ bounds the manifold $W$ corresponding
to $v$, which by definition of the JSJ decomposition is either Seifert
fibered or hyperbolic.  \qed

\begin{lemma} Suppose $M_i$ is non-hyperbolic and let $T$ be
the essential torus in $M_i$ given in Lemma 23.3(2).  Let $X, W$ be
the components of $M_i | T$, where $X \supset T_0$.  If $T$ is
compressible in $M_i(r_a)$ for some $a=1,2$, then both $X(r_a)$ and
$X(r_3)$ are solid tori, and $W$ is hyperbolic.  
\end{lemma}

\proof $T$ is compressible in $X(r_3)$ because $M_i(r_3)$ is a lens
space.  By assumption $T$ is also compressible in $M_i(r_a)$.  Since
$M_i(r_a)$ contains no lens space summand, by Lemma 23.2 either both
$X(r_a)$ and $X(r_3)$ are solid tori, or $X(r_a)$ is a solid torus and
$X(r_3) = (S^1 \times D^2) \# L(p,q)$ for some $p>1$.  We need to show
that the second case is impossible.

Let $m_j$ be the slope on $T$ which bounds a disk in $X(r_j)$, $j=a,
3$.  Then $M_i(r_3) = W(m_3) \# L(p,q)$.  Since $M_i(r_3)$ is a lens
space $L(p,q)$, we have $W(m_3) = S^3$, so $W$ is the exterior of a
knot in $S^3$.  If $W$ is Seifert fibered then it is the exterior of a
torus knot, so $M_i(r_a) = W(m_a)$ is obtained by Dehn surgery on a
torus knot in $S^3$ and hence is atoroidal, which contradicts Lemma
23.1.  Since by definition $W$ is Seifert fibered or hyperbolic, this
implies that $W$ is hyperbolic.  Note that by Lemma 22.2 $p \geq 3$
for all $i \in \{6,...,13\}$.  By Lemma 23.2(2) we have $\Delta(m_a,
m_3) \geq p \geq 3$.  Since $W(m_a) = M_i(r_a)$ is toroidal and
$W(m_3) = S^3$, this is a contradiction to [GLu, Theorem 1.1], which
says that only integral or half integral surgeries on hyperbolic knots
in $S^3$ can produce toroidal manifolds.  This completes the proof
that both $X(r_a)$ and $X(r_3)$ are solid tori.

If $W$ is not hyperbolic then by definition $W$ is Seifert fibered.
By the above $X(r_3)$ is a solid torus.  Let $m_3$ be a meridian slope
of $X(r_3)$.  Then $M_i(r_3) = W(m_3)$, so $M_i(r_3)$ being a lens
space implies that the orbifold of $W$ is a disk with two cone points,
in which case $M_i(r_a) = W(m_a)$ is either a connected sum of two
lens spaces or a Seifert fiber space with orbifold a sphere with at
most three cone points.  In the first case $M_i(r_a)$ contains no
essential torus, while in the second case the only possible essential
torus in $M_i(r_a)$ is a horizontal torus cutting the manifold into a
$T^2 \times I$, or two twisted $I$-bundles over the Klein bottle.
This is a contradiction because by Lemma 23.1 $M_i(r_a)$ contains an
essential torus cutting it into either a Seifert fiber space with
orbifold an annulus with a cone point of index 2, or two Seifert fiber
spaces, at least one of which is not a twisted $I$-bundle over the
Klein bottle.  \qed

\begin{lemma} The torus $T$ in Lemma 23.3(2) is incompressible
in both $M_i(r_1)$ and $M_i(r_2)$.  \end{lemma}

\proof First assume that $T$ is compressible in both $X(r_1)$ and
$X(r_2)$.  By Lemma 23.4 $X(r_j)$ is a solid torus for $j=1,2,3$.
Since $\Delta(r_1,r_2)>1$, by Lemma 23.2 we see that $X$ is the
exterior of a $(p,q)$ knot in a solid torus.  Since $T$ is not
boundary parallel, $p>1$.  Let $r$ be the cabling slope on $T_0$.
Since $X(r_j)$ is a solid torus, we have $r_j \neq r$.  Therefore by
[CGLS, Theorem 2.4.3] we must have $\Delta(r, r_j) = 1$ for $i=1,2,3$.
By Lemma 23.2 one can check that $\Delta(r_1, r_3) = 1$ and
$\Delta(r_2,r_3) \leq 2$.  Since $\Delta(r_1, r_2) \geq 4$, this is a
contradiction because any three slopes $r_1, r_2, r_3$ with distance
$1$ from a given slope $r$ have the property that $\Delta(r_a, r_b) +
\Delta(r_b, r_c) = \Delta(r_a, r_c)$ for some permutation $(a,b,c)$ of
$(1,2,3)$.

Now assume that $T$ is compressible in $M_i(r_1)$, say.  By Lemma 23.4
$W$ is hyperbolic.  On the other hand, by the above $T$ is
incompressible in $M_i(r_2)$, so $W$ is a submanifold in $M_i(r_2)$
bounded by an incompressible torus, hence by Lemma 23.1(3) it is
non-hyperbolic, which is a contradiction.  \qed

\begin{lemma} Let $T(a_1, b_1; a_2, b_2) = (S^3, L)$, where
$a_i, b_i \geq 2$.  If at least one of $a_1, b_1, a_2, b_2$ is
greater than $2$ then the exterior of $L$ is atoroidal, and there is
no M\"obius band $F$ in $S^3$ with $F \cap L$ a component of $L$.
\end{lemma}

\proof Denote by $T(a)$ a rational tangle with slope $1/a$, where $a$
is an integer.  Given a tangle $\tau=(B^3, \tau)$, denote $B^3 - \Int
N(\tau)$ by $E$ or $E(\tau)$, and call it the tangle space of $\tau$.
Since $T(a)$ is a trivial tangle in the sense that $\tau$ is rel
$\bdd$ isotopic to arcs on $\bdd B^3$, the tangle space $E(a)$ is
atoroidal, and any incompressible annulus in $B^3-\tau$ is trivial in
the sense that it is either parallel to an annulus on $\bdd
(B^3-\tau)$ or cuts off a $D^2 \times I$ in $B^3$ with $\tau \cap (D^2
\times I)$ a core arc.

The tangle space $E(r_1, r_2)$ of a Montesinos tangle $T(r_1, r_2)$ is
obtained by gluing $E(r_1), E(r_2)$ along a twice punctured disk $P =
E(r_1) \cap E(r_2)$.  The above implies that $E(r_1, r_2)$ is always
atoroidal.  If $A$ is an essential annulus in $E(r_1, r_2)$ with
minimal intersection with $P$, then an innermost circle outermost arc
argument shows that $A$ intersects $P$ in essential arcs or circles in
$A$.  If the intersection is a set of circles then each component of
$A \cap E(r_i)$ is a set of trivial annuli, which implies that $A$ is
also trivial.  If each component of $A \cap P$ is an essential arc
then each component of $A \cap E(r_i)$ is a bigon in the sense that it
is a disk intersecting $P$ in two arcs, which implies that $r_i = 2$
for $i=1,2$.  Therefore $E(r_1, r_2)$ contains no essential annulus
unless $r_1 = r_2 = 1/2$ mod $1$.

By definition $T(a_1, b_1; a_2, b_2)$ is the union of two Montesinos
tangles $T(a_i, b_i)$.  If the tangle space of $T(a_1, b_1; a_2, b_2)$
is toroidal then either one of the $T(a_i, b_i)$ is toroidal or they
are both annular.  By the above neither case is possible if at least
one of the $a_1, b_1, a_2, b_2$ is greater than $2$.

The proof for a M\"obius band is similar.  If $F$ is a M\"obius band
in $S^3$ bounded by a component of $L$ and has interior disjoint from
$L$ then after cutting along the surface $P_1 = E(a_1, b_1) \cap
E(a_2, b_2)$ it either lies in one of the $E(a_i, b_i)$ or intersects
each in bigons.  One can show that the first case is impossible, and
in the second case $a_i = b_i = 2$ for $i=1,2$.  \qed

\begin{lemma}  $M_i$ is hyperbolic for $i=6$ or $8\leq i \leq 13$.
\end{lemma}

\proof Let $T$ be the $\alpha$-equivariant essential torus in $M_i$
given in Lemma 23.3(2).  By Lemma 23.5 $T$ is incompressible in both
$M_i(r_a)$, $a=1,2$.  Since $T$ is $\alpha$-equivariant, its image $F$
in $Q_i = M_i/\alpha$ is a 2-dimensional orbifold with zero orbifold
Euler characteristic (see [Sct] for definition), and all the cone
points have indices $2$.  Hence it is $T^2$, $K^2$, $P^2(2,2)$,
$S^2(2,2,2,2)$, $A^2$, $M^2$, or $D^2(2,2)$, where the surfaces are
torus, Klein bottle, projective plane, sphere, annulus, M\"obius band
and disk, and the numbers indicate the indices of the cone points.
Note that in the last three cases the boundary of the surface is part
of the branch set of $\alpha$.  Since $T$ is incompressible in
$M_i(r_a)$, $F$ is incompressible in $Q_i(r_a)$ in the sense that if
some simple loop on $F$ bounds a disk in $Q_i(r_a)$ intersecting the
branch set at most once then it bounds such a disk on $F$.  We need to
show that for each type of surface above there is some $a=1,2$ such
that no such incompressible 2-dimensional orbifold exists in
$Q_i(r_a)$.

We have $Q_i = (B^3, K_i)$, where $B^3 = M_i/\alpha$ is a 3-ball and
$K_i$ is the branch set of $\alpha$.  Since $F$ lies in $B^3$, it
cannot be $K^2$ or $P^2$.  For all $i$ one can check that the branch
set $K_i$ of $\alpha$ in $Q_i$ contains at most one closed circle,
hence the case $A^2$ is also impossible.

By Lemma 22.2, $Q_i(r_1) = T(a_1, b_1; a_2, b_2)$ for some $a_j, b_j
\geq 2$, and $(a_j, b_j) \neq (2,2)$ for some $j$.  Therefore by Lemma
23.6 we have $F \neq T^2, M^2$ for $i=6,...,13$ because there is no
such surfaces in $Q_i(r_1)$.  It remains to show that $F \neq D^2(2,2)$
or $S^2(2,2,2,2)$.

For $i=13$, $Q_i(r_2) = (S^3, L)$, where $L$ consists of two parallel
copies of a trefoil knot $K$.  Since each component of $L$ is
non-trivial in $S^3$, $F \neq D^2(2,2)$ in this case.  Suppose $F =
S^2(2,2,2,2)$, and let $V$ be a regular neighborhood of the trefoil
knot containing $L$, intersecting $F$ minimally.  Then $F \cap V \neq
\emptyset$, and $F$ is not contained in $V$ as otherwise one can show
that $F-L$ would be compressible.  Therefore $F\cap V$ is a union of
two meridian disks, and $F \cap S^3 - \Int V$ is an essential annulus
in $S^3 - \Int V$.  Since $S^3 - \Int V$ contains no essential annulus
with the meridian of $V$ as boundary slope, this is a contradiction.

The proofs for the cases $i \in \{6,8,9,10,11,12\}$ are similar.  In
these cases $Q_i(r_2) = (S^3, L)$, and there is a torus $T'$ cutting
$S^3$ into two solid tori $V_1, V_2$, each containing some components
of $L$.  One can check that no component $L'$ of $L$ bounds a disk
intersecting $L-L'$ at two points, so $F \neq D^2(2,2)$.  If
$F=S^2(2,2,2,2)$ then either $F$ lies in one of the $V_j$, or it
intersects one of the $V_j$ in two meridional disks and the other
$V_k$ in an essential annulus with boundary slope the meridional slope
of $V_j$.  Neither case is possible for the $Q_i(r_2)$ listed in Lemma
22.2.  \qed

\begin{lemma}  $M_7$ is hyperbolic.
\end{lemma}

\proof By Lemma 22.2 we have $M_7(r_1) \in X(3,3;2,3)$ and $M_7(r_2)
\in X(2,2;2,3)$.  Consider the tangle decomposition sphere $P_a$ of
the orbifold $Q_7(r_a)$, $a=1,2$, which corresponds to a horizontal
plane in Figure 22.7(b), (d) respectively.  It lifts to an essential
torus $T_a$ in $M_7(r_a)$.

Each side of $P_a$ is a Montesinos tangle of type $T(r_1, r_2)$, which
is the sum of two rational tangles over a disk $D$.  The boundary of
$D$ determines the fibration of the double branched cover $X(r_1,
r_2)$ of $T(r_1, r_2)$, which has a unique Seifert fibration unless
$r_1 = r_2 = 2$, in which case the closed circle in the tangle is
isotopic (without crossing the arcs) to a curve on the punctured
sphere, which lifts to a fiber in the other fibration of $X(r_1,r_2)$.
It is easy to check from Figures 22.7(b) and (d) that the fiber curves
from the two sides of $P_a$ do not match, so $M_7(r_a)$ is not a
Seifert fiber space.  Since each side of $T_a$ is a small Seifert
fiber space with orbifold a disk with two cone points, it follows that
$M_7(r_a)$ contains no other essential torus.

Suppose $M_7$ is non-hyperbolic and let $T$ be the essential torus in
$M_7$ given by Lemma 23.3.  By Lemma 23.5 it is incompressible in both
$M_7(r_a)$, therefore by the uniqueness of $T_a$ above we see that $T
= T_a$ in $M_7(r_a)$ up to isotopy.  As before, denote by $W$ and $X$
the components of $M_7|T$, with $X \supset T_0$.  Then $W$ is the
manifold on one side of $T_a$ in $M_7(r_a)$.  Therefore we must have
$W = X(2,3)$, so $X(r_1) = X(3,3)$ and $X(r_2) = X(2,2)$.  We will
show that this is impossible.

Let $Y$ be the component of the JSJ decomposition of $X$ that contains
$T$.  Then $Y$ is either hyperbolic or Seifert fibered.  There are
three cases. 

Case 1. {\it $T_0 \subset \bdd Y$ and $Y$ is Seifert fibered.\/} By
Lemma 23.5 $T$ is incompressible in $Y(r_a)$ for $a=1,2$, so $r_a$ is
not the fiber slope on $T_0$.  Hence the Seifert fibration extends
over $Y(r_1)$ and $Y(r_2)$.  In this case $\bdd Y - T_0$ is
incompressible in $Y(r_a)$.  Since $X(r_a)$ is atoroidal, either
$Y(r_1) \cong Y(r_2) \cong T^2 \times I$, or $Y = X$.  In the first
case we have $X(r_1) = X(r_2)$, which is a contradiction because
$X(r_1) = X(3,3) \not \cong X(2,2) = X(r_2)$.  In the second case,
Since $X(r_1)$ has orbifold $D^2(3,3)$, the orbifold of $X$ must be
$A^2(3,3)$ or $A^2(3)$.  On the other hand, since $X(r_2)$ has
orbifold $D^2(2,2)$, the orbifold of $X$ must be $A^2(2,2)$ or
$A^2(2)$, which contradicts the fact that Seifert fibrations
for these manifolds are unique.

Case 2.  {\it $T_0 \subset \bdd Y$ and $Y$ is hyperbolic.\/} If $\bdd
Y$ has more than two boundary components then the fact that $X(r_a)$
is atoroidal implies that $Y(r_a)$ is either $\bdd$-reducible, or a
$T^2 \times I$.  If $\bdd Y = T \cup T_0$ then $Y=X$ and by assumption
both $Y(r_a) = X(r_a)$ are annular and atoroidal.  In either case
$Y(r_a)$ is either $\bdd$-reducible or annular and atoroidal.  Since
$\Delta(r_1, r_2) = 5$ and $Y$ is hyperbolic, this is a contradiction
to [Wu2, Theorem 1.1] if both $Y(r_a)$ are $\bdd$-reducible, to [GW2]
if one of them is $\bdd$-reducible and the other is annular, and to
[GW3] if both $Y(r_a)$ are annular and atoroidal.  (The main theorem
of [GW3] said that if $\Delta = 5$ then $Y$ and $r_1, r_2$ are listed
in one of the three possibilities in [GW3, Theorem 1.1], but in that
case both $Y(r_a)$ are toroidal.)

Case 3.  {\it $T_0 \not \subset \bdd Y$.\/} Let $X_1$ be the component
of $X|(\bdd Y)$ containing $T_0$, and let $X_2$ be the closure of $X -
X_1$, so $X = X_1 \cup X_2$.  Since $M_i$ contains no non-separating
torus (Lemma 23.3), $T_1 = X_1 \cap X_2$ is a single torus.  Since
$X(r_a)$ is atoroidal, $T_1$ must be compressible in $X_1(r_a)$ for
$a=1,2$.  Thus Lemma 23.2 and the fact that $X(r_a)$ contains no lens
space summand implies that $X_1(r_a)$ is a solid torus for $a=1,2$, as
in Lemma 23.2(1); moreover, since $\Delta(r_1, r_2) > 1$, by Lemma
23.2(1) $X_1$ is a $(p,q)$ cable space, and by [CGLS, Theorem 2.4.3]
$\Delta(r, r_a) = 1$ for $r$ the cabling slope.  It is easy to see
that the meridian slopes $m_a$ of $X_1(r_a)$ satisfies $\Delta(m_1,
m_2) = p \Delta(r_1,r_2) \geq 5$.  Now $X(r_a) = X_2(m_a)$, so by
Cases 1 and 2 above applied to $X_2$ we see that this case is also
impossible.  \qed

Denote by $c\cdot d$ the minimal intersection number between the two
isotopy classes of simple closed curves on a surface represented by
$c$ and $d$, respectively.  If $\varphi: F\to F$ is a homeomorphism
and $K$ a curve on the surface $F$, then $K$ is said to be {\it
  $\varphi$-full\/} if for any essential curve $c$ on $F$ there is
some $i$ such that $c \cdot \varphi^i(K) \neq 0$.

If $K$ is a knot in a 3-manifold $Y$ with a preferred
meridian-longitude, denote by $Y(K,p/q)$ the manifold obtained from
$Y$ by $p/q$ surgery on $K$.  Let $X = F \times I/\psi$ be an
$F$-bundle over $S^1$ with gluing map $\psi$, let $F_t = F \times t$,
$t \in I = [0,1]$, and let $K$ be an essential curve on $F_{1/2}$.
Then there is a preferred meridian-longitude pair $(m,l)$ on $\bdd
N(K)$, with $l$ the slope of $F_{1/2} \cap \bdd N(K)$.

\begin{lemma} Let $X = F \times I/\psi$.  Let $\eta: F\times I \to
F_0$ be the projection, $\varphi = \eta\circ \psi$, and $K$ an
essential curve on $F_{1/2}$.  If $\eta(K) \subset F_0$ is
$\varphi$-full and $q>1$, then $X(K, p/q)$ is hyperbolic.
\end{lemma}

\proof Let $A_i$ be an annulus in $X$ with $\bdd A_i = K \cup K_i$,
where $i=0,1$ and $K_i \subset F_i$.  Let $V_i$ be a regular
neighborhood of $A_i$.  Put $Y = F \times I$.  Then $Y = V_1 \cup W$,
where $W$ is homeomorphic to $F \times I$, and $V_1 \cap W$ is an
annulus $A'$.  After $p/q$ surgery on $K$ we have $Y(K, p/q) = V_1(K,
p/q) \cup W$.  Note that $V(K, p/q)$ is a solid torus with $A'$ an
annulus on $\bdd V_1(K, p/q)$ running $q$ times along the longitude.
By an innermost circle outermost arc argument one can show that $Y(K,
p/q)$ is irreducible, $\bdd$-irreducible, atoroidal, and any essential
annulus $A_2$ can be isotoped to be disjoint from $K_1$, i.e.\ $\bdd
A_2 \cdot K_1 = 0$.  Moreover, if $A_2$ has at least one boundary
component on $F_1$ then $A_2$ is either vertical in the sense
that it is isotopic to $c \times I \subset F\times I$ for some curve
$c \subset F$, or isotopic to $A'$ and hence has both boundary curves
parallel to $K_1$.  Similarly, using $A_0$ and $V_0$ one can show that
$\bdd A_2 \cdot K_0 = 0$, and if $A_2$ is not vertical and $A_2 \cap
F_0 \neq \emptyset$ then it has both boundary curves parallel to
$K_0$.

The above facts imply that $X(K, p/q) = Y(K, p/q)/\psi$ is
irreducible.  Since the non-separating surface $F_0$ cuts $X(K, p/q)$
into $Y(K, p/q)$, which is not an $I$-bundle, we see that $X(K, p/q)$
is not Seifert fibered.  It remains to show that $X(K, p/q)$ is
atoroidal.

If $T$ is an essential torus in $X(K, p/q)$ then it can be isotoped so
that $T \cap Y(K, p/q) = Q$ is a set of essential annuli.  Let
$C_i = Q \cap F_i$.  We claim that for any curve $c \subset
C_0$, $\varphi(c)$ is isotopic to a curve in $C_0$.

We have $\psi(C_0) = C_1$, so $\psi(c) \subset C_1$.  If $\psi(c)$
belongs to a vertical annulus $Q_1$ then $\varphi(c) = \eta(\psi(c))
\cong Q_1 \cap F_0 \subset C_0$.  If $\psi(c)$ belongs to a
non-vertical annulus then by the property proved above, $\psi(c)$ is
isotopic to $K_1$, so $\varphi(c) \cong \eta(K_1) = K_0$.  Note that
if $Q$ has a non-vertical component with boundary on $F_1$ then the
fact that $C_0, C_1$ have the same number of components implies that
there is also a non-vertical component $Q_0$ with boundary on $F_0$,
and we have shown that each component $c'$ of $\bdd Q_0$ is isotopic
to $K_0$, so $\varphi(c) \cong c' \subset \bdd Q_0 \subset C_0$.  This
completes the proof of the claim.

Let $c$ be a component of $C_0$.  We have shown above that $c \cdot
K_0 = 0$ for any $c \subset C_0$.  Applying the above to
$\varphi^{-1}$ we see that there is a curve $c' \subset C_0$ such that
$\varphi(c') \cong c$.  By induction we have $c \cdot \varphi^i(K_0) =
c' \cdot \varphi^{i-1}(K_0) = 0$ for all $i$, which is a contradiction
to the assumption that $K_0 \cong \eta(K)$ is $\varphi$-full and hence
$c \cdot \varphi^i(K_0) \neq 0$ for some $i$.  \qed

\begin{lemma}  The manifold $M_5$ is hyperbolic.
\end{lemma}

\proof Let $W = M_5 | F_b$, let $F_+, F_-$ be the two copies of $F_b$
in $W$, and let $A$ be the annulus $T_0 | \bdd F_b$.  Then $W$ is
obtained from $Y = F_+ \cup F_- \cup A$ by attaching faces of $\ga$
and then some 3-cells.  

Two faces of $\ga$ are {\it parallel\/} if their boundary curves are
parallel on $Y$.  Since parallel faces cobound a 3-cell in $W$, we
need only attach one such face among a set of parallel faces.  From
Figure 11.10 one can check that the four bigons are parallel faces,
and the two 6-gons are parallel to each other.  Therefore $W$ is
obtained from $Y$ by attaching one bigon $\sigma_1$, one 6-gon
$\sigma_2$ and then a 3-cell.  Let $\sigma_1$ be the bigon on Figure
11.10 between the edges $B$ and $G$, and assume that the edge $B
\subset F_+$.

Cutting $W$ along $\sigma_1$, we obtain a manifold $W_1$ with boundary
a torus, and it contains the $6$-gon $\sigma_2$.  Therefore it is a
solid torus such that the remnant of $F_+$, denoted by $F'_+$, runs
along the longitude three times.  If we replace $\sigma_2$ and the
attached 3-cell by a solid torus $J$ with meridian intersecting $F'_+$
in one essential arc then $W_1$ becomes a $F'_+ \times I$ and $W$
becomes $X = F_+ \times I$.  Therefore $W = X(K, p/q)$, where $K$ is
the core of $J$, and $q=3$.  Let $\psi: F_- \to F_+$ be the gluing
map, $\eta: F_+ \times I \to F_-$ the projection, and $\varphi = \eta
\circ \psi$.  By Lemma 23.9 we need only show that the curve $K_-$ on
$F_-$ isotopic to $K$ is $\varphi$-full.

In $M_5$ the bigon $\sigma_1$ has boundary edges $B \cup G$ on $F_b$,
as shown in Figure 11.10(b).  Suppose $B \subset F_+$ and $G\subset
F_-$ when we consider $\sigma_1$ as a bigon in $F_+\times I$.  Then
$\psi$ maps the curve $B$ on $F_-$ to the curve $B$ on $F_+$, which is
mapped to $G$ on $F_-$ by $\eta$.  Therefore $\varphi: F_- \to F_-$
maps $B$ to $G$.  Since $F_- | B$ is an annulus and $B$ is disjoint
from the curve $K_-$ above, this determines $K_-$.  Also
$\varphi(K_-)$ is the curve on $F_-$ disjoint from $\varphi(B) = G$,
so $K_-$ intersects $\varphi(K_-)$ transversely at a single point,
cutting $F_-$ into an annulus.  Therefore $K_-$ is $\varphi$-full, and
the result follows.  \qed

\begin{lemma}  The manifold $M_4$ is hyperbolic.
\end{lemma}

\proof The proof is similar to that of Lemma 23.10.  In this case $W =
M_4 | F_b$ is obtained from $F_+ \cup F_- \cup A_1 \cup A_2$ by
attaching two bigons $\sigma_1, \sigma_2$ and one $4$-gon $\sigma_3$,
so $W = X(K,p/q)$ with $q=2$, where $X = F_+ \times I$ and $K$ is
disjoint from $\sigma_1, \sigma_2$.  Choose $\sigma_1, \sigma_2$ to be
the bigons in Figure 11.9(a) bounded by $H\cup E$ and $E\cup N$,
respectively.  Then $F_-$ can be identified with $F_b$, and $\sigma_1,
\sigma_2$ intersects $F_-$ in the edges $E$ and $N$, respectively.
These cut $F_-$ into an annulus containing the curve $K_-$ isotopic to
the knot $K$ in $X = F_+ \times I$.  The map $\varphi: F_- \to F_-$
maps the edges $E$ and $N$ in Figure 11.9(b) to $H$ and $E$,
respectively, so $\varphi(K_-)$ is the curve in the annulus $F_- |
(H\cup E)$.  The curves $K_-$ and $\varphi(K_-)$ intersect transversely
at a single point, cutting $F_-$ into a neighborhood of $\bdd F_-$,
hence $K_-$ is $\varphi$-full, and $M_4$ is hyperbolic by Lemma 23.9.
\qed

\begin{lemma} Let $F$ be a closed orientable surface of genus 2, and
let $\alpha, \beta$ be two non-separating simple closed curves on
$F$, intersecting minimally, cutting $F$ into disks.  Let $X$ be
obtained from $F\times I$ ($I=[0,1]$) by attaching a 2-handle along
$\alpha \times \{0\}$.  Identify $F$ with $F \times \{1\} \subset
F\times I$, and let $T = \bdd X - F$.  Then

(1) a compressing disk $D$ of $F$ intersects $\beta$ at least 3
times; and 

(2) an incompressible annulus $A$ in $X$ with $\bdd A \subset F$ and
$\bdd A \cap \beta = \emptyset$ is boundary parallel.
\end{lemma}

\proof (1) Let $E$ be the disk in $X$ bounded by $\alpha \times 1$,
cutting $X$ into $X' = T \times I$.  Note that $X$ is a compression
body, and $\{E \}$ is the unique (up to isotopy) complete disk system
for $X$.

By assumption $\alpha \cup \beta$ cuts $F$ into disks, hence $|\alpha
\cap \beta| \geq 3$.  Since $\alpha$ intersects $\beta$ minimally, we
may choose a hyperbolic structure on the surface $F$ so that $\alpha,
\beta$ are geodesics.  Let $D$ be a compressing disk for $F$ in $X$.
Up to isotopy we may assume that $\gamma = \bdd D$ is a geodesic or a
slight push off of a geodesic if it is isotopic to $\alpha$ or
$\beta$.  Then both $|\gamma \cap \alpha|$ and $|\gamma \cap \beta|$
are minimal up to isotopy.

We may assume that $D \cap E$ consists of arcs.  If $D \cap E \neq
\emptyset$, by taking an arc that is outermost on $D$, surgering $E$
along the corresponding outermost disk in $D$, and discarding one of
the resulting components, we get a new disk $E'$ having fewer
intersections with $D$, such that $\{E'\}$ is a complete disk system
for $X$.  Since $\{E\}$ is the unique complete disk system for $X$,
$E'$ is isotopic to $E$.  Since $|\bdd E' \cap \bdd D| < |\bdd E \cap
\bdd D| = |\alpha \cap \gamma|$ and $\bdd E'$ is isotopic to $\bdd E$,
this is a contradiction to the fact that $|\gamma \cap \alpha|$ is
minimal up to isotopy.  Therefore $D\cap E = \emptyset$.  Hence $D$
either (a) is parallel to $E$, or (b) cuts off a solid torus
containing $E$.

In case (a) $\bdd D$ is a parallel copy of $\alpha$, so $|\bdd D \cap
\beta| = |\alpha \cap \beta| \geq 3$ and we are done.  In case (b),
let $F_1$ be the punctured torus on $F$ bounded by $\gamma$ which does
not contain $\alpha$.  If $|\gamma \cap \beta| \leq 2$ then $F_1$
contains at most one arc of $\beta$, so it contains an essential loop
disjoint from $\alpha \cup \beta$, which is a contradiction to the
assumption that $\alpha \cup \beta$ cuts $F$ into disks.

(2) Let $A$ be an incompressible annulus in $X$ with $\bdd A \subset F
- \beta$.  We may assume that $\alpha, \beta$ are hyperbolic
geodesics, and each component of $\bdd A$ is either a geodesic, or a
slight push off of a geodesic if it is parallel to $\alpha$, $\beta$
or another component of $\bdd A$.  Thus both $|\bdd A \cap \alpha|$
and $|\bdd A \cap \beta|$ are minimal; in particular, $\bdd A \cap
\beta = \emptyset$.  As in (1), this implies that $A \cap E$ consists
of essential arcs on $A$.  If $A \cap E = \emptyset$ then  $\bdd A$
lies in $F |(\alpha \cup \beta)$, but since $A$ is incompressible
while each component of $F |(\alpha \cup \beta)$ is a disk, this is
impossible.  Therefore we may assume that $A \cap E$ is a non-empty
set $C$ of essential arcs on $A$.

Let $B_i$ be a component of $A|C$.  Then $B_i$ is a disk in $X' =
T\times I$, so $\bdd B_i$ is a trivial loop on $T'$, bounding a disk
$B'_i$ on $T'$.  Let $E_1, E_2$ be the two copies of $E$ on $T'$.  If
$B'_i \cap (E_1 \cup E_2)$ is a single disk then one can use a disk
component of $B'_i \cap F$ to isotope $A$ to reduce $|\bdd A \cap \bdd
E| = |\bdd A \cap \alpha|$, which is a contradiction to the minimality
of $|\bdd A \cap \alpha|$.  Therefore $B'_i \cap (E_1 \cup E_2)$
consists of two disks, and $B' \cap F$ is a single disk $P_i$.  One
can check that $\cup P_i$ is an annulus on $F$ parallel to $A$.  \qed

\begin{lemma}
The manifold $M_{14}$ is hyperbolic.
\end{lemma}

\proof Cutting $M_{14}$ along the surface $F_b$, we obtain two
manifolds $X_1, X_2$, where $X_1$ is the one containing the four bigon
faces of $F_a$, and $X_2$ contains the two 4-gon faces of $F_a$.  Let
$\sigma_1, \sigma_2$ be the bigons on $F_a$ bounded by the edges
$E\cup F$ and $B\cup Y$ respectively in Figure 20.6, and let
$\sigma_3$ be the 4-gon bounded by the edges $B \cup C \cup Y \cup F$.
Note that any other face of $F_a$ is parallel in $X_i$ to one of
these.

Let $A_i = X_i \cap T_0$.  Then $X_1$ is obtained from the genus 2
surface $F_b \cup A_1$ by attaching $\sigma_1, \sigma_2$ and then a
3-cell, hence it is a handlebody of genus 2 because $\bdd \sigma_1,
\bdd \sigma_2$ are disjoint nonparallel nonseparating curves on
$F_b\cup A_1$.  The core of $A_1$ is a curve on $\bdd X_1$ such that
after attaching a 2-handle to $X_1$ along $A_1$ we get the manifold on
the side of $\hat F_b$ which contains no torus boundary component,
hence from Figure 22.14(b) we see that it is the double branched cover
of a Montesinos tangle $T(2,2)$, which is a twisted $I$-bundle over
the Klein bottle.  This implies that the surface $F_b = \bdd X_1 -
A_1$ is incompressible in $X_1$.

Now consider $X_2$.  Let $F$ be the genus 2 surface $F_b \cup A_2$,
$\alpha$ the boundary of $\sigma_3$, and $\beta$ the core of $A_2$.
Then $\alpha$ intersects $\beta$ minimally at four points.  From
Figure 20.6(b) we see that the edges $B,C,Y,F$ cut the surface $F_b$
into two disks, hence $\alpha \cup \beta$ cuts $F$ into disks.  $X_2$
is obtained from $F\times I$ by attaching a 2-handle along the curve
$\alpha \times \{0\}$.  Therefore it satisfies the conditions of Lemma
23.12.  In particular, $F_b$ is incompressible in $X_2$.

Since $X_1$ is a handlebody and $X_2$ is a compression body, they are
irreducible and atoroidal.  Since $M_{14}$ is obtained by gluing $X_1,
X_2$ along the incompressible surface $F_b$, $M_{14}$ is also
irreducible.  It is well known that an incompressible surface in a
Seifert fiber space is either vertical, and therefore an annulus or
torus, or horizontal, in which case it intersects all boundary
components.  Since the surface $F_b$ satisfies neither condition, we
see that $M_{14}$ is not Seifert fibered.  It remains to show that
$M_{14}$ is atoroidal.

Assume $M_{14}$ is toroidal and let $T_1$ be an essential torus in
$M_{14}$ intersecting $F_b$ minimally.  Since $X_i$ is atoroidal,
$T_1$ intersects $X_i$ in incompressible annuli.  A component $A'_2$
of $T_1 \cap X_2$ is an incompressible annulus in $X_2$ disjoint from
$\beta$, hence by Lemma 23.12 it is parallel to an annulus $A''$ on
$\bdd X_2$.  If $A'' \subset F_b$ then $T_1$ can be isotoped to reduce
$|T_1 \cap F_b|$, which is a contradiction to the minimality
assumption.  Therefore $A'' \supset \beta$ and hence $A'' \supset
A_2$, so each component of $\bdd A''$ is parallel to a component of
$\bdd F_b$.  Since this is true for all components of $T_1 \cap X_2$,
we see that each component of $T_1 \cap F_b$ is parallel to a
component of $\bdd F_b$.

Now let $A'_1$ be a component of $T_1 \cap X_1$.  By the above, the
two boundary components of $A'_1$ are parallel on $\bdd X_1$.  Since
$X_1$ is a handlebody, $A'_1$ is parallel to an annulus $A''_2$ on
$\bdd X_1$.  For the same reason as above, it must contain the annulus
$A_1$.  This is true for all components of $T_1 \cap X_1$. Let $A'_i$
be a component of $T_1 \cap X_i$ which is closest to $A_i$.  Then
$\bdd A'_1 = \bdd A'_2$, hence $T_1 = A'_1 \cup A'_2$.  It follows
that $T_1$ is parallel to $T_0$, contradicting the assumption that
$T_1$ is essential in $M_{14}$.  \qed

\begin{thm} The manifolds $M_i$ in Definition 21.3 are all
hyperbolic.  \end{thm}

\proof This follows from [GW1, Theorem 1.1] for $i=1,2,3$, and from
Lemmas 23.7, 23.8, 23.10, 23.11 and 23.13 for $i>3$.  \qed

\section {Toroidal surgery on knots in $S^3$}

Recall that each of the manifolds $M_1, M_2, M_3$ admits two toroidal
Dehn fillings $r'_i, r''_i$ on a torus boundary component $T_0$ with
distance $4$ or $5$.  These are the exteriors of the links $L_1, L_2,
L_3$ in Figure 24.1.  Let $L_i = K'_i \cup K''_i$, where $K'_i$ is the
left component of $L_i$.  Let $T_1 = \bdd N(K'_i)$, and $T_0 = \bdd
N(K''_i)$.

\bigskip
\leavevmode

\centerline{\epsfbox{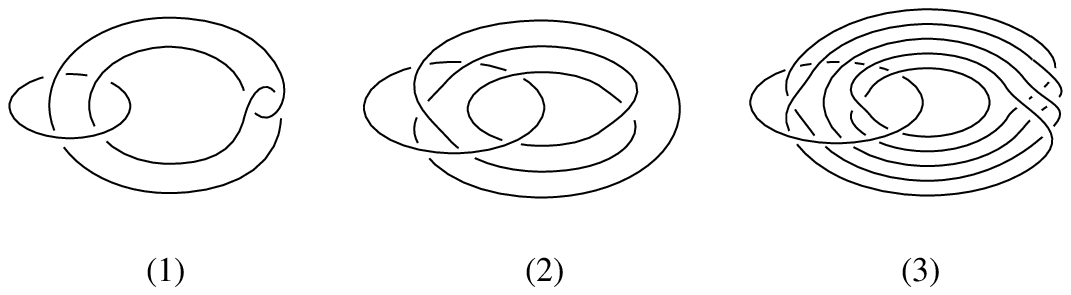}}
\bigskip
\centerline{Figure 24.1}
\bigskip

Each $M_i$ has a pair of toroidal slopes $r'_i, r''_i$ on $T_0$.
These are given in [GW1, Theorem 7.5] and shown in Figures 7.2, 7.4
and 7.5 of [GW1].

\begin{lemma}
  With respect to the preferred meridian-longitude pair of $K''_i$,
  the slopes $r'_i, r''_i$ are given as follows, up to relabeling.

(1) $r'_1 = 0$ and $r''_1 = 4$.

(2) $r'_2 = -2$ and $r''_2 = 2$.

(3) $r'_3 = -9$ and $r''_3 = -23/2$.
\end{lemma}

\proof (1) This is basically proved in [GW1, Lemma 7.1].  It was shown
that $M_1$ is the double branched cover of the tangle $Q_1$ in Figure
7.2(c) of [GW1].  Let $m$ be the meridian, $l$ the preferred
longitude, and $l'$ the blackboard longitude of the diagram of $K''_1$
in [GW1, Figure 7.2(a)].  Calculating the linking number of $l'$ with
$K''_1$ in Figure 7.2(a) we see that $l' = 2m + l$.  Let $\eta: M_1
\to Q_1$ be the branched covering map.  If $r$ is a slope on $T_0$
then $\eta(r)$ is a curve of a certain slope on the inside boundary
sphere, which will be denoted by a number in $\Bbb Q_0 = \Bbb Q \cup
\{\infty\}$.  One can check that $\eta(m) = 0/1$, and $\eta(l') =
1/0$.  The two toroidal slopes $r'_1, r''_1$ map to slopes $-1/2$ and
$1/2$, as shown in Figure 7.2(d) and (e) of [GW1].  We have
$\varphi(-2m+l') = (-2\times 0 + 1)/(-2 \times 1 + 0) = -1/2$, and
$\varphi(2m+l') = (2\times 0 + 1)/(2 \times 1 + 0) = 1/2$.  Therefore
$r'_1 = -2m + l' = -2m + (2m + l) = l$ and $r''_1 = 2m + l' = 4m + l$.

(2) This is similar to (1), using [GW1, Figure 7.4] instead.  We have
$\eta(m) = 0/1$, $\eta(l') = 1/2$, $\eta(r'_2) = 1/0$, $\eta(r''_2) =
1/4$, and $l = l'$.  Therefore $r'_2 = -2m + l$ and $r''_2 = 2m + l$.

(3) Use 7.5(k) to denote [GW1, Figure 7.5(k)].  7.5(a) shows that
$l' = l - 6m$.  A careful tracking of $l'$ during the modification
from 7.5(a) to 7.5(b) then to 7.5(c) shows that $\eta(m) = 0/1$ and
$\eta(l') = 1/3$.  From 7.5(c) and 7.5(e) we see that $\eta(r'_3) =
1/0$ and $\eta(r''_3) = 2/5$.  Therefore $r'_3 = -3m + l' = -9m + l$,
and $r''_3 = -m + 2l' = -m + 2(l-6m) = -13m + 2l$.
\qed

\begin{lemma}
  Suppose $K$ is a hyperbolic knot in $S^3$ admitting two toroidal
  Dehn surgeries $K(r_1), K(r_2)$ with $\Delta(r_1, r_2) = 4$ or $5$.
  Then there is an $i \in \{1,2,3\}$ and a  slope $s$ on $T_1$ of
  $M_i$ such that $(E(K), r_1, r_2) \cong (M_i(s), r'_i, r''_i)$.
\end{lemma}

\proof Let $F_a$ be an essential punctured torus in $M_K = S^3 - \Int
N(K)$ such that $\hat F_a$ is an essential torus in $K(r_a)$, chosen
so that $|\bdd F_a|$ is minimal.  By Theorem 21.4 the triple $(E(K),
r_1, r_2)$ is equivalent to either $(M_i, r'_i, r''_i)$ with $1\leq
i\leq 14$, or to $(M_i(s), r'_i, r''_i)$ for some $i=1,2,3,14$.
Therefore we need only show that the manifold $M_i$ ($i=4,...,14$) is
not the exterior of a knot or link in $S^3$.

When $i=4$, the surface $F_b$ has two boundary circles on $T_0$ with
the same orientation.  Let $A$ be an annulus on $T_0$ connecting these
two boundary components.  Then $F_b \cup A$ is a non-orientable closed
surface in $M_4$.  It follows that $M_4$ cannot be the exterior of a
knot in $S^3$ because $S^3$ contains no embedded non-orientable
surface.

For $i=5$, let $V_b$ be the Dehn filling solid torus of $M_5(r_b)$.
Then the $\Bbb Z_2$ homology group $H$ of $V_b \cup F_b$ is generated
by $\alpha$, $x$ and $y$, where $\alpha$ is the core of $V_b$, and $x,
y$ are represented by the edges $A$ and $E$ in Figure 11.10(b),
respectively.  A bigon in Figure 11.10(a) gives the relation $x=y$.
Consider the quotient group $H'$ obtained from $H$ by identifying $x$
with $y$.  Then $H' = \Bbb Z_2 \oplus \Bbb Z_2$ is generated by $\alpha$
and $x$.  Each corner of Figure 11.10(a) represents the element
$\alpha$, and each edge represents $x$ in $H'$.  Since each face in
Figure 11.10(a) has an even number of edges and an even number of
corners on its boundary, it represents $0$ in $H'$.  Therefore
$$H_1(M_5(r_b), \Bbb Z_2) = H_1 (V_b \cup F_b \cup F_a, \Bbb Z_2) = H' =
\Bbb Z_2 \oplus \Bbb Z_2.$$  
Since the $\Bbb Z_2$ homology of any manifold obtained by Dehn surgery
on a knot in $S^3$ is either trivial or $\Bbb Z_2$, it follows that
$M_5$ is not a knot exterior.

Now assume that $M_i$ is the exterior of a knot $K$ in $S^3$ for some
$i\geq 6$.  By Theorem 23.14 $K$ is hyperbolic.  Put $\Bbb Q_0 = \Bbb
Q \cup \{\infty \}$.  A number in $\Bbb Q_0$ is represented by $p/q$,
where $p,q$ are coprime integers, and $q\geq 0$.  Given a
meridian-longitude pair $(m,l)$ and $r = p/q \in \Bbb Q_0$, denote by
$K(r)$ the manifold obtained by surgery on $K$ along the slope
$pm+ql$.  There is a one to one correspondence $\eta: \Bbb Q_0 \to
\Bbb Q_0$ such that $\Delta(\eta(r), \eta(s)) = \Delta(r, s)$,
and $K(\eta(r))$ is the double branched cover of $Q_i(r)$, which is
the manifold $X_i(r)$ given in Lemma 22.2.  Since $K(\eta(r_3))$ is a
lens space, by the Cyclic Surgery Theorem [CGLS, p.237] the slope
$\eta(r_3)$ is an integer slope with respect to the preferred
meridian-longitude of $K$.  To simplify the calculation, let $l =
\eta(r_3)$.

By [GLu, Theorem 1.1] $\eta(r_1)$ and $\eta(r_2)$ are integer or half
integer slopes.  Suppose $\eta(r_i) = p_i/q_i$.  Then $p_3/q_3 = 0/1$.
By the above we have $q_i = 1$ or $2$ for $i=1,2$.  By Lemma 22.2,
$|p_1| = \Delta(r_1, r_3) = 1$, and $|p_2| = \Delta(r_2, r_3) \leq 2$.

If $|p_2| = 1$ then $4 \leq \Delta(r_1, r_2) = \Delta(\eta(r_1),
\eta(r_2)) = |p_1q_2 - p_2q_1|$ implies that $q_1 = q_2 = 2$.  This is
a contradiction to [GWZ, Theorem 1], which says that a hyperbolic knot
in $S^3$ admits at most one non-integral toroidal surgery.

We now have $|p_2|=2$, so $\eta(r_2) = p_2/q_2 = \pm 2/1$.  Since
$\Delta(r_1, r_2) = |p_1q_2 - p_2q_1| = |\pm 1 - (\pm 2)q_1| \geq 4$,
we must have $p_1/q_1 = \mp 1/2$, and $\Delta = 5$.  From Lemma 22.2
we see that for $i\in \{6,...,13\}$, the only $M_i$ satisfying
$\Delta(r_2, r_3) = 2$ and $\Delta(r_1, r_2)=5$ are the ones with
$i=7$, $10$ or $11$.  

Consider the case $i=10$.  Let $r_0$ be the slope such that
$\eta(r_0)$ is the meridian slope $1/0$.  Then we have $\Delta(r_0,
r_i) = \Delta(\eta(r_0), \eta(r_i)) = \Delta(1/0, p_i/q_i) = q_i$.
Therefore by the above we have $\Delta(r_0, r_i) = q_i = 2, 1, 1$ for
$i=1,2,3$, respectively.  By Lemma 22.2 we have $r_1 = 0/1$, $r_2 =
-5/2$ and $r_3 = 1/0$.  Let $r_0 = p'/q'$.  Then we have
\begin{eqnarray*}
& & \Delta(r_0, r_1) = |p'| = 2 \\
& & \Delta(r_0, r_2) = |2p'+5q'| = 1 \\
& & \Delta(r_0, r_3) = q' = 1 
\end{eqnarray*}
These equations have a unique solution $r_0 = -2/1$.  One can check
that $Q_{10}(-2/1)$ is the 2-bridge knot $K_{2/7}$, so its double
branched cover is $L(7,2) \neq S^3$, which is a contradiction.

The tangles $Q_7$, $Q_{11}$ and $Q_{14}$ in Figures 22.7, 22.11 and
22.14 have a circle component.  If $M_i$ is a knot exterior in
$S^3$ then there is a slope $r$ such that the double branched cover of
$Q_i(r)$ is $S^3$.  Since each of $Q_7(r)$ and $Q_{11}(r)$ has at
least two components, its double branched cover has nontrivial $\Bbb
Z_2$ homology [Sa, Sublemma 15.4], so $M_7$ and $M_{11}$ are not knot
exteriors in $S^3$.  Similarly $M_{14}$ is not the exterior of a link
in $S^3$.  
\qed

\begin{lemma} (1) Let $i\in \{1,2,3\}$.  If $r_1, r_2$ are toroidal
  slopes of $M_i$ on $T_0$ with $\Delta(r_1, r_2) \geq 4$, then
  $\{r_1, r_2\} = \{r'_i, r''_i\}$. 

(2) The slope $-7$ is a solid torus filling slope on $T_0$ of $M_3$,
and there is an orientation preserving homeomorphism of $M_3$ which
interchanges the two solid torus filling slopes $\{1/0, -7/1\}$ and
the two toroidal slopes $\{-9, -13/2\}$.  
\end{lemma}

\proof Since $M_1, M_2, M_3, M_{14}$ are the only ones in Definition
21.3 with two boundary components, by Theorem 21.4 $(M_i, r_1, r_2)$
is equivalent to one of the $(M_j, r'_j, r''_j)$ with $j=1,2,3,14$.
Since $M_1,M_2,M_3$ are link complements in $S^3$ and by Lemma 24.2
$M_{14}$ is not, we have $j\neq 14$.  Computing $H_1(M_j, T_1)$ shows
$H_1(M_1, T_0) = \Bbb Z$, $H_1(M_2, T_0) = \Bbb Z_3$, and $H_1(M_1,
T_0) = \Bbb Z_5$, hence we must have $j=i$.  

By definition there is a homeomorphism $\varphi: (M_i, r_1, r_2) \to
(M_i, r'_i, r''_i)$, up to relabeling of $r_1, r_2$.  For $i=1,2$, by
[Ga1] and [Be] the knot $K''_i$ has no nontrivial solid torus surgery,
hence $\varphi(m) = m$, where $m$ is a meridian of $K''_i$.  Since
$\Delta(m, r'_i) = \Delta(m, r''_i) = 1$, by the homeomorphism we also
have $\Delta(m, r_1) = \Delta(m, r_2) = 1$, so $r_1, r_2$ are also
integer slopes.  It follows that if $\{r_1, r_2\} \neq \{r'_1,
r''_1\}$ then there is a pair of toroidal slopes with distance at
least 5.  Since $M_1, M_2$ is not homeomorphic to $M_3$, this is a
contradiction to Theorem 21.4 and [Go].

Now suppose $i=3$.  By an isotopy one can deform the tangle in [GW1,
Figure 7.5(c)], which is shown in Figure 24.2(a), to the one in Figure
24.2(b), which is invariant under the $\pi$ rotation $\psi$ along the
forward slash diagonal.  The $1/0$ slope in Figure 24.2(b)
corresponds to the $1/2$ slope in Figure 24.2(a), which, by the
proof of Lemma 24.1(3), lifts to the slope $-7m + l$ on $T_0$.  The
two toroidal slopes $1/0$ and $5/2$ for the tangle in Figure 22.2(a)
correspond to the slopes $-1/2$ and $2$ in Figure 24.2(b), which are
interchanged by $\psi$.  It follows that $\psi$ lifts to an
orientation preserving homeomorphism $\psi': M_3 \to M_3$, which
interchanges the two solid torus filling slopes $\{1/0, -7/1\}$ and
the two toroidal slopes $\{-9, -13/2\}$.  In fact, $\psi'$ is
represented by the matrix 
$$A = \left(
\array{rr}
7 & 1 \\
-1 & 0
\endarray
\right)
$$
in the sense that if $A(p,q)^t = (p',q')^t$ (where $B^t$ denotes the
transpose of the matrix $B$) then $\psi'(pm + ql) = p'm + q'l$.  

Solid torus surgeries on knots in a solid torus have been completely
classified by Gabai [Ga1] and Berge [Be].  It was shown that there is
only one knot admitting two nontrivial solid torus surgeries, which is
a 7-braid.  Since $K''_3$ is a 5-braid, we see that $m = 1/0$ and $m'
= -7/1$ are the only solid torus filling slopes on $T_0$.  Therefore
the homeomorphism $\varphi: (M_3, r'_3, r''_3) \to (M_3, r_1, r_2)$
must map the set of two curves $\{m, m'\}$ to itself, possibly with
the orientation of one or both of the curves reversed.  If $\varphi$
preserves the orientation of $m'$ and reverses the orientation of $m$
then $\varphi(r'_3) = 9/1$ would also be a toroidal slope, which is a
contradiction to [Go] because $\Delta(-9/1, 9/1) = 18 > 8$.  Similarly
$\varphi$ cannot preserve the orientation of $m$ while reversing the
orientation of $m'$.  Therefore $\varphi$ is orientation preserving
and its induced map on the set of slopes on $T_0$ is either the
identity map, which fixes $\{r'_3, r''_3\}$, or the same as that
induced by $\psi'$ above, which interchanges $\{r'_3, r''_3\}$.  
\qed

\bigskip
\leavevmode

\centerline{\epsfbox{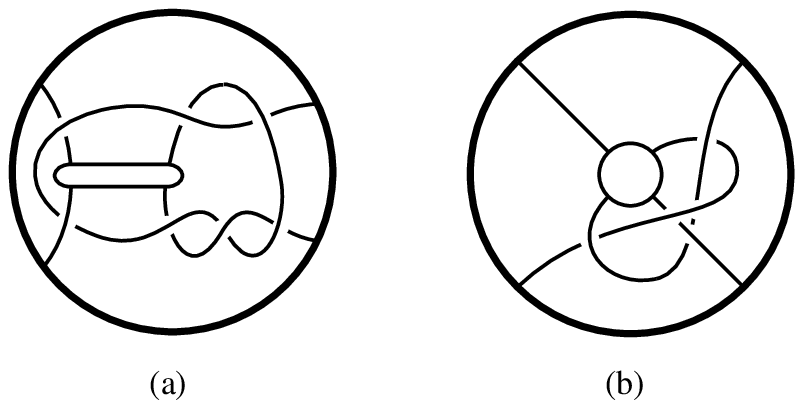}}
\bigskip
\centerline{Figure 24.2}
\bigskip

\begin{thm} Suppose $K$ is a hyperbolic knot in $S^3$ admitting two
toroidal surgeries $K(r_1), K(r_2)$ with $\Delta(r_1, r_2) \geq 4$.
Then $(K, r_1, r_2)$ is equivalent to one of the following, where $n$
is an integer.

(1) $K = L_1(n)$, $r_1 = 0$, $r_2 = 4$.

(2) $K = L_2(n)$, $r_1 = 2-9n$, $r_2 = -2-9n$.

(3) $K = L_3(n)$, $r_1 = -9 - 25n$, $r_2 = -(13/2) - 25n$.

(4) $K$ is the Figure 8 knot, $r_1 = 4$, $r_2 = -4$.
\end{thm}

\proof By [Go] the Figure 8 knot $L_1(-1)$ is the only hyperbolic knot
in $S^3$ admitting two toroidal surgeries of distance at least 6, so
we assume $\Delta = 4$ or $5$.  By Lemma 24.2 there is a homeomorphism
$\varphi: (M_i(s), \{r'_i, r''_i\}) \cong (E(K), \{r_1, r_2\})$ for some
$i=1,2,3$ and $s \subset T_1$.  We only need to show that $s = 1/n$
because the slopes $r_i$ can then be calculated using Lemma 24.1 and
the Kirby calculus [Ro, p.267].

If Dehn filling on $T_1$ of $\bdd M_i$ along slope $s$ produces a knot
exterior $E(K) = M_i(s)$, then the meridian-longitude of $K$ may be
different from that of $K''_i$ on $T_0$.  We use $(m'',l'')$ (resp.\
$(m,l)$) to denote a meridian-longitude pair of $K''_i$ (resp.\ $K$) in
$S^3$.

\medskip

\noindent
Claim 1.  {\it If $E(K) = M_1(s)$ for some $s$ on $T_1$ of $\bdd M_1$
then $s=1/n$.}

\medskip

Since the linking number between the two components of $L_1$ is $0$, a
$p/q$ Dehn filling on $T_1$ produces a manifold $M_1(p/q)$ with
$H_1(M_1(p/q), \Bbb Z) = \Bbb Z \oplus \Bbb Z_p$, hence $M_1(p/q)$ is
a knot complement only if $|p| = 1$.  It follows that $K = L_1(n)$,
where $n=qp$.  

\medskip

\noindent
Claim 2.  {\it If $E(K) = M_2(s)$ for some $s$ on $T_1$ of $\bdd M_2$
then $s=1/n$ for some $n$.}

\medskip

As before, let $L_2 = K'_i \cup K''_i$.  Let $M = E(K)$.  Assume $s
= p/q$ and $|p|>1$.  We have $K(m'') = L_2(s,m'') = K'_i(s) = L(p,q)$.
Therefore by the Cyclic Surgery Theorem [CGLS], $m''$ is an integer
slope with respect to $(m,l)$, say $m'' = am + l$.  By [GLu] the
toroidal slopes $r_1, r_2$ of $K$ are integer or half integer slopes
with respect to $(m,l)$.  Recall that $\varphi(r'_2) = r_1$ and
$\varphi(r''_2) = r_2$.  Since $m''$ is an integer slope, $r_1, r_2$
cannot both be integer slopes, otherwise $4 = \Delta(r_1, r_2) \leq
\Delta(r_1, m'') + \Delta(m'', r_2) = 2$, which is a contradiction.
Also by [GWZ] they cannot both be half integer slopes.

Now assume $r_1$ is an integer slope and $r_2$ is a half integer
slope with respect to $(m,l)$.  Since $m''$ is an integer slope, we
may choose $l=m''$.  Then $r_1 = p_1 m + l$ and $r_2 = p_2 m + 2l$,
so $\Delta(r_1, m'') = \Delta(r_2, m'') = 1$ implies $p_1, p_2 = \pm
1$.  But then $\Delta(r_1, r_2) = |2p_1 - p_2| \leq 3$, a
contradiction.  

\medskip
\noindent
Claim 3.  {\it If $E(K) = M_3(s)$ then there is an integer $n$ and a
  homeomorphism $\eta: (M_3(1/n), \{r'_3, r''_3\}) \to (E(K), \{r_1,
  r_2\})$.}

\medskip

By Lemma 24.3 $M_3(-7)$ is a solid torus, and the meridian
slope $m''$ and the slope $r=-7$ are the only solid torus filling
slopes on $T_0$.  If $\varphi(m'') = m$ then $S^3 = K(m'') = K'_i(s)$
implies that $s=1/n$ for some $n$, so $\eta = \varphi$ is the required
map.  If $\varphi(r) = m$, let $\psi$ be the orientation preserving
homeomorphism of $M_3$ given in Lemma 24.3, which maps $m''$ to $r$.  By
Lemma 24.3 $\psi$ interchanges the slopes $r'_3, r''_3$.  Let $s' =
\psi^{-1}(s)$.  Then $\varphi \circ \psi :(M_3(s'), r''_3, r'_3) \cong
(E(K), r_1, r_2)$ maps $m''$ to $m$.  As above this implies that $s' =
1/n$, hence $\eta = \varphi\circ \psi$ is the required map.

We now assume that $\varphi(m'') \neq m$ and $\varphi(r) \neq m$.
Note that $K(m'')$ and $K(r)$ can be obtained from the solid tori
$M_3(m'')$ and $M_3(r)$ by $s$ filling on $T_1$, so they have cyclic
$\pi_1$, hence by [CGLS] $r, m''$ are integer slopes of $K$.  Choose
$l = m''$.  Since $\Delta(m'', r)=1$, we may assume $r = 1/1$ up to
rechoosing the orientation of $l$.  The toroidal slopes $r'_3 = -9$
satisfy $\Delta(r'_3, m'') = 1$ and $\Delta(r'_3, r) = 2$, which
implies $r'_3 = -1/1$ or $1/3$ with respect to $(m,l)$.  The second is
impossible by [GLu].  Similarly the fact that $\Delta(r''_3, m'') = 2$
and $\Delta(r''_3, r) = 1$ implies that $r''_3 = 2$ with respect to
$(m,l)$.  But then we have $5 = \Delta(r'_3, r''_3) = \Delta(-1, 2) =
3$, a contradiction.  
\qed

\begin{cor}
  A hyperbolic knot $K$ in $S^3$ has at most four toroidal surgeries.
  If there are four, then they are consecutive integers.
\end{cor}

\proof By [GLu] a toroidal slope of $K$ must be integer or half
integer, and if it is a half integer then $K$ is a Eudave-Mu\~noz
knot.  By [T1, Corollary 1.2], if $K$ is a Eudave-Mu\~noz knot then it
has at most three toroidal slopes, hence the result is true if $K$ has
a half integer toroidal slope.  Therefore we may assume that all
toroidal slopes of $K$ are integer slopes.  The result follows if
$\Delta(r,s) < 4$ for all pairs of toroidal slopes $(r,s)$ of $K$.
Therefore by Theorem 24.4 we need only show that if $K$ is either
$L_1(n)$ or $L_2(n)$ for some $n$ then $K$ has at most three integer
toroidal slopes.

If $K$ is the knot $L_1(n)$ in Theorem 24.4(1) then by [BW] it has
exactly two toroidal slopes unless it is the Figure 8 knot, which has
three toroidal slopes.

Now consider a knot $K = L_2(n)$ in Theorem 24.4(2) and let $r$ be an
integral toroidal slope of $K$ other than $r_1, r_2$ in the Theorem.
Since $\Delta(r_i, r)\leq 4$, $r$ must be between $r_1$ and $r_2$.
Denote by $M_2(p/q)$ the $p/q$ filling on $T_0$ with respect to the
preferred meridian-longitude pair of $L_2$.  By the proof of Lemma
24.1(2), $M_2(-1)$ is the double branched cover of $Q_2(1)$.  Using
the tangle in [GW1, Figure 7.4(c)] one can check that $Q_2(1)$ is a
Montesinos tangle $T(1/2, -2/5)$, therefore $M_2(-1)$ is a small
Seifert fiber space with orbifold $D^2(2,5)$.  Since $L_2(n)(-1-9n)$
is obtained from $M_2(-1)$ by Dehn filling on $T_1$ and contains no
non-separating surface, it is atoroidal.  Because of symmetry ($L_2$
is amphicheiral), $M_2(1)$ is homeomorphic to $M_2(-1)$, so
$L_2(n)(1-9n)$ is also atoroidal.  It follows that the only possible
integer toroidal slopes of $L_2(n)$ are $j - 9n$ for $j=-2, 0, 2$.
This completes the proof.  (Actually it can be shown that $-9n$ is not
a toroidal slope of $L_2(n)$ either, so it has at most two integer
toroidal slopes.)  \qed

The following corollary is an immediate consequence of Theorem 24.4.

\begin{cor} Let $K$ be a hyperbolic knot in $S^3$ which admits two
  toroidal surgeries along slopes $r_1, r_2$, and $\Delta =
  \Delta(r_1, r_2) \geq 4$.  Then one of the $r_i$ is an integer, and
  the other one is an integer if $\Delta \neq 5$, and a half integer
  if $\Delta = 5$.
\end{cor}

Although there are infinitely many hyperbolic 3-manifolds $M$ with
toroidal fillings $M(r), M(s)$ at distance 4 or 5, we have shown that
they all come from finitely many cores $X(r,s)$ as defined in Section
21.

\begin{qtn}
  Are there only finitely many cores $X(r,s)$ of toroidal Dehn
  fillings on hyperbolic 3-manifolds with $\Delta(r,s) = 3$?
  $\Delta(r,s) = 2$?
\end{qtn}

We observe that the answer to Question 24.7 in the case
$\Delta(r,s)=1$ is almost certainly `no'; here is an outline of an
argument. Let $N$ be a closed irreducible $3$-manifold with a unique
incompressible torus $T$ up to isotopy. Let $F$ be a once-punctured
torus, regarded as a disk with two bands. It is intuitively clear
that, for any positive integer $n$, by tangling the bands in a
sufficiently complicated fashion we can construct an embedding $F_n$
of $F$ in $N$ so that if $K_n = \bdd F_n$, then $N-K_n$ is hyperbolic,
and $K_n$ cannot be isotoped to meet $T$ in fewer than $n$ points. Let
$M_n = N-\Int N(K_n)$, and let $r,s$ on $\bdd M_n$ be the meridian of
$K_n$ and the longitudinal slope defined by $F_n$, respectively. Then
$\Delta(r,s)=1$, $M_n(r) = N$ is toroidal by definition, and $M_n(s)$
contains the non-separating torus $\hat F_n = F_n \cup D$, where $D$
is a meridian disk of $V_s$. Hence, if we make sure that $M_n(s)$ does
not contain a non-separating sphere, then $M_n(s)$ is also toroidal.
Since the number of intersections of $K_n$ with $T$ is at least $n$,
the triples $(M_n,r,s)$ cannot all come from only finitely many cores.

\bigskip

\noindent
Dept.\ Mathematics, University of Texas at Austin, Austin, TX 78712
\\
Email: {\it gordon@math.utexas.edu}

\bigskip

\noindent
Department of Mathematics,  University of Iowa,  Iowa City, IA 52242
\\
Email: {\it wu@math.uiowa.edu}


\begin{thebibliography}{CGLS} 

\bibitem[Be]{Be} J.~Berge, {\em The knots in $D^2\times S^1$ which have
nontrivial Dehn surgeries that yield $D^2\times S^1$}, Topology
Appl.\ {\bf 38} (1991), 1--19.

\bibitem[BW]{BW} M.~Brittenham and Y-Q.~Wu, {\em The classification of
  exceptional Dehn surgeries on 2-bridge knots}, Comm.\ Anal.\
  Geom.\ {\bf 9} (2001), 97--113.

\bibitem[CGLS]{CGLS} M.\ Culler, C.\ Gordon, J.\ Luecke and P.\ Shalen, 
{\em  Dehn surgery on knots},  Annals Math.\ {\bf 125} (1987),
237--300.

\bibitem[Eu]{Eu} M.~Eudave-Mu\~noz, {\em Non-hyperbolic manifolds obtained
  by Dehn surgery on hyperbolic knots}, Geometric Topology
  (Athens, GA, 1993), AMS/IP Stud.\ Adv.\ Math., 2.1,
  Amer. Math. Soc., Providence, RI, 1997, pp.\ 35--61.

\bibitem[Ga1]{Ga1} D. Gabai, {\em Surgery on knots in solid tori},
Topology {\bf 28} (1989), 1-6.

\bibitem[Ga2]{Ga2} ------, {\em 1-bridge braids in solid tori}, Top.\
Appl.\ {\bf 37} (1990), 221-235.

\bibitem[Go]{Go} C.\ Gordon, {\em Boundary slopes of punctured tori in
3-manifolds},  Trans.\ Amer.\ Math.\ Soc.\ {\bf 350}  (1998),
1713--1790.

\bibitem[GLi]{Gli} C.\ Gordon and R.\ Litherland, {\em Incompressible
planar surfaces in 3-manifolds},  Topology Appl.\  {\bf 18} (1984),
121-144.

\bibitem[GLu]{Glu} C.\ Gordon and J.\ Luecke, {\em Non-integral toroidal
Dehn surgeries},  Comm.\ Anal.\ Geom.\  {\bf 12}  (2004),
417--485.

\bibitem[GW1]{GW1} C.\ Gordon and Y-Q.\ Wu,  {\em Toroidal and annular
Dehn fillings},  Proc.\ London Math.\ Soc.\ {\bf 78}  (1999),
662-700.

\bibitem[GW2]{GW2} ------,  {\em Annular and boundary reducing Dehn
surgery},  Topology {\bf 39}  (2000), 531-548.

\bibitem[GW3]{GW3} ------,  {\em Annular Dehn fillings},  Comment.\
Math.\ Helv.\ {\bf 75} (2000), 430--456.

\bibitem[GWZ]{GWZ} C.\ Gordon, Y-Q.\ Wu and X.\ Zhang,  {\em Non-integral
toroidal surgery on hyperbolic knots in $S^3$},  Proc.\ Amer.\ Math.\
Soc.\  {\bf 128}  (2000), 1869-1879.

\bibitem[HM]{HM} C.\ Hayashi and K.\ Motegi, {\em Only single twists on
unknots can produce composite knots}, Trans.\ Amer.\ Math.\ Soc.\
{\bf 349} (1997), 4465--4479.

\bibitem[L1]{L1} S.\ Lee, {\em Exceptional Dehn fillings on hyperbolic
3-manifolds with at least two boundary components}, preprint.

\bibitem[L2]{L2} ------, {\em Dehn fillings yielding Klein bottles},
  preprint.

\bibitem[MaS]{Mas} D.~Matignon and N.~Sayari, {\it Klein slopes on
    hyperbolic 3-manifolds}, preprint.

\bibitem[MeS]{MS} W.\ Meeks III and P.\ Scott, {\em Finite group
actions on 3-manifolds}, Invent.\ Math.\ {\bf 86} (1986), 287-346.

\bibitem[Oh]{Oh} S.\ Oh,  {\em Reducible and toroidal manifolds
obtained by Dehn filling},  Topology Appl.\ {\bf 75} (1997),
93--104.

\bibitem[Ro]{Ro} D.~Rolfsen, {\em Knots and Links} Publish or Perish,
1990.

\bibitem[Sa]{Sa} M.\ Sakuma,  {\em Homology of abelian coverings of
links and spatial graphs},  Canadian J.\ Math.\  {\bf 47}  (1995),
201-224.

\bibitem[Sch]{Sch}  M.\ Scharlemann,  {\em Producing reducible
$3$-manifolds by surgery on a knot},  Topology {\bf 29}  (1990),
481--500.

\bibitem[Sct]{Sct} P.\ Scott, {\em The geometry of 3-manifolds},
Bull.\ London Math.\ Soc.\ {\bf 15} (1983), 401--487.

\bibitem[T1]{T1} M.~Teragaito, {\em Distance between toroidal
  surgeries on hyperbolic knots in the 3-sphere}, preprint,
  arXiv:math.GT/0312201.

\bibitem[T2]{T2} ------, {\em Toroidal Dehn fillings on large
  hyperbolic 3-manifolds}, preprint, arXiv:math.GT/0508250.

\bibitem[Wu1]{Wu1} Y-Q.\ Wu,  {\em Dehn fillings producing reducible
manifolds and toroidal manifolds},  Topology {\bf 37}  (1998),
95--108.

\bibitem[Wu2]{Wu2} ------, {\em Incompressibility of surfaces in surgered
3-manifolds},  Topology {\bf 31}  (1992), 271--279.

\bibitem[Wu3]{Wu3} ------, {\em 
Sutured manifold hierarchies, essential laminations, and Dehn
surgery}, J.\ Diff.\ Geom., {\bf 48} (1998), 407--437.

\end{thebibliography}
\enddocument